\documentclass[10pt]{article}
\topmargin - 0.7 in \overfullrule = 0pt
\usepackage{comment}
\usepackage{graphics}
\usepackage{graphicx}
\usepackage{color}
\usepackage{amsfonts}
\usepackage{dsfont}
\usepackage{amsmath}
\usepackage{amstext}
\usepackage{amsopn}
\usepackage{amsbsy}
\usepackage{amscd}
\usepackage{amsxtra}
\usepackage{amsthm}
\usepackage{amssymb}
\usepackage{esint}
\usepackage{abstract,lipsum,mathrsfs}
\usepackage[new]{old-arrows}
\usepackage{cases}
\usepackage{tocloft}
\usepackage{ragged2e} 
\usepackage{enumitem}
\usepackage[utf8]{inputenc}

\newcommand{\h}{\hspace}
\newcommand{\p}{\partial}
\numberwithin{equation}{section}

\title{Pattern Formation in Landau--de Gennes Theory}
\author{Ho--Man Tai \hspace{6pt} and \hspace{6pt} Yong Yu}
\date{}

\newtheorem{thm}{Theorem}[section]
\newtheorem{defn}[thm]{Definition}
\newtheorem{lem}[thm]{Lemma}
\newtheorem{prop}[thm]{Proposition}

\newtheorem{rmk}[thm]{Remark}

\def\Xint#1{\mathchoice
  {\XXint\displaystyle\textstyle{#1}}%
  {\XXint\textstyle\scriptstyle{#1}}%
  {\XXint\scriptstyle\scriptscriptstyle{#1}}%
  {\XXint\scriptscriptstyle\scriptscriptstyle{#1}}%
  \!\int}
\def\XXint#1#2#3{{\setbox0=\hbox{$#1{#2#3}{\int}$}
  \vcenter{\hbox{$#2#3$}}\kern-.5\wd0}}

\def\dashint{\Xint-}

\newcounter{alphasect}
\def\alphainsection{0}

\let\oldsection=\section
\def\section{%
  \ifnum\alphainsection=1%
    \addtocounter{alphasect}{1}
  \fi%
\oldsection}%

\renewcommand\thesection{%
  \ifnum\alphainsection=1%
    \Alph{alphasect}%
  \else
    \arabic{section}%
  \fi%
}%
\newenvironment{alphasection}{%
  \ifnum\alphainsection=1%
    \errhelp={Let other blocks end at the beginning of the next block.}
    \errmessage{Nested Alpha section not allowed}
  \fi%
  \setcounter{alphasect}{0}
  \def\alphainsection{1}
}{%
  \setcounter{alphasect}{0}
  \def\alphainsection{0}
}%

\begin{document}
\maketitle
\begin{abstract}
\noindent \textbf{\normalsize Abstract:}  We study the spherical droplet problem in 3D--Landau de Gennes theory with finite temperature. By rigorously constructing the biaxial--ring solutions and split--core--segment solutions, we theoretically confirm the numerical results of Gartland--Mkaddem in \cite{GM00}. The structures of disclinations are also addressed.
\end{abstract}
\tableofcontents

\section{Introduction}\label{sec_1}\vspace{0.5pc}

The order parameter $\mathscr{Q}$ in the Landau--de Gennes (LdG for short) theory takes values in $\mathbb{S}_0$. Here $\mathbb{S}_0$ is the $5$--dimensional linear vector space consisting of all real $3 \times 3$ symmetric traceless matrices. Given a point $x$ in the domain, we let $\lambda_1(x) \leq \lambda_2(x) \leq \lambda_3(x)$ be the three eigenvalues of $\mathscr{Q}(x)$.  The state of liquid crystal can then be classified into four types according to the quantitative relationships of these three eigenvalues. 
\begin{enumerate}
\item[$\mathrm{(1)}$.] $\mathscr{Q}$ is called isotropic at $x$ if $\lambda_1(x) = \lambda_2(x) = \lambda_3(x) = 0$;
\item[$\mathrm{(2)}$.] $\mathscr{Q}$ is called negative uniaxial at $x$ if  $\lambda_1(x) < \lambda_2(x) = \lambda_3(x)$;
\item[$\mathrm{(3)}$.] $\mathscr{Q}$ is called positive uniaxial at $x$ if  $\lambda_1(x) = \lambda_2(x) < \lambda_3(x)$;
\item[$\mathrm{(4)}$.] $\mathscr{Q}$ is called biaxial at $x$ if $\lambda_1(x) < \lambda_2(x) < \lambda_3(x)$.
\end{enumerate}In the LdG theory, the director field of a liquid crystal is defined to be the normalized eigenvector associated with the largest eigenvalue of $\mathscr{Q}$. At a continuous point of $\mathscr{Q}$, the director field can be locally oriented if $\mathscr{Q}$ is biaxial or positive uniaxial at this point. Moreover, the oriented director field is continuous at this point. However, at the negative uniaxial or isotropic location, the director field might lose its orientability and continuity. It would be difficult for us to extend the definition of the director field to the negative uniaxial or isotropic locations continuously. In the liquid crystal theory, the locations where the oriented director field is misfit are called "disclinations" of the liquid crystal material. 

\subsection{Spherical droplet problem and some existing works}
In this article, we consider the so--called spherical droplet problem. Throughout the remaining arguments, $B_R(x)$ is the open ball in $\mathbb{R}^3$ with center $x$ and radius $R$. The ball $B_R(0)$ is simply denoted by $B_R$. For a $\mathbb{S}_0$--valued order parameter $\mathscr{Q}$ on $ B_R$, its LdG energy functional in the one--constant limit is read as follows:
\begin{align}
\int_{B_R} \dfrac{1}{2}\h{0.5pt}\big|\h{0.5pt}\nabla \mathscr{Q}\h{0.5pt}\big|^2 - \dfrac{a^2}{2}\big|\h{0.5pt}\mathscr{Q}\h{0.5pt}\big|^2-\sqrt{6} \h{1pt} \text{tr}\left(\mathscr{Q}^3\right)+\dfrac{1}{2}\h{0.5pt}\big|\h{0.5pt}\mathscr{Q}\h{0.5pt}\big|^4.
\label{Landau-de Gennes energy}
\end{align}
Here $-a^2$ is the reduced temperature. For any matrix $A \in \mathbb{S}_0$, the norm of $A$ is defined by $|\h{1pt} A \h{1pt}|^2 := \mathrm{tr}\big(A^2\big)$. Let $\mathrm{I}_3$ be  the $3 \times 3$ identity matrix. The Euler--Lagrange equation associated with  $(\ref{Landau-de Gennes energy})$ is 
\begin{align}
- \Delta \mathscr{Q} = a^2\mathscr{Q} +3\sqrt{6}\left(\mathscr{Q}^2 - \dfrac{1}{3}\h{1.5pt}|\mathscr{Q}|^2\h{1.5pt} \mathrm{I}_3 \right) - 2 |\mathscr{Q}|^2 \mathscr{Q} \qquad \text{in } B_R.
\label{el-eq of D}
\end{align}In the spherical droplet problem,  (\ref{el-eq of D}) is supplied with the following strong anchoring condition: \begin{align}\label{bdy cond of Q} \mathscr{Q} = \dfrac{\sqrt{3}}{2} \h{1pt}a \h{1pt}H_a \left( e_r \otimes e_r - \dfrac{1}{3} \mathrm{I}_3 \right) \h{15pt}\text{on $\p B_R$.}
\end{align}Note that $e_r$ is the radial direction in $\mathbb{R}^3$.  $H_a$ is the  constant given below: \begin{align}\label{H_a} H_a := \dfrac{3 + \sqrt{9 + 8a^2}}{2 \sqrt{2} \h{1pt}a}.
\end{align}

In 1988, a radial hedgehog solution to (\ref{el-eq of D})--(\ref{bdy cond of Q})  was considered  in \cite{SS88} by Schopohl--Sluckin. The solution   can be represented by $f(r)\big( 3 e_r \otimes e_r - \mathrm{I}_3\big)$ with $f(r)$ solving an ODE induced from (\ref{el-eq of D}). Here $r$ denotes the radial variable in $\mathbb{R}^3$. The radial hedgehog solution has an isotropic core at the origin. Right after \cite{SS88}, in 1989, Penzenstadler--Trebin \cite{PT89} discovered  that there may have a solution to (\ref{el-eq of D})--(\ref{bdy cond of Q}) with biaxial--ring disclination. In some parameter regime, hedgehog solution is not stable. The isotropic core of the hedgehog solution can be broadened to a  disclination ring with topological charge $1/2$.  It was until 2000 that the split--core--segment disclination   was numerically found by Gartland and Mkaddem in \cite{GM00}. Besides being broadened to a disclination ring, the isotropic core of the hedgehog solution can also be splitted into a segment disclination with strength $1$. Up to now, the core structure of the hedgehog solution is well understood. If we replace the domain $B_R$ with the whole space $\mathbb{R}^3$, then the asymptotic behavior of the entire hedgehog solution at spatial infinity is also known. See \cite{AF97, FG00, G97, MP10, L13}. To our surprise, there are few theoretical studies on the core structures of the biaxial--ring disclination and the split--core--segment disclination. A first attempt was made  in \cite{Y20}. It shows that there are two families of solutions to (\ref{el-eq of D})--(\ref{bdy cond of Q}) which can be suitably rescaled so that in the low--temperature limit ($a \to \infty$), one family of the rescaled solutions  converges to a limiting state with biaxial--ring disclination, while another family of the rescaled solutions  converges to a limiting state with split--core--segment disclination. For the two limiting states, the asymptotic behaviors of their director fields near associated disclinations are   explicitly calculated for the first time. Note that solutions to the LdG equation with ring--like disclinations have also been considered in \cite{DMP211, DMP212} when the order parameter $\mathscr{Q}$ satisfies the Lyuksyutov constraint. Interested readers may also refer to \cite{BZ11, BPP12, INSZ13, C15, INSZ15, ABL16, FRSZ16, INSZ16, C17, CL17, HMP17, ACS21} for various recent studies on solutions to the LdG equation with disclinations.  
 

\subsection{Axially symmetric formulation of LdG equation}

Before we give our main results, let us introduce an axially symmetric formulation of the LdG equation. This formulation can help us reduce the degrees of freedom of $(\ref{el-eq of D})$  from five to three. Firstly we define some notations.
Let $x=(x_1,x_2,z)$ denote a point in $\mathbb{R}^3$ and $(r,\phi,\theta)$ be the spherical coordinates of $\mathbb{R}^3$. Here $\phi \in [\h{0.5pt}0,\pi\h{0.5pt}]$ is the polar angle, while $\theta \in [\h{1pt}0,2\pi\h{0.5pt})$ is the azimuthal angle. Moreover, we denote by $\rho$ the radial variable in the $(x_1, \h{.5pt} x_2)$--plane and hence  $(\rho, z, \theta)$ are the cylindrical coordinates of $\mathbb{R}^3$. As for the linear vector space $\mathbb{S}_0$, it is spanned by the following five matrices: \begin{align*} &L_1 = \dfrac{1}{\sqrt{2}} \begin{pmatrix}\h{0.5pt}0 & 0 & 1 \\
0 & 0 & 0 \\
1 & 0 & 0 \h{0.5pt}\end{pmatrix}, \h{5pt} L_2 = \dfrac{1}{\sqrt{2}} \begin{pmatrix}\h{0.5pt}0 & 1 & 0 \\
1 & 0 & 0 \\
0 & 0 & 0 \h{0.5pt}\end{pmatrix},  \\[2mm]
&L_3 = \dfrac{1}{\sqrt{2}} \begin{pmatrix}\h{0.5pt}0 & 0 & 0 \\
0 & 0 & 1 \\
0 & 1 & 0 \h{0.5pt}\end{pmatrix}, \h{5pt} L_4 = \dfrac{1}{\sqrt{6}} \begin{pmatrix}\h{0.5pt}-1 & 0 & 0 \\
0 & -1 & 0 \\
0 & 0 & 2 \h{0.5pt}\end{pmatrix},  \h{5pt} L_5 = \dfrac{1}{\sqrt{2}} \begin{pmatrix}\h{0.5pt}1 & 0 & 0 \\
0 & -1 & 0 \\
0 & 0 & 0 \h{0.5pt}\end{pmatrix}.
\end{align*} Using the notations above, we put the unknown order parameter $\mathscr{Q}$ into the ansatz:
\begin{align}\label{ansatz}
&\mathscr{Q} = \dfrac{a}{\sqrt{2}} \Big\{ v_1 \big( \cos 2 \theta \h{1pt} L_5 + \sin 2 \theta \h{1pt} L_2 \big) + v_2 \h{1pt}L_4 + v_3 \big( \cos \theta \h{1pt} L_1 + \sin \theta \h{1pt} L_3 \big) \Big\}, \\[2mm] & \h{12pt}\text{where for $j = 1, 2, 3$, \,} v_j = v_j (\rho,z) \h{3pt}\text{are real--valued unknown functions.}  \nonumber
\end{align}
Meanwhile, we define
\begin{align}\label{scaled var}
u(x)=v(Rx), \qquad \text{for any\,\,} x \in B_1 \end{align} and let $\mathscr{L}$ be the augmented operator given as follows: \begin{align}\label{aug operator}\mathscr{L}\left[ V\right] := \left(V_1 \cos 2\theta, V_1 \sin 2\theta, V_2, V_3 \cos \theta, V_3 \sin \theta\right)^\top, \h{15pt}\text{for any $V = \left(V_1, V_2, V_3\right)^\top \in \mathbb{R}^3$.}\end{align}
Hence, the $\mathscr{Q}$--variable in (\ref{ansatz}) solves $(\ref{el-eq of D})$ if and only if $w := \mathscr{L}\left[u\right]$ satisfies
\begin{equation}
- \mu^{-1} \Delta w  = \dfrac{3}{\sqrt{2}}\h{1pt}\nabla_w S\left[w\right]  -a\h{0.5pt}\big(\h{0.2pt}|\h{0.5pt}w\h{0.5pt}|^2-1 \h{0.2pt}\big)\h{0.5pt}w \qquad \text{in } B_1.
\label{el-eq of u_a, intro}
\end{equation}Here and throughout the article, $\mu = a R^2$ is a fixed positive constant. For any $w = (w_1, w_2, w_3, w_4, w_5)^\top$, the degree--3 homogeneous polynomial $S\left[w\right]$ is defined by \begin{align*}S\left[w\right] := - w_3\left( w_1^2 + w_2^2\right) + \sqrt{3}\h{1pt}w_2w_4w_5 + \dfrac{1}{2}w_3\left(w_4^2 + w_5^2\right)+ \dfrac{1}{3} w_3^3  + \dfrac{\sqrt{3}}{2}w_1\left(w_4^2 - w_5^2\right).
\end{align*}Note that in terms of the variable $u = (u_1, u_2, u_3)$ in (\ref{scaled var}), the three eigenvalues of the matrix $a^{-1}\mathscr{Q}\big(R x\big)$, where $\mathscr{Q}$ is given in (\ref{ansatz}), can be explicitly calculated as follows: \begin{align}
\left\{\begin{aligned}
\lambda_1 &= - \dfrac{1}{2} \left(u_{1}+\dfrac{1}{\sqrt{3}} u_{2}\right) ; \\[2mm]
\lambda_2 &= \dfrac{1}{4} \left(u_{1}+\dfrac{1}{\sqrt{3}} \h{1pt} u_{2}\right)
- \dfrac{1}{4}\sqrt{\big(u_{1}-\sqrt{3} \h{1pt} u_{2}\big)^2+4\h{0.5pt}u_{3}^2}; \\[2mm]
\lambda_3 &= \dfrac{1}{4} \left(u_{1}+\dfrac{1}{\sqrt{3}} \h{1pt} u_{2}\right)
+ \dfrac{1}{4}\sqrt{\big(u_{1} - \sqrt{3} \h{1pt} u_{2}\big)^2+4\h{0.5pt}u_{3}^2}.
\end{aligned}\right.
\label{evalue of D with u= u^+_a,b}
\end{align}In light of the boundary condition for $\mathscr{Q}$ in (\ref{bdy cond of Q}), the system (\ref{el-eq of u_a, intro}) is subjected to the Dirichlet boundary condition: \begin{align}\label{boundary condition of u_a, intro} w = \mathscr{L}\left[\h{1pt} U_a^*\h{1pt}\right] \h{15pt}\text{on $\p B_1$,}\end{align}where with the  $H_a$ defined in (\ref{H_a}),
\begin{equation*}
U_a^* := H_a U^* := H_a\left( \dfrac{\sqrt{3}}{2} \sin^2\phi,\,  \dfrac{3}{2}\left(\cos^2 \phi-\dfrac{1}{3}\right) ,\, \sqrt{3}\sin \phi \cos \phi\right)^\top.
\end{equation*}It can be shown that  (\ref{el-eq of u_a, intro}) is the Euler--Lagrange equation of the energy functional:
\begin{align}\label{def of mathcal e a mu}
\mathcal{E}_{a, \mu}\left[w\right] := \int_{B_1} f_{a, \mu}\left(w\right), \h{15pt}\text{where $f_{a, \mu}\left(w\right) := \big| \nabla w \big|^2 + \mu\left[D_a-3\sqrt{2}\h{0.5pt}S\left[ w \right]+\dfrac{a}{2}\h{0.5pt}\big(\h{0.5pt}|\h{0.2pt}w\h{0.2pt}|^2-1\big)^2\right]$.}
\end{align}Notice that the constant $D_a$ is given by \begin{align}\label{def of Da} D_a :=\dfrac{27}{16\h{0.2pt}a^3} \left[ 1+\dfrac{4\h{0.2pt}a^2}{3}+\left(1+\dfrac{8\h{0.2pt}a^2}{9}\right)^{3/2} \right].
\end{align} It is chosen so that the minimum value of $2D_a- 6\sqrt{2}\h{0.5pt}S\left[w\right] +a\h{0.5pt}\big(\h{0.5pt}|\h{0.2pt}w\h{0.2pt}|^2-1\big)^2$ equals  $0$ when this polynomial is restricted on the set $\big\{ \mathscr{L}\left[\h{0.5pt}x\h{0.5pt}\right] : x \in \mathbb{R}^3\big\}$.   \vspace{0.2pc}

We particularly focus on a special class of axially symmetric solutions  to the boundary value problem (\ref{el-eq of u_a, intro}) and (\ref{boundary condition of u_a, intro}). These solutions are axially symmetric and meanwhile satisfy some reflective symmetry with respect to the $(x_1, x_2)$--plane.   \begin{defn} A $5$--vector field $w$ is $\mathscr{R}$--axially symmetric on some ball $B_r$ if it satisfies the following three conditions on $B_r$:\begin{enumerate}\item[$\mathrm{(1).}$] $w = \mathscr{L}\left[ \h{1pt}u\h{1pt}\right]$ on $B_r$. Here $u$ is a $3$--vector field  depending only on the $(\rho,z)$--variables; 
\item[$\mathrm{(2).}$]  $u_1$ and $u_2$ are even with respect to the $z$--variable; 
\item[$\mathrm{(3).}$] $u_3$ is odd with respect to the $z$--variable.\end{enumerate}
\end{defn}

 
\subsection{\normalsize Main results}\label{sec_1_2}\vspace{0.5pc} 

In this section, we introduce the main results of this article. Note that for both biaxial--ring solutions and split--core solutions discussed below, their director fields might coexist the biaxial--ring and split--core disclinations. To simplify the expositions of our main theorems, for the biaxial--ring solutions, we focus on the half--degree ring structure of their disclinations. For the split--core solutions, we focus on the split--core structure of their disclinations. The biaxial--ring disclinations in the the split--core solutions can be similarly studied as the biaxial--ring solutions. The split--core  disclinations in biaxial--ring solutions can also be  similarly considered as the split--core solutions.\vspace{0.2pc}

Firstly, we discuss the biaxial--ring solutions. \begin{thm}[\bf Biaxial--ring solutions and their ring disclinations] There exists a  constant $a_0 > 0$ so that for all $a > a_0$, the followings hold: \begin{enumerate} \item[$\mathrm{(1)}$.] There exists a $\mathscr{R}$--axially symmetric solution, denoted by $w_{a, +} = \mathscr{L}[\h{1pt}u_{a, +}\h{1pt}]$, to the boundary value problem (\ref{el-eq of u_a, intro}) and (\ref{boundary condition of u_a, intro}).  The origin is not zero of $w_{a, +}$. If $w_{a, +}$ has zeros, then all zeros of $w_{a, +}$  must be   on the $z$--axis. There must have even number (might be $0$) of zeros of $w_{a, +}$ on the set $\big\{ \big(0, 0, z\big) : 0 < z < 1\big\}$;
\item[$\mathrm{(2)}$.] Let $\lambda^+_{a; j}$ ($j = 1, 2, 3$) be the three eigenvalues in (\ref{evalue of D with u= u^+_a,b}) computed with $u = u_{a, +}$ there. Recall the  $\mathscr{Q}$ in (\ref{ansatz}) and denote by $ \mathscr{Q}^+_a$ the tensor field $a^{-1}\mathscr{Q}\big(Rx\big)$ with $v(y) = u_{a, +}\big(R^{-1}y\big)$. Then $\lambda^+_{a; j}$ ($j = 1, 2, 3$) are the three eigenvalues of $\mathscr{Q}^+_a$. There exist a $\delta_0 \in \big(0, 1/2\big)$ independent of $a$ and  a $\rho_a \in \big(\delta_0, 1 - \delta_0\big)$ so that  $\mathscr{Q}^+_a$ is negative uniaxial on  $\mathscr{C}_a := \big\{ \big(x_1, x_2, 0\big) : x_1^2 + x_2^2 = \rho_a^2\big\}$ with $\lambda^+_{a; 1} < 0 < \lambda^+_{a; 2} = \lambda^+_{a; 3}$ on the circle $\mathscr{C}_a$. Fix an $\epsilon > 0$ and denote by $\mathscr{T}_{a, \epsilon}$ the torus $\big\{ x \in \mathbb{R}^3 : \mathrm{dist}\big(x, \mathscr{C}_a\big) \leq \epsilon \big\}$. There exists a small $\epsilon$ depending on $a$ so that $\mathscr{Q}^+_a$ is biaxial on $\mathscr{T}_{a, \epsilon} \setminus \mathscr{C}_a$ with $\lambda^+_{a; 1} < \lambda^+_{a; 2} < \lambda^+_{a; 3}$ on $\mathscr{T}_{a, \epsilon} \setminus \mathscr{C}_a$; 
\item[$\mathrm{(3)}$.] Given a $3$--vector field $u$, we define the following vector field with unit length:\begin{align}
\begin{split}
 \kappa \h{0.5pt}[\h{.5pt} u \h{.5pt}] \h{2pt}:=\h{2pt} &\dfrac{\sqrt{2}}{2}\left(1+\dfrac{u_{1}-\sqrt{3}u_{2}}{\sqrt{(u_{1}-\sqrt{3}u_{2})^2 +4 u_{3}^2}}\right)^{1/2} e_\rho\\ &+\dfrac{\sqrt{2}u_{3}}{\sqrt{(u_{1}-\sqrt{3}u_{2})^2 +4 u_{3}^2}}\left(1+\dfrac{u_{1}-\sqrt{3}u_{2}}{\sqrt{(u_{1}-\sqrt{3}u_{2})^2 +4 u_{3}^2}}\right)^{-1/2}e_z.
\end{split}
\label{hat kappa_3, intro}
\end{align}Here $e_\rho := \left(\frac{x_1}{\rho}, \frac{x_2}{\rho}, 0\right)^\top$ and $e_z := (0, 0, 1)^\top$. Then the director field of   $\mathscr{Q}_a^+$ on   $\mathscr{T}_{a, \epsilon} \setminus \mathscr{C}_a$ can be oriented and expressed by  $\kappa\h{0.5pt}[\h{0.5pt}u_{a, +}\h{0.5pt}]$. It is the normalized eigenvector of $\mathscr{Q}_a^+$ associated with the eigenvalue $\lambda^+_{a; 3}$;
\item[$\mathrm{(4)}$.] The circle $\mathscr{C}_a$ is a ring disclination of the director field $\kappa\h{0.5pt}[\h{0.5pt}u_{a, +}\h{0.5pt}]$. In terms of the $(\rho,z)$--variables, its structure is described as follows. Let $x_a = (\rho_a, 0)$ be a point on the $(\rho, z)$--plane and denote by $D_r(x_a)$  the open disk on the $(\rho, z)$--plane with center $x_a$ and radius $r$. Here $\rho_a$ is given in the item (2). Fix an arbitrary $r \in (0, \epsilon)$. When we approach $x_{a, r} := (\rho_a - r,0)$ along the semi--circle $\p^- D_r (x_a) := \p D_r (x_a) \cap \big\{z \leq 0 \big\}$, the director field $\kappa\h{.5pt}[\h{.5pt} u_{a, +} \h{.5pt}]$ tends to $-e_z$. When we approach $x_{a, r}$ along $\p^+ D_r (x_a) :=  \p D_r (x_a) \cap \big\{z \geq 0 \big\}$, $\kappa\h{0.5pt}[\h{.5pt} u_{a, +} \h{.5pt}]$ converges to $e_z$. If we start from $x_{a, r}$ and rotate counter--clockwisely along the circle $\p D_r (x_a)$ back to $x_{a, r}$, $\kappa\h{0.5pt}[\h{.5pt} u_{a, +} \h{.5pt}]$ continuously varies from $- e_z$ to $e_z$. In addition, except at $x_{a, r}$, the image of $\kappa\h{0.5pt}[\h{0.5pt}u_{a, +}\h{0.5pt}]$  keeps strictly on the right--half part of the $(\rho,z)$--plane.  The angle of the director field $\kappa\h{0.5pt}[\h{.5pt} u_{a, +} \h{.5pt}]$ is totally changed by $\pi$ during this process. Note that the positive direction of the horizontal (vertical resp.) axis in the $(\rho, z)$--plane is given by $e_\rho$ ($e_z$ resp.);
\item[$\mathrm{(5)}$.] Let $\varphi'$ be an angular variable ranging from $[-\pi,\pi]$. Fixing an arbitrary $r \in (0, \epsilon)$, we define the value of $u_{a, +}\big( x_a + r(\cos\varphi',\sin \varphi') \big)$ at $\varphi' = - \pi$ ($\pi$ resp.) to be $- e_z$ ($e_z$ resp.). Then  the pointwise limit of $u_{a, +}\big( x_a + r(\cos\varphi',\sin \varphi') \big)$ as $r \to 0^+$ is given by
\begin{align*}
\mbox{\fontsize{8}{8}\selectfont\(
\left\{ \h{5pt}
\begin{aligned}
&-e_z&
\h{10pt} &\text{if $\varphi' = - \pi$;}\\
&\dfrac{\sqrt{2}}{2} \left( 1 - \dfrac{\varkappa_a \h{2pt}\textup{ctan} \varphi'}{\sqrt{ 4 + \varkappa_a^2 \h{2pt}\textup{ctan}^2 \varphi'}} \right)^{1/2} e_\rho
-  \sqrt{ \dfrac{2}{ 4 + \varkappa_a^2\h{2pt}\textup{ctan}^2 \varphi'} }
\left( 1 - \dfrac{\varkappa_a\h{2pt}\textup{ctan} \varphi'}{\sqrt{4+\varkappa_a^2\h{2pt}\textup{ctan}^2 \varphi'}} \right)^{-1/2} e_z&
\h{10pt} &\text{if $\varphi' \in (-\pi,0)$;}\\
&e_\rho&
\h{10pt} &\text{if $\varphi' =0$;}\\
& \dfrac{\sqrt{2}}{2} \left( 1+\dfrac{\varkappa_a \h{2pt}\textup{ctan} \varphi'}{\sqrt{ 4 + \varkappa_a^2 \h{2pt}\textup{ctan}^2 \varphi'}} \right)^{1/2} e_\rho
+ \sqrt{ \dfrac{2}{ 4 + \varkappa_a^2\h{2pt}\textup{ctan}^2 \varphi'} }
\left( 1+\dfrac{\varkappa_a\h{2pt}\textup{ctan} \varphi'}{\sqrt{4+\varkappa_a^2\h{2pt}\textup{ctan}^2 \varphi'}} \right)^{-1/2} e_z&
\h{10pt} &\text{if $\varphi' \in (0,\pi)$;}\\
&e_z&
\h{10pt} &\text{if $\varphi' =\pi$.}\\
\end{aligned}
\right.\)}
\end{align*}With $u_{a, +} = \big(u_{a, +; 1}, u_{a, +; 2}, u_{a, +; 3}\big)^\top$, the  constant  $\varkappa_a$ equals  $\dfrac{ \p_\rho u_{a, +;1}(x_a) - \sqrt{3}\h{1pt} \p_\rho u_{a, +;2}(x_a)}{ \p_z u_{a, +;3}(x_a)} \geq 0$.
\end{enumerate}\label{biaxial--ring solution}
\end{thm}

\noindent Near the disclination ring, the structure of the biaxial--ring solution and the distribution of its director field are illustrated in Figure 1 on the $(x_1, z)$--plane. \vspace{0.2pc}\begin{figure}[h]
\centering
    \includegraphics[scale=0.32]{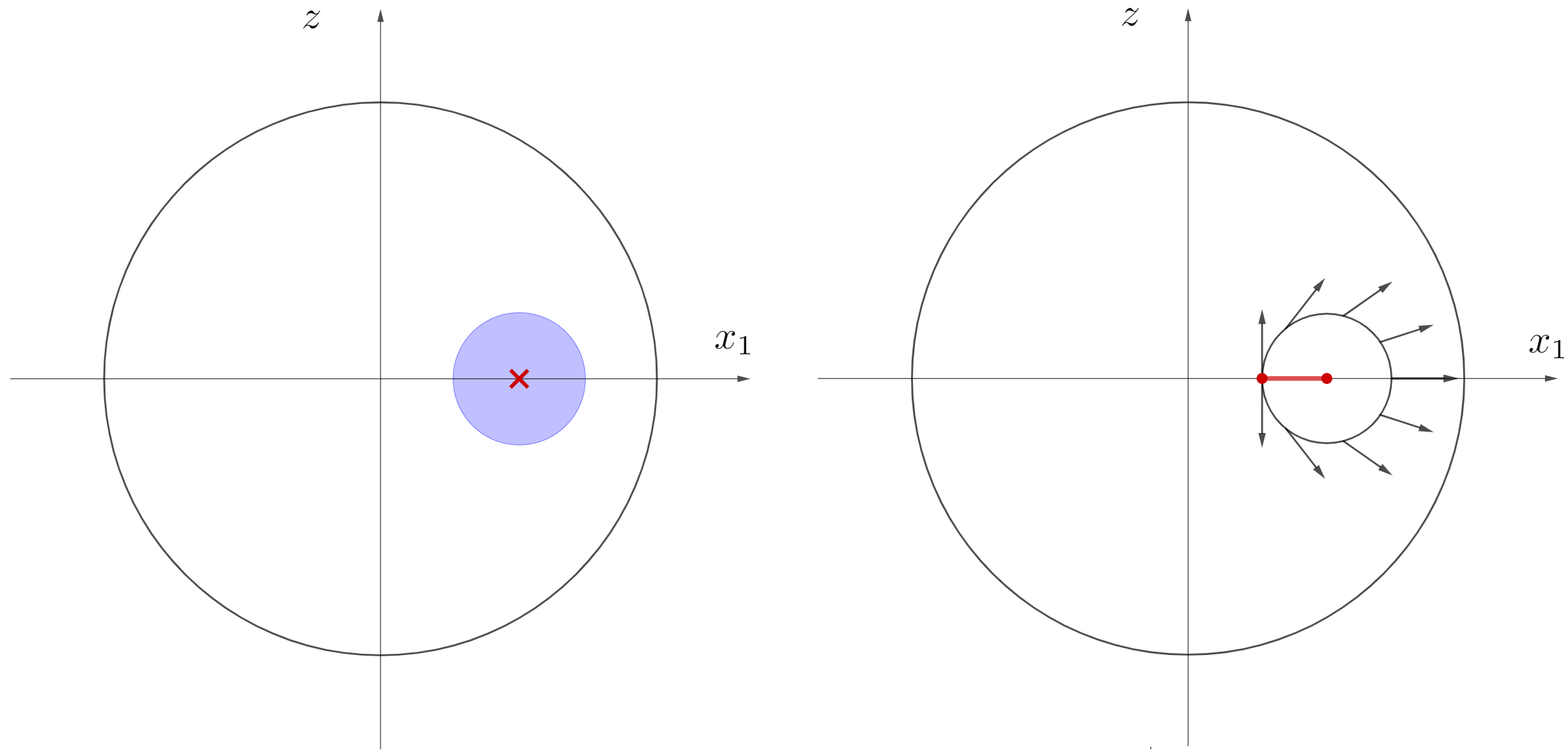}
\begin{tabular}{p{15cm}}
\begin{small}Figure 1. Biaxial--ring solution. The graph on the left indicates that the biaxial--ring solution is negative uniaxial at the point with red cross, while it is biaxial on the punctured disk (blue region). The graph on the right shows the distribution of the director field near the disclination ring. \end{small}
\end{tabular}
 \end{figure}

To discuss split--core solutions and their structures, it is better to define a notation for dumbbell first. \vspace{0.4pc}
\begin{defn}
Let $z^+_a = (0,0,z_a)$ be a point on the positive part of the $z$--axis. Its  symmetric point  with respect to the $(x_1, x_2)$--plane is denoted by $z_a^-$. Let $\epsilon > 0$ and $r>0$ be two constants with $\epsilon < r/2$. In the $(x_1,z)$--plane, the horizontal line $z = z_a - r + \epsilon$ has two intersections with the circle $\p D_{r}(z^+_a)$. Here we also use $D_\rho (x)$ to denote an open disk in the $(x_1,z)$--plane with centre $x$ and radius $\rho$. The intersection point with positive
$x_1$--coordinate is denoted by $x^+_1$, while the intersection point with negative $x_1$--coordinate is denoted by $x^+_2$. Similarly the horizontal line $z = - z_a + r- \epsilon$ also has two intersections with $\p D_{r}(z_a^-)$. Within these two intersections, the one with positive $x_1$--coordinate is denoted by $x_1^-$, while another intersection is denoted by $x_2^-$. \vspace{0.2pc}

The contour $\mathscr{C}_{r, \epsilon}\big(z_a^+, z_a^-\big)$ in the $(x_1,z)$--plane is then defined as follows: firstly we start from $x^+_1$ and rotate counter--clockwisely along $\p D_{r}(z_a^+)$ to $x^+_2$. Then we connect $x^+_2$ and $x^-_2$ by the straight segment between them. From $x^-_2$, we rotate counter--clockwisely along $\p
D_ {r} (z^-_a)$ to arrive at $x^-_1$. Finally we connect $x^-_1$ and $x^+_1$ by the straight segment between them. The dumbbell, denoted by $D_{r,\h{1pt}\epsilon}\big(z^+_a, z_a^-\big)$, refers to the region in the $(x_1,z)$--plane  enclosed by the contour $\mathscr{C}_{r, \epsilon}\big(z_a^+, z_a^-\big)$.
\label{dumbbell}
\end{defn}

Now we discuss our main results on the split--core solutions.
\begin{thm}[\bf Split--core solutions and their split--core disclinations] There exists a  constant $a_0 > 0$ so that for all $a > a_0$, the followings hold: \begin{enumerate} \item[$\mathrm{(1)}$.] There exists a $\mathscr{R}$--axially symmetric solution, denoted by $w_{a, -} = \mathscr{L}[\h{1pt}u_{a, -}\h{1pt}]$, to the boundary value problem (\ref{el-eq of u_a, intro}) and (\ref{boundary condition of u_a, intro}).   All zeros of $w_{a, -}$  must be   on the $z$--axis and different from $0$. There must have odd number  of zeros of $w_{a, -}$ on the set $\big\{ \big(0, 0, z\big) : 0 < z < 1\big\}$;
\item[$\mathrm{(2)}$.] Let $\lambda^-_{a; j}$ ($j = 1, 2, 3$) be the three eigenvalues in (\ref{evalue of D with u= u^+_a,b}) computed in terms of the $3$--vector field $u_{a, -}$. Recall the tensor field $\mathscr{Q}$ in (\ref{ansatz}) and denote by $ \mathscr{Q}^-_a$ the tensor field $a^{-1}\mathscr{Q}\big(Rx\big)$ with $v(y) = u_{a, -}\big(R^{-1}y\big)$. Then $\lambda^-_{a; j}$ ($j = 1, 2, 3$) are the three eigenvalues of $\mathscr{Q}^-_a$. Let $z_a^+$ be the lowest zero of $w_{a, -}$ on the positive part of the $z$--axis. $z_a^-$ is its symmetric point with respect to the $(x_1, x_2)$--plane. Then there exist $\epsilon > 0$ and $\epsilon_1 > 0$ with $\epsilon_1 < \epsilon / 2$ so that \begin{enumerate} \item[$\mathrm{(2.1)}$.] The tensor field $\mathscr{Q}_a^-$ is biaxial on $D_{\epsilon, \epsilon_1}\big(z_a^+, z_a^-\big) \setminus l_z$, where $l_z$ is the $z$--axis.  $D_{\epsilon, \epsilon_1}\big(z_a^+, z_a^-\big)$ is the dumbbell introduced in Definition \ref{dumbbell}. More precisely, there holds $\lambda_{a; 2}^- < \lambda_{a; 1}^- < \lambda_{a; 3}^-$ on $D_{\epsilon, \epsilon_1}\big(z_a^+, z_a^-\big) \setminus l_z$;
\item[$\mathrm{(2.2)}$.] The tensor field $\mathscr{Q}_a^-$ is isotropic at $z_a^+$ and $z_a^-$;
\item[$\mathrm{(2.3)}$.] Given two points $Z$ and $W$, we use $(Z, W)$ to denote the open segment connecting $Z$ and $W$. Then on $\big(z_a^+, z_a^+ + \epsilon \h{1pt}e_z \big)\h{1pt}\cup\h{1pt}\big(z_a^-, z_a^- - \epsilon \h{1pt}e_z\big)$, the tensor field $\mathscr{Q}_a^-$ is positive uniaxial. More precisely, it holds $\lambda^-_{a; 2} = \lambda^-_{a; 1} < \lambda^-_{a; 3}$ on the set $\big(z_a^+, z_a^+ + \epsilon \h{1pt}e_z \big)\h{1pt}\cup\h{1pt}\big(z_a^-, z_a^- - \epsilon \h{1pt}e_z\big)$;
\item[$\mathrm{(2.4)}$.] The tensor field $\mathscr{Q}_a^-$ is negative uniaxial on $(z_a^+, z_a^-)$.  More precisely, it holds $\lambda^-_{a; 2} < \lambda^-_{a; 1} = \lambda^-_{a; 3}$ on $(z_a^+, z_a^-)$. 
\end{enumerate}
The constants $\epsilon$ and $\epsilon_1$ are suitably small and independent of $a$;
\item[$\mathrm{(3)}$.]  The director field of   $\mathscr{Q}_a^-$ on   $D_{\epsilon, \epsilon_1}\big(z_a^+, z_a^-\big) \setminus l_z$ can be oriented and represented by  $\kappa\h{0.5pt}[\h{0.5pt}u_{a, -}\h{0.5pt}]$ (see (\ref{hat kappa_3, intro})). It is the normalized eigenvector of $\mathscr{Q}_a^-$ associated with the eigenvalue $\lambda^-_{a; 3}$;
\item[$\mathrm{(4)}$.] The closed segment connecting $z_a^+$ and $z_a^-$, denoted by $\big[z_a^+, z_a^-\big]$, is the split--core--segment disclination of $\mathscr{Q}_a^-$. Its structure is described as follows: \begin{enumerate} 
\item[$\mathrm{(4.1)}$.] Let $D_{r, r_1}\big(z_a^+, z_a^-\big)$ be an arbitrary dumbbell contained in $D_{\epsilon, \epsilon_1}\big(z_a^+, z_a^-\big)$ and denote by $\mathscr{C}_{r, r_1}\big(z_a^+, z_a^-\big)$ the boundary contour of $D_{r, r_1}\big(z_a^+, z_a^-\big)$. In the $(x_1, z)$--plane, when we move along $\mathscr{C}_{r, r_1}\big(z_a^+, z_a^-\big)$ clockwisely from $z_a^+ + r\h{1pt}e_z$ to the point $\left(\sqrt{r^2 - (r - r_1)^2}, 0, 0\right)$, then the dirctor field $\kappa\h{0.5pt}[\h{0.5pt}u_{a, -}\h{0.5pt}]$ varies continously from $e_z$ to $e_\rho$. When we continue to move from  the point $\left(\sqrt{r^2 - (r - r_1)^2}, 0, 0\right)$ to the point $z_a^- - r\h{1pt}e_z$, then $\kappa\h{0.5pt}[\h{0.5pt}u_{a, -}\h{0.5pt}]$ varies continously from  $e_\rho$ to $- e_z$. The image of $\kappa\h{0.5pt}[\h{0.5pt}u_{a, -}\h{0.5pt}]$ restricted on $\mathscr{C}_{r, r_1}\big(z_a^+, z_a^-\big) \setminus \big\{ z_a^+ + r\h{1pt}e_z, z_a^- - r\h{1pt}e_z\big\}$ keeps strictly on the right--half part of the $(\rho,z)$--plane;
\item[$\mathrm{(4.2)}$.] On $\big(z_a^+, z_a^-\big)$, the director field $\kappa\h{0.5pt}[\h{0.5pt}u_{a, -}\h{0.5pt}]$ equivalently equals  $e_\rho$. Note that this result yields that all points on $\big(z_a^+, z_a^-\big)$ are disclinations of $\mathscr{Q}_a^-$ with strength $1$ since in $\mathbb{R}^3$, $e_\rho$ depends on the azimuthal angle and is the radial direction in the $(x_1, x_2)$--plane; 
\item[$\mathrm{(4.3)}$.] At $z_a^+$, the director field $\kappa\h{0.5pt}[\h{0.5pt}u_{a, -}\h{0.5pt}]$ satisfies \begin{align*}
\lim_{(a^{-1}, \h{0.5pt}r) \h{1pt}\to\h{1pt} (0, 0)} \h{1.5pt}
\Big\lVert \h{2pt}\kappa\h{0.5pt}[\h{0.5pt}u_{a, -}\h{0.5pt}]
- \kappa^+ \big(\cdot - z^+_a\big)
\h{1.5pt} \Big\rVert_{\infty;\h{1pt} \p B_{r}  (z^+_a ) } = 0.
\end{align*}
Here $\| \cdot \|_{\infty; S}$ denotes the $L^\infty$--norm on some set $S$. $\kappa^+ (x)$  is given by
\begin{align*}\mbox{\fontsize{9}{9}\selectfont\(
\kappa^+ (x)
:=\left\{ \h{5pt}
\begin{aligned}
&e_z&
\h{5pt} &\text{if $\phi = 0$;}\\
&\dfrac{\sqrt{2}}{2}
\left( 1 - \dfrac{\sqrt{3} \h{2pt}\cos\phi}{\displaystyle\sqrt{3+\sin^2 \phi}} \right)^{1/2} e_\rho
+\sqrt{ \dfrac{2 \sin^2 \phi }{ 3 + \sin^2 \phi} }
\left( 1 - \dfrac{\displaystyle\sqrt{3} \h{2pt}\cos\phi}{\displaystyle\sqrt{3+\sin^2 \phi}} \right)^{-1/2} e_z&
\h{5pt} &\text{if $\phi \in (\h{0.5pt}0,\pi\h{0.5pt}]$.}
\end{aligned}
\right.\)}
\end{align*}The asymptotic behavior of $\kappa\h{0.5pt}[\h{0.5pt}u_{a, -}\h{0.5pt}]$ near $z_a^-$ is given as follows: \begin{align*}
\lim_{(a^{-1}, \h{0.5pt}r) \h{1pt}\to\h{1pt} (0, 0)} \h{1.5pt}
\Big\lVert \h{1.5pt} \kappa\h{0.5pt}[\h{0.5pt}u_{a, -}\h{0.5pt}]
- \kappa^-\big(\cdot - z^-_a\big)
\h{1.5pt} \Big\rVert_{\infty;\h{1pt} \p B_{r}  (z^-_a ) } = 0.
\end{align*}
Here $\kappa^- (x)$ is defined by
\begin{align*}\mbox{\fontsize{9}{9}\selectfont\(
\kappa^- (x)
:=\left\{ \h{5pt}
\begin{aligned}
&\dfrac{\sqrt{2}}{2}
\left( 1 + \dfrac{\sqrt{3} \h{2pt}\cos\phi}{\displaystyle\sqrt{3+\sin^2 \phi}} \right)^{1/2} e_\rho
-\sqrt{ \dfrac{2 \sin^2 \phi }{ 3 + \sin^2 \phi} }
\left( 1 + \dfrac{\sqrt{3} \h{2pt}\cos\phi}{\displaystyle\sqrt{3+\sin^2 \phi}} \right)^{-1/2} e_z&
\h{5pt} &\text{if $\phi \in [\h{0.5pt}0,\pi\h{0.5pt})$;}\\
&-e_z&
\h{5pt} &\text{if $\phi = \pi$.}
\end{aligned}
\right.\)}
\end{align*}
\end{enumerate}
\end{enumerate}\label{split--core solution}
\end{thm}
\noindent Near the split core, the structure of the split--core solution and the distribution of its director field are illustrated in Figure 2 on the $(x_1, z)$--plane. \vspace{0.2pc}\begin{figure}[h]
\centering
    \includegraphics[scale=0.12]{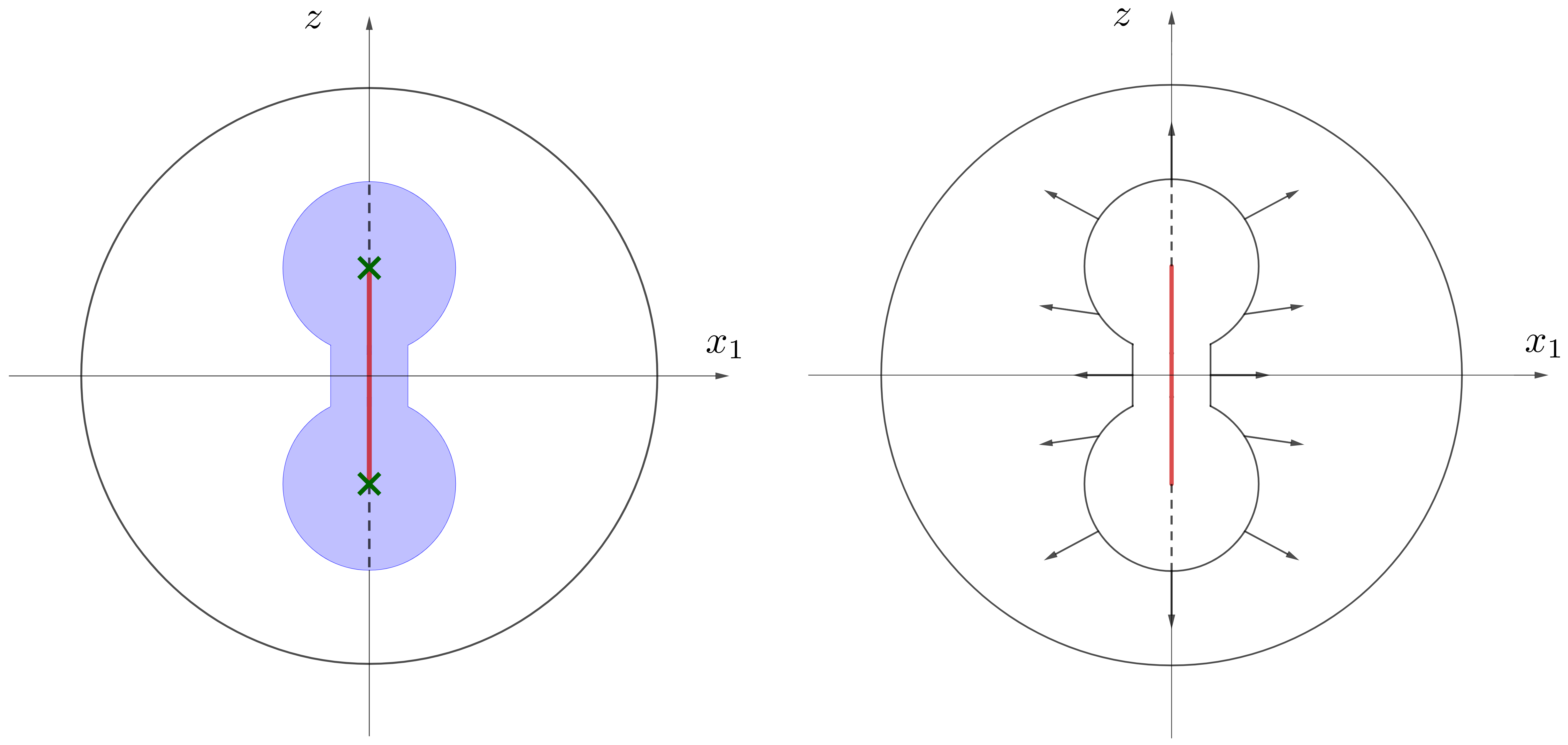}
\begin{tabular}{p{15cm}}
\begin{small}Figure 2. Split--core solution. The graph on the left indicates that the split--core solution is negative uniaxial on the red bold segment (end--points not included). It is isotropic at the two points with green cross. On the two dashed black segment (isotropic points not included), the solution is positive uniaxial. It is biaxial at the  points of the dumbbell off the $z$--axis (blue region). The graph on the right shows the distribution of the director field near the split core. \end{small}
\end{tabular}
 \end{figure} 

\subsection{\normalsize Main ideas and methodology}\label{sec_1_3}\vspace{0.5pc} 

We briefly discuss main ideas used to prove Theorems \ref{biaxial--ring solution} and \ref{split--core solution}.

\subsubsection{Vectorial Signorini problem}\vspace{0.5pc}
For any $5$--vector field  $w$ on $B_1$, we define its $\mathcal{E}_\mu$--energy by  \begin{align}\label{def of E mu} \mathcal{E}_\mu\left[w\right] := \int_{B_1} \big| \nabla w \big|^2  + \sqrt{2} \h{0.2pt}\mu \left( 1 - 3 S\left[ w\right]\right).\end{align}In the low--temperature limit $a \rightarrow \infty$ (possibly up to a subsequence), minimizers of $\mathcal{E}_{a, \mu}$ (see (\ref{def of mathcal e a mu})) within the configuration space: \begin{align}\label{def of F_a, mu}\mathscr{F}_{a} := \Big\{ w = \mathscr{L}\left[ u\right] : u = u\left(\rho, z\right) \in \mathbb{R}^3, \h{1.5pt}\mathcal{E}_{a, \mu} \left[ w \right] < \infty \h{2.5pt}\text{and $w$ satisfies (\ref{boundary condition of u_a, intro})} \Big\}\end{align}converges strongly in $H^1\left(B_1\right)$ to a minimizer of the $\mathcal{E}_\mu$--energy within the configuration space: \begin{align}\label{def of F mu} \mathscr{F} := \Big\{ w = \mathscr{L}\left[ u \right] : u = u \left(\rho, z\right) \in \mathbb{S}^2, \h{1.5pt} \mathcal{E}_\mu\left[w\right] < \infty \h{5pt}\text{and $w = \mathscr{L}\left[U^*\right]$ on $\p B_1$} \Big\}.
\end{align}Here $\mathbb{S}^2$ is the unit sphere in $\mathbb{R}^3$ with the center located at $0$. In other words, the $\mathcal{E}_{a, \mu}$--energy functional defined on the $\mathscr{F}_{a}$--space $\Gamma$--converges to the $\mathcal{E}_\mu$--energy functional defined on the $\mathscr{F}$--space. In \cite{Y20}, the second author studies a class of critical points of $\mathcal{E}_\mu$ in $\mathscr{F}$ with the $\mathscr{R}$--axial symmetry.
\noindent The following results on the  $\mathscr{R}$--axially symmetric critical points of $\mathcal{E}_\mu$ in $\mathscr{F}$ are shown in \cite{Y20}: \begin{thm}\label{results for limiting case}
Denote by $T$ the flat boundary of $B^+_1$ where $B_1^+$ is the upper--half part of $B_1$, and $\mathscr{F}^s$ the configuration space consisting of all $\mathscr{R}$--axially symmetric vector fields in $\mathscr{F}$. For any $b \in \mathrm{I}_-$ and $c \in \mathrm{I}_+$ where  $\mathrm{I}_- := \left(-1, - 1/2 \h{1pt}\right]$ and $\mathrm{I}_+ := \left[ - 1/2, 1\h{1pt}\right)$, we denote by $w_b^+$ and $w_c^-$ the minimizers of the following Signorini--type problems, respectively: \begin{align}\label{obstacle limit} \mathrm{Min} \Big\{ \mathcal{E}_\mu\left[ w\right] : w \in \mathscr{F}_b^+\Big\} \h{15pt}\text{and}\h{15pt}\mathrm{Min} \Big\{ \mathcal{E}_\mu\left[ w\right] : w \in \mathscr{F}_c^- \Big\}.
\end{align}With $w_j$ denoting the $j$--th component of a vector field $w$, $\mathscr{F}_b^+$ and $\mathscr{F}_c^-$ are configuration spaces given by  \begin{align*} \mathscr{F}_b^+ := \Big\{ w \in \mathscr{F}^s : w_3 \geq b \h{4pt}\text{on \h{1pt}$T$}\Big\} \h{15pt}\text{and}\h{15pt}\mathscr{F}_c^- := \Big\{ w \in \mathscr{F}^s : w_3 \leq c \h{4pt}\text{on \h{1pt}$T$}\Big\}.
\end{align*}Then we have \begin{enumerate}
\item[$\mathrm{(1)}$.] For any $b \in \mathrm{I}_-$ and $c \in \mathrm{I}_+$, there exist $u_b^+ : \mathbb{D} \longrightarrow \mathbb{S}^2$ and $u_c^- : \mathbb{D} \longrightarrow \mathbb{S}^2$, where \begin{align*}   \mathbb{D} := \Big\{ \left(\rho, z\right) : \rho > 0 \h{5pt}\text{and}\h{5pt}\rho^2 + z^2 < 1\Big\},
\end{align*}so that $w_b^+ = \mathscr{L}\big[ u_b^+\big]$ and $w_c^- = \mathscr{L}\big[ u_c^-\big]$. Moreover, $u_b^+$ and $u_c^-$ satisfy \begin{align*} - \dfrac{1}{\rho} D \cdot \left( \rho Du \right)  + \dfrac{1}{\rho^2}\begin{pmatrix}4u_1\\[2mm]0\\[2mm]u_3\end{pmatrix} - \dfrac{3\mu}{\sqrt{2}}\h{1pt}\nabla_u P\left[ u \right] = \left\{ | Du |^2 + \dfrac{1}{\rho^2}\left( 4 u_1^2 + u_3^2 \right) - \dfrac{9 \mu}{\sqrt{2}}\h{1pt} P\left[u\right]\right\} u \qquad \text{in } \mathbb{D}^+.
\end{align*}In the above, $D = \left(\p_\rho, \p_z\right)$ is the gradient operator on $(\rho,z)$--plane. $P\left[u\right]$ is defined by \begin{align}\label{def of P}P\left[u\right] :=-u_1^2u_2+\dfrac{\sqrt{3}}{2}u_1u_3^2+\dfrac{1}{3}u_2^3+\dfrac{1}{2}u_2u_3^2, \h{15pt}\text{for any $u \in \mathbb{R}^3$.}\end{align} $\mathbb{D}^+$ is the subset $\big\{ \left(\rho, z\right) \in \mathbb{D} : z > 0\big\}$.  In addition, it satisfies \begin{align*}  
u_b^+ = u_c^- = U^* \h{20pt}\text{on $\Big\{ \left(\rho, z\right) : \rho \geq 0 \h{3pt}\text{and}\h{3pt} \rho^2 + z^2 = 1\Big\}$.}
\end{align*}In light of the above boundary conditions and the equations satisfied by $u_{b; 1}^+$, $u_{b; 3}^+$, $u_{c; 1}^-$ and $u_{c; 3}^-$ in $\mathbb{D}^+$, where $u_{b; j}^+$ and $u_{c; j}^-$ are $j$--th components of $u_b^+$ and $u_c^-$ respectively, we can apply strong maximum principle to obtain \begin{align*}u_{b; 1}^+ > 0, \h{15pt}u_{b; 3}^+ > 0, \h{15pt}u_{c; 1}^- > 0, \h{15pt}u_{c; 3}^- > 0 \h{20pt}\text{in $\mathbb{D}^+$;}
\end{align*}
\item[$\mathrm{(2)}$.] For any $b \in \mathrm{I}_-$ and $c \in \mathrm{I}_+\setminus \big\{0\big\}$, $w_{b}^+$ and $w_{c}^-$ are weak solutions to the following Dirichlet boundary value problem: \begin{align} \label{eqn and bdy cond of w}\left\{ \begin{array}{lcl} - \Delta w - \dfrac{3 \mu}{\sqrt{2}}\h{1pt}\nabla_w S \left[ w \right] = \left\{ \big| \nabla w \big|^2 - \dfrac{9\mu}{\sqrt{2}}\h{1pt}S \left[w\right]\right\} w \h{20pt}&\text{$\mathrm{in}$ $B_1$;}\\[3mm]
\h{81pt}w = \mathscr{L}\left[\h{0.5pt} U^*\h{0.5pt}\right] &\text{\h{7pt}$\mathrm{on}$ $\p B_1$.}\end{array}\right.
\end{align}Moreover, $w_b^+$ and $w_c^-$ are smooth in $B_1$ up to the boundary $\p B_1$, except possibly at finitely many singularities. All the singularities of $w_{b}^+$ and $w_{c}^-$ must be on $l_z$, but different from $0$ and two poles;
\item[$\mathrm{(3)}$.] There exist $b_\star \in \mathrm{I}_-$, $c_\star \in \left(0, 1\right)$, $w_{b_\star}^+ \in \mathscr{F}_{b_\star}^+$ and $w_{c_\star}^- \in \mathscr{F}_{c_\star}^-$ such that \begin{align*} \mathcal{E}_\mu\left[ w_{b_\star}^+\right] = \mathcal{E}_\mu\left[ w_{c_\star}^-\right]  = \mathrm{Min}\Big\{ \mathcal{E}_\mu \left[w\right] : w \in \mathscr{F}^s \Big\};
\end{align*}
\item[$\mathrm{(4)}.$] There exist a $b_0 \in ( -1, - 1/2\h{0.3pt})$ and $c_0 \in (0, 1)$ so that  \begin{align*}  \inf  \h{2pt} \Big\{w_{b; \h{0.5pt}3}^+(x) :  x\h{0.5pt}\in\h{0.5pt}T,\h{4pt}b \in \left(-1, \h{0.5pt}b_\star\right)\h{4pt}\text{and}\h{4pt}w_b^+ \h{2pt}\text{is a minimizer of $\mathcal{E}_\mu$ in $\mathscr{F}_b^+$} \Big\}\h{2pt} \geq \h{2pt} b_0\end{align*}and\begin{align*}\sup  \h{2pt} \Big\{w_{c; \h{1pt}3}^-(x) :  x\h{0.5pt}\in\h{0.5pt}T, \h{4pt}c \in \left(c_\star, \h{0.5pt}1\right)\h{4pt}\text{and}\h{4pt}w_c^-\h{2pt}\text{is a minimizer of $\mathcal{E}_\mu$ in $\mathscr{F}_c^-$} \Big\}\h{2pt} \leq \h{2pt} c_0.
\end{align*}Here and in what follows, $w_{b; \h{0.5pt}j}^+$ and $w_{c;\h{0.5pt}j}^-$ are the $j$--components of $w_b^+$ and $w_c^-$, respectively.\end{enumerate}
\end{thm}

For our problems in this article where the reduced temperature $a$ is finite, we introduce two vectorial Signorini--type problems for $\mathcal{E}_{a, \mu}$--energy, similarly as the two minimization problems in (\ref{obstacle limit}) for the limiting energy $\mathcal{E}_\mu$. Firstly we define $\mathscr{F}_a^s$ to be the configuration space which consists of all $\mathscr{R}$--axially symmetric vector fields in $\mathscr{F}_a$. See the definition of $\mathscr{F}_a$ in (\ref{def of F_a, mu}). For any $b \in \mathrm{I}_-$ and $c \in \mathrm{I}_+$, we let $$\mathscr{F}_{a,b}^+ :=\bigg\{w \in \mathscr{F}^s_a :\h{1pt} w_3 \geq  b H_a \h{3pt}\text{on $T$}   \bigg\} \h{15pt}\text{ and }\h{15pt} \mathscr{F}_{a,c}^- :=\bigg\{w \in \mathscr{F}^s_a :\h{1pt} w_3  \leq  c H_a \h{3pt}\text{on $T$}   \bigg\}.$$ Associated with $\mathscr{F}_{a, b}^+$ and $\mathscr{F}_{a, c}^-$, we consider the  two Signorini--type minimization problems given as follows: \begin{align}\label{obstacle finite}\mathrm{Min}\Big\{\mathcal{E}_{a, \mu}[\h{.5pt}w\h{.5pt}] : w \in \mathscr{F}_{a, b}^+\Big\} \h{20pt}\text{and}\h{20pt} \mathrm{Min}\Big\{\mathcal{E}_{a, \mu}[\h{.5pt}w\h{.5pt}] : w \in \mathscr{F}_{a, c}^-\Big\}.\end{align}By the direct method of calculus of variations, it holds \begin{prop}For any $b \in \mathrm{I}_-$ and $c \in \mathrm{I}_+$, there exist $w^+_{a,b} \in \mathscr{F}^+_{a,b}$ and $w^-_{a,c} \in \mathscr{F}^-_{a,c}$ such that  \begin{align*}\mathcal{E}_{a, \mu} \big[ w_{a, b}^+\big] =\mathrm{Min}\Big\{\mathcal{E}_{a, \mu}[\h{.5pt}w\h{.5pt}] : w \in \mathscr{F}_{a, b}^+\Big\} \h{20pt}\text{and}\h{20pt} \mathcal{E}_{a, \mu} \big[ w_{a, c}^-\big] = \mathrm{Min}\Big\{\mathcal{E}_{a, \mu}[\h{.5pt}w\h{.5pt}] : w \in \mathscr{F}_{a, c}^-\Big\}.\end{align*} Similarly as in item (1) of Theorem \ref{results for limiting case},  there exist $u_{a, b}^+ : \mathbb{D} \longrightarrow \mathbb{R}^3$ and $u_{a, c}^- : \mathbb{D} \longrightarrow \mathbb{R}^3$  so that $$w_{a, b}^+ = \mathscr{L}\big[ u_{a, b}^+\big] \h{15pt}\text{and}\h{15pt}w_{a, c}^- = \mathscr{L}\big[ u_{a, c}^-\big].$$ In addition, $u_{a, b}^+$ and $u_{a, c}^-$ satisfy \begin{align}\label{eqn of ua, bin D plus}
- \dfrac{1}{\rho} D \cdot \left( \rho Du \right)   + \dfrac{1}{\rho^2}\begin{pmatrix}4u_1\\[2mm]0\\[2mm]u_3\end{pmatrix}= \mu\left\{ \dfrac{3 }{\sqrt{2}}\h{1pt}\nabla_u P\left[ u\right]-a\h{0.5pt} \Big(\h{0.2pt}|\h{0.5pt}u\h{0.5pt}|^2-1 \h{0.2pt}\Big)\h{0.5pt}u\right\} \qquad \text{in } \mathbb{D}^+.
\end{align}By the boundary condition in  (\ref{boundary condition of u_a, intro}), there also hold\begin{align}\label{bdy cond of u_a, b plus}u_{a, b}^+ = u_{a, c}^- = U_a^*\h{20pt}\text{on $\Big\{\left(\rho, z\right) : \rho \geq 0\h{3pt}\text{and}\h{3pt}\rho^2 + z^2 = 1\Big\}$}.\end{align} \label{existence of mini finite}
\end{prop}\begin{rmk}\label{rmk on sign and gamma conv} We would like to point out: \begin{enumerate} \item[$\mathrm{(1).}$] Suppose that $u = (u_1, u_2, u_3)^\top$ denotes either $u_{a, b}^+$ or $u_{a, c}^-$ in Proposition \ref{existence of mini finite} and let $\mathscr{L}\left[ u^\star\right]$ be the $\mathscr{R}$--axially symmetric vector field with $u^\star  = (\h{1pt}| u_1 |, u_2, | u_3 |\h{1pt})^\top$ on $\mathbb{D}^+$. If $\mathscr{L}\left[u\right]$ minimizes the $\mathcal{E}_{a, \mu}$--energy in either $\mathscr{F}_{a, b}^+$ or $\mathscr{F}_{a, c}^-$, then $\mathscr{L}\left[ u^\star\right]$ also minimizes the $\mathcal{E}_{a, \mu}$--energy in the same configuration space as $\mathscr{L}\left[u\right]$. In light of the first and third equations in (\ref{eqn of ua, bin D plus}) and the boundary condition in (\ref{bdy cond of u_a, b plus}), we can apply strong maximum principle to obtain the strict positivity of $u^\star_1$ and $u^\star_3$ in $\mathbb{D}^+$, which in turn induces the strict positivity of $u_1$ and $u_3$ in $\mathbb{D}^+$;
\item[$\mathrm{(2).}$] Denote by $w$ either $w_{a, b}^+$ or $w_{a, c}^-$ in Proposition \ref{existence of mini finite}. In addition, we let \begin{align*} \widetilde{w} := \left\{ \begin{array}{lcl} 
w, \h{20pt}&\text{if $\big| w \big| \leq H_a$;}\\[2mm]
H_a \widehat{w}, &\text{if $\big| w\big| > H_a$.}
\end{array}\right.
\end{align*}Here $\widehat{w}$ is the normalized vector field of $w$. For any $b \in \mathrm{I}_-$ and $c \in (0, 1)$, the vector field $\widetilde{w}$ lies in the same configuration space as $w$. As for the $\mathcal{E}_{a, \mu}$--energy of $\widetilde{w}$, firstly we have \begin{align*}\int_{\left|\h{0.5pt} w \h{0.5pt}\right| \h{1pt}>\h{1pt} H_a} \big|\h{0.2pt} \nabla w \h{0.2pt}\big|^2 = \int_{\left|\h{0.5pt} w \h{0.5pt}\right| \h{1pt}>\h{1pt} H_a} \big|\h{0.3pt} \nabla \big|\h{0.2pt} w \h{0.2pt}\big| \h{0.3pt}\big|^2 + \big|\h{0.2pt} w \h{0.2pt}\big|^2 \big| \nabla \widehat{w} \h{0.3pt}\big|^2 \geq  \int_{\left|\h{0.5pt} w \h{0.5pt}\right| \h{1pt}>\h{1pt} H_a} \big|\h{0.2pt} \nabla \widetilde{w} \h{0.2pt}\big|^2.
\end{align*}On the other hand, it holds \begin{align*} 2D_a - 6 \sqrt{2} S\left[w\right] + a \left( \big| w \big|^2 - 1 \right)^2 &= 2 D_a - 6 \sqrt{2}\h{2pt} \big| w \big|^3 S\big[\h{1pt}\widehat{w}\h{1pt}\big] + a \left( \big| w \big|^2 - 1 \right)^2\\[2mm]
&= 2 \sqrt{2}\h{2pt} \big| w\big|^3 \left( 1 - 3 S\big[\h{1pt}\widehat{w}\h{1pt}\big] \right) + 2D_a - 2\sqrt{2}\h{1pt}\big| w \big|^3 + a \left( \big| w \big|^2 - 1 \right)^2.
\end{align*}Note that $S \left[ w \right] \leq 1\big/3$ for any $w = \mathscr{L}\left[x \right]$ with the unit length. The polynomial $$2D_a - 2\sqrt{2}\h{1pt}h^3 + a \left( h^2 - 1 \right)^2$$ achieves its global minimum value $0$ at $ h = H_a$. It then turns out from the above arguments that \begin{align*}2 D_a - 6 \sqrt{2} \h{1pt}S\left[w\right] + a \left( \big| w \big|^2 - 1 \right)^2 >  2D_a - 6 \sqrt{2} \h{1pt}S\big[\h{1pt}\widetilde{w}\h{1pt}\big] + a \left( \big| \widetilde{w} \big|^2 - 1 \right)^2 \h{10pt}\text{if $\big| w \big| > H_a$.}
\end{align*}If the Lebesgue measure of $\Big\{ \big| w\big| > H_a\Big\}$ is strictly positive, then $\widetilde{w}$ has strictly smaller $\mathcal{E}_{a, \mu}$--energy than $w$ in the corresponding configuration space. This is a contradiction to the fact that $w$ saturates minimum $\mathcal{E}_{a, \mu}$--energy in its associated configuration space. Hence we have \begin{align*} \big| w_{a, b}^+ \big| \leq H_a \h{15pt}\text{and}\h{15pt}\big| w_{a, c}^- \big| \leq H_a \h{20pt}\text{a.\h{1pt}e. on $B_1$, for all $b \in \mathrm{I}_-$ and $c \in (0, 1)$;}\end{align*}Then the smoothness of $w_{a, b}^+$ and $w_{a,c}^-$ on $B_1^+$ infers \begin{align*} \big| w_{a, b}^+\left(x\right) \big| \leq H_a \h{15pt}\text{and}\h{15pt}\big| w_{a, c}^-\left(x\right) \big| \leq H_a, \h{20pt}\text{for any  $x \in B^+_1$,  $b \in \mathrm{I}_-$ and $c \in (0, 1)$;}\end{align*}
\item[$\mathrm{(3).}$] Fix arbitrary $b \in \mathrm{I}_-$ and $c \in \mathrm{I}_+$.
For any sequence $\big\{ a_n \big\}$ with $a_n \rightarrow \infty$ as $n \rightarrow \infty$, there are  two vector fields $w_b^+ \in \mathscr{F}_{b}^+$ and $w_c^- \in \mathscr{F}_c^-$  so that up to a subsequence which we still denote by $\big\{ a_n\big\}$, $$w_{a_n,b}^+ \longrightarrow w_b^+ \h{15pt}\text{ and }\h{15pt} w_{a_n,c}^- \longrightarrow w_c^-, \h{20pt}\text{  strongly in $H^1(B_1)$  as $n \to \infty$.}$$  In addition,  $$\int_{B_1} a_n\Big[\h{1pt}\big| w_{a_n, b}^+\big|^2-1 \h{1pt}\Big]^2 \longrightarrow 0 \h{10pt}\text{and} \h{10pt}\int_{B_1} a_n\Big[\h{1pt}\big| w_{a_n, c}^-\big|^2-1 \h{1pt}\Big]^2  \longrightarrow  0 \qquad \text{as } n \to \infty. $$ The mappings $w_b^+$ and $w_c^-$ are minimizers of the two minimization problems in (\ref{obstacle limit}), respectively. 
\label{strong H1 conv. of u_a}
\end{enumerate}\label{positivity and gamma conv}
\end{rmk}
In light of the definitions of $H_a$ and $D_a$ in (\ref{H_a}) and (\ref{def of Da}), the proof of item (3) in Remark \ref{positivity and gamma conv} is standard. We omit it here.

\subsubsection{Multiple $\mathscr{R}$--axially symmetric solutions to the spherical droplet problem}\vspace{0.5pc}
In light of item (3) in Remark \ref{rmk on sign and gamma conv}, we expect that $w_{a, b}^+$ and $w_{a, c}^-$  can give us two different solutions to (\ref{el-eq of D}), at least for large $a$. Here $w_{a, b}^+$ and $w_{a, c}^-$ are minimizers of the two problems in (\ref{obstacle finite}), respectively.  Therefore, we need to prove the smoothness of $w_{a, b}^+$ and $w_{a, c}^-$ on $T$. As is known, the scalar Signorini  problem is a thin obstacle problem where an obstacle condition is supplied on a thin set of codimension 1. It is known that the $C^{1, \frac{1}{2}}$--regularity is the optimal regularity of solutions to the scalar Signorini problem on the thin set.  See Chapter 9 in \cite{PSU12}. To our surprise,  item (2) in Theorem \ref{results for limiting case} tells us that except possibly the case when $c = 0$, solutions to the two Signorini--type problems in (\ref{obstacle limit}) are smooth on the thin set  $T$.
The main reason is that for the scalar Signorini problem, solutions might not solve the Euler--Lagrange equation of the associated energy functional on the whole domain containing the thin set weakly. However, the minimization problems in (\ref{obstacle limit}) are different. The obstacle conditions in (\ref{obstacle limit}) are only supplied on the third components of vector fields in $\mathscr{F}_b^+$ or $\mathscr{F}_c^-$. Therefore, the solutions to (\ref{obstacle limit}) satisfy weakly all the equations in (\ref{eqn and bdy cond of w}) except possibly the equation for the third component. Notice that the minimization problems in (\ref{obstacle limit}) are for $\mathbb{S}^4$--valued mappings. Under the circumstance that the sign of the third components of the solutions can be determined, the unit length condition of the solutions to (\ref{obstacle limit}) allows us to represent the third components of the solutions in terms of their remaining components. Due to the equations satisfied by the remaining components, it is possible for us to verify that the third components of the solutions to (\ref{obstacle limit}) satisfy in the weak sense the third equation in (\ref{eqn and bdy cond of w}). Therefore, solutions to the two minimization problems in (\ref{obstacle limit}) can solve all the equations in (\ref{eqn and bdy cond of w}) weakly, at least for all $b \in \mathrm{I}_-$ and $c \in \mathrm{I}_+ \setminus \{0\}$. As a consequence,  we can apply Schoen--Uhlenbeck's arguments for harmonic maps (see \cite{SU82, SU83}) to get the smoothness of the solutions to (\ref{obstacle limit}) on $T$. Readers may refer to \cite{Y20} for details.\vspace{0.2pc}

Different from the two configuration spaces $\mathscr{F}_b^+$ and $\mathscr{F}_c^-$, there are no unit length condition for vector fields in the configuration spaces $\mathscr{F}_{a, b}^+$ and $\mathscr{F}_{a, c}^-$. It is not quite straightforward to prove, for all $b \in \mathrm{I}_-$, $c \in \mathrm{I}_+$ and $a > 0$, the smoothness of $w_{a, b}^+$ and $w_{a, c}^-$ on the whole $\overline{B_1}$, particularly on $T$. In the next, with appropriate assumptions on the parameters $a$, $b$, $c$, we discuss our methodology of studying the regularity of $w_{a, b}^+$ and $w_{a, c}^-$ on $T$. We only focus on $w_{a, b}^+$. The arguments for $w_{a, c}^-$ are similar if we assume  $c \in (0, 1)$. \vspace{0.4pc}

\noindent \textbf{\normalsize Step 1. Reduction to an interior uniform convergence on $T$.} To show the smoothness of $w_{a, b}^+$ on $T$, we need \begin{align}\label{drop signo cond} \Big\{ x \in T :  w_{a, b; 3}^+ \left(x\right) = b \Big\} = \emptyset, \h{20pt}\text{for some $a$, $b$ suitably chosen}.
\end{align}Here $ w_{a, b; j}^+ $ denotes the $j$--th component of $w_{a, b}^+$. Recalling the constants $b_0$ and $b_\star$ in Theorem \ref{results for limiting case}, we fix a parameter $b$ so that \begin{align}\label{cond of b parameter} -1 < b < \min \Big\{b_0, b_\star\Big\}.\end{align}If there exists $a_0 = a_0\left(b\right) > 0$ so that $ w_{a, b; 3}^+ \geq \left(b_0 + b\right)\big/2 $ on $T$ for any $a > a_0$, then (\ref{drop signo cond}) follows for $b$ satisfying (\ref{cond of b parameter}) and $a > a_0$. Now we assume on the contrary that there exist $a_n \rightarrow \infty$ as $n \rightarrow \infty$ and a sequence of points $\big\{x_n\big\} \subset T$ so that \begin{align}\label{contradiction to violate}  w_{a_n, b;3}^+ \left(x_n\right) < \big(b_0 + b\big)\big/2, \h{15pt}\text{for all $n \in \mathbb{N}$.}
\end{align}In light of item (3) in Remark \ref{rmk on sign and gamma conv}, we can find a $w_b^+$ solving the first problem in (\ref{obstacle limit}) so that up to a subsequence $w_{a_n, b}^+$ converges strongly in $H^1\left(B_1\right)$ to $w_b^+$. By (\ref{cond of b parameter}) and item (4) in Theorem \ref{results for limiting case}, it turns out $w_{b; 3}^+ \geq b_0 > b$ on $T$. Therefore, if $w_{a_n, b}^+$ converges to $w_b^+$ uniformly on $T$, then when $n$ is large, it turns out $ w_{a_n, b; 3}^+  > \left(b_0 + b\right)\big/2$ on $T$, which gives us a contradiction to (\ref{contradiction to violate}). By the Dirichlet boundary condition satisfied by $w_{a_n, b}^+$ on $\p B_1$, we only need the uniform convergence of $w_{a_n, b}^+$ on $T \cap \overline{B_{1 - \delta_0}}$ for some constant $\delta_0$ suitably small. In fact, for fixed $\delta_0$ sufficiently small and $a_n$ sufficiently large, any $x_n \in T$ satisfying (\ref{contradiction to violate}) must be contained in $T \cap \overline{B_{1 - \delta_0}}$. See Step 4 in the proof of Proposition \ref{prop 4.1}.   \vspace{0.4pc}

\noindent \textbf{\normalsize Step 2. Energy--decay estimates.}  $w_b^+$ is smooth on $T$. For any $\epsilon_0 > 0$, there is  $r_{\epsilon_0} \in (0, \epsilon_0)$ so that \begin{align*} r^{-1}\int_{B_{r}\left(x\right)\h{1pt}\cap\h{1pt} B_1} \big| \nabla w_b^+ \big|^2 + \sqrt{2} \mu \big[\h{1pt} 1 - 3 \h{1pt}S\big[ w_b^+\big]\h{1pt}\big] < \epsilon_0,\h{15pt}\text{for all $x \in T$ and $r \in \big(\h{0.5pt}0, r_{\epsilon_0}\h{0.5pt}\big]$.} \end{align*} In light of items (2) and (3) in Remark \ref{positivity and gamma conv}, there is $N_{r, \epsilon_0, b} \in \mathbb{N}$ depending on $r$, $\epsilon_0$  and $b$  such that \begin{align}\label{energy density small approx mapping}r^{-1} \int_{B_{r}\left(x\right) \h{1pt}\cap\h{1pt} B_1}f_{a_n, \mu}\big(w_{a_n, b}^+\big) < \epsilon_0, \h{15pt}\text{for any $x \in T$, $r \in \big(\h{0.5pt}0, r_{\epsilon_0}\h{0.5pt}\big]$ and $n \geq N_{r, \epsilon_0, b}$.}
\end{align} Here $f_{a_n, \mu}$ is the energy density function given in (\ref{def of mathcal e a mu}). With the small energy condition in (\ref{energy density small approx mapping}), we can derive some energy--decay estimates related to $w_{a_n, b}^+$. These energy--decay estimates imply the interior uniform convergence of $w_{a_n, b}^+$ on $T$. See the proof of Proposition \ref{prop 4.1}. Due to different centers of balls in our energy--decay estimates, we divide the following arguments into two cases. \vspace{0.4pc}

\noindent \textbf{Step 2.1. Energy--decay estimate on $B_r$.} If the localized energy (\ref{energy density small approx mapping}) is evaluated on $B_r \subset B_1$, then we have  \begin{prop}\label{small energy implies energy decay} Fix   $b \in \mathrm{I}_-$. There exist three positive constants $a_0$, $\epsilon_1$ and $\nu_0$ with $\nu_0 \in (0, 1/2)$, such that for any  $a > a_0$, if it satisfies \begin{align}\label{small energy cond}\mathcal{E}_{\h{0.2pt}a, \h{0.2pt}\mu; \h{1pt}0, r}\big[w_{a, b}^+\big] := r^{-1} \int_{B_r}f_{a, \mu}\big(w_{a, b}^+\big) < \epsilon_1,\end{align} then either one of the followings holds: \begin{align*} \mathrm{(1).}\h{5pt}\mathcal{E}_{\h{0.2pt}a,\h{0.2pt} \mu;  \h{1pt}0, \nu_0 \h{0.1pt}r}\big[w_{a, b}^+\big] \h{2pt}\leq\h{2pt} r^{3/2}; \h{20pt}\mathrm{(2).} \h{5pt}\mathcal{E}_{\h{0.2pt}a,\h{0.2pt} \mu; \h{1pt}0, \nu_0 \h{0.1pt}r}\big[w_{a, b}^+\big] \h{2pt}\leq\h{2pt} \dfrac{1}{2} \h{2pt} \mathcal{E}_{\h{0.2pt}a,\h{0.2pt} \mu; \h{1pt}0, r}\big[w_{a, b}^+\big].
\end{align*} The constants $a_0$, $\epsilon_1$ and $\nu_0$ depend on the parameters $b$ and $\mu$. 
\end{prop} \noindent This result will be shown in Section 2. \vspace{0.4pc}

\noindent \textbf{Step 2.2. Energy--decay estimate on balls in the family $\mathscr{J}$.} For our convenience, we denote by $\mathscr{J}$ the family of balls given as follows: \begin{align}\label{def of family J} \mathscr{J} := \Big\{ B_r\left(x\right) \subset B_1 \setminus l_z : x \in T\Big\}.
\end{align}In light of the $\mathscr{R}$--axial symmetry of $w_{a, b}^+$, our energy--decay estimate for $B_r(x) \in \mathscr{J}$ can be reduced from $3$D to $2$D. Firstly, we fix $x \in T$ and let $B_r\left(x\right) $ be a ball so that $B_r\left(x\right) \cap l_z = \emptyset$. Moreover, we define $\rho_x$ to be the radial coordinate of $x$ in the $\left(x_1, x_2\right)$--plane. Obviously, $\rho_x > r$ since $B_r\left(x\right)\cap l_z = \emptyset$. Now we revolve the disk in the $\left(x_1, z\right)$--plane with the center $\left(\rho_x, 0, 0\right)$ and radius $r/4$ around the $z$--axis. The obtained solid torus is denoted by $T_x$. We can use finitely many balls with radius $r$ to cover $T_x$. Meanwhile, the centers of these covering balls are located on $\mathfrak{C}_x := \Big\{\left(y_1, y_2, 0\right) : y_1^2 + y_2^2 = \rho_x^2\Big\}$. Note that the least number of the covering balls used to cover $T_x$, denoted by $N = N\left(r, x\right)$, can be bounded from above as follows: \begin{align}\label{upper bound on number of convering balls} N = N\left(r, x\right)\h{2pt}\leq\h{2pt}\dfrac{4\h{0.2pt}\pi \h{0.2pt} \rho_x}{r}.
\end{align}Now we let $\Big\{ \mathfrak{B}_1, ..., \mathfrak{B}_N\Big\}$ be $N$ balls covering $T_x$. The centers of $\mathfrak{B}_j$ ($j = 1, ..., N$) are located on $\mathfrak{C}_x$, while the radii of these balls equal  $r$. It then follows \begin{align*}  \int_{T_x \h{0.5pt}\cap \h{0.5pt}B_1} f_{a, \mu} \big( w_{a, b}^+\big) \h{2pt}\leq\h{2pt} \int_{  \cup_{j = 1}^N  \big(\mathfrak{B}_j \h{0.5pt}\cap\h{0.5pt}B_1\big)}f_{a, \mu} \big(w_{a, b}^+\big) \h{2pt}\leq\h{2pt}\sum_{j = 1}^N   \int_{\mathfrak{B}_j \h{0.5pt}\cap\h{0.5pt}B_1}f_{a, \mu} \big(w_{a, b}^+\big).
\end{align*}Due to the $\mathscr{R}$--axial symmetry of $w_{a, b}^+$, for any $j \in \big\{1, ..., N\big\}$, the integration of $f_{a, \mu}\big(w_{a, b}^+\big)$ over $\mathfrak{B}_j \h{0.5pt}\cap\h{0.5pt}B_1$ equal  the integration of $f_{a, \mu}\big(w_{a, b}^+\big)$ over $B_r\left(x\right)\h{0.5pt}\cap\h{0.5pt}B_1$. Combining this fact with (\ref{upper bound on number of convering balls}), we can reduce the last estimate to \begin{align}\label{decay on torus} \int_{T_x\h{0.5pt}\cap\h{0.5pt}B_1} f_{a, \mu} \big( w_{a, b}^+\big) \h{2pt}\leq\h{2pt} N   \int_{B_r\left(x\right)\h{0.5pt}\cap\h{0.5pt}B_1}f_{a, \mu} \big(w_{a, b}^+\big) \h{2pt}\lesssim\h{2pt} \dfrac{\rho_x}{r} \int_{B_r\left(x\right)\h{0.5pt}\cap\h{0.5pt}B_1}f_{a, \mu} \big(w_{a, b}^+\big).
\end{align}Define \begin{align}\label{small energy cond, off 0}\mathcal{E}_{\h{0.2pt}a, \h{0.2pt}\mu; \h{0.2pt}x, \h{0.2pt}r}\big[w_{a, b}^+\big] := r^{-1} \int_{B_r(x) \h{0.5pt}\cap\h{0.5pt}B_1}f_{a, \mu}\big(w_{a, b}^+\big).\end{align} Then the  estimate in (\ref{decay on torus}) induces \begin{align*}\rho_x^{-1}\int_{T_x \h{0.5pt}\cap\h{0.5pt}B_1} f_{a, \mu} \big( w_{a, b}^+\big) \h{2pt}\lesssim\h{2pt}\mathcal{E}_{\h{0.2pt}a, \h{0.2pt}\mu; \h{0.2pt}x, \h{0.2pt}r}\big[w_{a, b}^+\big].
\end{align*} Utilizing the notations in Theorem \ref{results for limiting case} and Proposition \ref{existence of mini finite} and noticing that $w_{a, b}^+ = \mathscr{L}\big[ u_{a, b}^+\big]$, we further rewrite the last estimate as follows: \begin{align}\label{small energy condition rho z coordinate}
\rho_x^{-1} \int_{D_{r/4}\left(\rho_x, 0\right)\h{0.5pt}\cap\h{1pt}\mathbb{D}} \rho \h{1pt}e_{a, \mu}\big[u_{a, b}^+\big] \h{2pt}\lesssim\h{2pt}\mathcal{E}_{\h{0.2pt}a, \h{0.2pt}\mu; \h{0.2pt}x, \h{0.2pt}r}\big[w_{a, b}^+\big].
\end{align}Here $D_r\left(\rho_0, z_0\right)$ denotes the disk in the $\left(\rho, z\right)$--plane with radius $r$ and center $\left(\rho_0, z_0\right)$.  $\mathbb{D}$ is given in (1) of Theorem \ref{results for limiting case}. For any vector field $u : \mathbb{D} \to \mathbb{R}^3$, the energy density $e_{a, \mu}\left[u\right]$  is read as:   \begin{align}\label{def of G_a} &e_{a, \mu}\left[u\right] := \big| D u \big|^2 + G_{a, \mu}\left(\rho, u\right), \nonumber\\[1.5mm]
&\h{30pt}\text{where}\h{5pt}G_{a, \mu}\left(\rho, u\right) := \dfrac{4 u_1^2 + u_3^2}{\rho^2}  + \mu \left[ D_a - 3 \sqrt{2} P\left[u\right] + \dfrac{a}{2} \left( | \h{0.3pt}u \h{0.3pt}|^2 - 1 \right)^2 \right].
\end{align}  Since for any $\left(\rho, z\right) \in D_{r/4}\left(\rho_x, 0\right)$, it satisfies $\rho > 3 \h{0.2pt}\rho_x/4$, we then have from (\ref{small energy condition rho z coordinate}) that \begin{align}\label{small energy condition rho z coordinate, simplified}\int_{D_{r/4}\left(\rho_x, 0\right) \h{0.5pt}\cap\h{1pt}\mathbb{D}} e_{a, \mu}\big[u_{a, b}^+\big] \h{2pt}\lesssim\h{2pt}\mathcal{E}_{\h{0.2pt}a, \h{0.2pt}\mu; \h{0.2pt}x, \h{0.2pt}r}\big[w_{a, b}^+\big].
\end{align}Note that this estimate holds for any $B_r(x)$ satisfying $x \in T$  and $B_r(x) \cap l_z = \emptyset$.\vspace{0.3pc} 

Motivated by the above arguments, we introduce a localized energy functional:
\begin{eqnarray}E_{a, \mu; \h{1pt}x, \h{0.5pt} \delta}\h{0.5pt}[\h{0.5pt}u\h{0.5pt}] :=\displaystyle\int_{D_\delta\left(\rho_x, 0\right)}  e_{a, \mu}[\h{0.5pt}u\h{0.5pt}], \h{15pt}\text{for any $B_r\left(x\right) \in \mathscr{J}$ and $\delta \in \left(0, \frac{r}{4}\h{1.5pt}\right]$.} \label{localized energy on disk}\end{eqnarray}
Then we have 
\begin{prop}
There exist four positive constants $a_0$, $\epsilon_1$, $\lambda$ and $\theta_0$, where $\lambda$ and $\theta_0 $ are less than $1/4$, such that for any $a > a_0$ and
$B_{4r}\left(x\right) \in \mathscr{J}$, if  $ E_{a, \mu;\h{1pt}x, r}\h{1pt}\big[\h{.5pt}u_{a, b}^+\h{.5pt}\big] < \epsilon_1$, then either one of the followings holds:
\begin{align*}\mathrm{(1). }\h{5pt} E_{a, \mu;\h{1pt}x, \lambda\h{0.5pt} \theta_0 \h{0.5pt}r}\h{1pt}\big[\h{.5pt}u_{a, b}^+\h{.5pt}\big] \h{2pt}\leq\h{2pt} r^{3/2}
\h{2pt};\h{25pt}
\mathrm{ (2). }\h{5pt} E_{a, \mu;\h{1pt}x,  \lambda\h{0.5pt}\theta_0\h{0.5pt} r}\h{1pt}\big[\h{.5pt}u_{a, b}^+\h{.5pt}\big] \h{2pt}\leq\h{2pt} \dfrac{1}{2}\h{1.5pt}E_{a, \mu;\h{1pt}x, \lambda\h{0.5pt}r}\h{1pt}\big[\h{.5pt}u_{a, b}^+\h{.5pt}\big].\end{align*}
\label{variant decay lemma}
\end{prop} 
\noindent This result will be shown in Section 3. \vspace{0.4pc}

\noindent \textbf{Step 3. Contradiction to (\ref{contradiction to violate}).} In Section 4, we obtain the  H\"{o}lder estimate of $w_{a_n, b}^+$ on the interior of $T$. Strictly away from $\p B_1$, the estimate is uniform in $n$. This H\"{o}lder estimate relies on  Propositions \ref{small energy implies energy decay} and \ref{variant decay lemma},  together with a trace argument and a Campanato--Morrey type estimate. Still by Proposition \ref{variant decay lemma} and  the Dirichlet boundary condition of $w_{a_n, b}^+$ on $\p B_1$, we can show that any $x_n$ satisfying (\ref{contradiction to violate}) must be in $T \cap \overline{B_{1 - \delta_0}}$, provided that $\delta_0$ is sufficiently small and $a_n$ is sufficiently large. A contradiction to (\ref{contradiction to violate}) is then obtained by  the interior H\"{o}lder estimate of $w_{a_n, b}^+$ on $T$, Arzel\`{a}--Ascoli theorem and the fact that $ w_{b; 3}^+  \geq b_0 > b$ on $T$. See the proof of Proposition \ref{prop 4.1}.

\begin{rmk} We would like to put more words on the proofs of Propositions \ref{small energy implies energy decay} and \ref{variant decay lemma}. Similarly as in \cite{Y20}, our proofs rely on some Luckhaus--type arguments. That is to study the limiting map of a blow--up sequence. However, in the current work, our temperature is finite. During the blow--up process, we should have $a_n \to \infty$ and $\overline{r}_n \to 0$ as $n \to \infty$. Here $a_n$ is a sequence of parameter $a$ and $\overline{r}_n$ is a sequence of radii of blow--up balls/disks. Different limits of $a_n \overline{r}_n^{\h{0.5pt}2}$ as $n \to \infty$ lead  to different energy minimization problems satisfied by the limiting map of the blow--up sequence. In the following, $(a_n, \overline{r}_n)$ is called in the small--scale, intermediate--scale and large--scale regimes if $a_n\overline{r}_n^{\h{0.5pt}2}$  converges to $0$, some finite positive number $L $ and $\infty$, respectively. If the centers of the blow--up locations are at $0$, then we have the following energy minimization problems satisfied by the limiting map of the blow--up sequence: \begin{table}[h] 
\centering
\bgroup
\def\arraystretch{1.35}
\setlength\tabcolsep{5pt}
\begin{tabular}{|c|c|c|c|c|}
\hline
 &\textbf{Signorini problem}
 & \textbf{Non--Signorini problem} \\ \hline

  Small--scale regime
  &   \vtop{
  \hbox{\strut Dirichlet energy in $\overline{M}_k$} }
  &  \vtop{
  \hbox{\strut  Dirichlet energy in $M_k$}  } \\
\hline

  Intermediate--scale regime
  & \vtop{
\hbox{\strut N.A.}}
  &  $\mathcal{E}^{sc}_L$--energy in $M_k$ \\ \hline

 Large--scale regime
 & \vtop{
\hbox{\strut N.A.}}
 & \vtop{
  \hbox{\strut \h{33pt} Dirichlet energy in $M_k$}
\hbox{\strut (one component is a constant function)}
 } \\ \hline
\end{tabular}
\egroup
\caption{Energy minimization for limit of blow--up sequence (balls centering  at the origin)}
\end{table}

\noindent Note that $M_k$ and $\overline{M}_k$ are two configuration spaces given in (\ref{configuration space M_k, origin}). The energy $\mathcal{E}_L^{sc}$ is defined in Lemma \ref{strong H1 convergence in intermediate case, origin}. Due to the axial symmetry of $w_{a_n, b}^+$, when the center of the blow--up location is at $0$, the limit of the blow--up sequence does not satisfy any Signorini--type obstacle problem in the intermediate and large--scale regimes. However, when the centers of the blow--up locations are different from $0$, the Signorini--type obstacle problems might occur  in all three regimes. See Table 2 below: \begin{table}[h] 
\centering
\bgroup
\def\arraystretch{1.35}
\setlength\tabcolsep{5pt}
\begin{tabular}{|c|c|c|c|c|}
\hline
 &\textbf{Signorini problem}
 & \textbf{Non--Signorini problem} \\ \hline

  Small--scale regime
  &   \vtop{
  \hbox{\strut Dirichlet energy in $\overline{\mathfrak{M}}_k$} }
  &  \vtop{
  \hbox{\strut  Dirichlet energy in $\mathfrak{M}_k$}  } \\
\hline

  Intermediate--scale regime
  & \vtop{
\hbox{\strut $E_{L, h}$--energy in $\overline{\mathfrak{M}}_k$}}
  &  $E_{L, h}$--energy in $\mathfrak{M}_k$ \\ \hline

 Large--scale regime
 & \vtop{
\hbox{\strut Dirichlet energy in $\overline{N}_k$}}
 & \vtop{
  \hbox{\strut  Dirichlet energy in $N_k$}
 } \\ \hline
\end{tabular}
\egroup
\caption{Energy minimization for limit of blow--up sequence (balls in $\mathscr{J}$)}
\end{table}

\noindent Note that in Table 2, $\mathfrak{M}_k$, $\overline{\mathfrak{M}}_k$, $N_k$ and $\overline{N}_K$ are configurations spaces defined in (\ref{configuration space M_k}) and (\ref{N_k and bar N_k}), respectively. The energy $E_{L, h}$ is given in  Lemma \ref{strong H1 convergence in intermediate case}. It is the three possible regimes associated with the blow--up sequence that make our analysis more complicated than the harmonic map case studied in \cite{Y20}.
\end{rmk}
\subsubsection{Biaxial--ring disclination}\vspace{0.5pc}

Let  $w_{a, b}^+ = \mathscr{L} \big[ u_{a, b}^+\big]$ be a biaxial--ring solution with $b \in \mathrm{I}_-$. Now, we discuss the reason why $w_{a, b}^+$ induces ring disclinations when $a$ is large.  In the following, when we discuss $u_{a, b}^+$,  the notation $T$ refers to the set $\big\{ \big(\rho, 0 \big) : \rho \in [\h{1pt}0, 1\h{0.5pt}]\h{1pt}\big\}$. When we discuss $w_{a, b}^+$,  $T$ is the flat boundary of $B_1^+$. Notice that by the strict positivity of $u^+_{a, b; 1}$ on $\mathbb{D}$ (see Lemma \ref{sign of Un1}), $w_{a, b}^+$ cannot yield any isotropic point on $T$. Here $u_{a, b; j}^+$ denotes the $j$--th component of the vector field  $u_{a, b}^+$. If $w_{a, b}^+$ admits  disclination on $T$, then the disclination must be negative uniaxial. By the $\mathscr{R}$--axial  symmetry of the biaxial--ring solutions, we have $u_{a, b; 3}^+ = 0$ on $T$. It therefore can be shown from (\ref{evalue of D with u= u^+_a,b}) that the points on $T$ where $u^+_{a, b; 1} = \sqrt{3}\h{1pt}u^+_{a, b; 2}$ must be  negative uniaxial locations. For the solutions $u_{a, b}^+$, when $a$ is suitably large, $u_{a, b; 2}^+(0, 0) > 0$  in that $u_{a, b}^+$ converges to some $u_b^+$ uniformly near the origin as $a \to \infty$, and $u_b^+(0, 0) = (0, 1, 0)^\top$. It then turns out that $u_{a, b; 1}^+ - \sqrt{3}\h{1pt}u^+_{a, b; 2} < 0$  at the origin when $a$ is large. Meanwhile, the boundary condition (\ref{boundary condition of u_a, intro}) induces that $u^+_{a, b; 1} - \sqrt{3}\h{1pt}u^+_{a, b; 2} = \sqrt{3}\h{1pt}H_a > 0$ at the right--end point $(1, 0)$. Hence, one can use the continuity of  $u_{a, b}^+$ to prove the existence of points on $T$ on which $u_{a, b; 1}^+ - \sqrt{3}\h{1pt}u_{a, b; 2}^+ = 0$. By the analyticity of the solutions $u_{a, b}^+$, the number of these points is finite. There must exist a point on $T$ so that in a small neighborhood of this point on $T$, the value of $u_{a, b; 1}^+ - \sqrt{3}\h{1pt}u_{a, b; 2}^+$ varies from negative to positive, as $\rho$ increases. Moreover, $u^+_{a, b; 1} - \sqrt{3}\h{1pt}u^+_{a, b; 2}$ vanishes at this point. This location gives a biaxial--ring disclination. Readers may refer to Section 6 for the details.
\subsubsection{Split--core--segment disclination}\vspace{0.5pc}

Let  $w_{a, c}^- = \mathscr{L} \big[ u_{a, c}^-\big]$ be a split--core solution with $c \in (0, 1)$ and denote by $u_{a, c; j}^-$ the $j$--th component of $u_{a, c}^-$. For large $a$, we have $u_{a, c; 2}^- < 0$ at the origin in that $u_{a, c}^-$ converges to some $u_c^-$ uniformly near the origin as $a \to \infty$, and $u_c^-(0, 0) = (0, -1, 0)^\top$. The boundary condition (\ref{boundary condition of u_a, intro}) infers that $u_{a, c; 2}^- = H_a > 0$ at the north pole. For large $a$, the solution $w_{a, c}^-$ must admit at least one zero on $l_z^+$, the positive part of the $z$--axis. As $a \to \infty$, these zero locations usually converge, at least up to a subsequence, to some singularity of the limiting map (see Lemma \ref{accu of zero}). Therefore, near the zeros of $w_{a, c}^-$, the amplitude of $w_{a, c}^-$ decays to $0$ sharply from values close to $1$. So far, most convergence results in the Landau--de Gennes theory are valid only on the places strictly away from the zeros with a positive distance independent of $a$. Compared with the size of the core regions  which is   approximately of the order $\mathrm{O}\big(a^{-1/2}\big)$, these places where we have  the uniform convergence of $w_{a, c}^-$  are far away from the core regions. We are  lack of nice uniform convergence of $w_{a, c}^-$ as $a \to \infty$ near the zeros of $w_{a, c}^-$. It is this reason that yields the major difficulty in our studies of the split--core--segment disclination, particularly in the core regions. \vspace{0.2pc}

However, in light of the three eigenvalues given in (\ref{evalue of D with u= u^+_a,b}), the amplitude of $w_{a, c}^-$, equivalently $u_{a, c}^-$, is not important. Denote by $\lambda_{a, c; j}^-$ ($j = 1, 2, 3$) the three eigenvalues in (\ref{evalue of D with u= u^+_a,b}) computed in terms of $u_{a, c}^-$. To compare relative quantitative relationships of the values $ \lambda_{a, c; j}^- $ ($j = 1, 2, 3$) in the core regions is equivalent to compare the quantitative relationships of their scaled values $ \frac{\lambda_{a, c; j}^-}{ | u_{a, c}^-|}$ ($j = 1, 2, 3$), provided that $u_{a, c}^-$ has single zero in each core region. Here $ \frac{\lambda_{a, c; j}^-}{ | u_{a, c}^-|}$ ($j = 1, 2, 3$) depend only on the normalized vector field of $u_{a, c}^-$. Based on this consideration, the mutual distances of the zeros of $w_{a, c}^-$ are studied in the item (3) of Proposition \ref{strict isolation}. More precisely, in the item (3) of Proposition \ref{strict isolation}, the zeros of $w_{a, c}^-$ are shown to be well--apart from each other in the sense that  their mutual distances have a strictly positive lower bound independent of $a$.  Hence, the scaled values  $ \frac{\lambda_{a, c; j}^-}{ | u_{a, c}^-|}$ are indeed well--defined in each core region except at the associated zero. The well--apartness result of the zeros of $w_{a, c}^-$ is a consequence of the non--degeneracy result in the item (2) of Proposition \ref{strict isolation}. To prove the non--degeneracy result, we systematically apply the division trick  of Mironescu \cite{M96} on the Ginzburg--Landau equation and the blow--down analysis of Lin--Wang \cite{LW02}. Our proof is also motivated by the work of  Millot--Pisante \cite{MP10} on the three dimensional Ginzburg--Landau functional. Now we briefly discuss the key ideas used in the proof of the non--degeneracy result. Let $\big\{w_{a_n, c}^-\big\}$ be a sequence of split--core solutions with a zero $z_n$ on the $z$--axis. Without loss of generality, we can assume $z_n$ is on the positive part of the $z$--axis. Moreover, we pick up a sequence of radii, denoted by $\big\{r_n\big\}$, which converges to $0$ as $n \to \infty$. It is crucial to understand the limits of the blow--up sequence $w^{(n)} (\zeta) := w_{a_n, c}^-\big(z_n + r_n\h{1pt}\zeta\big)$ as $n \to \infty$. Here we have three limiting regimes to study. If $a_n r_n^2 \to 0$ as $n \to \infty$, in Lemma \ref{lim map in interior core}, the limit map of $w^{(n)}$ is shown to be $0$. If $a_n r_n^2 \to \infty$ as $n \to \infty$, in Lemma \ref{w^infty is homogeneous zero}, the limiting map is shown to be $0$--homogeneous. In light of the results in \cite{Y20}, the limiting map equals $\Lambda_+$ or $\Lambda_-$ in the large--scale regime. See  (\ref{def of lambda plus and minus}). The most interesting regime is the intermediate regime where $a_n r_n^2 \to L$ for some $L > 0$. In this regime, $r_n$ is comparable with the size of the core regions. Let $w^\infty_\star$ be the limit of $w^{(n)}$ in the intermediate regime. Note that $w^\infty_\star$ is globally defined on $\mathbb{R}^3$. In Lemma \ref{lim map in interior core}, we consider the energy--minimal property of $w_\star^\infty$. Furthermore, in Lemma \ref{characterization of wstar infty}, we rigorously characterize the limit map $w_\star^\infty$ by showing that $w_\star^\infty = f\big( \h{1pt}\sqrt{L\mu} \h{0.5pt}|\h{0.5pt}\zeta\h{0.5pt}|\h{0.5pt}\big) \Lambda$. Here $\Lambda = \Lambda_+$ or $\Lambda_-$. $f$ is a radial function satisfying the ODE problem in (2) of Proposition \ref{strict isolation}. It is in the proof of Lemma \ref{characterization of wstar infty} where the division trick and the blow--down analysis come into play. Lemma \ref{characterization of wstar infty} implies that asymptotically near the zeros, the amplitude of $w_{a, c}^-$ is approximately homogeneous with respect to the angular variables. The phase mapping $\frac{w_{a, c}^-}{| \h{0.5pt}w_{a, c}^- \h{0.5pt}|}$ indeed  approximately equals  $\Lambda_+$ or $\Lambda_-$, the limit map of $w^{(n)}$ in the large--scale regime $a_n r_n^2 \to \infty$. With the convergences of $w^{(n)}$ in different regimes, we prove the uniform convergence of $\frac{w_{a, c}^-}{| \h{0.5pt}w_{a, c}^- \h{0.5pt}|}$ near zeros and as $a \to \infty$. See Proposition \ref{num and asy of singu}. The convergence in Proposition \ref{num and asy of singu} is sufficient for us to compare the quantitative relationships of the values $ \frac{\lambda_{a, c; j}^-}{ | u_{a, c}^-|}$ ($j = 1, 2, 3$) in the core regions. Moreover, asymptotic behavior of the director field   near the core regions can also be studied. Note that by (\ref{hat kappa_3, intro}), the director field $\kappa\h{0.5pt}[\h{0.5pt}u_{a, c}^-\h{0.5pt}]$ depends only on the normalized mapping $\frac{u_{a, c}^-}{| \h{0.5pt}u_{a, c}^- \h{0.5pt}|}$ as well.


\subsection{Notations\vspace{0.5pc}}

Most notations will be given at the first places when they will be used. Here we give some notations that will be frequently used in the following sections. \begin{itemize}\item For a vector field $w_{a_1, ..., a_n}$ where $a_1$, ..., $a_n$ are some notations or parameters, we use $w_{a_1, ..., a_n ;j}$ to denote its $j$--th component. Sometimes, we also use $[\h{1pt} w_{a_1, ..., a_n} \h{1pt}]_j$ to denote the $j$--th component of $w_{a_1, ..., a_n}$ interchangeably. For a vector field denoted by a simple notation $w$, we directly use $w_j$ to represent its $j$--th component;
\item Given $k \in \mathbb{N}$, $p \in [\h{0.5pt}1, \infty\h{0.5pt}]$ and a set $\Omega$, the notation $\| \cdot \|_{k, p; \Omega}$ denotes the norm in  $W^{k, p}(\Omega)$. Moreover, we use $\| \cdot \|_{p; \Omega}$ to denote the norm in $L^p(\Omega)$;
\item To compare two quantities $A$ and $B$, we use $A \lesssim_{c_1, ..., c_n} B$ to denote $A \leq c B$ with $c$ depending on $c_1$, ..., $c_n$. Here $c$ can also depend on the parameter $\mu$. Without parameters, we also use $A \lesssim B$ to denote $A \leq c B$ with $c$ a constant depending probably only on the parameter $\mu$;
\item Letting $a$ and $b$ be two quantities, we define $a \vee b := \mathrm{max}\big\{a, b \big\}$;
\item With  $\delta_{jk}$ denoting the standard Kronecker delta, $e_j$ is the unit vector in $\mathbb{R}^5$ whose $k$--th component $e_{j; k} = \delta_{jk}$. Here $j, k = 1, ..., 5$. $e^*_j$ is the unit vector in $\mathbb{R}^3$ whose $k$--th component $e^*_{j; k} = \delta_{jk}$. Here $j, k = 1, ..., 3$;
\item For a vector field $w$, we use $\widehat{w}$ to denote its normalized vector field $w \big/ |w|$. Letting $w$ be a non--zero $n+1$--vector and $\mathbb{S}^n$ be the standard unit sphere in $\mathbb{R}^{n+1}$ with center $0$, we also use $\Pi_{\mathbb{S}^n}[\h{0.5pt}w\h{0.5pt}]$ to denote the normalized vector of $w$ interchangeably;
\item Letting $f \in L^1\left( \Omega; \mathrm{d}\nu\right)$, where $\mathrm{d} \nu$ is a measure on some set $\Omega$, we use $\displaystyle \dashint_\Omega f \h{1pt}\mathrm{d} \nu$ to denote the average of $f$ on $\Omega$ with respect to the measure $\mathrm{d} \nu$;
\item Given a set $\Omega$ in some Euclidean space, $\Omega^+$ ($\Omega^-$ resp.) contains all points in $\Omega$ whose last component are positive (negative resp.).
\end{itemize} 

\cftaddtitleline{toc}{section}{\normalsize Part I: Existence of multiple $\mathscr{R}$--axially symmetric solutions \vspace{0.5pc}}{}


\section{Energy--decay estimate on $B_r$} 

We use a Luckhaus--type argument to prove Proposition \ref{small energy implies energy decay}. In this section, we consider the case in which the blow--up location is at $0$. Due to the radial symmetry of the balls $B_r$, this case is easier to be dealt with than the case in which the blow--up balls are in the family $\mathscr{J}$.
\subsection{Blow--up sequence and some preliminary results}
Fix a  constant $\nu_0 \in \big(0, 1/2\big)$ which depends only on $\mu$ and will be determined later in the proof. Suppose  that Proposition \ref{small energy implies energy decay} fails. There exist $a_n$ and $\epsilon_n$ with
\begin{equation}
a_n \longrightarrow \infty \h{15pt} \text{and}  \h{15pt} \epsilon_n \longrightarrow 0 \qquad \text{as\,\,}n \to \infty
    \label{convergence of parameters, origin}
\end{equation}
so that for any $n \in \mathbb{N}$, we can find a radius $r_n$ with which the followings hold by the mapping $w_n := w_{a_n, b}^+$:
\begin{equation}
\textup{(i). } \mathcal{E}_{\h{0.2pt}a_n, \h{0.2pt}r_n}\left[ w_n\right] \h{2pt}<\h{2pt}\epsilon_n;\quad \textup{(ii). }   \mathcal{E}_{\h{0.2pt}a_n, \h{0.2pt}\nu_0\h{0.1pt}r_n}\left[ w_n\right] \h{2pt}>\h{2pt}r_n^{3/2}; \quad \textup{(iii). } \mathcal{E}_{\h{0.2pt}a_n,  \h{0.2pt}\nu_0\h{0.1pt}r_n}\left[ w_n\right] \h{2pt}>\h{2pt}\dfrac{1}{2}\h{1pt}\mathcal{E}_{\h{0.2pt}a_n, \h{0.2pt}r_n}\left[ w_n\right].
\label{basic assumptions, origin}
\end{equation}
Here we have dropped the parameter $\mu$ and the origin $0$ from the subscripts and simply use $\mathcal{E}_{\h{0.2pt}a_n, \h{0.2pt}r}\left[ w_n\right]$ to denote $\mathcal{E}_{a_n, \mu; \h{1pt}0, r}\left[ w_n\right]$. Moreover, we assume $b \in \mathrm{I}_-$.  Define
\begin{equation}
s_n^2 :=\mathcal{E}_{\h{0.2pt}a_n, \h{0.2pt}r_n}\left[ w_n\right], \h{5pt} y_n :=\displaystyle\dashint_{B_{r_n}} w_n, \h{5pt} \mathcal{W}_n\big(\zeta\big) : = w_n\big(r_n\zeta\big), \h{5pt} \mathcal{W}_n^{sc}\big(\zeta\big):=\dfrac{\mathcal{W}_n\big(\zeta\big)-y_n}{s_n}, \h{10pt} \text{where } \zeta \in B_1.
\label{s_n and y_n,  origin}
\end{equation} Then by  Poincar\'{e}'s inequality,  the scaled mappings $\Big\{\mathcal{W}_n^{sc}\Big\}$ is uniformly bounded in $H^1\left(B_1\right)$. Hence, there is a subsequence, which we still denote by $\Big\{ \mathcal{W}_n^{sc} \Big\}$, so that as $n \rightarrow \infty$,
\begin{equation}
\mathcal{W}_n^{sc} \longrightarrow \mathcal{W}_\infty^{sc} \h{20pt} \text{ weakly in } H^1\left(B_1\right), \text{ strongly in }L^2\left(B_1\right)\h{5pt}\text{and strongly in } L^2\left( T\right).
\label{conv. of u^L_n}
\end{equation}
Recall that $T$ is the flat boundary of $B_1^+$. Let \begin{align}\label{def of F_n}F_n \left[w\right] := \mu\left[D_{a_n}-3\sqrt{2}\h{0.5pt}S\left[ w \right]+\dfrac{a_n}{2}\h{0.5pt}\big(\h{0.5pt}|\h{0.2pt}w\h{0.2pt}|^2-1\big)^2\right].\end{align} (iii) in (\ref{basic assumptions, origin}) induces
\begin{equation}
\nu_0^{-1} \int_{B_{\nu_0}} \big|\nabla \mathcal{W}^{sc}_n \big|^2 + \left(\dfrac{r_n}{s_n}\right)^2F_n\big(\mathcal{W}_n\big) > \dfrac{1}{2},\h{15pt}\text{for all $n$.}
    \label{limit of contradictory inequality, origin}
\end{equation}
To estimate the potential term in the above inequality, we need
\begin{lem}
$s_n + \dfrac{r_n}{s_n} \longrightarrow 0 $ as $n \to \infty$.
\label{s_n,r_n/s_n tend to 0,w}
\end{lem}
\begin{proof}[\bf Proof]
From the condition (ii) in $(\ref{basic assumptions, origin})$, it turns out that
$$r_n\h{0.2pt}s_n^2  \h{2pt}\geq\h{2pt} \int_{B_{\nu_0 \h{0.2pt}r_n}} \big| \nabla w_n \big|^2 + F_n\big(w_n\big) \h{2pt}>\h{2pt} \nu_0\h{0.2pt}r_n^{5/2}.$$ This lemma then follows since  by (i) in (\ref{basic assumptions, origin}) and the convergence of $\epsilon_n$ in (\ref{convergence of parameters, origin}), $s_n \longrightarrow 0$ as $n \rightarrow \infty$. 
\end{proof}

In the following, we discuss some facts related to  $\big\{y_n\big\}$ in (\ref{s_n and y_n,  origin}). By the $\mathscr{R}$--axial symmetry of $w_n$, \begin{align}\label{yn in origin case} y_n = \big(0, 0, y_{n; 3}, 0, 0\big)^\top, \h{15pt}\text{where $y_{n; 3}$ is a constant in $\mathbb{R}$}.\end{align}
Owing to (i) in (\ref{basic assumptions, origin}), the convergence of $\epsilon_n$ in (\ref{convergence of parameters, origin}) and the uniform boundedness in item (2) of Remark \ref{rmk on sign and gamma conv}, up to a subsequence, $\mathcal{W}_n$ converges strongly in $H^1(B_1)$ to a constant vector $y_*$. Meanwhile, the $y_n$ defined in (\ref{s_n and y_n,  origin}) (see also (\ref{yn in origin case})) converges to $y_*$ as well when $n \rightarrow \infty$. Obviously for some constant $y_{*; 3} \in \mathbb{R}$, we have $y_* = \big(0, 0, y_{*; 3}, 0, 0\big)^\top$. If samely as before we use $T$ to denote $$\Big\{ \big(\zeta_1, \zeta_2, 0 \big) : \zeta_1^2 + \zeta_2^2 \leq 1\Big\},$$ then trace theorem infers that $\mathcal{W}_n$ converges to $y_*$ strongly in $L^2\left(T\right)$. Notice that   $\mathcal{W}_{n; 3} \geq H_{a_n}\h{0.5pt}b$ on $T$ in the sense of trace. Taking $n \rightarrow \infty$, we  then obtain  
\begin{equation}
y_n \longrightarrow y_* = \big(0,0,y_{*;3},0,0\big)^\top, \h{5pt} \text{where $y_{*;3}$ is a constant satisfying $y_{*; 3} \geq b$.}
\label{y_*=(0,0,y_*,3,0,0)}
\end{equation}If in addition it holds \begin{eqnarray*}
\liminf_{n \rightarrow \infty} \h{2pt} \left|  \dfrac{H_{a_n}\h{0.2pt}b - y_{n; 3}}{s_n} \right| \h{2pt}< \h{2pt}\infty,
\end{eqnarray*}then there exists a constant $w_* \in \mathbb{R}$ so that up to a subsequence,  \begin{eqnarray} \lim_{n \rightarrow \infty} \h{2pt} \dfrac{H_{a_n} b - y_{n; 3}}{s_n} =  w_*.
\label{pi_L convergence, origin}
\end{eqnarray}In this case, we have
\begin{lem}
If (\ref{pi_L convergence, origin}) holds, then $\mathcal{W}^{sc}_{\infty; 3} \h{1pt}\geq\h{1pt} w_{*}$ on $T$ in the sense of trace. 
\label{ratio convergence and limit boundary condition, origin}
\end{lem}
\begin{proof}[\bf Proof] The third component of $\mathcal{W}^{sc}_n$ can be decomposed into
$$\mathcal{W}^{sc}_{n; 3}=\dfrac{ \mathcal{W}_{n; 3} -  H_{a_n}\h{0.2pt}b }{s_n} + \dfrac{H_{a_n}\h{0.2pt}b  - y_{n; 3} }{s_n} \quad \text{on\,\,}T.$$
The lemma then follows by the Signorini obstacle boundary condition satisfied by $\mathcal{W}_n$ on $T$.
\end{proof}
By (2) in Remark \ref{rmk on sign and gamma conv}, Lemma \ref{ratio convergence and limit boundary condition, origin} and Fatou's lemma, the following results hold:
\begin{lem}
There exist an increasing positive sequence $\big\{\sigma_k\big\}$ which tends to $1$ as $k \to \infty$, a sequence of positive numbers $\big\{b_k\big\}$ and a subsequence of $\big\{\mathcal{W}_n\big\}$, still denoted by $\big\{\mathcal{W}_n\big\}$, such that 
\begin{enumerate}
  \item[$\mathrm{(1).}$] For any $k \in \mathbb{N}$, the  mappings $\mathcal{W}^{sc}_n$ and their weak $H^1\left(B_1\right)$--limit $\mathcal{W}^{sc}_\infty$ satisfy
  $$  \sup_{n \h{0.5pt}\in\h{0.5pt} \mathbb{N}\h{0.7pt}\cup\h{0.5pt} \{\infty\}} \int_{\p B_{\sigma_k}} \big|\h{0.2pt} \mathcal{W}^{sc}_n \h{0.2pt}\big|^2 +  \big| \nabla \mathcal{W}^{sc}_n \big|^2 \h{2pt} \leq \h{2pt} b_k \h{1pt}.$$Here $\mathcal{W}^{sc}_n$ is the scaled mapping of $\mathcal{W}_n$ given in (\ref{s_n and y_n, origin});
  \item[$\mathrm{(2).}$] For any $k \in \mathbb{N}$, the sequence $\mathcal{W}^{sc}_n$ converges to $\mathcal{W}^{sc}_\infty$ strongly in $L^2\left(\p B_{\sigma_k}\right)$ as $n \rightarrow \infty$;
  \item[$\mathrm{(3).}$] For any $n, k \in \mathbb{N}$, it holds $\mathcal{W}_{n; 3} \geq H_{a_n}\h{0.2pt}b$ on $\p B_{\sigma_k} \cap T$. Moreover, $\mathcal{W}_n$ satisfies $\big| \mathcal{W}_n\big| \leq H_{a_n}$ on $\p B_{\sigma_k}$;
\item[$\mathrm{(4).}$]  If (\ref{pi_L convergence, origin}) holds, then for any $k \in \mathbb{N}$, we have $\mathcal{W}^{sc}_{\infty;3} \geq w_*$ on $\p B_{\sigma_k}  \cap T$.
\end{enumerate}
\label{uniform convergence of u^L_n in S_sigma_k}
\end{lem}  \noindent Using $\{\sigma_k\}$ obtained in Lemma \ref{uniform convergence of u^L_n in S_sigma_k}, we introduce two configuration spaces:
\begin{equation}
\begin{split}
& M_k :=\Big\{w\in H^1(B_{\sigma_k};\mathbb{R}^5)\h{1pt}:\h{1pt}\text{$w $ is $\mathscr{R}$--axially symmetric in $B_{\sigma_k}$ and}\h{2pt} w=\mathcal{W}^{sc}_\infty \text{ on }\p B_{\sigma_k}\Big\}
;\\[2mm]
& \overline{M}_{\h{0.5pt}k} :=  \Big\{ w \in M_{\h{0.5pt}k} : w_3 \h{2pt}\geq\h{2pt}w_{*} \h{4pt}\text{on $T_{\sigma_k} := B_{\sigma_k} \cap T$} \Big\}.
\end{split}
\label{configuration space M_k, origin}
\end{equation}
These spaces will be used in Sections 2.2--2.4 as configuration spaces of $\mathcal{W}_\infty^{sc}$ in different limiting regimes.

\subsection{Energy--decay estimate in small--scale regime}
The minimization problem satisfied by the limiting map $\mathcal{W}^{sc}_\infty$ is different if $a_n r_n^2$ converges in different regime. In this section we suppose that $a_n r_n^2 \longrightarrow 0$ as $n \to \infty$. \vspace{0.2pc} We now prove  the following minimizing property of $\mathcal{W}^{sc}_\infty$ in the small--scale regime.
\begin{lem}
Suppose that $a_n r_n^2 \rightarrow 0$ as $n \rightarrow \infty$. For any natural number $k$, if it satisfies
\begin{eqnarray}
\liminf_{n \rightarrow \infty} \h{2pt} \left|  \dfrac{ H_{a_n}\h{0.2pt}b - y_{n; 3}}{s_n} \right| \h{2pt}= \h{2pt}\infty,
\label{infty case, origin}
\end{eqnarray}
then $\mathcal{W}^{sc}_\infty$ minimizes the Dirichlet energy within the configuration space $M_k$. If (\ref{pi_L convergence, origin}) holds, then $\mathcal{W}^{sc}_\infty$ minimizes the Dirichlet energy within the configuration space $\overline{M}_{\h{0.5pt}k}$. In both cases, $\mathcal{W}^{sc}_n$ converges to $\mathcal{W}^{sc}_{\infty}$ strongly in $H_{\mathrm{loc}}^1\big(B_1\big)$ as $n \rightarrow \infty$.
\label{strong H1 convergence, small, origin}
\end{lem}
\begin{proof}[\bf Proof]
\textbf{Step 1. Comparison map}\\
Suppose that $w$ is an arbitrary map in $ M_k $. Then we define
\begin{equation}
M_{n, R}\h{0.5pt} [\h{0.5pt}w\h{0.5pt}] :=\left\{ \begin{array}{lcl} y_n+R\h{1pt}s_n \h{1pt}\dfrac{w}{|\h{1pt}w\h{1pt}|\vee R}, \h{40pt}&&\text{if (\ref{infty case, origin}) holds;}\vspace{0.8pc}\\
  y_n^*+R\h{1pt}s_n \h{1pt}\dfrac{w - w_* \h{0.2pt}e_3 }{|\h{1pt}w - w_* \h{0.2pt}e_3\h{1pt}|\vee R}, &&\text{if (\ref{pi_L convergence, origin}) holds.}
\end{array}\right.
\label{definition of F_n^s, origin}
\end{equation}
In this definition, $R>0$ is a positive constant. The constant vector $y_n^*$   equals  $\left(H_{a_n}\h{0.2pt}b \right)\h{0.2pt}e_3$. Now we fix an arbitrary $s \in (0,1)$ and  introduce
\begin{equation}
v_{n, s, R}\big(\zeta\big):=\left\{
\begin{aligned}
&M_{n, R}[\h{0.5pt}w\h{0.5pt}]\left(\dfrac{\zeta}{1-s}\right) &\text{if }& \zeta \in B_{(1-s)\sigma_k}; \\[2mm]
&\dfrac{\sigma_k-|\h{1pt}\zeta\h{1pt}|}{s\sigma_k}M_{n, R}\big[\h{.5pt}\mathcal{W}^{sc}_\infty\h{.5pt}\big]\big(\h{1pt}\sigma_k \widehat{\zeta}\h{1pt}\big) + \dfrac{|\h{1pt}\zeta\h{1pt}|-(1-s)\sigma_k}{s\sigma_k}\mathcal{W}_n\big(\h{1pt}\sigma_k\widehat{\zeta}\h{1pt}\big)&\text{if }& \zeta \in B_{\sigma_k} \setminus B_{(1-s)\sigma_k}.
\end{aligned}
\right.
\label{definition of v_n, origin}
\end{equation} The map $v_{n, s, R}$ is our comparison map.  \vspace{0.4pc}\\
\textbf{Step 2. Upper bound}\vspace{0.2pc}\\
Notice that (1) in Lemma \ref{uniform convergence of u^L_n in S_sigma_k} infers the absolute continuity of $\mathcal{W}^{sc}_\infty$ and $\mathcal{W}_n$ on $\p B_{\sigma_k}$ near $\p B_{\sigma_k} \cap T$, with respect to the polar angle $\phi$. Therefore, $\mathcal{W}^{sc}_{\infty; 4} = \mathcal{W}^{sc}_{\infty; 5} = \mathcal{W}_{n; 4} = \mathcal{W}_{n; 5} = 0$ on $\p B_{\sigma_k} \cap T$. Combined this result with the fact that $w$ is $\mathscr{R}$--axially symmetric in $B_{\sigma_k}$, the comparison map $v_{n, s, R}$ is $\mathscr{R}$--axially symmetric in $B_{\sigma_k}$ as well. Here we have also used the definition of $y_n^*$ and  $y_n$ given in (\ref{yn in origin case}).  If (\ref{pi_L convergence, origin}) holds, then by (3)--(4) in Lemma \ref{uniform convergence of u^L_n in S_sigma_k} and the assumption that $w \in \overline{M}_k$, the third component of $v_{n, s, R}$, denoted by $ [ v_{n, s, R} ]_3$,  satisfies $ [ v_{n, s, R} ]_3  \geq H_{a_n}\h{0.2pt}b$ on $T_{\sigma_k}$. If (\ref{infty case, origin}) holds, then we have \begin{align}\label{ratio conv to negative infinity} \dfrac{{H_{a_n}\h{0.2pt}b - y_{n; 3}}}{s_n} \longrightarrow - \infty \h{15pt}\text{as $n \rightarrow \infty$,}
\end{align}since by Signorini obstacle condition, it satisfies \begin{align*}  \dfrac{{H_{a_n}\h{0.2pt}b - y_{n; 3}}}{s_n} \leq \dfrac{{ \mathcal{W}_{n, 3} - y_{n; 3}}}{s_n} \h{15pt}\text{on $T$, for all $n \in \mathbb{N}$.}
\end{align*}The right--hand side above indeed converges  almost everywhere to the trace of $\mathcal{W}^{sc}_{\infty; 3}$  on $T$  as $ n\rightarrow \infty$. So the limit in (\ref{ratio conv to negative infinity}) must diverge to $- \infty$ instead of $\infty$. Owing to (\ref{ratio conv to negative infinity}) and (\ref{definition of F_n^s, origin}), we still have $[v_{n, s, R}]_3 \geq H_{a_n}\h{0.2pt}b$ on $T_{\sigma_k}$ if (\ref{infty case, origin}) holds. Note that we should take $n$ suitably large with the largeness of $n$ depending on $R$. Now we apply the energy minimizing property of $\mathcal{W}_n$. It then turns out 
\begin{equation}
\int_{B_{\sigma_k}} \big| \nabla  \mathcal{W}^{sc}_n \big|^2 \h{1pt} \leq \h{1pt} \int_{B_{\sigma_k}}  s^{-2}_n  \big| \nabla  v_{n, s, R} \big|^2 +\left(\dfrac{r_n}{s_n}\right)^2F_n\left(\cdot\right)\bigg|^{v_{n, s, R}}_{\mathcal{W}_n}.
    \label{small scale basic inequality, origin}
\end{equation}Here and in what follows, $F_n(\cdot)\h{1pt}\Big|^{v_{n, s, R}}_{\mathcal{W}_n} := F_n\big(v_{n, s, R}\big) - F_n\big(\mathcal{W}_n\big)$. \vspace{0.3pc}\\
\noindent \textbf{Step 3. Estimate of the Dirichlet energy of $v_{n, s, R}$}\vspace{0.2pc}\\
Direct calculations yield 
\begin{equation}
    \int_{B_{\sigma_k}} \big| \nabla v_{n, s, R} \big|^2 =  \left(1 - s \right) \int_{B_{\sigma_k}} \Big|\h{0.5pt}\nabla M_{n, R}\h{0.5pt}[\h{0.5pt}w\h{0.5pt}]\h{0.5pt}\Big|^2 + \int_{B_{\sigma_k}\setminus B_{(1-s)\sigma_k}} \big| \nabla  v_{n, s, R} \big|^2.
    \label{decomp. of grad v_n in B_k, origin}
\end{equation}
Using the fact $w \in H^1\left(B_{\sigma_k};\mathbb{R}^5\right)$, we obtain from the definition of $M_{n, R}\left[w\right]$ in (\ref{definition of F_n^s, origin}) that
\begin{equation}
       s_n^{-2} \int_{B_{\sigma_k}} \Big|\h{0.5pt}\nabla M_{n, R}\h{0.5pt}[\h{0.5pt}w\h{0.5pt}]\h{0.5pt}\Big|^2 \h{2pt}\text{is independent of $n$ and converges to $\int_{B_{\sigma_k}} \big|\nabla w\big|^2$ as $R\rightarrow \infty$.}
    \label{conv. main part, origin}
\end{equation}
To deal with the second term on the right--hand side of (\ref{decomp. of grad v_n in B_k, origin}), we denote by $\big(\tau,\Phi,\Theta\big)$ the spherical coordinates in the $\zeta$--space. Here $\tau$ is the radial variable, $\Phi$ is the polar angle and $\Theta$ is the azimuthal angle. Then the Dirichlet energy of $v_{n, s, R}$ on $B_{\sigma_k} \setminus B_{(1 - s) \sigma_k}$ can be expressed by
\begin{equation}
    \int_{B_{\sigma_k}\setminus B_{(1-s)\sigma_k}} \big|\nabla  v_{n, s, R}\big|^2 = \int_{B_{\sigma_k}\setminus B_{(1-s)\sigma_k}}\big|\p_\tau v_{n, s, R}\big|^2 + \dfrac{1}{\tau^2}\big|\p_\Phi v_{n, s, R}\big|^2 + \dfrac{1}{\tau^2 \sin^2\Phi}\big|\p_\Theta v_{n, s, R}\big|^2.
    \label{decomp of grad v_n in B_k / B_(1-s), origin}
\end{equation}
The definition of $v_{n, s, R}$ on $ B_{\sigma_k} \setminus B_{(1-s)\sigma_k }$ (see (\ref{definition of v_n, origin})) induces
$$ \int_{B_{\sigma_k}\setminus B_{(1-s)\sigma_k}}\big|\p_\tau v_{n, s, R}\big|^2 \h{2pt}\leq\h{2pt} \big(\h{0.2pt} s \h{0.2pt} \sigma_k^3\h{0.2pt}\big)^{-1}\int_{\p B_{\sigma_k} }\big|M_{n, R}\h{0.5pt}\big[\mathcal{W}^{sc}_\infty\big] - \mathcal{W}_n \big|^2.$$
By the definition of $M_{n, R}$ in (\ref{definition of F_n^s, origin}), if (\ref{infty case, origin}) holds, then  we control $\p_\tau v_{n, s, R}$ as follows:
\begin{align*}
s_n^{-2}\int_{B_{\sigma_k}\setminus B_{(1-s)\sigma_k}}\big|\p_\tau v_{n, s, R}\big|^2 \h{2pt}\lesssim\h{2pt}\left( s\sigma_k^3 \right)^{-1} \left[ \int_{\p B_{\sigma_k} }\big| \mathcal{W}^{sc}_n - \mathcal{W}^{sc}_\infty\big|^2 + \int_{\p B_{\sigma_k}} \left|\h{1pt}  R\h{0.5pt}\dfrac{\mathcal{W}_\infty^{sc}}{ \big| \mathcal{W}_\infty^{sc}\big| \vee R } - \mathcal{W}_\infty^{sc} \right|^2\h{2pt}\right].
\end{align*}If (\ref{pi_L convergence, origin}) holds, then we use \begin{align*}
s_n^{-2}\int_{B_{\sigma_k}\setminus B_{(1-s)\sigma_k}}\big|\p_\tau v_{n, s, R}\big|^2
\h{2pt}&\lesssim\h{2pt} \left( s\sigma_k^3 \right)^{-1}\int_{\p B_{\sigma_k}}\left| \h{1pt}\big( \mathcal{W}^{sc}_n - \mathcal{W}^{sc}_\infty\big) - \left[ \dfrac{y_n^* - y_n}{s_n} - w_* \h{1pt}e_3 \right]\h{1.5pt} \right|^2\\[2mm]
&+\h{2pt}\left( s\sigma_k^3 \right)^{-1}\int_{\p B_{\sigma_k}}   \left|\h{1pt} R\h{1pt} \dfrac{\mathcal{W}_\infty^{sc} - w_* e_3}{ \big| \mathcal{W}_\infty^{sc} - w_* e_3\big| \vee R } - \big(\mathcal{W}_\infty^{sc} - w_* e_3 \big) \right|^2.
\end{align*}
In both cases, we can apply (2) in Lemma \ref{uniform convergence of u^L_n in S_sigma_k} to get 
\begin{equation}
 \limsup_{n \to \infty} s_n^{-2}\int_{B_{\sigma_k}\setminus B_{(1-s)\sigma_k}}\big|\p_\tau v_{n, s, R}\big|^2 \h{2pt}\lesssim\h{2pt}\left( s\sigma_k^3 \right)^{-1}\int_{\p B_{\sigma_k}}   \left|\h{1pt} R\h{1pt} \dfrac{\mathcal{W}_\infty^{sc} - Y_*}{ \big| \mathcal{W}_\infty^{sc} - Y_*\big| \vee R } - \big(\mathcal{W}_\infty^{sc} - Y_*\big) \right|^2.
\label{conv. of radial derivative of v_n, origin}
\end{equation} Here for our convenience, $Y_* = 0$ if (\ref{infty case, origin}) holds. If (\ref{pi_L convergence, origin}) holds, then $Y_* = w_* e_3$. \vspace{0.2pc}

Still by the definition of $v_{n, s, R}$ on $B_{\sigma_k}\setminus B_{(1-s)\sigma_k}$, we observe  that
\begin{align*}
\int_{B_{\sigma_k}\setminus B_{(1-s)\sigma_k}} \dfrac{1}{\tau^2}\left[ \big|\p_\Phi v_{n, s, R}\big|^2 + \left(\dfrac{\big|\p_\Theta v_{n, s, R}\big|}{ \sin \Phi}\right)^2\right] \h{2pt}\lesssim\h{2pt}  \int_{B_{\sigma_k}\setminus B_{(1-s)\sigma_k}} \dfrac{1}{\tau^2}\bigg[\h{1.5pt}\Big| \nabla M_{n, R}\big[ \mathcal{W}^{sc}_\infty \big]\Big|^2+\Big| \nabla \mathcal{W}_n \Big|^2\bigg]\big(\sigma_k \widehat{\zeta}\h{1pt}\big).
\end{align*}Therefore, (1) in Lemma \ref{uniform convergence of u^L_n in S_sigma_k} induces\begin{align}
    s_n^{-2} \int_{B_{\sigma_k}\setminus B_{(1-s)\sigma_k}} \dfrac{1}{\tau^2}\left[ \big|\p_\Phi v_{n, s, R}\big|^2 + \left(\dfrac{\big|\p_\Theta v_{n, s, R}\big|}{ \sin \Phi}\right)^2\right] \h{2pt}&\lesssim\h{2pt} \dfrac{s}{\sigma_k} \int_{\p B_{\sigma_k}} \Big| \nabla \mathcal{W}_n^{sc}\Big|^2+\left| R\h{1pt} \nabla \dfrac{\mathcal{W}_\infty^{sc} - Y_*}{ \big| \mathcal{W}_\infty^{sc} - Y_*\big| \vee R } \right|^2 \nonumber\\[2mm]
&\lesssim\h{2pt}\dfrac{s}{\sigma_k} \left[ b_k+ \int_{\p B_{\sigma_k}} \left| R\h{1pt} \nabla \dfrac{\mathcal{W}_\infty^{sc} - Y_*}{ \big| \mathcal{W}_\infty^{sc} - Y_*\big| \vee R } \right|^2\right].
     \label{ineq of phi derivative of v_n, origin}
\end{align}
By (\ref{decomp of grad v_n in B_k / B_(1-s), origin}) and the estimates in (\ref{conv. of radial derivative of v_n, origin})--(\ref{ineq of phi derivative of v_n, origin}), it holds
$$ \lim_{R \to \infty} \limsup_{n \to \infty}  s_n^{-2} \int_{B_{\sigma_k}\setminus B_{(1-s)\sigma_k}}\big|\nabla  v_{n, s, R} \big|^2 \h{2pt}\lesssim\h{2pt}\dfrac{s}{\sigma_k} \left[ b_k+ \int_{\p B_{\sigma_k}} \big|   \nabla \mathcal{W}_\infty^{sc}  \big|^2\right] \h{2pt}\lesssim\h{2pt}\dfrac{s}{\sigma_k}\h{1pt}b_k.$$
The last estimate above has also used (1) in Lemma \ref{uniform convergence of u^L_n in S_sigma_k}. Now we divide $s_n^2$ from both sides of (\ref{decomp. of grad v_n in B_k, origin}). In light of (\ref{conv. main part, origin}) and the last limit, it then follows
\begin{equation}
  \lim_{s \to 0} \lim_{R \to \infty} \lim_{n \to \infty}   s_n^{-2}\int_{B_{\sigma_k}} \big|\nabla  v_{n, s, R}\big|^2  = \int_{B_{\sigma_k}}\big|\nabla w\big|^2.
    \label{conv. of grad v_n, origin}
\end{equation}
\textbf{Step 4. Estimate of potential term}\vspace{0.2pc}\\
For the potential term, we notice that
\begin{equation}
\int_{B_{\sigma_k}}F_{n}\left(\cdot\right)\Big|_{\mathcal{W}_n}^{v_{n, s, R}} = I^s_1+I^s_2,
\label{I_1 to I_2, origin}
\end{equation}
where the terms on the right--hand side above are defined and estimated as follows.\vspace{0.2pc}\\
\textbf{Estimate of $I^s_1$.} $I^s_1$ is defined by
\begin{align}\label{expression of Is1}I^s_1 := - 3\sqrt{2}\h{0.2pt}\mu\int_{B_{\sigma_k}}\big(v_{n, s, R} - \mathcal{W}_n\big)\cdot \int_0^1 \nabla_w S \h{1pt}\Big|_{w\h{0.2pt}=\h{0.2pt}t\h{0.1pt}v_{n, s, R}+(1-t)\h{0.1pt}\mathcal{W}_n} \mathrm{d} \h{1pt} t.\end{align} If (\ref{infty case, origin}) holds, then the definition (\ref{definition of F_n^s, origin}) yields $\big| M_{n, R}\left[ w \right] - y_n \big| \leq s_n\h{0.2pt}R$. If (\ref{pi_L convergence, origin}) holds, then we can take $n$ large enough depending on $w_*$ so that $\big| M_{n, R}\left[ w \right] - y_n \big| \leq 3\h{0.2pt}s_n\h{0.2pt}R$ for all $R > | w_* |$. In both cases, we have \begin{align}\label{diff between Mnr and yn} \big| M_{n, R}\left[ w \right] - y_n \big| \leq 3\h{0.2pt}s_n\h{0.2pt}R \h{15pt}\text{on $B_{\sigma_k}$, for large $n$ and $R$.}
\end{align}As for $\mathcal{W}_n$, by (1) in Lemma \ref{uniform convergence of u^L_n in S_sigma_k}, we have \begin{align*} \int_{B_{\sigma_k} \setminus B_{(1 - s) \sigma_k}}\left| \dfrac{\mathcal{W}_n\big(\sigma_k \widehat{\zeta}\h{0.2pt}\big) - y_n}{s_n}\right|^2 \leq \dfrac{s}{\sigma_k} \int_{\p B_{\sigma_k}} \big| \mathcal{W}_n^{sc} \big|^2 \leq \dfrac{s}{\sigma_k}\h{0.2pt}b_k.
\end{align*}Utilizing this estimate, (\ref{diff between Mnr and yn}) and the definition of $v_{n, s, R}$ in (\ref{definition of v_n, origin}) then yield \begin{align}\label{estimate of vnsr and yn, origin} s_n^{-2} \int_{B_{\sigma_k}} \big| \h{1pt}v_{n, s, R} - y_n \h{1pt}\big|^2 \h{2pt}\lesssim\h{2pt} R^2 + \dfrac{s}{\sigma_k}\h{0.2pt}b_k, \h{15pt}\text{for large $n$ and $R$.}
\end{align}The $\nabla_w S$ in (\ref{expression of Is1}) is quadratic in terms of the variable $w$. By H\"{o}lder's inequality, Sobolev's imbedding and (\ref{estimate of vnsr and yn, origin}), it turns out \begin{align}\label{bounds of I1s} \big|\h{1pt} I_1^s\h{1pt}\big| \h{2pt}&\lesssim\h{2pt}s_n  \left[\h{1pt} \big\| \mathcal{W}_n \big\|_{4; B_{\sigma_k}} + \big\| v_{n, s, R} \big\|_{1, 2; B_{\sigma_k}}\right]^2 \left[ \left( s_n^{-2} \int_{B_{\sigma_k}} \big| v_{n, s, R} - y_n \big|^2 \right)^{1/2} + \left( \int_{B_{\sigma_k}} \big| \mathcal{W}_n^{sc} \big|^2 \right)^{1/2} \right] \nonumber\\[2mm]
&\lesssim \h{2pt}s_n  \left( R^2 + \dfrac{s}{\sigma_k}\h{0.2pt}b_k \right)^{1/2}\left[\h{1pt} \big\| \mathcal{W}_n \big\|_{4; B_{\sigma_k}} + \big\| v_{n, s, R} \big\|_{1, 2; B_{\sigma_k}}\right]^2.
\end{align}Due to the above estimate, the $L^\infty$--bound of $\mathcal{W}_n$ on $B_{\sigma_k}$, (\ref{decomp. of grad v_n in B_k, origin})--(\ref{ineq of phi derivative of v_n, origin}),  (\ref{estimate of vnsr and yn, origin}) and Lemma \ref{s_n,r_n/s_n tend to 0,w}, we obtain
\begin{equation}
 \lim_{n \rightarrow \infty}\left(\dfrac{r_n}{s_n}\right)^2\big| \h{0.5pt} I^s_1 \h{0.5pt} \big| = 0.
\label{I_1^s, origin}
\end{equation}
\textbf{Estimate of $I^s_2$.} $I^s_2$ is defined by
\begin{align} I^s_2 := 2 \h{1pt}a_n \mu \int_{ B_{\sigma_k} }\big( v_{n, s, R} - \mathcal{W}_n \big) \cdot \int_0^1 \Big( \mathcal{W}_n + t \big( v_{n, s, R} - \mathcal{W}_n \big) \Big) \left( \Big| \mathcal{W}_n + t \big( v_{n, s, R} - \mathcal{W}_n \big) \Big|^2 - 1 \right)  \mathrm{d} \h{1pt}  t. \label{def of I2s}\end{align}
Using the $L^\infty$--bounds of $v_{n, s, R}$ and $\mathcal{W}_n$, for $R$ suitably large, we can estimate $I_2^s$ as follows:
\begin{align*}
\big| \h{0.5pt} I^s_2 \h{0.5pt} \big| \h{2pt}& \lesssim \h{2pt}a_n R\int_{ B_{\sigma_k} }  \Big| v_{n, s, R} - \mathcal{W}_n \Big|  \int_0^1  \left|\h{1pt} \Big|\h{1pt} \mathcal{W}_n + t \big( v_{n, s, R} - \mathcal{W}_n \big) \Big|^2 - 1 \right| \h{2pt} \mathrm{d} \h{1pt}  t\\[2mm]
&  \lesssim \h{2pt} a_n R^2  \int_{ B_{\sigma_k} }  \Big| v_{n, s, R} - \mathcal{W}_n \Big|^2  + a_n R \left(\int_{ B_{\sigma_k} }  \Big| v_{n, s, R} - \mathcal{W}_n \Big|^2\right)^{1/2}\left(\int_{B_{\sigma_k}} \left|\h{1pt}  \big| \mathcal{W}_n  \big|^2 - 1 \h{1pt}\right|^2\right)^{1/2}.
\end{align*}
By (\ref{estimate of vnsr and yn, origin}) and (\ref{s_n and y_n,  origin}), the last estimate is reduced to \begin{align*}\big| \h{0.5pt} I^s_2 \h{0.5pt} \big| \h{2pt}& \lesssim \h{2pt}a_n \h{0.5pt}s_n^2 R^2 \left[ R^2 + \dfrac{s}{\sigma_k}\h{0.2pt}b_k\right] + \sqrt{a_n} \h{0.5pt}s^2_n \h{0.5pt}r_n^{-1} R  \left[ R^2 + \dfrac{s}{\sigma_k}\h{0.2pt}b_k\right]^{1/2}.
\end{align*}Therefore, in the small--scale regime, we obtain
\begin{equation}\label{est of I2s, origin}  \left( \dfrac{r_n}{s_n} \right)^2 \big|\h{0.2pt} I_2^s \h{0.2pt}\big| \h{2pt} \lesssim \h{2pt} a_n \h{0.2pt}r_n^2 R^2 \left[ R^2 + \dfrac{s}{\sigma_k}\h{0.2pt}b_k\right] + \sqrt{a_n} \h{0.4pt} r_n  R  \left[ R^2 + \dfrac{s}{\sigma_k}\h{0.2pt}b_k\right]^{1/2}
\longrightarrow 0 \quad \text{as\,\,} n \to \infty.\end{equation}
\noindent \textbf{Step 5. Completion of proof} \vspace{0.2pc}\\
Finally, we apply (\ref{conv. of grad v_n, origin}), (\ref{I_1 to I_2, origin}), (\ref{I_1^s, origin}) and (\ref{est of I2s, origin})   to the right--hand side of (\ref{small scale basic inequality, origin}). By the lower semi--continuity, it follows
\begin{align*}
\int_{ B_{\sigma_k} } \Big| \nabla \mathcal{W}^{sc}_\infty \Big|^2
&\h{2pt}\leq\h{2pt} \liminf_{ n \to \infty } \int_{ B_{\sigma_k} } \Big| \nabla \mathcal{W}^{sc}_n \Big|^2 \h{2pt}\leq\h{2pt}
 \limsup_{ n \to \infty } \int_{ B_{\sigma_k} } \Big| \nabla \mathcal{W}^{sc}_n \Big|^2 \\[2mm]
&\h{2pt}\leq\h{2pt}  \lim_{s \to 0} \lim_{R \to \infty} \lim_{n \to \infty}  \int_{ B_{\sigma_k} }  s_n^{-2} \h{0.2pt}\big| \nabla  v_{n, s, R} \big|^2 + \left( \dfrac{r_n}{s_n} \right)^2 F_n ( \cdot ) \h{1pt}\Big|_{ \mathcal{W}_n }^{v_{n, s, R}}
= \int_{ B_{\sigma_k} } \big|\nabla w \big|^2.
\end{align*}The proof is completed.
\end{proof}
\noindent\textbf{Proof of Proposition \ref{small energy implies energy decay} in small--scale regime.}\vspace{0.2pc}\\
We will take $n \rightarrow \infty$ in (\ref{limit of contradictory inequality, origin}). Firstly, by the uniform boundedness of $\mathcal{W}_n$ and $y_n$, we have \begin{align}\label{splitting y_n2 - 1, origin}\int_{B_1} \left( \big| y_n \big|^2 - 1 \right)^2 &= \int_{B_1} \left( \big| y_n - \mathcal{W}_n \big|^2 + 2\h{1pt} \mathcal{W}_n \cdot \big(y_n - \mathcal{W}_n\big) + \big| \mathcal{W}_n\big|^2 - 1 \right)^2 \nonumber\\[2mm]
&\lesssim \int_{B_1} \big| y_n - \mathcal{W}_n\big|^2 + \left|\h{1pt} \big| \mathcal{W}_n\big|^2 - 1\h{1pt}\right|^2.
\end{align}By using notations in (\ref{s_n and y_n,  origin}) and the uniform upper--bound of $\mathcal{W}_n^{sc}$ in $L^2\left(B_1\right)$, this estimate infers \begin{align*}\int_{B_1} \left( \big| y_n \big|^2 - 1 \right)^2 \h{2pt}\lesssim\h{2pt}s_n^2 + a_n^{-1} + a_n^{-1} \left( \dfrac{s_n}{r_n}\right)^2,
\end{align*}which furthermore gives us \begin{align*} \left(\dfrac{r_n}{s_n}\right)^2 a_n \left( \big| y_n \big|^2 - 1 \right)^2 \h{2pt}\lesssim\h{2pt}a_n r_n^2 + \left(\dfrac{r_n}{s_n}\right)^2 + 1.
\end{align*}Recalling Lemma \ref{s_n,r_n/s_n tend to 0,w}, in the small--scale regime, we can take $n \to \infty$ in the last estimate and get \begin{align}\label{lim of yn2 - 1, small, origin}\limsup_{n \rightarrow\infty}\left(\dfrac{r_n}{s_n}\right)^2 a_n \left( \big| y_n \big|^2 - 1 \right)^2 \h{2pt}\lesssim\h{2pt} 1.
\end{align}

Due to the strong $H^1_{\mathrm{loc}}$--convergence of $\mathcal{W}_n^{sc}$ shown in Lemma \ref{strong H1 convergence, small, origin} and the uniform boundedness of $\mathcal{W}_n$,  we can take $n \rightarrow \infty$ in (\ref{limit of contradictory inequality, origin}) and obtain \begin{align}\label{contra, taking limit, small, origin}\nu_0^{-1}\int_{B_{\nu_0}}\big| \nabla \mathcal{W}_\infty^{sc}\big|^2 + \dfrac{ \mu}{2\h{1pt} \nu_0} \liminf_{n \rightarrow \infty} a_n \left(\dfrac{r_n}{s_n}\right)^2 \int_{B_{\nu_0}} \left( \big| \mathcal{W}_n\big|^2 - 1 \right)^2 \geq \dfrac{1}{2}.
\end{align}Switching $\mathcal{W}_n$ and $y_n$ in (\ref{splitting y_n2 - 1, origin}) and changing the integration domain from $B_1$ to $B_{\nu_0}$, we have \begin{align*}\int_{B_{\nu_0}} \left( \big| \mathcal{W}_n \big|^2 - 1 \right)^2 \h{2pt}\lesssim\h{2pt} \int_{B_{\nu_0}} \big| \mathcal{W}_n - y_n \big|^2 + \left( \big| y_n \big|^2 - 1 \right)^2 \h{2pt}\lesssim\h{2pt} s_n^2 + \nu_0^3 \left( \big| y_n \big|^2 - 1 \right)^2.
\end{align*} By applying this estimate to the second term on the left--hand side of (\ref{contra, taking limit, small, origin}) and recalling (\ref{lim of yn2 - 1, small, origin}), in the small--scale regime, it holds  
\begin{align}\label{contra to energy decay of minimizer, small, origin}\dfrac{1}{2} \h{2pt}\leq\h{2pt} K \nu_0^2+ \nu_0^{-1}\int_{B_{\nu_0}} \Big| \nabla \mathcal{W}^{sc}_\infty \Big|^2.\end{align} Here and in what follow, $K$ is a universal constant depending only on $\mu$. Recall the minimization problem satisfied by $\mathcal{W}_\infty^{sc}$ in Lemma \ref{strong H1 convergence, small, origin}. If (\ref{infty case, origin}) holds, then all the components of $\mathcal{W}^{sc}_\infty$ are harmonic in $B_{1}$.  Standard elliptic estimate yields \begin{align}\label{unif bound of nabla lim scla, origin} \big\| \nabla \mathcal{W}^{sc}_\infty\big\|_{\infty; \h{0.3pt}B_{1/2}} \h{2pt}\lesssim\h{2pt} \big\| \mathcal{W}^{sc}_\infty\big\|_{2; \h{0.3pt}B_1} \h{2pt}\lesssim\h{2pt}1.
\end{align}Applying this estimate to the right--hand side of (\ref{contra to energy decay of minimizer, small, origin}), we then can take $\nu_0 \in \big(0, 1/2\big)$ suitably small so that (\ref{contra to energy decay of minimizer, small, origin}) fails. If (\ref{pi_L convergence, origin}) holds, then $\mathcal{W}^{sc}_{\infty; \h{0.2pt}j}$ are also harmonic in $B_1$, where $j = 1, 2, 4, 5$. Hence,  $\nabla \mathcal{W}^{sc}_{\infty; \h{0.2pt}j}$ are $L^\infty$--bounded on $B_{1/2}$ as well by a universal constant, for $j = 1, 2, 4, 5$. In the remaining of the proof, we consider $\mathcal{W}^{sc}_{\infty; \h{0.2pt}3}$. Firstly, we introduce a standard Signorini problem on $B_{\sigma_k}^+$ as follows:  \begin{align}\label{standard signorini, origin}
\mathrm{Min} \left\{ \int_{B^+_{\sigma_k}} \big| \nabla u \big|^2 : u \in \mathfrak{R}\right\}, \h{5pt}\text{where $\mathfrak{R} := \left\{ u \in W^{1,2}\left(B_{\sigma_k}^+\right): \h{4pt}u \h{2pt}\Big|_{T_{\sigma_k}} \geq w_*, \h{4pt}u \h{2pt}\Big|_{\left(\p B_{\sigma_k}\right)^+} = \mathcal{W}^{sc}_{\infty; \h{0.2pt}3} \right\}$.}\end{align}As a convention, $T_{\sigma_k}$ and $\left(\p B_{\sigma_k}\right)^+$  in (\ref{standard signorini, origin}) are flat and spherical boundaries of $B_{\sigma_k}^+$, respectively. We claim that $\mathcal{W}^{sc}_{\infty, 3}$ saturates the minimal energy in the minimization problem (\ref{standard signorini, origin}). In fact, for  any $u \in \mathfrak{R}$, the function \begin{align*} u_\sharp := \dashint_0^{2 \pi} u\left(\rho, z, \theta\right)\mathrm{d} \theta
\end{align*} lies in $\mathfrak{R}$ as well. Moreover, the Dirichlet energy of $u_\sharp$ is bounded from above by the Dirichlet energy of $u$. Using this symmetrization and the fact that minimizer to the problem (\ref{standard signorini, origin}) is unique, we know that the minimizer to (\ref{standard signorini, origin}) must be axially symmetric. Let $u_{\star}$ be the unique minimizer to (\ref{standard signorini, origin}) and extend $u_\star$ to $B_{\sigma_k}$ so that the extension, still denoted by $u_\star$, is even with respect to the $z$--variable. The minimizing property of $\mathcal{W}^{sc}_\infty$ shown in Lemma \ref{strong H1 convergence, small, origin}   infers that \begin{align*} \int_{B_{\sigma_k}} \big| \nabla \mathcal{W}^{sc}_{\infty; 3} \big|^2 \leq \int_{B_{\sigma_k}} \big| \nabla u_\star\big|^2 \h{10pt}\iff\h{10pt}\int_{B^+_{\sigma_k}} \big| \nabla \mathcal{W}^{sc}_{\infty; 3} \big|^2 \leq \int_{B^+_{\sigma_k}} \big| \nabla u_\star\big|^2.
\end{align*}Since $\mathcal{W}^{sc}_{\infty; 3} \in \mathfrak{R}$, the reverse direction of the above inequality holds trivially. Hence $\mathcal{W}^{sc}_{\infty; 3}$ is the unique minimizer of (\ref{standard signorini, origin}). Applying Lemma 9.1 in \cite{PSU12}, we obtain \begin{align}\label{cons of book 1, origin} \rho^{-1} \int_{B_{\rho/2}} \Big| \nabla^2 \mathcal{W}^{sc}_{\infty; 3} \Big|^2 \h{2pt}\lesssim\h{2pt}\rho^{-3} \int_{B_\rho \setminus B_{\rho/2}} \big| \nabla \mathcal{W}^{sc}_{\infty; 3} \big|^2, \h{15pt}\text{for any $ \rho \in \big(0, 1/4\big)$.}
\end{align}On the other hand, Poincar\'{e}'s inequality induces that \begin{align}\label{poincare, small, origin}\int_{B_\rho \setminus B_{\rho/2}} \left|\h{1pt} \nabla \mathcal{W}^{sc}_{\infty; 3} \h{1pt} \right|^2 = \int_{B_\rho \setminus B_{\rho/2}} \left|\h{1pt} \nabla \mathcal{W}^{sc}_{\infty; 3} - \dashint_{B_\rho \setminus B_{\rho/2}} \nabla \mathcal{W}^{sc}_{\infty; 3} \h{1pt}\right|^2 \h{2pt}\lesssim\h{2pt}\rho^2 \int_{B_\rho \setminus B_{\rho/2}} \Big| \nabla^2 \mathcal{W}^{sc}_{\infty; 3} \h{1pt}\Big|^2.
\end{align}The first equality above holds due to the axial symmetry of $\mathcal{W}^{sc}_{\infty; 3}$ and the even symmetry of $\mathcal{W}^{sc}_{\infty; 3}$ with respect to the $z$--axis. Combining the estimates in (\ref{cons of book 1, origin})--(\ref{poincare, small, origin}) and utilizing the filling hole argument, we then get \begin{align*}
\int_{B_{\rho/2}} \Big| \nabla^2 \mathcal{W}^{sc}_{\infty; 3} \Big|^2 \h{2pt}\leq\h{2pt}\theta \int_{B_\rho} \Big| \nabla^2 \mathcal{W}^{sc}_{\infty; 3} \Big|^2, \h{15pt}\text{for any $\rho \in \big(0, 1/4\big)$.}
\end{align*}Here $\theta \in \big(0, 1\big)$ is a universal constant. By the above  estimate, standard iteration argument yields \begin{align*} \int_{B_\rho} \Big| \nabla^2 \mathcal{W}^{sc}_{\infty; 3} \Big|^2 \h{2pt}\lesssim\h{2pt} \rho^{\alpha}\int_{B_{1/16}} \Big| \nabla^2 \mathcal{W}^{sc}_{\infty; 3} \Big|^2, \h{10pt}\text{for some  $\alpha \in \big(0, 1\big)$ depending only on $\theta$ and any $\rho \in \big(0, 1/16\big)$.}
\end{align*}Taking $\rho = 1/8$ in (\ref{cons of book 1, origin}) and using the upper--bound of the $L^2$--norm of $\nabla \mathcal{W}^{sc}_{\infty; 3}$ on $B_1$, we then obtain from the above estimate that 
 \begin{align*} \int_{B_\rho} \Big| \nabla^2 \mathcal{W}^{sc}_{\infty; 3} \Big|^2 \h{2pt}\lesssim\h{2pt} \rho^{\alpha}, \h{10pt}\text{for some $\alpha \in \big(0, 1\big)$ depending only on $\theta$ and any $\rho \in \big(0, 1/16\big)$,}
\end{align*}
which furthermore infers \begin{align*}\int_{B_\rho } \left| \nabla \mathcal{W}^{sc}_{\infty; 3}  \right|^2 = \int_{B_\rho} \left| \nabla \mathcal{W}^{sc}_{\infty; 3} - \dashint_{B_\rho } \nabla \mathcal{W}^{sc}_{\infty; 3} \right|^2 \h{2pt}\lesssim\h{2pt}\rho^2  \int_{B_\rho} \Big| \nabla^2 \mathcal{W}^{sc}_{\infty; 3} \Big|^2 \h{2pt}\lesssim\h{2pt}\rho^{2 + \alpha}, \h{10pt}\text{for any $\rho \in \big(0, 1/16\big)$.}
\end{align*}In light of this estimate and the harmonicity of the remaining components in $\mathcal{W}^{sc}_\infty$, we can still take $\nu_0$ suitably small so that the estimate (\ref{contra to energy decay of minimizer, small, origin}) fails in the case when (\ref{pi_L convergence, origin}) holds. The proof is then finished.
\subsection{Energy--decay estimate in intermediate--scale regime}
In this section we suppose that $a_n r_n^2 \longrightarrow L$ as $n \to \infty$. Here $L \in \big(0,\infty\big)$ is constant.  
\begin{lem}\label{properties of y_n in intermediate scale regime}In the intermediate--scale regime, we can keep extracting a subsequence of $\big\{y_n\big\}$, still denoted by $\big\{y_n\big\}$, so that for some constant $c_1 \in \mathbb{R}$ depending on $L$, the following limit holds:
\begin{equation}
\dfrac{ \big| y_n \big| - 1 }{s_n} \longrightarrow c_1 \quad \text{as\,\,}n \to \infty.
\label{|y_n|-1/s_n to c_6, origin}
\end{equation}
Due to the above limit and (\ref{y_*=(0,0,y_*,3,0,0)}), 
\begin{equation}
y_n \longrightarrow y_* = e_3 \h{15pt}\text{as $n \rightarrow \infty$}.
\label{y_*=0,0,1,0,0, origin}
\end{equation}Therefore, in the intermediate--scale regime, the limit in (\ref{infty case, origin}) must hold. Moreover, except item (4), we still have all the first three items in Lemma \ref{uniform convergence of u^L_n in S_sigma_k} in the intermediate--scale regime.
\end{lem}
\begin{proof}[\bf Proof]
In light of the notations in (\ref{s_n and y_n,  origin}) and Lemma \ref{s_n,r_n/s_n tend to 0,w},
\begin{equation}
    \int_{B_1} \big| \nabla \mathcal{W}_n \big|^2+  a_n \h{0.3pt}r_n^2 \left( \big| \mathcal{W}_n \big|^2-1\right)^2 \h{2pt}\lesssim\h{2pt} s_n^2, \qquad \text{for $n$ suitably large.}
    \label{E_n < s_n^2, origin}
\end{equation} On the other hand, 
$$\dfrac{ \big| y_n \big|^2 - 1 }{s_n}
= \dfrac{ \big| \mathcal{W}_n - s_n \h{1pt} \mathcal{W}^{sc}_n \big|^2-1}{s_n}
= \dfrac{\big| \mathcal{W}_n \big|^2-1 - 2\h{0.2pt}s_n \h{1pt} \mathcal{W}_n \cdot \mathcal{W}^{sc}_n + s_n^2 \big| \mathcal{W}^{sc}_n \big|^2}{s_n}.$$
The lemma then follows by this decomposition, (\ref{E_n < s_n^2, origin}), $L^\infty$--boundedness  of $\mathcal{W}_n$ and the uniform $L^2$--boundedness  of $\mathcal{W}^{sc}_n$ on $B_1$. Here we also used the non--zero assumption on the limit of  $a_n r_n^2$  as $n \rightarrow \infty$.
\end{proof}In the next, we study the minimization problem satisfied by $\mathcal{W}^{sc}_\infty$ in the intermediate--scale regime.
\begin{lem}Recall $M_k$ defined in (\ref{configuration space M_k, origin}).  For any $k \in \mathbb{N}$, the mapping $\mathcal{W}^{sc}_\infty$ minimizes the $\mathcal{E}^{sc}_L$--energy on the space $M_k$. Here with the limit $c_1$ obtained in (\ref{|y_n|-1/s_n to c_6, origin}),
$$\mathcal{E}^{sc}_L[\h{0.5pt}w\h{0.5pt}]:=\int_{B_{\sigma_k}} \big| \nabla  w \big|^2 + 2 \h{0.5pt} L \h{0.5pt} \mu \left( w_3 + c_1 \right)^2, \qquad \text{for all\,\,} w \in M_k.$$
Moreover, $\mathcal{W}^{sc}_n$ converges to $\mathcal{W}^{sc}_\infty$ strongly in $H^1_{\textup{loc}}\left(B_1\right)$ as $n \rightarrow \infty$.
\label{strong H1 convergence in intermediate case, origin}
\end{lem}
\begin{proof}[\bf Proof] We use the same comparison map introduced in the proof of Lemma \ref{strong H1 convergence, small, origin}. Recalling (\ref{conv. of grad v_n, origin}), (\ref{estimate of vnsr and yn, origin}) and the definition of $v_{n, s, R}$ in (\ref{definition of v_n, origin}), we have, up to a subsequence, the following convergence: \begin{align}\label{L4 conv of vnsR, origin} \dfrac{v_{n, s, R} - y_n}{s_n} \longrightarrow w, \h{15pt}\text{strongly in $L^4\left(B_{\sigma_k}\right)$, as $n \rightarrow \infty$, $R \to \infty$ and $s \to 0$, successively.}\end{align}By Sobolev embedding, it satisfies \begin{align}\label{L4 conv of W-nsc, origin}\mathcal{W}^{sc}_n \longrightarrow \mathcal{W}^{sc}_\infty, \h{15pt}\text{strongly in $L^4\left(B_{1}\right)$, as $n \rightarrow \infty$.}\end{align} With the limits (\ref{|y_n|-1/s_n to c_6, origin})--(\ref{y_*=0,0,1,0,0, origin}),  (\ref{L4 conv of vnsR, origin})--(\ref{L4 conv of W-nsc, origin}) and the assumption that $a_n r_n^2$ converges to $L$ as $n \rightarrow \infty$, the following two convergences hold:
\begin{align*}
\left( \dfrac{r_n}{s_n} \right)^2 a_n \int_{B_{\sigma_k}} \left( \big| v_{n, s, R} \big|^2 - 1 \right)^2 &=  a_n r_n^2 \int_{B_{\sigma_k}} \left( s_n \left| \dfrac{v_{n, s, R} - y_n}{s_n}\right|^2 + 2 \h{1pt}y_n \cdot \left( \dfrac{v_{n, s, R} - y_n }{s_n}\right) +  \dfrac{\big|y_n \big|^2 - 1}{s_n} \right)^2 \\[2mm]
 & \longrightarrow 4\h{0.5pt}L  \int_{B_{\sigma_k}} \big( w_3 + c_1 \big)^2, \h{15pt} \text{as $n \rightarrow \infty$, $R \rightarrow \infty$ and $s \rightarrow 0$, successively}
\end{align*}and
\begin{align}\label{conv of | wn |^2 - 1, origin}\left( \dfrac{r_n}{s_n} \right)^2 a_n \int_{B_r} \left( \big| \mathcal{W}_n \big|^2 - 1 \right)^2 &= \left( \dfrac{r_n}{s_n} \right)^2 a_n \int_{B_r} \left( s_n^2 \h{1pt}\big| \mathcal{W}_n^{sc}\big|^2 + 2 \h{1pt}s_n\h{1pt}y_n \cdot \mathcal{W}^{sc}_n + \big|y_n \big|^2 - 1 \right)^2 \nonumber\\[2mm]
 & \longrightarrow 4\h{0.5pt}L  \int_{B_r} \big( \mathcal{W}^{sc}_{\infty; 3} + c_1 \big)^2, \h{20pt} \text{as $n \rightarrow \infty$, for any $r \in [\h{0.2pt}0, 1 \h{0.2pt}]$.}
\end{align}Owing to the last two limits, in the intermediate--scale regime, \begin{align*} \left(\dfrac{r_n}{s_n}\right)^2 I_2^s \longrightarrow 2 \h{1pt}L\h{1pt}\mu   \int_{B_{\sigma_k}} \big( w_3 + c_1 \big)^2 -   \big( \mathcal{W}^{sc}_{\infty; 3} + c_1 \big)^2, \h{5pt}\text{as $n \to \infty$, $R \to \infty$ and $s \to 0$, successively.}
\end{align*}Here $I_2^s$ is given in (\ref{def of I2s}). By applying the above limit,  (\ref{conv. of grad v_n, origin}), (\ref{I_1 to I_2, origin}) and (\ref{I_1^s, origin})  to the right--hand side of (\ref{small scale basic inequality, origin}), it follows
\begin{align*}
\int_{ B_{\sigma_k} } \Big| \nabla \mathcal{W}^{sc}_\infty \Big|^2 + 2 \h{0.5pt}L\h{0.5pt}\mu \big( \mathcal{W}^{sc}_{\infty; 3} + c_1 \big)^2
&\h{2pt}\leq\h{2pt}   \int_{ B_{\sigma_k} } \big|\nabla w \big|^2 + 2 \h{0.5pt}L\h{0.5pt}\mu \big( w_3 + c_1 \big)^2, \h{5pt}\text{for any $w \in M_k$.}
\end{align*} The proof is then finished.
\end{proof}
\begin{proof}[\bf Proof of Proposition \ref{small energy implies energy decay}  in the intermediate--scale regime.]\
\\[3mm]
Recall (\ref{E_n < s_n^2, origin}).  It holds \begin{align*} \int_{B_1} \big| \nabla \mathcal{W}^{sc}_n\big|^2 + \left(\dfrac{r_n}{s_n}\right)^2 a_n \left( \big| \mathcal{W}_n\big|^2 -1\right)^2 \h{2pt}\lesssim\h{2pt}1.
\end{align*}Noticing the convergence in (\ref{conv of | wn |^2 - 1, origin}), we can take $n \rightarrow \infty$ in the above estimate and obtain \begin{align}\label{bound, small, origin, L problem} \int_{B_1} \big| \nabla \mathcal{W}^{sc}_\infty\big|^2 + L \left(\mathcal{W}^{sc}_{\infty; 3} + c_1 \right)^2 \h{2pt}\lesssim\h{2pt}1.
\end{align}The components $\mathcal{W}^{sc}_{\infty; j}$ ($j = 1, 2, 4, 5$) satisfy the harmonic equation on $B_1$. Therefore, by (\ref{bound, small, origin, L problem}), \begin{align}\label{Linfty bound of nabla w sc infty j, origin} \big\| \nabla \mathcal{W}^{sc}_{\infty; j}\big\|_{\infty; \h{0.3pt}B_{1/2}} \h{2pt}\lesssim\h{2pt} \big\| \mathcal{W}^{sc}_{\infty; j}\big\|_{2; \h{0.3pt}B_1} \h{2pt}\lesssim\h{2pt} \big\| \nabla \mathcal{W}^{sc}_{\infty; j}\big\|_{2; \h{0.3pt}B_1} \h{2pt}\lesssim\h{2pt}1, \h{15pt}j = 1, 2, 4, 5.
\end{align}The second estimate above uses Neumann--Poincar\'{e} inequality. As for $\mathcal{W}^{sc}_{\infty; 3}$, by Lemma \ref{strong H1 convergence in intermediate case, origin}, it satisfies the equation:\begin{align*} \Delta \mathcal{W}^{sc}_{\infty; 3} = 2 \h{0.2pt}L\h{0.2pt}\mu \left( \mathcal{W}^{sc}_{\infty; 3} + c_1 \right) \h{15pt}\text{in $B_1$.}
\end{align*}Hence, the function $\big| \nabla \mathcal{W}^{sc}_{\infty; 3} \big|^2 + L\h{0.2pt}\mu\left(\mathcal{W}^{sc}_{\infty; 3} + c_1 \right)^2$ is subharmonic on $B_1$. Using (\ref{bound, small, origin, L problem}) and the standard local boundedness result for the subharmonic functions (see Theorem 4.1 in \cite{HL11}), we have  \begin{align}\label{est Wsd inty 3, origin}  \big| \nabla \mathcal{W}^{sc}_{\infty; 3} \big|^2 + L\h{0.2pt}\mu\left(\mathcal{W}^{sc}_{\infty; 3} + c_1 \right)^2  \h{2pt}\lesssim\h{2pt} \int_{B_1}\big| \nabla \mathcal{W}^{sc}_{\infty; 3} \big|^2 + L\h{0.2pt}\mu\left(\mathcal{W}^{sc}_{\infty; 3} + c_1 \right)^2 \h{2pt}\lesssim\h{2pt}1 \h{15pt}\text{on $B_{1/2}$.}
\end{align}Utilizing the strong $H^1$--convergence of $\mathcal{W}^{sc}_n$ in Lemma \ref{strong H1 convergence in intermediate case, origin} and (\ref{conv of | wn |^2 - 1, origin}), we take $n \rightarrow \infty$ in (\ref{limit of contradictory inequality, origin}) and get \begin{align*} \nu_0^{-1}\int_{B_{\nu_0}} \big| \nabla \mathcal{W}^{sc}_{\infty}\big|^2 + 2\h{0.2pt}L\h{0.2pt}\mu \left( \mathcal{W}^{sc}_{\infty; 3} + c_1 \right)^2 \h{2pt}\geq \h{2pt} \dfrac{1}{2}.
\end{align*}However, in light of (\ref{Linfty bound of nabla w sc infty j, origin})--(\ref{est Wsd inty 3, origin}), the above estimate fails for $\nu_0$ suitably small.  The proof is completed.
\end{proof}
\subsection{Energy--decay estimate in large--scale regime}
In this section we suppose that  $a_n\h{0.5pt} r_n^2 \longrightarrow \infty$ as $n \to \infty$. \begin{lem}\label{properties of sequence in large scale, origin} In the large--scale regime, the followings hold up to a subsequence:
\begin{enumerate}
\item[$\mathrm{(1).}$] In light of (\ref{E_n < s_n^2, origin}), we have \begin{eqnarray}\label{conv of scaled |w_n|^2 - 1, origin, large scale}
s_n^{-2}  \int_{B_{1}}  \left(  \big| \mathcal{W}_n \big|^2 - 1 \right)^2 \h{3pt}\lesssim\h{3pt} \left(a_n\h{0.5pt} r_n^2\h{1pt}\right)^{-1} \longrightarrow 0, \h{20pt}\text{as $n \rightarrow \infty$.}
\end{eqnarray}Moreover, the limit (\ref{y_*=0,0,1,0,0, origin}) is still satisfied by the sequence $\big\{y_n\big\}$. Therefore, in the large--scale regime, the limit in (\ref{infty case, origin}) must hold;
\item[$\mathrm{(2).}$] By (\ref{conv of scaled |w_n|^2 - 1, origin, large scale}) and similar arguments as in the proof of Lemma \ref{properties of y_n in intermediate scale regime}, there exists a $c_2 \in \mathbb{R}$ so that \begin{equation}
\dfrac{ \big| y_n \big| - 1 }{s_n} \longrightarrow c_2 \quad \text{as\,\,}n \to \infty.
\label{|y_n|-1/s_n to c_6, large scale, origin}\end{equation}Moreover, $\mathcal{W}^{sc}_{\infty; 3} + c_2 \equiv 0$ on $B_1$;
\item[$\mathrm{(3).}$] Except item (4), we have all the first three items in Lemma \ref{uniform convergence of u^L_n in S_sigma_k} in the large--scale regime. By Fatou's lemma, we can in addition have \begin{align}\label{large scale, |w_n|^ - 1, sphere, origin} \sup_{n \h{0.5pt}\in\h{0.5pt}\mathbb{N}} \h{2pt} \left( \dfrac{ r_n }{ s_n } \right)^2 a_n  \int_{\p B_{\sigma_k}} \left(  \h{1pt} \big| \h{1pt} \mathcal{W}_n \h{1pt} \big|^2 - 1 \right)^2 \h{2pt} \leq \h{2pt} b_k, \h{15pt}\text{for any $k \in \mathbb{N}$}. \end{align}
\end{enumerate}
\end{lem}
Before we prove the strong convergence of $\mathcal{W}^{sc}_n$, we need the following lemma:
\begin{lem}
Recall the sequence $\{\sigma_k\}$ obtained in Lemma \ref{uniform convergence of u^L_n in S_sigma_k}. For any $k \in \mathbb{N}$, we have
$$ \mathcal{W}_n \longrightarrow e_3 \h{10pt} \text{ in $C^0\big(\p B_{\sigma_k}\big)$, as $n \to \infty$.}  $$
This convergence is obtained up to a subsequence.
\label{uniform conv. of |W_n| on bdy of B_sigma}
\end{lem}
\begin{proof}[\bf Proof] Most part of the proof has been contained in \cite{MZ10} already. We just sketch the ideas and point out the minor differences between our case and \cite{MZ10}. 
Note that the $s_n$ defined in (\ref{s_n and y_n,  origin}) converges to $0$ as $n \rightarrow \infty$. Therefore, on $B_{\frac{1- \sigma_k}{8}r_n}(r_nq_k)$, where $q_k$ is the north pole of the ball $B_{\sigma_k}$, we can follow  the same arguments used in the proof of Proposition 4 in \cite{MZ10}. It then turns out  $\big| w_n \big| \geq 1/2$ on $B_{ \frac{1- \sigma_k}{8}r_n}(r_nq_k)$ for $n$ suitably large depending on $k$. Here we take $\sigma_k$ close to $1$. With this lower bound, we can apply Lemma \ref{energy density e^L_a[w_a] bdd} to get \begin{align}\label{uni bound gradien}  \left(\dfrac{1- \sigma_k}{8}\right)^2 r_n^2 \sup_{B_{\frac{1 - \sigma_k}{16}r_n}(r_nq_k)} f_{a_n, \mu}\big(w_n\big) \lesssim 1. \end{align}Here we also have used the convergence: $$r_n^{-1} \int_{B_{\frac{1 - \sigma_k}{8}r_n}(r_n q_k)} f_{a_n, \mu}\big(w_n\big) \longrightarrow 0, \h{15pt}\text{as $n \rightarrow \infty$.}$$It is a consequence of (\ref{E_n < s_n^2, origin}). In light of the estimate in (\ref{uni bound gradien}), we then obtain the following uniform boundedness of the gradient of $\mathcal{W}_n$: $$ \sup_{B_{\frac{1 - \sigma_k}{16}}(q_k)} \big| \nabla \mathcal{W}_n \big| \h{2pt}\lesssim \h{2pt}\dfrac{8}{1 - \sigma_k} . $$
The above inequality shows that $\mathcal{W}_n$ is equi--continuous on the closure of $ B_{\frac{1 - \sigma_k}{16}}(q_k)$. By the $\mathscr{R}$--axial symmetry, $\mathcal{W}_n$ is also equi--continuous near $- q_k$.\vspace{0.2pc}

Now we discuss the points on $\p B_{\sigma_k}$ away from $\pm q_k$. Let $\Phi_0$ be the polar angle of the points on $ \p B_{\sigma_k}  \cap  \p B_{\frac{1-\sigma_k}{16}}(q_k) $ and suppose that $\mathcal{W}_n = \mathscr{L}\left[\h{0.5pt}v_n\h{0.5pt}\right]$ for some $3$--vector field $v_n = v_n\left(\rho, z\right)$. For any polar angles $\Phi_1, \Phi_2$ satisfying $\Phi_0 \leq \Phi_1 < \Phi_2 \leq \pi-\Phi_0$, it holds
\begin{eqnarray*}
\big| \h{1pt} v_n(\sigma_k,\Phi_2) -  v_n(\sigma_k,\Phi_1) \h{1pt} \big|
&=& \left|\h{1pt} \int_{\Phi_1}^{\Phi_2} \p_\Phi \h{1pt} v_n(\sigma_k,\Phi) \h{1pt}\mathrm{d} \h{1pt} \Phi \h{1pt} \right|  \h{2pt}\leq\h{2pt} \big|\h{1pt} \Phi_2 - \Phi_1 \h{1pt}\big|^{1/2} \left( \int_{\Phi_1}^{\Phi_2} \big| \p_\Phi \h{1pt} v_n(\sigma_k,\Phi) \big|^2 \h{1pt}\mathrm{d} \h{1pt} \Phi \right)^{1/2} \\[2mm]
&\leq& \left(\dfrac{\big|\h{1pt} \Phi_2 - \Phi_1 \h{1pt}\big|}{\sin \Phi_0}\right)^{1/2}\left( \int_0^{2\pi}\int_{\Phi_1}^{\Phi_2} \big| \p_\Phi \h{1pt} v_n(\sigma_k,\Phi) \big|^2 \h{2pt}\sigma_k^2 \h{2pt}\sin\Phi \h{2pt}\mathrm{d} \h{1pt} \Phi \h{1pt}\mathrm{d} \h{1pt} \Theta \right)^{1/2}.
\end{eqnarray*}
Applying the fact that $\displaystyle\int_{\p B_{\sigma_k}} \big| \p_{\Phi} \h{1pt} v_n \big|^2 \leq \displaystyle 2 \int_{\p B_{\sigma_k}} \big| \nabla \h{1pt} \mathcal{W}_n \big|^2$ and (1) in Lemma \ref{uniform convergence of u^L_n in S_sigma_k}, we have
$$
\big| \h{1pt} v_n(\sigma_k,\Phi_2) -  v_n(\sigma_k,\Phi_1) \h{1pt} \big|
\lesssim s_n \left(\dfrac{b_k}{ \sin\Phi_0}\right)^{1/2}\big|\h{1pt} \Phi_1 - \Phi_2 \h{1pt}\big|^{1/2}.
$$

In light of the above arguments and the relationship between $v_n$ and $\mathcal{W}_n$, we know that $\mathcal{W}_n$ is equi-continuous on  $\p B_{\sigma_k}$. Since $\mathcal{W}_n$ is uniform bounded in $B_1$, we conclude by Arzel\`{a}--Ascoli theorem that up to a subsequence, $\mathcal{W}_n$ converges in $C^0\left(\p B_{\sigma_k}\right)$ as $n \rightarrow \infty$. The lemma then follows by (2) of Lemma \ref{uniform convergence of u^L_n in S_sigma_k} and (\ref{y_*=0,0,1,0,0, origin}).
\end{proof}
With the aid of this lemma, the mapping $\mathcal{W}_{\infty}^{sc}$ satisfies
\begin{lem}
For any natural number $k$,  the mapping $\mathcal{W}^{sc}_\infty$ minimizes the Dirichlet energy over $ M_k $. Moreover, $\mathcal{W}^{sc}_n$ converges to $\mathcal{W}^{sc}_\infty$ strongly in $ H_{ \mathrm{loc} }^1 \big(B_1 \big) $ as $n \rightarrow \infty$. In the large--scale regime, it also holds
$$ a_n \h{1pt} \left( \dfrac{ r_n } { s_n } \right)^2 \int_{ B_{\sigma_k}} \Big( \h{1pt} \big| \h{1pt} \mathcal{W}_n \h{1pt} \big|^2 - 1 \Big)^2 \longrightarrow 0, \h{15pt}\text{as $n\rightarrow \infty$.}$$
\label{strong H1 convergence, large, origin}
\end{lem}
\begin{proof}[\bf Proof] Suppose that $w$ is an arbitrary map in $ M_{ \h{0.5pt} k }$. Then we define the same $v_{n, s, R}$--mapping as in (\ref{definition of v_n, origin}). Note that here $M_{n, R}[ w ]$ is defined by the expression in (\ref{definition of F_n^s, origin}) for the (\ref{infty case, origin}) case. In terms of $v_{n, s, R}$, we define our comparison map as follows:  
\begin{eqnarray}
\widetilde{v}_{n,s,R} \big(\zeta\big) :=\left\{
\begin{aligned}
&    \Pi_{ \h{0.5pt} \mathbb{S}^4 } \left[  v_{n,s,R} \left( \dfrac{ \zeta }{ 1-s }\right) \right] \quad &\text{if }& \zeta \in B_{ ( 1 - s) \h{1pt} \sigma_k }; \\[4mm]
&\dfrac{ \sigma_k - | \h{1pt} \zeta \h{1pt} | } { s \h{0.5pt} \sigma_k } \h{2pt} \Pi_{ \h{0.5pt} \mathbb{S}^4 } \left[ \mathcal{W}_n\h{.5pt} \big(\h{1pt}\sigma_k \widehat{\zeta}\h{1pt}\big) \right]
+ \dfrac{ | \h{1pt} \zeta \h{1pt} | - ( 1 - s ) \h{1pt} \sigma_k } {s \h{0.5pt}\sigma_k}\h{1pt}\mathcal{W}_n\h{.5pt} \big(\h{1pt}\sigma_k \widehat{\zeta}\h{1pt}\big) \quad &\text{if }& \zeta \in B_{ \sigma_k } \setminus B_{ (1-s) \sigma_k }.
\end{aligned}
\right.
\label{definition of v_n in large scale, origin}
\end{eqnarray}By Lemma \ref{uniform conv. of |W_n| on bdy of B_sigma}, $\mathcal{W}_n$ converges uniformly to $e_3$ on $\p B_{\sigma_k}$ as $n \to \infty$. In light of this convergence and (\ref{y_*=0,0,1,0,0, origin}), for each $R$ fixed, we can take $n$ large enough such that \begin{equation}
\big|\h{0.5pt} \mathcal{W}_n \big| \geq 1/2 \h{5pt} \text{ on $\p B_{\sigma_k}$ } \h{10pt} \text{ and } \h{10pt}
\big|\h{0.5pt} v_{n,s,R} \big| \geq 1/2 \h{5pt} \text{ in $ B_{\sigma_k}$}.
\label{|W_n| and |h|>1/2}
\end{equation}Hence, the projections to $\mathbb{S}^4$ in the definition of   $\widetilde{v}_{n,s,R}$ are well--defined. Still by the convergences of   $M_{n, R}\big[\h{.5pt}\mathcal{W}^{sc}_\infty\h{.5pt}\big] $ and $\mathcal{W}_n$   on $\p B_{\sigma_k}$ and the limit in (\ref{y_*=0,0,1,0,0, origin}), when $n$ is  large,  $\big[ \widetilde{v}_{n, s, R}\big]_3$ satisfies the Signorini obstacle boundary condition: $
\big[ \widetilde{v}_{n,s,R} \big]_3 \h{2pt} \geq \h{2pt} H_{a_n} \h{0.5pt} b$  on $T_{\sigma_k}$. Due to the energy minimizing property of $\mathcal{W}_n$, it then turns out
\begin{equation}
\int_{B_{\sigma_k}} \big| \nabla  \mathcal{W}^{sc}_n \big|^2
\leq \int_{B_{\sigma_k}} \big| \nabla  \mathcal{W}^{sc}_n \big|^2 \h{1pt}
+ \left(\dfrac{r_n}{s_n}\right)^2F_n\left(\mathcal{W}_n\right)
\leq \h{1pt} \int_{B_{\sigma_k}}  s^{-2}_n  \big| \nabla  \widetilde{v}_{n, s, R} \big|^2 +\left(\dfrac{r_n}{s_n}\right)^2F_n\left(\widetilde{v}_{n, s, R}\right).
\label{large scale basic inequality, origin}
\end{equation}

The Dirichlet energy of $\widetilde{v}_{n,s,R}$ is computed as follows:
\begin{align}
\int_{B_{\sigma_k}}\Big|\h{.5pt} \nabla \widetilde{v}_{n,s,R} \h{.5pt}\Big|^2
&\h{2pt}=\h{2pt}(1-s)\int_{B_{\sigma_k}} \Big| \nabla  \h{1pt}\Pi_{ \h{0.5pt} \mathbb{S}^4 } \left[ v_{n,s,R}\h{.5pt} \big(\h{1pt}\zeta\h{1pt}\big) \right] \h{1pt}\Big|^2
+ \int_{B_{\sigma_k} \setminus B_{(1-s)\h{1pt}\sigma_k}} \big| \nabla \widetilde{v}_{n,s,R} \big|^2 \nonumber\\[2mm]
&\h{2pt}\leq\h{2pt} \int_{B_{ \sigma_k }} \dfrac{1-s}{\big| v_{n,s,R}\h{.5pt} \big|^{2}}
\Big| \nabla \h{.5pt}v_{n,s,R}\h{.5pt} \Big|^2 +\int_{B_{\sigma_k} \setminus B_{(1-s)\h{1pt}\sigma_k}} \big| \nabla \widetilde{v}_{n,s,R} \big|^2. \label{decomp. in shell in large scale, origin}
\end{align}Utilizing the uniform convergence of $v_{n, s, R}$ to $e_3$ on $B_{\sigma_k}$ and (\ref{conv. of grad v_n, origin}), we have
\begin{equation}
\lim_{s \to 0}\lim_{R \to \infty} \lim_{n \to \infty} \dfrac{1-s}{s_n^2}\int_{B_{ \sigma_k }} \big| v_{n,s,R}\h{.5pt} \big|^{-2}
\Big| \nabla \h{.5pt}v_{n,s,R}\h{.5pt} \Big|^2
= \int_{B_{\sigma_k}} |\nabla w|^2.
\label{limit of first term}
\end{equation}
To estimate the last integral in (\ref{decomp. in shell in large scale, origin}), we can follow exactly the same arguments as in Step 3 of the proof for Lemma \ref{strong H1 convergence, small, origin}. As a consequence, it holds \begin{equation}
\lim_{s \to 0}    \limsup_{n \to \infty}   s_n^{-2}\int_{B_{\sigma_k}\setminus B_{(1-s)\sigma_k}} \big|\nabla  \widetilde{v}_{n, s, R}\big|^2  =0.
    \label{conv. of grad v_n in annulus, origin}
\end{equation}
Note that to derive (\ref{conv. of grad v_n in annulus, origin}), we combine to use  the lower bound of $\mathcal{W}_n$ on $\p B_{\sigma_k}$ given in (\ref{|W_n| and |h|>1/2}), the limit $a_nr_n^2 \longrightarrow \infty$, (\ref{large scale, |w_n|^ - 1, sphere, origin})  and item (1) in Lemma \ref{uniform convergence of u^L_n in S_sigma_k}. Now we divide $s_n^2$ from both sides of (\ref{decomp. in shell in large scale, origin}). In light of (\ref{limit of first term})--(\ref{conv. of grad v_n in annulus, origin}), it then follows
\begin{equation}
  \limsup_{s \to 0} \h{1pt}\limsup_{R \to \infty} \h{1pt} \limsup_{n \to \infty}   s_n^{-2}\int_{B_{\sigma_k}} \big|\nabla  \widetilde{v}_{n, s, R}\big|^2  \leq \int_{B_{\sigma_k}}\big|\nabla w\big|^2.
    \label{conv. of grad v_n, large, origin}
\end{equation}

As for the estimate of the potential term, by the definition of $\widetilde{v}_{n,s,R}$ in (\ref{definition of v_n in large scale, origin}), we have
\begin{eqnarray*}
a_n \h{1pt} \left( \dfrac{ r_n } { s_n } \right)^2 \int_{ B_{ \sigma_k } } \Big( \h{1pt} \big| \h{1pt} \widetilde{v}_{n,s,R} \h{1pt} \big|^2 - 1 \Big)^2
&=& a_n \h{1pt} \left( \dfrac{ r_n } { s_n } \right)^2 \int_{ B_{ \sigma_k } \setminus B_{ (1 - s) \h{0.5pt} \sigma_k }} \Big( \h{1pt} \big| \h{1pt} \widetilde{v}_{n,s,R} \h{1pt} \big|^2 - 1 \Big)^2 \\[2mm]
& \lesssim & a_n \h{1pt} \left( \dfrac{  r_n } { s_n } \right)^2 \int_{ B_{ \sigma_k } \setminus B_{ (1 - s) \h{0.5pt} \sigma_k } } \left( \h{1pt} \left| \h{1pt} \mathcal{W}_n\h{.5pt} \big(\h{1pt}\sigma_k \widehat{\zeta}\h{1pt}\big) \h{1pt} \right|^2 - 1 \right)^2.
\end{eqnarray*}
It then turns out by (3) in Lemma \ref{properties of sequence in large scale, origin} that
$$
a_n \h{1pt} \left( \dfrac{ r_n } { s_n } \right)^2 \int_{ B_{ \sigma_k } } \Big( \h{1pt} \big| \h{1pt} \widetilde{v}_{n,s,R} \h{1pt} \big|^2 - 1 \Big)^2
\lesssim \h{3pt} s \h{1pt} \sigma_k \h{1pt} b_k \longrightarrow 0 \quad \text{ as $s \to 0$}.
$$
The uniform boundedness of $\widetilde{v}_{n,s,R}$, the limit of $r_n/s_n$ in Lemma \ref{s_n,r_n/s_n tend to 0,w} and the above limit yield
\begin{eqnarray}
 \limsup_{n \to \infty} \h{1pt} \left( \dfrac{ r_n } { s_n } \right)^2 \int_{ B_{ \sigma_k } } F_n\big( \widetilde{v}_{n,s,R} \big) \quad \text{ is independent of $R$ and converges to $0$ as $s \to 0$.}
\label{conv. of potential, large, origin}
\end{eqnarray}

Finally, we apply (\ref{conv. of grad v_n, large, origin})--(\ref{conv. of potential, large, origin}) to the right--hand side of (\ref{large scale basic inequality, origin}). By the lower semi--continuity, it follows
\begin{align*}
\int_{ B_{\sigma_k} } \Big| \nabla \mathcal{W}^{sc}_\infty \Big|^2
&\h{2pt}\leq\h{2pt} \liminf_{ n \to \infty } \int_{ B_{\sigma_k} } \Big| \nabla \mathcal{W}^{sc}_n \Big|^2 \h{2pt}\leq\h{2pt}
 \limsup_{ n \to \infty } \int_{ B_{\sigma_k} } \Big| \nabla \mathcal{W}^{sc}_n \Big|^2 + \left(\dfrac{r_n}{s_n}\right)^2F_n \big(\mathcal{W}_n\big) \\[2mm]
&\h{2pt}\leq\h{2pt}  \limsup_{s \to 0} \h{1pt} \limsup_{R \to \infty} \h{1pt} \limsup_{n \to \infty}  \int_{ B_{\sigma_k} }  s_n^{-2} \h{0.2pt}\big| \nabla  \widetilde{v}_{n, s, R} \big|^2 + \left( \dfrac{r_n}{s_n} \right)^2 F_n ( \widetilde{v}_{n, s, R} )
\leq \int_{ B_{\sigma_k} } \big|\nabla w \big|^2.
\end{align*}The proof is completed.
\end{proof}

\noindent\textbf{Proof of Proposition \ref{small energy implies energy decay} in the large--scale regime.}\vspace{0.2pc}\\
In light of Lemma \ref{strong H1 convergence, large, origin}, we can take $n\to \infty$ in (\ref{limit of contradictory inequality, origin}) and obtain
\begin{align}\label{contra to energy decay of minimizer, large, origin}\dfrac{1}{2} \h{2pt}\leq\h{2pt}  \nu_0^{-1}\int_{B_{\nu_0}} \Big| \nabla \mathcal{W}^{sc}_\infty \Big|^2.\end{align}
Since all the components of $\mathcal{W}^{sc}_\infty$ are harmonic in $B_{1}$, the estimates in (\ref{unif bound of nabla lim scla, origin}) still hold. Applying the estimates in (\ref{unif bound of nabla lim scla, origin}) to the right--hand side of (\ref{contra to energy decay of minimizer, large, origin}), we then can take $\nu_0 \in \big(0, 1/2\big)$ suitably small so that (\ref{contra to energy decay of minimizer, large, origin}) fails. The proof is then finished.


\section{Energy--decay estimate on balls in $\mathscr{J}$}

\noindent In this section, we prove Proposition  \ref{variant decay lemma}. 
\subsection{Blow--up sequence and some preliminary results}
The constant $\theta_0$ in Proposition \ref{variant decay lemma} is a universal constant. It will be determined during the course of the proof. Suppose on the contrary that Proposition \ref{variant decay lemma} fails. There exist $a_n$, $\lambda_n$ and $\epsilon_n$ with \begin{eqnarray} a_n \longrightarrow \infty,  \h{15pt} \lambda_n \longrightarrow 0\h{15pt}\text{and}\h{15pt} \epsilon_n \longrightarrow 0  \h{15pt}\text{as $n \rightarrow \infty$,}\label{convergence of parameters}\end{eqnarray}so that for any $n \in \mathbb{N}$, we can find a $B_{4 r_n}\left(x_n\right) \in \mathscr{J}$ with which the mapping $u_n := u_{a_n, b}^+$ satisfies  \begin{eqnarray}\text{(i).} \h{5pt}  E_n(r_n) < \epsilon_n; \h{15pt} \text{(ii).}\h{5pt} E_n\big(\lambda_n\h{1pt}\theta_0 \h{1pt} r_n\big) > r_n^{3/2}; \h{15pt} \text{(iii).} \h{5pt}E_n\big( \lambda_n\h{1pt}\theta_0 \h{1pt} r_n\big) > \dfrac{1}{2} \h{1pt}E_n(\lambda_n\h{1pt}r_n).
\label{basic assumptions}
\end{eqnarray}Here and in what follows, $E_n(r) := E_{a_n, \mu; \h{1pt}x_n, r}\h{1pt}[ u_n ]$. See (\ref{localized energy on disk}). Since now, we define
\begin{equation}
s_n^2 :=   E_n \big( \lambda_n\h{1pt}r_n\big) \h{20pt}\text{and}\h{20pt} y_n := \dashint_{D_{ \lambda_n\h{0.5pt}r_n}(\rho_{x_n}, 0)}u_n.
\label{s_n and y_n}
\end{equation}
As a convention, $\rho_x $ is the $\rho$--coordinate of $x$ in the cylindrical coordinate system. Moreover, for a new coordinate system $\big(\xi_1, \xi_2\big)$, the notation $D_r(\xi)$ is still used to represent the disk in the $\xi$\h{1pt}--\h{1pt}plane with center $\xi$ and radius $r$. If $\xi = 0$, then  $D_r(0)$ is simply denoted by $D_r$. With these notations, we let\begin{eqnarray}\label{blow sequence, j2} U_n(\xi)  :=  u_n\big(\left(\rho_{x_n}, 0\right)+  \lambda_n\h{0.2pt}r_n\xi \h{0.3pt} \big)\h{15pt}\text{and}\h{15pt}U^{sc}_n(\xi) :=\dfrac{U_n(\xi) - y_n}{s_n}, \h{20pt}\text{for any $\xi \in D_1$.}
\end{eqnarray}By this change of variables, (\ref{def of G_a}) and item (iii) in (\ref{basic assumptions}), it turns out   
\begin{eqnarray} \int_{D_{\theta_0}}  \big| D_{\xi} U^{sc}_n \big|^2 + \left( \dfrac{ \lambda_n\h{0.2pt}r_n}{s_n} \right)^2 G_{a_n, \mu} \big( \rho_{x_n} + \lambda_n\h{0.5pt} r_n \h{0.5pt} \xi_1, \h{1pt}U_n \big)  \h{2pt}>\h{2pt}\dfrac{1}{2}, \h{15pt}\text{for any $n \in \mathbb{N}$.}
\label{limit of contradictory inequality}
\end{eqnarray}
Here $D_\xi$ is the gradient operator with respect to the variable $\xi$. Utilizing   Poincar\'{e}'s inequality, $\Big\{U^{sc}_n\Big\}$ is uniformly bounded in $ H^1\big(D_1\big)$. Hence there is a subsequence, still denoted by $\Big\{U^{sc}_n\Big\}$, so that as $n \rightarrow \infty$,\begin{eqnarray}
U^{sc}_n \longrightarrow U^{sc}_{\infty} \h{20pt}  \text{ weakly in }H^1\big(D_1\big), \h{2pt}\text{strongly in } L^2\big(D_1\big)  \h{2pt} \text{and strongly in }L^2\big(T\big).
    \label{conv. of u_n}
\end{eqnarray}In (\ref{conv. of u_n}), $T$ is also used to denote the interval $\Big\{ \left(\xi_1, 0\right) : \xi_1 \in [-1, 1\h{0.3pt}]\h{2pt}\Big\}$ on the $\xi$--plane without ambiguity. \vspace{0.2pc}

 In the following, we introduce some preliminary results regarding  $r_n$, $s_n$ and $y_n$. Firstly for the sake of estimating the potential term in the integral on the left--hand side of (\ref{limit of contradictory inequality}), we need 
\begin{lem}\label{parameter estimate, j2}
$s_n + \dfrac{r_n}{s_n} \longrightarrow 0$ as $n\to \infty$. Moreover, it holds \begin{eqnarray}\label{control of yn1 over sn} 0 \h{2pt}\leq\h{2pt}\dfrac{ y_{n; 1} }{s_n} \h{3pt}\lesssim\h{3pt} \dfrac{ \rho_{x_n}}{ \lambda_n\h{0.2pt}r_n}, \h{15pt}\text{for  any $n \in \mathbb{N}$.}
\end{eqnarray}As for the third component, we have $y_{n; 3} = 0$, for any $n \in \mathbb{N}$. 
\label{estimates of s_n and y_n}
\end{lem}
\begin{proof}[\bf Proof] The convergence of $s_n$ follows by (i) in (\ref{basic assumptions}) and the convergence of $\epsilon_n$ in (\ref{convergence of parameters}). Moreover by (ii) in (\ref{basic assumptions}), it satisfies \begin{eqnarray*} s_n^2 \h{2pt}\geq\h{2pt} E_n\big( \lambda_n\h{0.5pt}\theta_0\h{0.5pt}r_n\big) \h{2pt}>\h{2pt}r_n^{3/2}, \h{10pt}\text{which infers $\dfrac{r_n}{s_n} \longrightarrow 0$, as $n \rightarrow \infty$.}
\end{eqnarray*}Due to H\"{o}lder's inequality, it turns out \begin{align*} y_{n; 1}  =  \dashint_{D_1 } U_{n; 1} \h{1.5pt}  \h{2pt}\lesssim\h{2pt}\left(\int_{D_1} U_{n; 1}^2\right)^{1/2}.
\end{align*}Since $4 r_n < \rho_{x_n}$, we then induce from the above estimate that  \begin{align*}  y_{n; 1}   \h{2pt}\lesssim\h{2pt}\rho_{x_n}\left(\int_{D_1} \left(\dfrac{U_{n; 1}}{\rho_{x_n} + \lambda_n\h{0.2pt}r_n \h{0.3pt}\xi_1}\right)^2\right)^{1/2} \h{2pt}\lesssim\h{2pt}s_n \h{1pt}\dfrac{\rho_{x_n}}{\lambda_n\h{0.2pt}r_n}.
\end{align*}The estimate of $y_{n; 1}$ in Lemma \ref{parameter estimate, j2} holds. In the end, $y_{n; 1} \geq 0$ due to item (1) in Remark \ref{rmk on sign and gamma conv} and  $y_{n; 3} = 0$ by the odd symmetry of $u_{n; 3}$  with respect to the $z$--variable. 
\end{proof}
 Owing to (i) in (\ref{basic assumptions}), the convergence of $\epsilon_n$ in (\ref{convergence of parameters}) and the uniform boundedness in item (2) of Remark \ref{rmk on sign and gamma conv}, up to a subsequence, $U_n$ converges strongly in $H^1(D_1)$ to a constant vector $y_*$.  Meanwhile, the $y_n$ defined in (\ref{s_n and y_n})  converges to $y_*$ as well when $n \rightarrow \infty$. Here one should not be confused with the $y_*$ used in Section 2, though we are using same notation to denote the limiting location of $y_n$ defined in (\ref{s_n and y_n}).  Due to $y_{n; 1} \geq 0$, we can infer $y_{*; 1} \geq 0$. However,  we cannot have in general $y_{*; 1} = 0$. It is the reason that makes our analysis for the current case complicated. As for $y_{*; 2}$, we notice that $U_{n; 2}$  satisfies $U_{n; 2} \geq H_{a_n}\h{0.5pt}b$ on $T$ in the sense of trace. Taking $n \rightarrow \infty$ and summarizing the above arguments, we  then obtain  
\begin{equation}
y_n \longrightarrow y_* = \big(y_{*;1} ,y_{*;2},0\big)^\top, \h{5pt} \text{where $y_{*; 1} $ and $y_{*; 2}$ are  constants satisfying $y_{*; 1} \geq 0$ and $y_{*; 2} \geq b$.}
\label{y_*=(y_*, 1y_*,2,0)}
\end{equation}If in addition it holds \begin{eqnarray*}
\liminf_{n \rightarrow \infty} \h{2pt} \left|  \dfrac{H_{a_n}\h{0.2pt}b - y_{n; 2}}{s_n} \right| \h{2pt}< \h{2pt}\infty,
\end{eqnarray*}then there exists a constant $\overline{w}_* \in \mathbb{R}$ so that up to a subsequence,  \begin{eqnarray} \lim_{n \rightarrow \infty} \h{2pt} \dfrac{H_{a_n} b - y_{n; 2}}{s_n} =  \overline{w}_*.
\label{pi_L conv}
\end{eqnarray}In this case, we have
\begin{lem}
If (\ref{pi_L conv}) holds, then $U^{sc}_{\infty; 2} \h{1pt}\geq\h{1pt} \overline{w}_{*}$ on $T$ in the sense of trace. 
\label{ratio convergence and limit boundary condition}
\end{lem}
\begin{proof}[\bf Proof] The second component of $U^{sc}_n$ can be decomposed into
$$U^{sc}_{n; 2}=\dfrac{ U_{n; 2} -  H_{a_n}\h{0.2pt}b }{s_n} + \dfrac{H_{a_n}\h{0.2pt}b  - y_{n; 2} }{s_n} \quad \text{on\,\,}T.$$
The lemma then follows by the Signorini obstacle boundary condition satisfied by $U_n$ on $T$.
\end{proof}

\subsection{Energy--decay estimate in small--scale regime}
In this section, we prove Proposition \ref{variant decay lemma} by supposing $a_n \big(\lambda_n \h{0.5pt} r_n\big)^2 \longrightarrow 0$ as $n \to \infty$. By Lemma \ref{ratio convergence and limit boundary condition} and Fatou's lemma, it holds  \begin{lem}
There exist an increasing positive sequence $\big\{\sigma_k\big\}$ which tends to $1$ as $k \to \infty$, a sequence of positive numbers $\big\{b_k\big\}$ and a subsequence of $\big\{U^{sc}_n\big\}$, still denoted by $\big\{U^{sc}_n\big\}$, so that for any $k$, we have \begin{enumerate}
\item[$\mathrm{(1).}$] The uniform upper bound: $$\sup_{n \h{0.2pt}\in\h{0.2pt}\mathbb{N}\h{0.4pt}\cup\h{0.2pt}\{\infty\}} \left\{\big\lVert  U^{sc}_n\big\rVert_{\infty;\h{1pt}\p D_{\sigma_k}} + \int_{\p D_{\sigma_k}} \big|\h{1pt}D_\xi U^{sc}_n\h{1pt}\big|^2 \right\}   \h{2pt}\leq\h{2pt}b_k;$$
\item[$\mathrm{(2).}$] The convergence $U^{sc}_n \longrightarrow U^{sc}_\infty$ in $C^0\big(\p D_{\sigma_k}\big)$ as $n \rightarrow \infty$;
\item[$\mathrm{(3).}$]  The second component of $U_n = y_n + s_n \h{0.4pt}U^{sc}_n$ satisfies $U_{n; 2}   \geq H_{a_n} b$ at  $\big(\pm \sigma_k, 0 \big)$; 
\item[$\mathrm{(4).}$] The third component of $U^{sc}_\infty$ satisfies $U^{sc}_{\infty; 3} = 0$ at  $\big( \pm \sigma_k, 0\h{1pt}\big)$; 
\item[$\mathrm{(5).}$] If (\ref{pi_L conv}) holds, then  $U^{sc}_{\infty; 2} \geq \overline{w}_*$ at  $\big( \pm \sigma_k, 0\h{1pt}\big)$.
\end{enumerate}
\label{uniform convergence of u_n in S_sigma_k}
\end{lem}
Using Lemma \ref{ratio convergence and limit boundary condition} and $\big\{\sigma_k\big\}$ obtained in Lemma \ref{uniform convergence of u_n in S_sigma_k}, we introduce two configuration spaces: \begin{align}   &\mathfrak{M}_{\h{0.5pt}k} := \Big\{ u \in H^1\left( D_{\sigma_k} ; \h{2pt}\mathbb{R}^3 \right) : u = U^{sc}_{\infty} \h{4pt}\text{on $\p D_{\sigma_k}$}, \text{$u_1$ and $u_2$ are even and $u_3$ is odd with respect to $\xi_2$}\Big\};\nonumber\\[2mm]
& \overline{\mathfrak{M}}_{\h{0.5pt}k} :=  \Big\{ u \in \mathfrak{M}_{\h{0.5pt}k} : u_2 \h{2pt}\geq\h{2pt}\overline{w}_{*} \h{4pt}\text{on $T_{\sigma_k} := \Big\{ \left(\xi_1, 0 \right) : \xi_1 \in \big[- \sigma_k, \sigma_k\big] \Big\}$} \Big\}.\label{configuration space M_k}
\end{align}We now prove the following energy--minimizing property of $U^{sc}_\infty$ in the small--scale regime.
\begin{lem}For any natural number $k$, if it satisfies \begin{eqnarray}
\liminf_{n \rightarrow \infty} \h{2pt} \left|  \dfrac{ H_{a_n} b - y_{n; 2}}{s_n} \right| \h{2pt}= \h{2pt}\infty,
\label{infty case}
\end{eqnarray}then $U^{sc}_\infty$ minimizes the Dirichlet energy within the space $\mathfrak{M}_{\h{0.5pt}k}$. Otherwise, if (\ref{pi_L conv}) holds, then $U^{sc}_\infty$ minimizes the Dirichlet energy within the space $\overline{\mathfrak{M}}_{\h{0.5pt}k}$. In both cases, $U^{sc}_n$ converges to $U^{sc}_{\infty}$ strongly in $H_{\mathrm{loc}}^1\big(D_1\big)$.
\label{strong H1 convergence in small case}
\end{lem}

\begin{proof}[\bf Proof] The proof is similar to the proof of Lemma \ref{strong H1 convergence, small, origin}. Here we just show the differences.  Suppose that $v$ is an arbitrary map in $\mathfrak{M}_{\h{0.5pt}k}$. Then we define
\begin{equation}
 F_{n, R}\h{0.5pt} [\h{0.5pt}v\h{0.5pt}] :=\left\{ \begin{array}{lcl}  y_n +R\h{1pt}s_n \h{1pt}\dfrac{v - Y_*}{|\h{1pt}v - Y_*\h{1pt}|\vee R}, \h{30pt}&&\text{if (\ref{infty case}) holds;}\vspace{0.8pc}\\
  y_n^*+R\h{1pt}s_n \h{1pt}\dfrac{v - Y_* }{|\h{1pt}v - Y_*\h{1pt}|\vee R}, &&\text{if (\ref{pi_L conv}) holds.}
\end{array}\right.
\label{definition of F^s_n}
\end{equation}
Here $R > 0$ is a positive constant. $y_n^* = \big(y_{n; 1}, H_{a_n} b, 0 \big)^{\top}$. The vector $Y_* = 0$ if (\ref{infty case}) holds. If (\ref{pi_L conv}) holds, then $Y_* = \overline{w}_* e^*_2$. In addition for any fixed $s \in (0, 1)$, we introduce a comparison map:
\begin{eqnarray}\mathscr{V}_{n, s, R}\big(\xi\big) :=\left\{
\begin{aligned}
&F_{n, R}\h{1pt}[\h{1pt}v\h{1pt}]\left(\dfrac{\xi}{1-s}\right) \quad &\text{if }& \xi \in D_{(1-s)\h{1pt}\sigma_k}; \\[2mm]
&\dfrac{\sigma_k-|\h{1pt}\xi\h{1pt}|}{s\h{0.5pt}\sigma_k} \h{1pt} F_{n, R} \h{1pt}\big[ U^{sc}_{\infty} \big] \big(\sigma_k\h{1pt}\widehat{\xi}\h{2pt}\big)+\dfrac{|\h{1pt}\xi\h{1pt}|-(1-s)\h{1pt}\sigma_k}{s\h{0.5pt}\sigma_k}\h{1pt}U_n\big(\sigma_k\h{1pt}\widehat{\xi}\h{2pt}\big) \quad &\text{if }& \xi \in D_{\sigma_k} \setminus D_{(1-s)\sigma_k}.
\end{aligned}
\right.
\label{definition of v_n}
\end{eqnarray}
It can be shown that $\big[ \mathscr{V}_{n, s, R} \big]_2 \geq H_{a_n} b$ on $T_{\sigma_k}$ if  (\ref{infty case})  or   (\ref{pi_L conv})  holds.   In light of the energy--minimizing property of $u_n$, it turns out
\begin{eqnarray}&&\int_{D_{\sigma_k}} \big|D_{\xi} U^{sc}_n\big|^2 \left(1 + \dfrac{\lambda_n \h{1pt}r_n}{\rho_{x_n}} \h{1pt}\xi_1\right)\nonumber\\ [2mm]
 &&\leq \int_{D_{\sigma_k}}\left\{ \dfrac{1}{ s_n^2 } \h{2pt} \big|D_{\xi} \mathscr{V}_{n, s, R}\big|^2 + \left(\dfrac{\lambda_n\h{1pt}r_n}{s_n}\right)^2 G_{a_n, \mu}\big(\rho_{x_n} + \lambda_n\h{1pt}r_n\h{1pt}\xi_1,  \h{2pt}\cdot \h{3pt} \big)  \h{2pt}\bigg|_{U_n}^{\mathscr{V}_{n, s, R}} \right\} \left(1 + \dfrac{\lambda_n \h{1pt}r_n}{\rho_{x_n}} \h{1pt}\xi_1\right).
\label{small scale basic inequality}
\end{eqnarray}Note that $G_{a_n, \mu}$ is given in (\ref{def of G_a}). Slightly modifying the proof for (\ref{conv. of grad v_n, origin}), we have  
\begin{eqnarray}\lim_{s \rightarrow 0}  \lim_{R \rightarrow \infty}  \lim_{n \rightarrow \infty} \dfrac{1}{s_n^2}\h{2pt}\int_{D_{\sigma_k}} \big|\h{1pt}D_\xi \mathscr{V}_{n, s, R}\h{1pt}\big|^2 \left( 1 + \dfrac{\lambda_n \h{1pt}r_n}{\rho_{x_n}} \h{1pt} \xi_1\right) = \int_{D_{\sigma_k}}\big|\h{1pt}D_\xi v \h{1pt}\big|^2. \label{main limit on right}\end{eqnarray}
It remains to study the potential term $G_{a_n, \mu}$. Notice that  \begin{eqnarray} \int_{D_{\sigma_k}} G_{a_n, \mu}\big(\rho_{x_n} + \lambda_n\h{1pt}r_n\h{1pt}\xi_1,  \h{2pt}\cdot \h{3pt} \big)  \h{2pt}\bigg|_{U_n}^{\mathscr{V}_{n, s, R}} \left( 1 + \dfrac{\lambda_n \h{1pt}r_n}{\rho_{x_n}} \h{1pt} \xi_1\right)   =  I^s_1 + I^s_2 + I^s_3 + I^s_4.
\label{I_1 to I_4}
\end{eqnarray}Now we define and estimate the terms on the right--hand side above. One should not be confused with the $I_1^s$ and $I_2^s$ introduced in (\ref{I_1 to I_2, origin}). The expressions of $I_1^s$ and $I_2^s$ used in this proof will be given as follows.\\
\\
\textbf{Estimate of $I^s_1$ and $I_2^s$.} Similarly as in the proof of Lemma \ref{strong H1 convergence, small, origin}, we define \begin{eqnarray*}I^s_1 :=- 3 \sqrt{2}\h{1pt} \mu\int_{D_{\sigma_k}} \big(\mathscr{V}_{n, s, R} - U_n\big) \cdot \int_0^1 \nabla_u P \h{2pt} \bigg|_{u \h{1pt}=\h{1pt} t \h{1pt}\mathscr{V}_{n, s, R} + ( 1 - t ) \h{1pt}U_n} \left( 1 + \dfrac{\lambda_n \h{1pt}r_n}{\rho_{x_n}} \h{1pt} \xi_1\right)
\end{eqnarray*}and \begin{eqnarray*}
I^s_2 := 2 \h{1pt}a_n \h{1pt}\mu  \int_{D_{\sigma_k}}\left( 1 + \dfrac{\lambda_n \h{1pt}r_n}{\rho_{x_n}} \h{1pt} \xi_1\right) \big(\mathscr{V}_{n, s, R} - U_n \big) \cdot
\int_0^1 \Big( \big| \h{1pt}U_n + t \h{1pt}\big( \mathscr{V}_{n, s, R} - U_n \big) \h{1pt}\big|^2 - 1 \Big) \h{2pt}\Big( U_n + t \h{1pt}\big( \mathscr{V}_{n, s, R} - U_n \big)\Big) .
\end{eqnarray*} Then similar arguments for (\ref{I_1^s, origin}) and (\ref{est of I2s, origin}) yield \begin{align}\label{conv of I1s and I2s, J2} \lim_{n \rightarrow \infty}\left(\dfrac{\lambda_n r_n}{s_n}\right)^2\big| I_1^s\big|  + \left(\dfrac{\lambda_n r_n}{s_n}\right)^2\big| I_2^s\big| \longrightarrow 0, \h{15pt}\text{as $n \rightarrow \infty$.}
\end{align}
\textbf{Estimate of $I^s_3$.} $I^s_3$ is defined by \begin{eqnarray*}I^s_3 := \int_{D_{\sigma_k}} \dfrac{4 \h{1pt}\big( [\mathscr{V}_{n, s, R}]_1 - U_{n; 1} \big)^2 + \big([\mathscr{V}_{n, s, R}]_3 - U_{n; 3} \big)^2}{\big( \rho_{x_n} + \lambda_n \h{1pt} r_n \h{1pt} \xi_1 \big)^2}  \left( 1 + \dfrac{\lambda_n \h{1pt}r_n}{\rho_{x_n}} \h{1pt} \xi_1\right) . \end{eqnarray*}In light of (1) in Lemma \ref{uniform convergence of u_n in S_sigma_k} and the definitions of $\mathscr{V}_{n, s, R}$, $U^{sc}_n$, for suitably large $R > 0$, it satisfies \begin{eqnarray}\left|\h{1pt}\dfrac{\mathscr{V}_{n, s, R} - U_n}{s_n} \h{1pt}\right| \h{2pt}\lesssim\h{2pt} b_k + R + |\h{1pt}U^{sc}_n\h{1pt}|  \h{20pt}\text{on $D_{\sigma_k}$.}
\label{v_n - U_n}
\end{eqnarray}Therefore we obtain \begin{eqnarray} \big|\h{1pt} I^s_3 \h{1pt}\big|  \h{2pt}\lesssim\h{2pt} \big[\h{1pt} b_k +  R \h{2pt}\big]^2 \left(\dfrac{s_n}{\rho_{x_n}}\right)^2.
\label{I_3^s}
\end{eqnarray}Here we also have used the uniform boundedness of $U_n^{sc}$ in $L^2\big(D_1\big)$.\vspace{0.4pc}

\noindent \textbf{Estimate of $I^s_4$.} $I^s_4$ is defined by \begin{eqnarray*}I^s_4 := \int_{D_{\sigma_k}} \dfrac{8 \h{1pt}U_{n; 1} \h{1pt}\big( [\mathscr{V}_{n, s, R}]_1 - U_{n; 1} \big) + 2 \h{1pt}U_{n; 3} \h{1pt} \big( [\mathscr{V}_{n, s, R}]_3 - U_{n; 3} \big)}{\big( \rho_{x_n} + \lambda_n \h{1pt} r_n \h{1pt} \xi_1 \big)^2} \left( 1 + \dfrac{\lambda_n \h{1pt}r_n}{\rho_{x_n}} \h{1pt} \xi_1\right).\end{eqnarray*}Utilizing (\ref{v_n - U_n}), the definition of $U^{sc}_n$ and $y_{n, 3} = 0$, we have \begin{eqnarray*}\big|\h{1pt}I^s_4\h{1pt}\big| \h{3pt}\lesssim\h{3pt}\big[ \h{1pt}b_k + R  \h{2pt}\big]\h{1pt}\left(\dfrac{s_n}{\rho_{x_n}}\right)^2  \left[  1 + \dfrac{ y_{n; 1} }{s_n}\right], \end{eqnarray*} which induces by (\ref{control of yn1 over sn}) in Lemma \ref{estimates of s_n and y_n} the following estimate: \begin{eqnarray}\big|\h{1pt}I^s_4\h{1pt}\big| \h{2pt}\lesssim\h{2pt}\big[ \h{1pt}b_k + R  \h{2pt}\big]  \dfrac{\rho_{x_n}}{\lambda_n\h{1pt}r_n} \left(\dfrac{s_n}{\rho_{x_n}}\right)^2.
\label{I_4^s}
\end{eqnarray} 
By (\ref{conv of I1s and I2s, J2}), (\ref{I_3^s}) and (\ref{I_4^s}), there follows \begin{eqnarray*} &&\left|\h{2pt}\int_{D_{\sigma_k}} \left(\dfrac{\lambda_n\h{1pt}r_n}{s_n}\right)^2 G_{a_n, \mu}\big(\rho_{x_n} + \lambda_n\h{1pt}r_n\h{1pt}\xi_1,  \h{2pt}\cdot \h{3pt} \big)  \h{2pt}\bigg|_{U_n}^{\mathscr{V}_{n, s, R}} \left(1 + \dfrac{\lambda_n \h{1pt}r_n}{\rho_{x_n}} \h{1pt}\xi_1\right) \h{2pt}\right|  \h{3pt}\longrightarrow\h{3pt} 0, \h{20pt}\text{as $n \rightarrow \infty$.}
\end{eqnarray*}Here we also have used the convergence of $\lambda_n$ in (\ref{convergence of parameters}).
Applying this limit and (\ref{main limit on right}) to the right--hand side of (\ref{small scale basic inequality}) yields
\begin{eqnarray*}
\int_{D_{\sigma_k}} \big|\h{1pt}D_{\xi} U^{sc}_{\infty} \h{1pt}\big|^2 \leq \int_{D_{\sigma_k}} \big|D_{\xi} v\big|^2.
\end{eqnarray*} The proof is completed.
\end{proof}
\begin{proof}[\bf Proof of Proposition \ref{variant decay lemma}  in small--scale regime]\
\\[2mm]
The proof is to find a universal constant $\theta_0$ so that (\ref{limit of contradictory inequality}) fails. Notice that for any $\theta \in (0, 1)$, \begin{eqnarray} \int_{D_{\theta}} G_{a_n, \mu} \big( \rho_{x_n} + \lambda_n \h{0.5pt} r_n \h{0.5pt} \xi_1, \h{1pt}U_n \big) = J^s_1 + J^s_2 + J^s_3,\label{decom of potential} \end{eqnarray}
where \begin{eqnarray*}
J^s_1 := \mu \int_{D_{\theta}}  D_{a_n}  - 3 \sqrt{2}\h{1pt}P(U_n), \h{15pt} J^s_2 :=  \int_{D_{\theta}}  \dfrac{4\h{1pt} U_{n; 1}^2 + U_{n; 3}^2}{\big( \rho_{x_n} + \lambda_n\h{0.5pt}r_n\h{0.5pt}\xi_1\big)^2}, \h{15pt}J^s_3 := \dfrac{a_n\h{0.5pt}\mu}{2} \int_{D_{\theta}}  \Big( \big| U_n \big|^2 - 1 \Big)^2.
\end{eqnarray*}In light of the uniform boundedness of $D_{a_n}$, $U_n$ and the limit of $r_n/s_n$ in Lemma \ref{estimates of s_n and y_n}, it turns out \begin{eqnarray} \left(\dfrac{\lambda_n\h{0.5pt}r_n}{s_n}\right)^2 \big|\h{1pt}J_1^s \h{1pt}\big| \h{2pt}\lesssim\h{2pt}\left(\dfrac{\lambda_n\h{0.5pt}r_n}{s_n}\right)^2 \longrightarrow 0,\h{20pt}\text{as $n \rightarrow \infty$.}
\label{limit of J_1^s}
\end{eqnarray}
By (\ref{control of yn1 over sn}) in Lemma \ref{estimates of s_n and y_n}, we can find a non--negative constant $c_3$ so that up to a subsequence, there holds  \begin{eqnarray*} \left(\dfrac{ \lambda_n\h{0.5pt}r_n \h{0.5pt} y_{n; 1}}{\rho_{x_n}  \h{0.5pt}s_n }\right)^2 \longrightarrow c_3, \h{20pt}\text{as $n \rightarrow \infty$.}
\end{eqnarray*}
On the other hand, it satisfies \begin{eqnarray*}\left(\dfrac{\lambda_n\h{0.5pt}r_n}{\rho_{x_n}}\right)^2 \int_{D_{\theta}} \big|\h{0.5pt}U^{sc}_n\h{0.5pt}\big|^2 \longrightarrow 0, \h{20pt}\text{as $n \rightarrow \infty$.} \end{eqnarray*}Utilizing the last two limits and the fact that $U_n = y_n + s_n\h{0.5pt}U^{sc}_n$, we obtain \begin{eqnarray}\left(\dfrac{\lambda_n\h{0.5pt}r_n}{s_n}\right)^2 J_2^s \longrightarrow 4 \h{0.5pt}\pi\h{0.5pt}c_3\h{0.5pt}\theta^{\h{0.5pt}2}, \h{20pt}\text{as $n \rightarrow \infty$.}
\label{limit of J_2^s}
\end{eqnarray}As for $J_3^s$ term, still by $U_n = y_n + s_n\h{0.5pt}U^{sc}_n$, we have \begin{eqnarray} \left(\dfrac{\lambda_n\h{0.5pt}r_n}{s_n}\right)^2J_3^s = \left(\dfrac{\lambda_n\h{0.5pt}r_n}{s_n}\right)^2\dfrac{a_n\h{0.5pt}\mu}{2} \int_{D_{\theta}}  \Big( \big| y_n \big|^2 - 1 + 2 \h{0.5pt}s_n\h{0.5pt}y_n \cdot U^{sc}_n + s_n^2 \big|U^{sc}_n\big|^2\Big)^2.
\label{decom of J_3^s}
\end{eqnarray}
 In light that $U^{sc}_n$ is uniformly bounded in $L^4\big(D_1\big)$, it satisfies \begin{eqnarray} \left(\dfrac{\lambda_n\h{0.5pt}r_n}{s_n}\right)^2\dfrac{a_n\h{0.5pt}\mu}{2} \int_{D_{\theta} }  \left|  2 \h{0.5pt}s_n\h{0.5pt}y_n \cdot U^{sc}_n + s_n^2 \big|U^{sc}_n\big|^2\right|^2 \h{3pt}\lesssim\h{3pt} a_n \big(\lambda_n r_n\big)^2.
\label{small term in J_3^s}
\end{eqnarray}
Utilizing the limits (\ref{limit of J_1^s})--(\ref{limit of J_2^s}) and the fact that \begin{eqnarray}  \left(\dfrac{\lambda_n\h{0.5pt}r_n}{s_n}\right)^2 \int_{D_{\theta}} G_{a_n, \mu} \big( \rho_{x_n} + \lambda_n \h{0.5pt} r_n \h{0.5pt} \xi_1, \h{1pt}U_n \big) \h{2pt}\leq 1,
\label{bound of potential}
\end{eqnarray}
we obtain the uniform boundedness of $\left(\dfrac{\lambda_n\h{0.5pt}r_n}{s_n}\right)^2 J_3^s$. By this uniform boundedness and (\ref{decom of J_3^s})--(\ref{small term in J_3^s}), there exists a universal non--negative constant $c_4$ so that up to a subsequence, it holds \begin{eqnarray}  \left(\dfrac{\lambda_n\h{0.5pt}r_n}{s_n}\right)^2 \dfrac{a_n\h{0.5pt}\mu}{2} \left(\big| y_n \big|^2 - 1\right)^2 \longrightarrow c_4, \h{20pt}\text{as $n \rightarrow \infty$.}
\label{limit of term in potential}
\end{eqnarray}Applying this limit and (\ref{small term in J_3^s}) to (\ref{decom of J_3^s}), in the small--scale regime, we have \begin{eqnarray*}\left(\dfrac{\lambda_n\h{0.5pt}r_n}{s_n}\right)^2J_3^s \longrightarrow  \pi \h{0.5pt}c_4\h{0.5pt}\theta^{\h{0.5pt}2}, \h{20pt}\text{as $n \rightarrow \infty$.}
\end{eqnarray*}By this limit and (\ref{limit of J_1^s})--(\ref{limit of J_2^s}), we can obtain from the decomposition (\ref{decom of potential}) the limit \begin{eqnarray}\left(\dfrac{\lambda_n\h{0.5pt}r_n}{s_n}\right)^2 \int_{D_{\theta}} G_{a_n, \mu} \big( \rho_{x_n} + \lambda_n \h{0.5pt} r_n \h{0.5pt} \xi_1, \h{1pt}U_n \big) \longrightarrow  \pi\h{0.5pt}\theta^{\h{0.5pt}2} \big( 4 c_3 + c_4\big), \h{20pt}\text{as $n \rightarrow \infty$.}
\label{limit of potential}
\end{eqnarray}The last limit and (\ref{bound of potential}) infer the bound \begin{align}\label{bound of c_3} \pi \h{0.5pt} \big[4 c_3 + c_4 \big] \leq  \theta^{\h{0.5pt}-2}, \h{10pt}\text{which induces \h{3pt} $ \pi \big[ 4 c_3 + c_4 \big] \leq 1$.}\end{align} Here we take $\theta \rightarrow 1^-$. By the above bound and (\ref{limit of potential}), for any $\theta \in (0, 1)$,  it holds \begin{eqnarray*} \lim_{n \rightarrow \infty} \left(\dfrac{\lambda_n\h{0.5pt}r_n}{s_n}\right)^2 \int_{D_{\theta}} G_{a_n, \mu} \big( \rho_{x_n} + \lambda_n \h{0.5pt} r_n \h{0.5pt} \xi_1, \h{1pt}U_n \big) \leq \theta^{\h{0.5pt}2}.
\end{eqnarray*} The above estimate and strong $H_{\h{0.5pt}\mathrm{loc}}^1$\h{0.5pt}--\h{0.5pt}convergence of $U^{sc}_n$ obtained in Lemma \ref{strong H1 convergence in small case} can be applied to the left--hand side of (\ref{limit of contradictory inequality}). Therefore, by taking $n \to \infty$ in (\ref{limit of contradictory inequality}), it follows \begin{eqnarray}\label{contra, j2, small} \dfrac{1}{2} \leq \theta_0^{\h{0.5pt}2} + \int_{D_{\theta_0}} \big|\h{0.5pt}D_{\xi} U^{sc}_\infty \h{0.5pt}\big|^2.  	
\end{eqnarray} For the components of $U_{\infty}^{sc}$ which are harmonic functions on $D_1$, we can apply uniform $H^1\left(D_1\right)$--boundedness of $U_\infty^{sc}$ to get the uniform $L^\infty$--boundedness of the gradient of these components on $D_{1/2}$.   If $U_{\infty; 2}^{sc}$ solves the Signorini obstacle problem, then we can have \begin{align*} \int_{D_{\theta_0}} \big|\h{0.5pt}D_{\xi} U^{sc}_{\infty; 2} \h{0.5pt}\big|^2 \lesssim  \theta_0 \left(\int_{D_{\theta_0}} \big|\h{0.5pt}D_{\xi} U^{sc}_{\infty; 2} \h{0.5pt}\big|^4\right)^{1/2} \leq  \theta_0 \left(\int_{D_{1/2}} \big|\h{0.5pt}D_{\xi} U^{sc}_{\infty; 2} \h{0.5pt}\big|^4\right)^{1/2}. 
\end{align*}By Sobolev inequality and Lemma 9.1 in \cite{PSU12}, we  have smallness of the right--hand side  above by choosing a small and universal $\theta_0$. Therefore (\ref{contra, j2, small}) fails, provided that $\theta_0$ is suitably small. The smallness is universal. The proof is completed.
\end{proof}
\subsection{Energy--decay estimate in intermediate--scale regime}
In this section we suppose that  $a_n \big(\lambda_n \h{0.5pt} r_n\big)^2 \longrightarrow L$ as $n \to \infty$. Here $L \in (0, \infty)$ is a constant.  \begin{lem}There exists a  $h \in \mathbb{R}$ so that up to a subsequence \begin{eqnarray} \dfrac{ |\h{0.5pt}y_n\h{0.5pt}| - 1 }{s_n} \longrightarrow h, \h{20pt}\text{as $n \rightarrow \infty$.}\label{ration convergence }\end{eqnarray} 
\label{y_*, interior, large scale}
\end{lem} The proof of this lemma follows by (\ref{limit of term in potential}) and the assumption that $a_n\big(\lambda_n r_n \big)^2 \longrightarrow L \neq 0$ as $n \rightarrow \infty$.  We are ready to characterize the energy--minimizing property of $U^{sc}_\infty$ in the intermediate--scale regime.
\begin{lem}
For any  $k \in \mathbb{N}$, if (\ref{infty case}) holds, then $U^{sc}_\infty$ minimizes the $E_{L, h}$--energy over $\mathfrak{M}_k$. Here
\begin{eqnarray*}
E_{L, h}[\h{0.5pt}u\h{0.5pt}]
:= \int_{D_{\sigma_k}} \big|\h{1pt}D_{\xi} u \h{1pt}\big|^2 + 2\h{0.5pt}L\h{0.5pt}\mu \h{1pt}\big(h + y_* \cdot u \h{1pt}\big)^2, \h{20pt}\text{for all $ u \in \mathfrak{M}_k$.}
\end{eqnarray*}The constant $h$ is obtained in (\ref{ration convergence }). $y_*$ is the limit of $y_n$ in (\ref{y_*=(y_*, 1y_*,2,0)}). If (\ref{pi_L conv}) holds, then $U^{sc}_\infty$ minimizes $E_{L, h}$--energy over $\overline{\mathfrak{M}}_{\h{0.5pt}k}$. In both cases, $U^{sc}_n$ converges to $U^{sc}_{\infty}$ strongly in $H_{\mathrm{loc}}^1\big(D_1\big)$.
\label{strong H1 convergence in intermediate case}
\end{lem}
\begin{proof}[\bf Proof] We use the same comparison mapping and notations as in the proof  of Lemma  \ref{strong H1 convergence in small case} and the proof of Proposition \ref{variant decay lemma} for the small--scale case. Recall (\ref{decom of potential}). In the intermediate--scale regime, we still have (\ref{limit of J_1^s})--(\ref{limit of J_2^s}). For $J_3^s$, we can take $n \rightarrow \infty$ in   (\ref{decom of J_3^s}). By (\ref{ration convergence }) and the strong $L^4\left(B_1\right)$--convergence of $U_n^{sc}$ to $U_\infty^{sc}$, it turns out  \begin{eqnarray*}\left(\dfrac{\lambda_n\h{0.5pt}r_n}{s_n}\right)^2J_3^s \longrightarrow 2 \h{0.5pt}L\h{0.5pt}\mu \int_{D_\theta}\big( h + y_* \cdot U^{sc}_\infty\big)^2  \h{20pt}\text{as $n \rightarrow \infty$, for any $\theta \in (0, 1)$.}
\end{eqnarray*}Taking $\theta = \sigma_k$ and utilizing the last limit and (\ref{decom of potential})--(\ref{limit of J_2^s}), we then get \begin{eqnarray}\left(\dfrac{\lambda_n\h{0.5pt}r_n}{s_n}\right)^2 \int_{D_{\sigma_k}} G_{a_n, \mu} \big( \rho_{x_n} + \lambda_n \h{0.5pt} r_n \h{0.5pt} \xi_1, \h{1pt}U_n \big) \longrightarrow 4 \h{0.5pt}\pi \h{0.5pt}c_3\h{0.5pt} \sigma_k^2  +  2 \h{0.5pt}L\h{0.5pt}\mu \int_{D_{\sigma_k}}\big( h + y_* \cdot U^{sc}_\infty\big)^2, \h{10pt}\text{as $n \rightarrow \infty$.}
\label{conv. of potential U_n in intermediate}
\end{eqnarray}

By the definition of $\mathscr{V}_{n, s, R}$ in (\ref{definition of v_n}), the limit in (\ref{pi_L conv}) and the uniform convergence in (2) of Lemma \ref{uniform convergence of u_n in S_sigma_k}, the comparison map  $\mathscr{V}_{n, s, R}$ equals $y_n + s_n \h{0.5pt}\omega^{sc}_n$ with  $$\omega^{sc}_n \longrightarrow \omega^{sc}_\infty \h{15pt}\text{uniformly in $D_{\sigma_k}$ as $n \rightarrow \infty$. } $$
Here $\omega^{sc}_\infty$ is defined as follows: $$\left\{\begin{array}{lcl} Y_* + R \h{1pt}\dfrac{v - Y_*}{\big|\h{1pt}v - Y_*\h{1pt}\big| \vee R}\h{2pt}\bigg|_{\frac{\xi}{1-s}} \h{15pt}&&\text{if $\xi \in D_{(1 - s)\h{0.5pt}\sigma_k}$;}\\[8mm]
\dfrac{\sigma_k-|\h{1pt}\xi\h{1pt}|}{s\h{0.5pt}\sigma_k} \h{1pt} \left[ Y_* + R \h{1pt}\dfrac{U^{sc}_\infty - Y_*}{\big|\h{1pt}U^{sc}_\infty - Y_*\h{1pt}\big| \vee R} \h{2pt}\bigg|_{\sigma_k \widehat{\xi} } \h{1pt}\right]+\dfrac{|\h{1pt}\xi\h{1pt}|-(1-s)\h{1pt}\sigma_k}{s\h{0.5pt}\sigma_k}\h{1pt}U^{sc}_\infty\left(\sigma_k\h{1pt}\widehat{\xi}\h{1pt}\right)  &&\text{if $\xi \in D_{\sigma_k} \setminus D_{(1-s)\sigma_k}$}.
\end{array}\right.$$ Applying the same derivations for (\ref{conv. of potential U_n in intermediate}), we get \begin{eqnarray*}\left(\dfrac{\lambda_n\h{0.5pt}r_n}{s_n}\right)^2 \int_{D_{\sigma_k}} G_{a_n, \mu} \big( \rho_{x_n} + \lambda_n \h{0.5pt} r_n \h{0.5pt} \xi_1, \h{1pt}\mathscr{V}_{n, s, R} \big) \longrightarrow 4\h{0.5pt}\pi\h{0.5pt}c_3\h{0.5pt}\sigma_k^2  +  2 \h{0.5pt}L\h{0.5pt}\mu \int_{D_{\sigma_k}}\big( h+ y_* \cdot \omega^{sc}_\infty\big)^2, \h{10pt}\text{as $n \rightarrow \infty$.}
\end{eqnarray*}In light of this limit, (\ref{conv. of potential U_n in intermediate}), (\ref{main limit on right}) and the weak convergence in (\ref{conv. of u_n}), we then can take $n \rightarrow \infty$, $R \rightarrow \infty$ and $s \rightarrow 0$ successively on both sides of (\ref{small scale basic inequality}) and obtain $E_{L, h} [ \h{1pt}U^{sc}_\infty\h{0.5pt}]\h{2pt}\leq\h{2pt}E_{L, h} [ \h{1pt}v\h{0.5pt}]$.
\end{proof}

By strong $H^1$--convergence in Lemma \ref{strong H1 convergence in intermediate case}, (\ref{conv. of potential U_n in intermediate}) and (\ref{bound of c_3}), we can take $n \rightarrow \infty$ in  (\ref{limit of contradictory inequality}) and obtain \begin{align}\label{con in interm scale} \int_{D_{\theta_0}}\big| D_\xi U_\infty^{sc} \big|^2 + 2 L \mu \big(h + y_* \cdot U_{\infty}^{sc}\big)^2 \h{2pt} \geq \h{2pt} \dfrac{1}{2} -  4\h{0.4pt}\pi\h{0.4pt}c_3\h{0.4pt}\theta_0^2 \h{2pt}\geq\h{2pt} \dfrac{1}{2} - \theta_0^2 \h{2pt}\geq\h{2pt} \dfrac{1}{4}.
\end{align}
Here $\theta_0$ is assumed to be in $ (0, 1/4)$. Now, we derive a contradiction to (\ref{con in interm scale}) with a small radius $\theta_0$ independent of $L$.  
\begin{proof}[\bf Proof of Proposition \ref{variant decay lemma} in intermediate--scale regime]\
\\[3mm]
If (\ref{infty case}) holds, then we can apply same arguments as in the proof of Proposition \ref{small energy implies energy decay} in the intermediate--scale regime. The reason is due to the sub--harmonicity of $\big|\h{1pt}D_{\xi} U_{\infty}^{sc} \h{1pt}\big|^2 + 2\h{0.5pt}L\h{0.5pt}\mu \h{1pt}\big(h + y_* \cdot U_\infty^{sc} \h{1pt}\big)^2$.   The remaining of this proof is deovted to studying the case when (\ref{pi_L conv}) holds. Due to (\ref{y_*=(y_*, 1y_*,2,0)}), (\ref{pi_L conv}) and (\ref{ration convergence }), we have $y_* = (y_{*; 1}, y_{*; 2}, 0)^\top$ with $y_{*; 1} = \sqrt{1-b^2}$ and $y_{*; 2} = b$ in the current case.  Notice that
\begin{align}
E_{L, h}\h{1pt} [\h{0.5pt}u\h{0.5pt}] \h{2pt}=\h{2pt} & E_{L, h}^{\star}\h{1pt}[ \h{0.5pt}u_1, u_2\h{0.5pt}] + \int_{D_{\sigma_k}}  \big| D_{\xi} u_3 \big|^2\h{2pt}, \nonumber\\[2mm]&\text{where \h{2pt} $E_{L, h}^{\star}[\h{0.5pt}u_1, u_2\h{0.5pt}] := \int_{D_{\sigma_k}} e_{L, h}^\star\left[u\right] := \int_{D_{\sigma_k}}  \sum_{j = 1}^2 \big| D_{\xi} u_j \big|^2  +   2\h{0.2pt}L\h{0.2pt}\mu\h{0.5pt}\Big( h + \sum_{j = 1}^2 y_{*; j}\h{1pt} u_j  \Big)^2.$}
\label{E_L^star}
\end{align}By Lemma \ref{strong H1 convergence in intermediate case}, the third component $U^{sc}_{\infty; 3}$ is harmonic on $D_{\sigma_k}$. Standard elliptic estimate yields \begin{eqnarray}   \big\| D_\xi U^{sc}_{\infty; 3} \big\|_{\infty; D_{1/2}}  \h{2pt}\lesssim\h{2pt}   \big\| U^{sc}_{\infty; 3} \big\|_{2; D_1} \h{2pt}\lesssim\h{2pt}1.\label{decay of 3rd component}
\end{eqnarray}Here and throughout the rest of the proof, we take $\sigma_k \in (7/8, 1)$. \vspace{0.2pc}

 Define the configuration space:
\begin{eqnarray*}
H_{k, \overline{w}_*} := \bigg\{ v = \big(v_1, v_2\big) \in H^1\big(D_{\sigma_k}; \mathbb{R}^2\big) : \text{$v\left(\xi_1, \xi_2\right) = v\left(\xi_1, - \xi_2\right)$ for any $\xi \in D_{\sigma_k}$  and }v_2 \geq \overline{w}_* \h{5pt}\text{on $T_{\sigma_k}$}\bigg\}.
\end{eqnarray*} Then $u$ is called a solution of the problem $S_{L, h, \overline{w}_*}$ on $D_{\sigma_k}$ if $u $ minimizes the energy $E_{L, h}^\star$ among all functions in 
\begin{eqnarray}  H_{k, \overline{w}_*, u} := \bigg\{ v \in H_{k, \overline{w}_*} : v = u \h{3pt}\text{on $\p D_{\sigma_k}$}\bigg\}.
\label{min. problem S_L}
\end{eqnarray} \begin{lem} Recalling the energy density $e_{L, h}^\star$ defined in (\ref{E_L^star}), we have \begin{enumerate}\item[$\mathrm{(1).}$] For any $b_* > 0$, there are two constants  $\nu \in (0, 1)$ and $L_0 > 0$ depending on $b_*$ so that for any pair $\big(h, \overline{w}_*\big)$ and any solution $u$ to the Problem $S_{L, h, \overline{w}_*}$ on $D_{\sigma_k}$, if it satisfies  $L > L_0$, $\| u \|_{1, 2; D_{\sigma_k}}  \leq b_*$ and $E_{L, h}^\star\left[u\right] \leq 1$, then either one of the followings holds: \begin{eqnarray}
\mathrm{(i).} \int_{D_{1/4} }
 e_{L, h}^\star[\h{0.2pt}u\h{0.2pt}] \h{2pt}\leq\h{2pt}\dfrac{1}{16}; \h{20pt}\mathrm{(ii).} \int_{D_{1/8}}
 e_{L, h}^\star[\h{0.2pt}u\h{0.2pt}] \h{2pt}\leq\h{2pt} \nu  \int_{D_{1/4}}
 e_{L, h}^\star[\h{0.2pt}u\h{0.2pt}]; \label{two possible cases}\end{eqnarray}
\item[$\mathrm{(2).}$] There is a positive universal constant $c $ so that for any pair $\big(h, \overline{w}_*\big)$ and any solution $u$  of the Problem  $S_{L, h, \overline{w}_*}$ on $D_{\sigma_k}$, if it satisfies $E_{L, h}^\star[\h{1pt}u\h{1pt}] \leq 1$, then   \begin{align} \int_{D_{R}} e_{L, h}^\star\left[ u \right] \h{2pt}\leq\h{2pt}c\h{1pt}\big(1 + L \big) R, \h{20pt}\text{for all $R \in (0, 1/2\h{1pt}]$.}
\label{linear growth}
\end{align}\end{enumerate}
\label{decay for large L}
\end{lem}\noindent In order not to interrupt the current proof, this lemma is postponed to be proved in Appendix. Now we finish the proof of Proposition \ref{variant decay lemma} in the intermediate--scale regime. By Lemma \ref{strong H1 convergence in intermediate case}, $U^{sc, \bot}_\infty := \big(U_{\infty; 1}^{sc}, U_{\infty, 2}^{sc}\big)$ solves the Problem ($S_{L, h, \overline{w}_*}$) on $D_{\sigma_k}$. Letting $C_{np}$ be the optimal Poincar\'{e}--constant for Neumann--Poincar\'{e} inequality on $D_1$ with $p = 2$, we then have \begin{align}\label{upp bd for uinfty} \big\|\h{1pt}U_\infty^{sc, \bot}\h{1pt}\big\|_{1, 2; D_1}^2 = \int_{D_1} \big|\h{1pt}U_{\infty}^{sc, \bot}\h{1pt}\big|^2 + \big|\h{1pt}D_\xi U_\infty^{sc, \bot}\h{1pt}\big|^2 \h{1.5pt}\leq\h{1.5pt}\left( C_{np}^2 + 1 \right)\int_{D_1}   \big|\h{1pt}D_\xi U_\infty^{sc, \bot}\h{1pt}\big|^2 \h{1.5pt}\leq\h{1.5pt}C_{np}^2 + 1.
\end{align}Here we have used the condition that the average of $U_\infty^{sc, \bot}$ over $D_1$ equals $0$. Therefore, if  $u$ in Lemma \ref{decay for large L} equals $U_\infty^{sc, \bot}$, then we can take the constant $b_*$ in Lemma \ref{decay for large L} to be $\big(C_{np}^2 + 1 \big)^{1/2}$. The constants $\nu$ and $L_0$ are then universal constants depending only on $C_{np}$. \vspace{0.3pc}

Suppose that $L \leq L_0$. By item (2) in Lemma \ref{decay for large L}, it turns out \begin{align*}\int_{D_{R}} e_{L, h}^\star\left[ U_\infty^{sc, \bot} \right] \h{2pt}\leq\h{2pt}c\h{1pt}\big(1 + L_0 \big) R, \h{15pt}\text{for all $R \in \big(0, 1/2\h{0.5pt}\big]$.}
\end{align*}Applying this result together with (\ref{decay of 3rd component}) to the left--hand side of (\ref{con in interm scale}) yields \begin{align*} \frac{1}{4}\h{1.5pt}\leq\h{1.5pt}K \theta_0^2 + c\h{1pt}\big(1 + L_0\big) \theta_0.
\end{align*}Here $K$ and $L_0$ are all positive universal constants. In this case, we can easily find a universal $\theta_0$ small enough so that the above inequality fails. \vspace{0.3pc}

In the next, we assume that $L > L_0$. If $\mathrm{(i)}$ in (\ref{two possible cases}) is satisfied by $u = U^{sc, \bot}_\infty$, then we can apply it together with (\ref{decay of 3rd component}) to the left--hand side of (\ref{con in interm scale}). It follows  \begin{eqnarray*} \dfrac{1}{4} \leq K \theta_0^2 + \dfrac{1}{16}, \h{15pt}\text{for any $\theta_0 \in (0, 1/4)$.} \end{eqnarray*} In this case, we can also find a  universal $\theta_0$ small enough so that the above inequality fails. \vspace{0.3pc}

In the remaining arguments, we assume that $L > L_0$ and  $\mathrm{(i)}$ in (\ref{two possible cases}) fails. Therefore,  $\mathrm{(ii)}$ in (\ref{two possible cases})  is satisfied by $u = U^{sc, \bot}_\infty$. Using $\nu$ in $\mathrm{(ii)}$ of (\ref{two possible cases}), we can find a natural number $l_0$ so that \begin{align}\label{ineq of nu}\nu^{l_0} < 1/8.\end{align} For any $i \in \mathbb{N}\h{1pt}\cup\h{1pt}\big\{0\big\}$, we define  $U^{(i)}_\infty(\xi) := U^{sc, \bot}_\infty\left(2^{- i} \xi\right)$. It then holds \begin{eqnarray}\label{energy bound unif}\int_{D_{\sigma_k}}  e_{4^{-i}L, \h{0.4pt}h}^\star\big[  U^{(i)}_\infty\big] = \int_{D_{2^{- i}\sigma_k}}  e_{L, \h{0.4pt}h}^\star \big[ U^{sc, \bot}_\infty \big]   \h{1.5pt}\leq\h{1.5pt} 1, \h{15pt}i = 0, 1, ....
\end{eqnarray}Therefore, $U^{(i)}_\infty$ solves the Problem $S_{4^{-i}L, \h{0.4pt}h, \h{0.4pt}\overline{w}_*}$ on $D_{\sigma_k}$. Using this result and (\ref{linear growth}), we can induce\begin{align}\int_{D_{2^{-l_0}R_0}}  e_{L, \h{0.4pt}h}^\star\big[ U^{sc, \bot}_\infty\big] \h{2pt}\leq\h{2pt} c\h{0.5pt}\left(1 + 4^{- l_0} L\right) R_0  \h{2pt}\leq\h{2pt} c \h{1pt} \big(1 + L_0\big) R_0 \h{2pt}<\h{2pt} 1/8, \h{10pt}\text{if $4^{- l_0} L \leq L_0$.}
\label{less case}
\end{align}Here $R_0$ is small enough so that $c \h{1pt} \big(1 + L_0\big) R_0 < 1/8$. Now we assume $4^{- l_0} L > L_0$. In addition, we let $$Y^{(i)} := \dashint_{D_{\sigma_k}} U_\infty^{(i)}  \h{15pt}\text{and}\h{15pt}V_\infty^{(i)} := U_{\infty}^{(i)} - Y^{(i)}.$$ Then  $V_{\infty}^{(i)}$ solves the problem $S_{4^{-i}L, \h{1pt}h^{(i)}, \h{1pt}\overline{w}_*^{(i)}}$ on $D_{\sigma_k}$. Here $$h^{(i)} := h + \sum_{j = 1}^2  y_{*; j} \h{1pt}Y^{(i)}_j \h{15pt}\text{ and } \h{15pt} \overline{w}_*^{(i)} := \overline{w}_* - Y^{(i)}_2.$$  Moreover, by Neumann--Poincar\'{e} inequality, it holds $$\big\|\h{1pt}V_\infty^{(i)}\h{1pt}\big\|_{1, 2; D_{\sigma_k}}^2 = \int_{D_{\sigma_k}} \big|\h{1pt}V_{\infty}^{(i)}\h{1pt}\big|^2 + \big|\h{1pt}D_\xi V_\infty^{(i)}\h{1pt}\big|^2 \h{1.5pt}\leq\h{1.5pt}\left( C_{np}^2 + 1 \right)\int_{D_{\sigma_k}}   \big|\h{1pt}D_\xi U_\infty^{(i)}\h{1pt}\big|^2 \h{1.5pt}\leq\h{1.5pt}C_{np}^2 + 1.$$ Notice that the upper--bound in this estimate is identical with the upper--bound in (\ref{upp bd for uinfty}). Therefore, when we apply (ii) in (\ref{two possible cases}) to $V_\infty^{(i)}$, the constant $\nu$ is the same as the one that we have used in (\ref{ineq of nu}).  If $$\int_{D_{1/4}}e^\star_{4^{-  (l_0 - l)}L, \h{1pt}h^{(l_0 - l)}} \left[ V^{(l_0 - l)}_\infty \right] =  \int_{D_{1/4}}e^\star_{4^{-  (l_0 - l)}L, \h{0.4pt}h} \left[ U^{(l_0 - l)}_\infty \right] \h{1.5pt}\leq\h{1.5pt} \dfrac{1}{16}, \h{15pt}\text{for some $l \in \Big\{0, 1, ..., l_0\Big\}$,}$$ then it follows \begin{align*} \int_{D_{2^{- 2 - l_0 }}} e^\star_{L, \h{0.4pt}h} \left[ U^{sc, \bot}_\infty \right] \leq \dfrac{1}{16}.
\end{align*}Otherwise,  utilizing (ii) in (\ref{two possible cases}), (\ref{ineq of nu}) and (\ref{energy bound unif}), we obtain   \begin{align*}\int_{D_{1/8}}  e^\star_{4^{- l_0}L, \h{1pt}h^{(l_0)}} \left[  V^{(l_0)}_\infty \right]  & \h{2pt}\leq\h{2pt} \nu \int_{D_{1/4}}  e_{4^{- l_0} L, \h{1pt}h^{(l_0)}}^\star \left[ V^{(l_0)}_\infty \right]  \h{2pt}= \h{2pt} \nu \int_{D_{1/8}}  e_{ 4^{- (l_0 - 1)}L, \h{1pt}h^{(l_0 - 1)}}^\star \left[ V^{(l_0 - 1)}_\infty\right]  \\[2mm]
\h{2pt}&\h{2pt}\leq\h{2pt}\cdot\cdot\cdot \leq\h{2pt} \nu^{l_0} \int_{D_{1/8}}  e^\star_{L, \h{0.4pt}h} \left[ U^{sc, \bot}_\infty\right] \h{2pt}<\h{2pt} 1/8.
\end{align*}  In any case, it turns out \begin{eqnarray}\int_{D_{2^{- 3 - l_0}}}  e^\star_{L, \h{0.4pt}h} \left[ U^{sc, \bot}_\infty\right] < 1/8, \h{10pt}\text{provided that $4^{- l_0} L > L_0$.}
\label{large case}
\end{eqnarray}In light of (\ref{decay of 3rd component}) and (\ref{less case})--(\ref{large case}), we can find a positive universal constant $\theta_0 < \mathrm{min}\Big\{2^{- l_0} R_0,  \h{2pt}2^{- 3 - l_0} \Big\}$ so that (\ref{con in interm scale}) fails. The proof finishes.
\end{proof}
\subsection{Energy--decay estimate in large--scale regime}
In this section we suppose that  $a_n \big(\lambda_n \h{0.5pt} r_n\big)^2 \rightarrow \infty$ as $n \to \infty$. Firstly, we recall and introduce some notations.  The notations defined in Section 3.1 will also be used in the following arguments. Throughout the section, $\Gamma$ denotes the equator of $\mathbb{S}^2$ which is formed by all points in $\mathbb{S}^2$ with the third coordinate $0$. The set $\Gamma_b$ contains all points in $\Gamma$ with the second coordinate greater than or equaling  $b$.  Moreover, we use  $\Pi_{\h{0.5pt}\Gamma}$ to denote the shortest--distance projection to  $\Gamma$. \vspace{0.2pc}

Recall $s_n^2$  and $U_n$ defined in (\ref{s_n and y_n}) and (\ref{blow sequence, j2}). It turns out
\begin{eqnarray}  s_n^{-2} \int_{D_{1}}  \left(  \big| U_n \big|^2 - 1  \right)^2 \h{3pt}\lesssim\h{3pt} \left[a_n\h{0.5pt}\big(\lambda_n r_n\big)^2\h{1pt}\right]^{-1} \longrightarrow 0, \h{20pt}\text{as $n \rightarrow \infty$.}
\label{basic limit in large scale}
\end{eqnarray}With the above convergence, we have \begin{lem}
The following results hold up to a subsequence: \begin{enumerate}
\item[$\mathrm{(1)}$.] There is a $y_* \in \Gamma_b$ so that $y_n \rightarrow y_*$ as $n \to \infty$. The projection $\Pi_{\h{0.5pt}\Gamma}\big(y_n\big)$ is well--defined for suitably large $n$. In addition, the sequences $\big\{y_n\big\}$ and  $\big\{ \Pi_{\h{0.5pt}\Gamma}\big(y_n\big) \big\}$ satisfy \begin{eqnarray} \lim_{n \rightarrow \infty} \h{2pt} \dfrac{\Pi_{\h{0.5pt}\Gamma}\big(y_n\big) - y_n}{s_n} =  v_*, \h{10pt}\text{for some $v_* \in \mathbb{R}^3$;}
\label{pi_gamma - y_n}
\end{eqnarray}
\item[$\mathrm{(2).}$] Let $X_n$ be the point $\left( \big(1 - H_{a_n}^2 b^2\big)^{\frac{1}{2}}, H_{a_n} b, 0 \right)^\top$. If it satisfies  \begin{eqnarray} \mathrm{(i).}\h{2pt}\big[\h{0.3pt}X_n - \Pi_{\h{0.5pt}\Gamma}\big(y_n\big) \h{0.3pt}\big]_2 \geq 0, \h{3pt}\text{for all  $n$ suitably large,} \h{5pt}\text{or}\h{8pt}\mathrm{(ii).}\h{2pt} \liminf_{n \rightarrow \infty} \h{2pt} \left|  \dfrac{X_n - \Pi_{\h{0.5pt}\Gamma}\big(y_n\big) }{s_n} \right| \h{2pt}< \h{2pt}\infty,
\label{positivity of second component}
\end{eqnarray}then up to a subsequence, it satisfies \begin{align} \lim_{n \rightarrow \infty} \h{2pt} \dfrac{X_n - \Pi_{\h{0.5pt}\Gamma}\big(y_n\big)}{s_n} =  \gamma_* \h{0.5pt}t_*, \h{10pt}\text{for some $\gamma_* \in \mathbb{R}$. }
\label{pi_L convergence in large scale}
\end{align}Here $t_* = e^*_3 \times y_*$. The above limit induces \begin{align}\label{y* at boundary case}y_* = \left(\h{1pt} \big(1 - b^2\big)^{\frac{1}{2}}, b, 0 \right)^\top \h{10pt}\text{ and } \h{10pt} t_* = \left( - b, \big(1 - b^2\big)^{\frac{1}{2}}, 0\right)^\top.\end{align}
 In this case, $U^{sc}_{\infty; 2} \h{1pt}\geq\h{1pt}v_{*; 2} + \gamma_{*} \h{1pt} y_{*; 1}$ on $T$ in the sense of trace.
\end{enumerate}
\label{ratio convergence and limit boundary condition in large scale}
\end{lem}
\begin{proof}[\bf Proof] Firstly, we consider the location of $y_*$. Here $y_*$ is the limit of the sequence $\big\{ y_n \big\}$. Still by  (\ref{limit of term in potential}), we have $y_* \in \mathbb{S}^2$. In light of (\ref{y_*=(y_*, 1y_*,2,0)}), it turns out $y_* \in \Gamma_b$. For suitably large $n$, the fact that $y_{n; 3} = 0$ then induces\begin{align}\label{bou of diff of pi gamma yn} \big|\h{1pt}\Pi_{\h{0.5pt}\Gamma}\big(y_n\big)- y_n \h{1pt}\big| \h{2pt}\leq\h{2pt} \left|\h{1pt} \Pi_{\h{0.5pt}\mathbb{S}^2}\big(U_n\big) - y_n \h{1pt}\right| &\h{2pt}\leq\h{2pt}  \left|\h{1pt} \Pi_{\h{0.5pt}\mathbb{S}^2}\big(U_n\big) - U_n\h{1pt} \right|+  \big|\h{1pt}U_n - y_n \h{1pt}\big| \nonumber\\[2mm] &\h{2pt}\leq\h{2pt} \left|\h{1pt}1 - \big|\h{0.5pt}U_n\h{0.5pt}\big|^2 \h{1pt}\right| +\big|\h{1pt}U_n - y_n \h{1pt}\big| \h{30pt}\text{on $D^+_1$.}  \end{align}The projection $\Pi_{\h{0.5pt}\mathbb{S}^2}\big(U_n\big)$ in (\ref{bou of diff of pi gamma yn}) is well--defined since $U_{n; 1} > 0 $ on $D_1^+$. See item (1) in Remark \ref{rmk on sign and gamma conv}. Integrating both sides of (\ref{bou of diff of pi gamma yn}) over $D_1^+$, by $D_1^+ \subset D_1$ and Poincar\'{e}'s inequality, we obtain \begin{eqnarray*}\big|\h{1pt}\Pi_{\h{0.5pt}\Gamma}\big(y_n\big)- y_n \h{1pt}\big|^2 \h{2pt}\lesssim\h{2pt} \int_{D_1} \big|\h{1pt}D_\xi U_n\h{1pt}\big|^2 + \int_{D_1} \left( 1 - \big|\h{1pt}U_n\h{1pt}\big|^2 \right)^2 \lesssim s_n^2.
\end{eqnarray*}Here we also have used (\ref{basic limit in large scale}) and the definition of $s_n$ in (\ref{s_n and y_n}).  The limit in (\ref{pi_gamma - y_n}) then follows by the last estimate.  Result (1) in the lemma is proved. \vspace{0.2pc}

If  (i) in (\ref{positivity of second component}) holds, then together with the boundary condition $U_{n; 2} \geq H_{a_n}\h{0.5pt}b$ on $T$, it turns out
\begin{eqnarray} \left[ \dfrac{U_n  - \Pi_{\h{0.5pt}\Gamma}\big(y_n\big)}{s_n} \right]_2 =  \left[\dfrac{U_n  - X_n }{s_n}\right]_2 +  \left[\dfrac{X_n  - \Pi_{\h{0.5pt}\Gamma}\big(y_n\big)}{s_n} \right]_2 \geq \left[\dfrac{X_n  - \Pi_{\h{0.5pt}\Gamma}\big(y_n\big)}{s_n} \right]_2 \geq 0 \h{20pt}\text{on $T$.} \label{result 2 lower bound}
\end{eqnarray}By trace theorem, $U^{sc}_n$ converges strongly in $L^2\big(T\big)$ to $U^{sc}_\infty$. This convergence together with (\ref{pi_gamma - y_n}) infers the almost everywhere convergence on $T$ of the quantity on the most left--hand side of (\ref{result 2 lower bound}). Hence, by (\ref{result 2 lower bound}), we have, up to a subsequence, that \begin{align}\label{con secon cop} \lim_{n \to \infty}\left[\dfrac{X_n  - \Pi_{\h{0.5pt}\Gamma}\big(y_n\big)}{s_n} \right]_2 = d_*, \h{15pt}\text{for some $d_* \geq 0$.}
\end{align} (\ref{y* at boundary case}) then follows by this limit. Now we write $$X_n = \big(\cos \alpha_n', \sin \alpha_n', 0\big)^\top \h{10pt}\text{ and } \h{10pt}  \Pi_{\h{0.5pt}\Gamma}\big(y_n\big) = \big(\cos \beta_n', \sin \beta_n', 0\big)^\top$$ with $\alpha_n'$ and $ \beta_n' $ converging to $\arcsin b $ as $n \to \infty$. It then follows  by (\ref{con secon cop}) that \begin{align*} \lim_{n \to \infty} \dfrac{\sin \alpha_n'  - \sin \beta_n'}{s_n}  = d_*,  \end{align*}which furthermore implies \begin{align*}\lim_{n \to \infty} \left[\dfrac{X_n  - \Pi_{\h{0.5pt}\Gamma}\big(y_n\big)}{s_n} \right]_1 =  \lim_{n \to \infty}\dfrac{\cos \alpha_n'  - \cos \beta_n'}{s_n} = - \dfrac{b}{\sqrt{1 - b^2}}\h{2pt}d_*.
\end{align*}(\ref{pi_L convergence in large scale}) holds by the above limit and (\ref{con secon cop}). \vspace{0.2pc}

If (ii) in (\ref{positivity of second component}) is satisfied, then we immediately have   (\ref{pi_L convergence in large scale}) up to a subsequence. \vspace{0.1pc}

In the end, we decompose  $U^{sc}_n$ as follows: \begin{eqnarray} U^{sc}_n = \dfrac{U_n - X_n}{s_n}  + \dfrac{X_n - \Pi_{\h{0.5pt}\Gamma}\big(y_n\big)}{s_n} + \dfrac{\Pi_{\h{0.5pt}\Gamma}\big(y_n\big) - y_n}{s_n}.
\label{decomposition of u_n in large scale case with free boundary}
\end{eqnarray}Utilizing (\ref{pi_gamma - y_n}), (\ref{pi_L convergence in large scale}) and the boundary condition $U_{n; 2} \geq H_{a_n}\h{0.5pt}b$ on $T$, we can take $n \to \infty$ on both sides of (\ref{decomposition of u_n in large scale case with free boundary}) and obtain $U^{sc}_{\infty; 2} \geq v_{*; 2} + \gamma_*\h{0.5pt}y_{*; 1}$ for almost every point on $T$. The proof is completed.
\end{proof}
 Before proceeding, let us consider the strict positivity of $U_{n; 1}$ in $D_1$. By showing $U_{n; 1} > 0 $ in $D_1$, we can define $\Pi_{\h{0.5pt}\mathbb{S}^2}\big(U_n\big)$ for every point in $D_1$. \begin{lem}For all $n \in \mathbb{N}$, we have $U_{n; 1} > 0$ in $D_1$.\label{sign of Un1}
\end{lem}\begin{proof}[\bf Proof] By the definition of $U_n$ in (\ref{blow sequence, j2}), the problem is reduced to proving $u_{n; 1} = \big[ u_{a_n, b}^+\big]_1 > 0$ in $\mathbb{D}$. Here $\mathbb{D}$ is defined in the item (1) of Theorem \ref{results for limiting case}. In light of the item (1) in Remark \ref{rmk on sign and gamma conv}, it follows $u_{n; 1} > 0$ in $\mathbb{D}^+$. We are left to show $u_{n; 1} > 0$ on $\mathbb{T} := \big\{\left(\rho, 0\right) : 0 < \rho < 1\big\}$. Since $\mathscr{L}\big[ u_n\big]$ minimizes the $\mathcal{E}_{a_n, \mu}$--energy in $\mathscr{F}^+_{a_n, b}$ (see Proposition \ref{existence of mini finite}), the first component $u_{n; 1}$ solves weakly the following equation: \begin{align}\label{eqn of un1}D \cdot \big( \rho \h{0.5pt} D \h{0.5pt} u_{n; 1}\big)   -  \mathcal{C} (\rho, u_n ) \rho\h{0.5pt} u_{n; 1} = - \dfrac{3 \sqrt{6}}{4} \h{1pt} \mu \h{1pt}\rho\h{1pt}u_{n; 3}^2 \h{15pt}\text{in $\mathbb{D}$.}
\end{align}Here $D = \big(\p_\rho, \p_z\big)$ and $$\mathcal{C}(\rho, u_n) := 4 \rho^{-2} + 3 \sqrt{2}\h{1pt}\mu\h{1pt}u_{n; 2} + a_n\h{1pt}\mu\h{1pt}\left( | u_n |^2 - 1 \right).$$ In light of the  item (2) in Remark \ref{rmk on sign and gamma conv}, we  have $u_{n; 1} \in C^{1, \alpha}\big(\mathbb{D}\big)$, for any $\alpha \in (0, 1)$. Here we have used Theorem 3.13 in \cite{HL11}. Moreover, $u_{n, 1}$ is smooth in $\mathbb{D}^+$. Suppose that there is a $\xi_0 \in \mathbb{T}$ so that $u_{n; 1}\left(\xi_0\right) = 0$. We then can find an open disk, denoted by $\mathscr{D}$, so that $\mathscr{D} \subset \mathbb{D}^+$ and $\p \mathscr{D}$ touches the set $\mathbb{T}$ at $\xi_0$. Note that $u_{n; 1} \geq 0$ in $\mathbb{D}$. By (\ref{eqn of un1}), it turns out  \begin{align*}\p_\rho^2 u_{n; 1} + \p_z^2 u_{n; 1} + \rho^{-1} \p_\rho u_{n; 1} - \mathcal{C}^+(\rho, u_n) u_{n; 1} \leq 0 \h{15pt}\text{in $\mathscr{D}$.}\end{align*}Here $\mathcal{C}^+(\rho, u_n)$ is the positive part of $\mathcal{C}(\rho, u_n)$. Note that $u_{n; 1} > 0$ in $\mathscr{D}$ and $u_{n; 1}\left(\xi_0\right) = 0$. We can then apply Hopf's lemma to get $$\dfrac{\p u_{n; 1}}{\p z} \Big|_{\xi_0} > 0.$$ However, this is impossible since by the $C^{1, \alpha}$--regularity of $u_{n; 1}$ in $\mathbb{D}$ and the even symmetry of $u_{n; 1}$ with respect to the $z$--variable, it must satisfy $$\dfrac{\p u_{n; 1}}{\p z} \Big|_{\xi_0} = 0.$$ The proof is completed.
\end{proof}
In the next, we are  concerned about the image  of the limiting map $U^{sc}_\infty$. \begin{lem}The following results hold for the limiting map $U^{sc}_\infty$. \begin{enumerate}\item[$\mathrm{(1).}$] The image of $U^{sc}_\infty$ lies in  $v_* + \mathrm{Tan}_{y_*}\mathbb{S}^2$ for almost all points in $D_1$;
\item[$\mathrm{(2).}$]  The image of $U^{sc}_\infty$ lies in   $v_* + \mathrm{Tan}_{y_*} \Gamma$ on $T$ in the sense of trace;
\item[$\mathrm{(3).}$] If (\ref{pi_L convergence in large scale}) holds, then on $T$, it satisfies $U^{sc}_\infty = v_* + w \h{1pt}t_*$ with $w \geq \gamma_*$ in the sense of trace. \end{enumerate}
In the results listed above,  $\gamma_*$, $v_*$ and $t_*$ are given in Lemma \ref{ratio convergence and limit boundary condition in large scale}. $\mathrm{Tan}_{y_*} \mathbb{S}^2$ contains all vectors in $\mathbb{R}^3$ which are perpendicular to $y_*$. $\mathrm{Tan}_{y_*} \Gamma$ contains all vectors in $\mathrm{Tan}_{y_*}\mathbb{S}^2$ with the third coordinate $0$.
\label{images of u_infty on boundary and interior}
\end{lem}
\begin{proof}[\bf Proof] The convergences and pointwise relationships in the proof are understood in the sense of almost everywhere, except otherwise stated.  Note that $\Pi_{\h{0.5pt}\mathbb{S}^2}\big(U_n\big)$ is well--defined for all points in $D_1^+$ by $U_{n; 1} > 0$ in $D_1^+$.  We then  decompose $U^{sc}_n$ as follows: \begin{eqnarray}
 U^{sc}_n = \dfrac{U_n - \Pi_{\h{0.5pt}\mathbb{S}^2}\big(U_n\big)}{s_n} + \dfrac{\Pi_{\h{0.5pt}\mathbb{S}^2}\big(U_n\big) - \Pi_{\h{0.5pt}\Gamma}\big(y_n\big)}{s_n}+ \dfrac{\Pi_{\h{0.5pt}\Gamma}\big(y_n\big) - y_n}{s_n} \h{15pt}\text{in $D^+_1$.}\label{deco of Uscn}
\end{eqnarray}
Due to (\ref{basic limit in large scale}), the first term on the right--hand side above converges to $0$ pointwisely in  $D_1^+$ as $n \to \infty$. In light of (\ref{pi_gamma - y_n}) and the pointwise convergence of $U_n^{sc}$, the second term on the right--hand side of (\ref{deco of Uscn}) also converges pointwisely as $n \to \infty$. Since $\Pi_{\h{0.5pt}\mathbb{S}^2}\big(U_n\big)$ converges to $y_*$ pointwisely on $D_1^+$ and by (\ref{pi_gamma - y_n}), $\Pi_{\h{0.5pt}\Gamma}\big(y_n\big)$ converges to $y_*$ as well as $n \to \infty$, then on $D_1^+$, the limit of the second term on the right--hand side of (\ref{deco of Uscn}) takes values in $\mathrm{Tan}_{y_*}\mathbb{S}^2$ as $n \to \infty$. Therefore, by (\ref{basic limit in large scale}), (\ref{pi_gamma - y_n}) and (\ref{deco of Uscn}), it follows $U^{sc}_\infty \in v_* + \mathrm{Tan}_{y_*}\mathbb{S}^2$ pointwisely in $D^+_1$. \vspace{0.2pc}

Utilizing trace theorem, we also have $U^{sc}_\infty \in v_* + \mathrm{Tan}_{y_*}\mathbb{S}^2$ pointwisely on $T$. Since \begin{eqnarray*}
 U^{sc}_n = \dfrac{U_n - \Pi_{\h{0.5pt}\Gamma}\big(y_n\big)}{s_n} + \dfrac{\Pi_{\h{0.5pt}\Gamma}\big(y_n\big) - y_n}{s_n} \h{15pt}\text{on $T$,} 
\end{eqnarray*}then  the limit of the first term on the right--hand side above lies in $\mathrm{Tan}_{y_*}\mathbb{S}^2$ pointwisely on $T$ as $n \to \infty$. In fact, this limit must be in $\mathrm{Tan}_{y_*} \Gamma$ pointwisely on $T$ in that $\left[ \Pi_{\h{0.5pt}\Gamma}\big(y_n\big)\right]_3 = 0$ and $U_{n; 3} = 0$ on $T$ in the sense of trace. If (\ref{pi_L convergence in large scale}) holds, by Lemma \ref{ratio convergence and limit boundary condition in large scale}, we have $U^{sc}_{\infty; 2} \geq v_{*; 2} + \gamma_*\h{1pt}y_{*; 1}$  on $T$ in the sense of trace. Letting $U^{sc}_\infty = v_* + w\h{1pt}t_*$ on $T$, we then obtain $w \geq \gamma_*$ on $T$ in the sense of trace.  Note that $y_{*; 1} = \big(1 - b^2 \big)^{\frac{1}{2}} \neq 0$ if the limit in (\ref{pi_L convergence in large scale}) holds.
\end{proof}Due to Lemma \ref{images of u_infty on boundary and interior} and the limit $a_n \big(\lambda_n\h{1pt}r_n\big)^2 \to \infty$ , we have the following modification of Lemma \ref{uniform convergence of u_n in S_sigma_k}:
 \begin{lem}
There exist an increasing positive sequence $\big\{\sigma_k\big\}$ tending to $1$ as $k \to \infty$, a sequence of positive numbers $\big\{b_k\big\}$ and a subsequence of $\big\{U_n\big\}$, still denoted by $\big\{U_n\big\}$, so that for any $k$, \begin{enumerate}
\item[$\mathrm{(1).}$] The mappings $U_n$, $U_n^{sc}$ and $U_\infty^{sc}$ satisfy the uniform boundedness given below: \begin{eqnarray*}\sup_{n \h{1pt}\in\h{1pt}\mathbb{N}\h{1pt}\cup\h{1pt}  \{ \infty \}}  \big\lVert  U^{sc}_n \big\rVert_{\infty;\h{1pt}\p D_{\sigma_k}} +  \int_{\p D_{\sigma_k}} \big|\h{1pt}D_\xi U^{sc}_\infty\h{1pt}\big|^2 +    \sup_n  \int_{\p D_{\sigma_k}} \big|\h{1pt}D_\xi U^{sc}_n\h{1pt}\big|^2 + a_n \left( \dfrac{\lambda_n\h{1pt}r_n}{s_n}\right)^2 \left[ \h{1pt}\big|\h{1pt}U_n\h{1pt}\big|^2 - 1 \right]^2 \leq b_k;\end{eqnarray*}
\item[$\mathrm{(2).}$] The sequence  $\big\{ U^{sc}_n \big\}$ converges to $U^{sc}_\infty$ in $C^0\big( \p D_{\sigma_k}\big)$ as $n \rightarrow \infty$;
\item[$\mathrm{(3).}$] The second component of $U_n$ satisfies $U_{n; 2} \geq H_{a_n} b$ at $\big(\pm \sigma_k, 0 \h{1pt}\big)$;
\item[$\mathrm{(4).}$] The third component  of $U^{sc}_\infty$ satisfies $ U^{sc}_{\infty; 3} = 0$ at $\big( \pm \sigma_k, 0\h{1pt}\big)$;
\item[$\mathrm{(5).}$] If (\ref{pi_L convergence in large scale}) holds, then $U^{sc}_\infty = v_* + w \h{1pt}t_*$ on $T_{\sigma_k}$ with $w \geq \gamma_*$ at  $\big( \pm \sigma_k, 0\h{1pt}\big)$;
\item[$\mathrm{(6).}$] The following uniform boundedness holds: \begin{eqnarray*} \sup_{n\h{1pt}\in\h{1pt}\mathbb{N}} \h{2pt}Q_n^2  \h{2pt}\leq\h{2pt}b_k  \h{ 10pt}\text{at $\big(\pm \sigma_k, 0\h{1.5pt}\big)$. Here $Q_n  := \sqrt{a_n} \h{2pt}\dfrac{\lambda_n\h{1pt}r_n}{s_n} \Big|\h{1pt}\big|\h{1pt}U_n \h{1pt}\big| - 1 \Big|^2.$}
\end{eqnarray*}
\end{enumerate}
\label{uniform convergence of u_n in S_sigma_k in large scale}
\end{lem}
\begin{proof}[\bf Proof] We only consider the item (6) in the lemma.  Firstly, we note that
\begin{eqnarray*}   \int_{D^+_1} Q_n^2  \h{2pt}\leq\h{2pt} a_n \left( \dfrac{\lambda_n\h{1pt}r_n}{s_n}\right)^2  \int_{D^+_1} \left|\h{1pt}\big|\h{1pt}U_n \h{1pt}\big|^2 - 1 \right|^2.
\end{eqnarray*}By H\"{o}lder's inequality, it holds \begin{eqnarray*} \int_{D_1^+}\big|\h{1pt} D_\xi  \h{1pt}Q_n \h{1pt}\big| \h{2pt}\lesssim\h{2pt}\left( a_n \left(\dfrac{\lambda_n\h{1pt}r_n}{s_n}\right)^2 \int_{D_1^+} \Big|\h{1pt}\big|\h{1pt}U_n\h{1pt}\big| - 1 \Big|^2 \right)^{1/2}\left( \int_{D_1^+}  \big|\h{1pt}D_\xi U_n\h{1pt}\big|^2  \right)^{1/2}.
\end{eqnarray*}According to (\ref{basic limit in large scale}) and the uniform boundedness of $U_n$ in $H^1\big(D_1\big)$, the last two estimates induce the uniform boundedness of the sequence $\big\{ Q_n \big\}$  in $W^{1, 1}\big(D_1^+\big)$. Since $W^{1,1}\big(D_1^+\big)$ is embedded into $L^1\big(T\big)$ continuously, the sequence  $\big\{ Q_n \big\}$ is uniformly bounded in $L^{1}\big(T\big)$. By Fatou's lemma,  it turns out
\begin{eqnarray*}
\int_{T }\liminf_{n \rightarrow \infty} Q_n
\h{2pt}\leq\h{2pt}\liminf_{n \rightarrow \infty}\int_{T } Q_n \h{2pt}<\h{2pt}\infty.
\end{eqnarray*}We therefore can assume that the value of $\displaystyle  \liminf_{n \rightarrow \infty} \h{1.5pt} Q_n$ is finite at $(\pm \sigma_k, 0)$. (6) in the lemma then follows.\end{proof}

Noticing Lemma \ref{images of u_infty on boundary and interior} and using the $\sigma_k$ obtained in Lemma \ref{uniform convergence of u_n in S_sigma_k in large scale}, we define  \begin{align}  &N_{\h{0.5pt}k} := \left\{ u \in H^1\left( D_{\sigma_k} ; \h{2pt}v_* + \mathrm{Tan}_{y_*} \mathbb{S}^2 \h{1pt}\right) \left| \begin{array}{lcl} u = U^{sc}_{\infty} \h{4pt}\text{on $\p D_{\sigma_k}$;} \h{16pt} u \in v_* + \mathrm{Tan}_{y_*} \Gamma \h{4pt}\text{on $T_{\sigma_k}$;} \\[2mm]
\text{$u_1$ and $u_2$ are even and $u_3$ is odd with respect to $\xi_2$--variable} \end{array} \right. \right\};\nonumber\\[2mm]
& \overline{N}_{\h{0.5pt}k} :=  \Big\{ u \in N_{\h{0.5pt}k} : u = v_* + w \h{1pt}t_*  \h{2pt}\text{on $T_{\sigma_k}$} \h{4pt} \text{with $w \geq \gamma_*$} \h{4pt}\text{on $T_{\sigma_k}$} \Big\}.
\label{N_k and bar N_k}
\end{align}
In the definition of $\overline{N}_{\h{0.5pt}k}$, the notions $t_*$ and $\gamma_*$ are given in (\ref{y* at boundary case}) and (\ref{pi_L convergence in large scale}), respectively. \vspace{0.2pc}

Now we show our main result in this section.
\begin{lem}
Fix a natural number $k$. If the following two conditions are satisfied 
\begin{eqnarray}
\mathrm{(i).}\h{3pt} \Big[ X_n - \Pi_{\h{0.5pt}\Gamma}\big(y_n\big)\Big]_2 < 0 \h{10pt}\text{for all $n$;}\h{20pt}   \mathrm{(ii).}\h{3pt} \liminf_{n \rightarrow \infty} \h{2pt} \left|  \dfrac{X_n - \Pi_{\h{0.5pt}\Gamma}\big(y_n\big) }{s_n} \right| \h{2pt}= \h{2pt}\infty,
\label{infinite case in large scale}
\end{eqnarray}
then $U^{sc}_\infty$ minimizes the Dirichlet energy over $N_k$. If one of the two conditions in (\ref{infinite case in large scale}) fails, then  (\ref{pi_L convergence in large scale}) holds up to a subsequence. In this case, $U^{sc}_\infty$ minimizes the Dirichlet energy over $\overline{N}_{\h{0.5pt}k}$. In all cases, $U^{sc}_n$ converges to $U^{sc}_{\infty}$ strongly in $H_{\mathrm{loc}}^1\big(D_1\big)$ as $n \to \infty$. Moreover, it satisfies$$a_n \h{1pt}\left( \dfrac{\lambda_n r_n}{s_n}\right)^2\int_{D_{\sigma_k}} \Big(\h{1pt}\big|\h{1pt}U_n\h{1pt}\big|^2 - 1 \Big)^2\left(1 + \dfrac{\lambda_n \h{1pt}r_n}{\rho_{x_n}} \h{1pt}\xi_1\right)  \longrightarrow 0,\h{10pt} \text{as $n \to \infty$.}
$$
\label{strong H1 convergence in large case}
\end{lem}
\begin{proof}[\bf Proof] We divide the proof into five steps. \vspace{0.3pc}\\
\textbf{Step 1. Construction of comparison map}\vspace{0.5pc}\\
Suppose that $v$ is an arbitrary map in $N_{\h{0.5pt}k}$. Then for any $R > 0$, we define
\begin{eqnarray}
 F^{\h{0.5pt}l}_{n, R}\h{0.5pt} [\h{0.5pt}v\h{0.5pt}] :=\left\{ \begin{array}{lcl} \Pi_{\h{0.5pt}\Gamma}\big(y_n\big)+R\h{1pt}s_n \h{1pt}\dfrac{v - Z_*}{|\h{1pt}v - Z_*\h{1pt}|\vee R}, \h{30pt}&&\text{if (i) and (ii) in (\ref{infinite case in large scale}) hold;}\vspace{0.8pc}\\
 \h{18pt}X_n + R\h{1pt}s_n \h{1pt}\dfrac{v - Z_* }{|\h{1pt}v - Z_*\h{1pt}|\vee R}, &&\text{if (\ref{pi_L convergence in large scale}) holds.}
\end{array}\right.
\label{definition of F_n^l}
\end{eqnarray}
Here  $Z_* = v_*$ if (i) and (ii) in (\ref{infinite case in large scale}) hold. If (\ref{pi_L convergence in large scale}) holds, then $Z_* = v_* + \gamma_*\h{1pt}t_*$. Since $Z_{*; 3} = X_{n; 3} = \big[\h{1pt} \Pi_{\h{0.5pt}\Gamma}\big(y_n\big)\h{1pt}\big]_3 =  0$,  the first two components of $F_{n, R}^l [ v]$ are even and the third component of $F_{n, R}^l [ v]$ is odd with respect to the $\xi_2$--variable. Now we let $$\Gamma_n :=  \Big\{u \in \Gamma : u_2 \geq H_{a_n}b\Big\}$$ and define, for any $\xi \in \p D_{\sigma_k}$, the mapping $J_n\big(\xi\big)$ as follows:
\begin{align}\label{definition of J_n} J_n(\xi) :=
\left\{ \begin{array}{lcl}
\Pi_{\h{0.5pt}\mathbb{S}^2} \big[\h{1pt}U_n\big(\xi\big) \h{1pt}\big] \h{40pt} &&\text{for Case A: \h{2pt} $\Pi_{\h{0.5pt}\mathbb{S}^2} \big[\h{1pt} U_n\big( \pm \sigma_k, 0 \h{1pt}\big) \h{1pt} \big] \in \Gamma_n$;}\\[4mm]
\mathrm{Rot}_{\alpha_1} \Pi_{\h{0.5pt}\mathbb{S}^2} \big[\h{1pt}U_n\big(\xi\big) \h{1pt}\big] &&\text{for Case B:  \h{2pt} $\Pi_{\h{0.5pt}\mathbb{S}^2} \big[\h{1pt} U_n\big(\sigma_k, 0 \h{1pt}\big) \h{1pt} \big] \not\in \Gamma_n$ } \\[4mm]
&&\h{15pt} \text{and $\left[\h{1pt} \Pi_{\h{0.5pt}\mathbb{S}^2} \big[\h{1pt} U_n\big(\sigma_k, 0 \h{1pt}\big) \h{1pt} \big] \h{1pt}\right]_2 \h{2pt}\leq\h{2pt} \left[\h{1pt} \Pi_{\h{0.5pt}\mathbb{S}^2} \big[\h{1pt} U_n\big(- \sigma_k, 0 \h{1pt}\big) \h{1pt} \big] \h{1pt}\right]_2$;}\\[4mm]
\mathrm{Rot}_{\beta_1} \Pi_{\h{0.5pt}\mathbb{S}^2} \big[\h{1pt}U_n\big(\xi\big) \h{1pt}\big] &&\text{for Case C: \h{2pt} $\Pi_{\h{0.5pt}\mathbb{S}^2} \big[\h{1pt} U_n\big(- \sigma_k, 0 \h{1pt}\big) \h{1pt} \big] \not\in \Gamma_n$ } \\[4mm]
&&\h{15pt} \text{and $\left[\h{1pt} \Pi_{\h{0.5pt}\mathbb{S}^2} \big[\h{1pt} U_n\big(- \sigma_k, 0 \h{1pt}\big) \h{1pt} \big] \h{1pt}\right]_2 \h{2pt}\leq\h{2pt} \left[\h{1pt} \Pi_{\h{0.5pt}\mathbb{S}^2} \big[\h{1pt} U_n\big( \sigma_k, 0 \h{1pt}\big) \h{1pt} \big] \h{1pt}\right]_2$.}
\end{array}\right.
\end{align}
Before proceeding, we explain with more details the definition of $J_n$ above. Firstly, $\Pi_{\h{0.5pt}\mathbb{S}^2}\big(U_n\big(\pm \sigma_k, 0\big)\big)$ is well--defined by Lemma \ref{sign of Un1}. Due to (1) in Lemma \ref{uniform convergence of u_n in S_sigma_k in large scale}, $U_{n; 3}$ is absolutely continuous on $\p D_{\sigma_k}$. Together with the odd symmetry of $U_{n; 3}$ with respect to the $\xi_2$--variable, we get $U_{n; 3}\big(\pm \sigma_k, 0\big) = 0$. This result infers  $\Pi_{\h{0.5pt}\mathbb{S}^2}\big(U_n\big(\pm \sigma_k, 0\big)\big) \in \Gamma$. Note that we have covered all the possibilities for the locations of $\Pi_{\h{0.5pt}\mathbb{S}^2}\big(U_n\big(\pm \sigma_k, 0\big)\big)$ when we define $J_n$ in (\ref{definition of J_n}). In fact, if $\Pi_{\h{0.5pt}\mathbb{S}^2} \big[\h{1pt} U_n\big(\sigma_k, 0 \h{1pt}\big) \h{1pt} \big] \not\in \Gamma_n$ and the inequality in Case B of (\ref{definition of J_n}) is not held, then we have \begin{align*}\left[\h{1pt} \Pi_{\h{0.5pt}\mathbb{S}^2} \big[\h{1pt} U_n\big(- \sigma_k, 0 \h{1pt}\big) \h{1pt} \big] \h{1pt}\right]_2 \h{2pt}<\h{2pt} \left[\h{1pt} \Pi_{\h{0.5pt}\mathbb{S}^2} \big[\h{1pt} U_n\big(\sigma_k, 0 \h{1pt}\big) \h{1pt} \big] \h{1pt}\right]_2 \h{2pt}<\h{2pt}H_{a_n}b,
\end{align*}which infers the conditions in Case C of (\ref{definition of J_n}). Similarly if $\Pi_{\h{0.5pt}\mathbb{S}^2} \big[\h{1pt} U_n\big(-\sigma_k, 0 \h{1pt}\big) \h{1pt} \big] \not\in \Gamma_n$ and the inequality in Case C of (\ref{definition of J_n}) is not satisfied, then we can infer the conditions in Case B of (\ref{definition of J_n}).  As for $\alpha_1$ in (\ref{definition of J_n}), it lies in $(0, \pi)$ and is the angle between $\Pi_{\h{0.5pt}\mathbb{S}^2} \big[\h{1pt} U_n\big(\sigma_k, 0 \h{1pt}\big) \h{1pt}\big]$ and $X_n$. $\beta_1$ lies in $(0, \pi)$ and is the angle between $\Pi_{\h{0.5pt}\mathbb{S}^2} \big[\h{1pt} U_n\big(- \sigma_k, 0 \h{1pt}\big)\h{1pt}\big] $ and $X_n$. If Case B and Case C are satisfied simultaneously, then $\alpha_1 = \beta_1$ since now $ \Pi_{\h{0.5pt}\mathbb{S}^2} \big[\h{1pt} U_n\big(- \sigma_k, 0 \h{1pt}\big) \h{1pt} \big] = \Pi_{\h{0.5pt}\mathbb{S}^2} \big[\h{1pt} U_n\big( \sigma_k, 0 \h{1pt}\big) \h{1pt} \big]$. For any angle $\alpha$, $\mathrm{Rot}_{\alpha}$ denotes the following rotation matrix: \begin{eqnarray*} \left(\begin{array}{lcl} \cos \alpha &- \sin \alpha &0 \\
 \sin \alpha &\h{8pt}\cos \alpha &0 \\
\h{8pt} 0 &\h{2pt}0 &1 \end{array}\right).
\end{eqnarray*}We  emphasize that the rotation matrices $\mathrm{Rot}_{\alpha_1}$ and $\mathrm{Rot}_{\beta_1}$ in (\ref{definition of J_n}) are used in order to obtain \begin{align}\label{low boun of Jn pm sigmak}\big[ J_n\left(\pm \sigma_k, 0\right)\big]_2 \geq H_{a_n}b.\end{align} In fact, by (1) in Lemma \ref{uniform convergence of u_n in S_sigma_k in large scale}, we have $U_n$ converges to $y_*$ in $C^0\big(\p D_{\sigma_k}\big)$. If Case B in (\ref{definition of J_n}) happens, then \begin{align*}\left[\h{1pt}\Pi_{\h{0.5pt}\mathbb{S}^2} \big[\h{1pt} U_n\big(\sigma_k, 0 \h{1pt}\big) \h{1pt} \big]\h{1pt}\right]_2 < H_{a_n} b.
\end{align*}Taking $n \to \infty$ on both sides above yields $y_{*; 2} \leq b$. Noticing (\ref{y_*=(y_*, 1y_*,2,0)}), we then obtain $y_* = \left(\big(1 - b^2\big)^{\frac{1}{2}}, b, 0\right)^{\top}$. This result on $y_*$ and the inequality in Case B of (\ref{definition of J_n}) induce \begin{align*}H_{a_n} b \h{2pt}=\h{2pt} \left[\h{1pt} \mathrm{Rot}_{\alpha_1} \Pi_{\h{0.5pt}\mathbb{S}^2} \big[\h{1pt} U_n\big(\sigma_k, 0 \h{1pt}\big) \h{1pt} \big] \h{1pt}\right]_2 \h{2pt}\leq\h{2pt} \left[\h{1pt} \mathrm{Rot}_{\alpha_1}\Pi_{\h{0.5pt}\mathbb{S}^2} \big[\h{1pt} U_n\big(- \sigma_k, 0 \h{1pt}\big) \h{1pt} \big] \h{1pt}\right]_2.
\end{align*}(\ref{low boun of Jn pm sigmak}) follows if Case B in (\ref{definition of J_n}) holds. Similar arguments can also be applied to Case C in (\ref{definition of J_n}). As for the symmetry of $J_n$ with respect to the $\xi_2$--variable, we firstly note that $U_{n; 1}$ and $U_{n; 2}$ are even and $U_{n; 3}$ is odd with respect to $\xi_2$--variable. Moreover, (1) in Lemma \ref{uniform convergence of u_n in S_sigma_k in large scale} infers the absolutly continuity of $U_n$ on $\p D_{\sigma_k}$. Therefore, $U_{n; 1}$ and $U_{n; 2}$ are even and $U_{n; 3}$ is odd with respect to the $\xi_2$--variable when they are restricted on $\p D_{\sigma_k}$. Since the rotation matrix $\mathrm{Rot}_\alpha$ is planar, $\mathrm{Rot}_\alpha \Pi_{\h{0.5pt}\mathbb{S}^2}\big[ U_n\big(\xi\big)\big]$ has the same third component as     $\Pi_{\h{0.5pt}\mathbb{S}^2}\big[ U_n\big(\xi\big)\big]$. All these arguments induce that the first two components of $J_n\big(\xi\big)$ are even and the third component of $J_n\big(\xi\big)$ is odd with respect to the $\xi_2$--variable. Particularly, \begin{align}\label{0 of Jn3} \big[ \h{1pt}J_n\big(\pm \sigma_k, 0\big)\h{1pt}\big]_3 = U_{n; 3}\big(\pm \sigma_k, 0 \big) = 0.
\end{align} 

With $F^l_{n, R}$ defined in (\ref{definition of F_n^l}), now we  fix  a $s \in (0, 1)$ and introduce
\begin{eqnarray}h_{n, s, R}\big(\xi\big) :=\left\{
\begin{aligned}
&  F^{\h{0.5pt}l}_{n, R}\h{1pt}[\h{1pt}v\h{1pt}]\left(\dfrac{\xi}{1-s}\right) \quad &\text{if }& \xi \in D_{(1-s)\h{1pt}\sigma_k}; \\[4mm]
&\dfrac{\sigma_k-|\h{1pt}\xi\h{1pt}|}{s\h{0.5pt}\sigma_k} \h{1pt}   F^{\h{0.5pt}l}_{n, R} \h{1pt}\big[ U^{sc}_{\infty} \big] \big(\sigma_k\h{1pt}\widehat{\xi}\h{2pt}\big)   +\dfrac{|\h{1pt}\xi\h{1pt}|-(1-s)\h{1pt}\sigma_k}{s\h{0.5pt}\sigma_k}\h{1pt}J_n\big(\sigma_k\h{1pt}\widehat{\xi}\h{2pt}\big) \quad &\text{if }& \xi \in D_{\sigma_k} \setminus D_{(1-s)\sigma_k}.
\end{aligned}
\right.
\label{definition of h_n in large scale}
\end{eqnarray}Still by (1) in Lemma \ref{uniform convergence of u_n in S_sigma_k in large scale}, it follows \begin{align}\label{unif conv of hnsr} h_{n, s, R} \longrightarrow y_* \h{10pt}\text{in $C^0\big(\overline{D_{\sigma_k}}\big)$ as $n \to \infty$.}
\end{align}
Our comparison map $\overline{v}_{n, s, R}$ in the large--scale regime is then defined by \begin{eqnarray}\overline{v}_{n, s, R} \big(\xi\big) :=\left\{
\begin{aligned}
&  \Pi_{\h{0.5pt}\mathbb{S}^2} \left[ h_{n, s, R}\left(\dfrac{\xi}{1-s}\right) \right] \quad &\text{if }& \xi \in D_{(1-s)\h{1pt}\sigma_k}; \\[4mm]
&\dfrac{\sigma_k-|\h{1pt}\xi\h{1pt}|}{s\h{0.5pt}\sigma_k} \h{2pt}   J_n  \big(\sigma_k\h{1pt}\widehat{\xi}\h{1.5pt}\big)   +\dfrac{|\h{1pt}\xi\h{1pt}|-(1-s)\h{1pt}\sigma_k}{s\h{0.5pt}\sigma_k}\h{1pt}U_n\big(\sigma_k\h{1pt}\widehat{\xi}\h{1.5pt}\big) \quad &\text{if }& \xi \in D_{\sigma_k} \setminus D_{(1-s)\sigma_k}.
\end{aligned}
\right.
\label{definition of v_n in large scale}
\end{eqnarray}
By our definition of $\overline{v}_{n, s, R}$ in (\ref{definition of v_n in large scale}), it turns out $\overline{v}_{n, s, R} = U_n$ on $\p D_{\sigma_k}$. In addition, from the symmetry obeyed by $h_{n, s, R}$, $J_n$, $U_n$ and the fact (\ref{0 of Jn3}), the first two components of $\overline{v}_{n, s, R}$ are even and the third component of $\overline{v}_{n, s, R}$ is odd with respect to the $\xi_2$--variable. To show that $\overline{v}_{n, s, R}$ is an appropriate comparison map, we are left to verify the Signorini obstacle condition satisfied by the second component of $\overline{v}_{n, s, R}$. \vspace{0.4pc}

\noindent\textbf{Step 2. Signorini obstacle condition of $\overline{v}_{n, s, R}$} \vspace{0.2pc}

We claim that  \begin{eqnarray}\big[\h{1pt}\overline{v}_{n, s, R} \h{1pt}\big]_2 \h{2pt}\geq\h{2pt} H_{a_n}\h{0.5pt}b \h{20pt}\text{ on $T_{\sigma_k}$.} \label{obstacle condition of v_n}\end{eqnarray}
The proof of (\ref{obstacle condition of v_n}) is divided into two cases. In the following,  $ \big<\cdot, \cdot \big>$ is the standard inner product in $\mathbb{R}^3$.\vspace{0.4pc}\\
\noindent\textbf{Case I:} \h{3pt} Suppose that (\ref{pi_L convergence in large scale}) holds. Then  any vector field $v \in \overline{N}_k$  can be represented by $v = v_* + w \h{1pt}t_*$ on $T_{\sigma_k}$, where $w$ is some function  satisfying \begin{align}\label{low boun w case I}w \geq \gamma_* \h{15pt}\text{ on $T_{\sigma_k}$.}\end{align} This representation of $v$ and (\ref{definition of F_n^l}) infer  \begin{eqnarray*}  F^{\h{0.5pt}l}_{n, R}\h{0.5pt} [\h{0.5pt}v\h{0.5pt}]  = X_n + \dfrac{R\h{1pt}s_n \h{0.5pt}\big(w - \gamma_*\big)}{|\h{1pt}w - \gamma_*\h{1pt}|\vee R}\h{2pt}t_* \h{20pt}\text{on $T_{\sigma_k}$.}
\end{eqnarray*}Define $t_n := e^*_3 \times X_n$ and write  $t_* = \big< t_*, t_n\big> \h{1pt}t_n + \big< t_*, X_n\big> \h{1pt}X_n$. The last equality can then be rewritten by \begin{eqnarray}  F^{\h{0.5pt}l}_{n, R}\h{0.5pt} [\h{0.5pt}v\h{0.5pt}] =  \dfrac{R\h{1pt}s_n \h{0.5pt}\big(w - \gamma_*\big)}{|\h{1pt}w - \gamma_*\h{1pt}|\vee R} \h{2pt}\h{1pt}\big<t_*, t_n\big>\h{2pt}t_n + \left[ 1 + \dfrac{R\h{1pt}s_n \h{0.5pt}\big(w - \gamma_*\big)}{|\h{1pt}w - \gamma_*\h{1pt}|\vee R} \h{2pt}\big<t_*, X_n\big>\right]\h{2pt}X_n \h{20pt}\text{on $T_{\sigma_k}$.}\label{representation of F_n^l[v]}
\end{eqnarray} Notice (\ref{low boun w case I}) and the fact that $\big< t_*, t_n\big> > 0$. It holds \begin{align}\label{sign of tn coor of F} \left[\h{1pt} F^{\h{0.5pt}l}_{n, R}\h{0.5pt} [\h{0.5pt}v\h{0.5pt}] \h{1pt}\right]_{t_n} \geq 0 \h{15pt}\text{  on $T_{\sigma_k}$.}\end{align}Here  $[\h{1pt}X\h{1pt}]_{t_n}$ and $[\h{1pt}X\h{1pt}]_{X_n}$   denote the $t_n$ and $X_n$ coordinates of a vector  $X$, respectively. Moreover, \begin{align}\label{sing of xn flnr} \left[\h{1pt} F^{\h{0.5pt}l}_{n, R}\h{0.5pt} [\h{0.5pt}v\h{0.5pt}] \h{1pt}\right]_{X_n} \longrightarrow 1 \h{15pt}\text{as $n \to \infty$.} \end{align} In light of  (\ref{representation of F_n^l[v]})--(\ref{sing of xn flnr}), we have \begin{eqnarray} \left[\h{1pt}\Pi_{\h{0.5pt}\mathbb{S}^2} \Big[ \h{1pt}F^{\h{0.5pt}l}_{n, R}\h{0.5pt} [\h{0.5pt}v\h{0.5pt}] \h{1pt}\Big] \h{2pt}\right]_2 \h{2pt}\geq\h{2pt} H_{a_n}\h{1pt}b \h{20pt}\text{on $T_{\sigma_k}$, if $n$ is suitably large.}
\label{boundary condition on piT_sigma_k}
\end{eqnarray}By (5) in Lemma \ref{uniform convergence of u_n in S_sigma_k in large scale}, we can apply the same arguments for deriving (\ref{sign of tn coor of F}) and (\ref{sing of xn flnr}) to get  \begin{eqnarray}\left[\h{1pt} F^{\h{0.5pt}l}_{n, R}\h{0.5pt} \big[\h{0.5pt}U^{sc}_\infty\h{0.5pt}\big]  \h{1pt}\right]_{t_n} \h{2pt}\geq\h{2pt} 0 \h{10pt}\text{and} \h{10pt}\left[ \h{1pt}F^{\h{0.5pt}l}_{n, R}\h{0.5pt} \big[\h{0.5pt}U^{sc}_\infty\h{0.5pt}\big] \h{1pt}\right]_{X_n} \h{2pt}\longrightarrow\h{2pt} 1  \h{10pt}\text{as $n \rightarrow \infty$ \h{5pt} at $(\pm \sigma_k, 0 )$.}\label{sign of F_u_infty}\end{eqnarray} In the current case,  the distance between  $X_n$ and $J_n\big(\pm \sigma_k, 0 \big) $ tends to $0$ as $n \to \infty$. By (\ref{low boun of Jn pm sigmak}), it follows \begin{eqnarray}\label{sing of Jn at sig} \left[\h{1pt} J_n\big(\pm \sigma_k, 0 \big) \h{1pt}\right]_{t_n} \h{2pt}\geq \h{2pt} 0 \h{15pt}\text{and}\h{15pt} \left[\h{1pt} J_n\big(\pm \sigma_k, 0 \big) \h{1pt}\right]_{X_n} \h{2pt}\geq\h{2pt} 0. \label{J_n at sigma_k}\end{eqnarray} Owing to (\ref{boundary condition on piT_sigma_k})--(\ref{sing of Jn at sig}), we obtain from the definition of $h_{n, s, R}$ in (\ref{definition of h_n in large scale}) that $$\big[\h{1pt}\Pi_{\h{0.5pt}\mathbb{S}^2}\big[\h{1pt} h_{n, s, R}\h{1pt}\big] \h{1pt}\big]_2 \h{2pt}\geq\h{2pt} H_{a_n}b \h{15pt}\text{ on $T_{\sigma_k} $.}$$ By this inequality, (\ref{low boun of Jn pm sigmak}) and (3) in Lemma \ref{uniform convergence of u_n in S_sigma_k in large scale}, (\ref{obstacle condition of v_n}) then follows by the definition of $\overline{v}_{n, s, R}$ in  (\ref{definition of v_n in large scale}). \\[3mm]
\textbf{Case II:} \h{3pt}Suppose that (i) and (ii) in (\ref{infinite case in large scale}) hold. If $y_* \in \Gamma_b \setminus \p \h{1pt}\Gamma_b$, then (\ref{obstacle condition of v_n}) follows by (\ref{unif conv of hnsr}), (\ref{low boun of Jn pm sigmak}) and (3) in Lemma \ref{uniform convergence of u_n in S_sigma_k in large scale},   provided that $n$ is suitably large. We are left to  consider the case in which $y_*$ and $t_*$ satisfy (\ref{y* at boundary case}). Define $\gamma_n := \big\langle\,\Pi_{\h{0.5pt}\Gamma}\big(y_n\big),X_n\,\big\rangle\, X_n$. Since $X_n$ and $\Pi_{\h{0.5pt}\Gamma}\big(y_n\big)$ converge to $y_*$ as $n \to \infty$, it holds
\begin{eqnarray*}2\h{1pt}\big|\h{1pt}\Pi_{\h{0.5pt}\Gamma}\big(y_n\big)-\gamma_n\h{1pt}\big| \h{2pt}\geq\h{2pt} \big|\h{1pt}\Pi_{\h{0.5pt}\Gamma}\big(y_n\big) - X_n\h{1pt}\big| \h{10pt} \text{ for large } n.\end{eqnarray*}
This estimate together with  (ii) in (\ref{infinite case in large scale}) induce that \begin{equation}
    \liminf_{n \rightarrow \infty}\left|\h{1pt}\dfrac{\Pi_{\h{0.5pt}\Gamma}\big(y_n\big)-\gamma_n}{s_n}\h{1pt}\right|=\infty.
    \label{gamma_n to infty III}
\end{equation}
Utilizing $\gamma_n$, we decompose $F^{\h{0.5pt}l}_{n, R}\h{0.5pt} [\h{0.5pt}v\h{0.5pt}]$ as follows
\begin{eqnarray}
F^{\h{0.5pt}l}_{n, R}\h{0.5pt} [\h{0.5pt}v\h{0.5pt}] \h{2pt}=\h{2pt}\gamma_n+s_n\left[\h{1pt}\dfrac{\Pi_{\h{0.5pt}\Gamma}\big(y_n\big) - \gamma_n}{s_n}+R\h{1pt}\dfrac{v - v_*}{|v-v_*|\vee R}\right].\label{decompo in case III}
\end{eqnarray}Recalling the $t_n$ introduced in Case I, we have $\Pi_{\h{0.5pt}\Gamma}\big(y_n\big) - \gamma_n = \big<\h{1.5pt} \Pi_{\h{0.5pt}\Gamma}\big(y_n\big), t_n \big> \h{2pt}t_n$. Moreover, $\big<\h{1.5pt} \Pi_{\h{0.5pt}\Gamma}\big(y_n\big), t_n \big> > 0$ by (i) in (\ref{infinite case in large scale}). (\ref{decompo in case III}) can then be rewritten as  \begin{eqnarray}
F^{\h{0.5pt}l}_{n, R}\h{0.5pt} [\h{0.5pt}v\h{0.5pt}] \h{2pt}=\h{2pt}\gamma_n+s_n\left[\h{2pt}\left|\h{2pt}\dfrac{\Pi_{\h{0.5pt}\Gamma}\big(y_n\big) - \gamma_n}{s_n} \h{2pt}\right| \h{2pt}t_n +R\h{1pt}\dfrac{v - v_*}{|v-v_*|\vee R}\right].\label{decompo in case III, 2}
\end{eqnarray} We still represent $v = v_* + w \h{1pt}t_*$ on $T_{\sigma_k}$.  (\ref{decompo in case III, 2}) then induces\begin{eqnarray*}
F^{\h{0.5pt}l}_{n, R}\h{0.5pt} [\h{0.5pt}v\h{0.5pt}] = s_n\left[\h{3pt}\left|\dfrac{\Pi_{\h{0.5pt}\Gamma}\big(y_n\big)-\gamma_n}{s_n}\right| + R\h{1pt}\dfrac{w \h{1pt}\big<t_*, t_n\big> }{| w |\vee R} \right] \h{1pt}t_n + \left[ \big<\Pi_{\h{0.5pt}\Gamma} \big(y_n\big), X_n\big> + R\h{1pt}s_n\h{1pt}\dfrac{w \h{1pt}\big<t_*, X_n\big> }{| w |\vee R}\right]\h{1pt} X_n \h{10pt}\text{on $T_{\sigma_k}$.}
\end{eqnarray*}
Firstly, (\ref{gamma_n to infty III}) induces $\left[ \h{1pt}F^{\h{0.5pt}l}_{n, R}\h{0.5pt} [\h{0.5pt}v\h{0.5pt}]\h{1pt}\right]_{t_n} > 0$ on $T_{\sigma_k}$ for large $n$. Moreover, $\left[ \h{1pt}F^{\h{0.5pt}l}_{n, R}\h{0.5pt} [\h{0.5pt}v\h{0.5pt}]\h{1pt}\right]_{X_n} > 0$ on $T_{\sigma_k}$ if $n$ is large in that $ \Pi_{\h{0.5pt}\Gamma}\big(y_n\big)$ and $X_n$ converge to $y_*$ as $n \rightarrow \infty$. Hence, (\ref{boundary condition on piT_sigma_k}) still holds in the current case. Following the same arguments as in Case I, we obtain (\ref{obstacle condition of v_n}) as well in the current case.\vspace{0.6pc}\\
\textbf{Step 3. Convergence of potential energy}\vspace{0.5pc}\\
We claim \begin{eqnarray}\limsup_{n \to \infty} \h{1pt}a_n \h{1pt}\left( \dfrac{\lambda_n r_n}{s_n}\right)^2\int_{D_{\sigma_k}} \Big(\h{1pt}\big|\h{1pt}\overline{v}_{n, s, R}\h{1pt}\big|^2 - 1 \Big)^2   \h{10pt}\text{is independent of $R$ and converges to $0$ as $s \rightarrow 0$.}
\label{convergence of ginzburg-landau term}
\end{eqnarray}The norm of $\overline{v}_{n, s, R}$ identically equals $1$ on $D_{(1 - s)\h{0.5pt}\sigma_k}$. It turns out  \begin{eqnarray}\label{identity of overline}
  a_n \h{1pt}\left( \dfrac{\lambda_n r_n}{s_n}\right)^2\int_{D_{\sigma_k}} \Big(\h{1pt}\big|\h{1pt}\overline{v}_{n, s, R}\h{1pt}\big|^2 - 1 \Big)^2 \h{2pt}=\h{2pt} a_n \h{1pt}\left( \dfrac{\lambda_n r_n}{s_n}\right)^2\int_{D_{\sigma_k} \setminus D_{(1 - s)\h{0.5pt}\sigma_k}} \Big(\h{1pt}\big|\h{1pt}\overline{v}_{n, s, R}\h{1pt}\big|^2 - 1 \Big)^2.
\end{eqnarray}Therefore, the left--hand side above is indepdent of $R$. \vspace{0.3pc}

If Case A in (\ref{definition of J_n}) holds, then by (\ref{definition of v_n in large scale})
\begin{eqnarray*}\overline{v}_{n, s, R}(\xi) =  \dfrac{\sigma_k - |\h{1pt}\xi \h{1pt}|}{s \h{0.5pt} \sigma_k} \h{1pt} \Pi_{\h{0.5pt}\mathbb{S}^2}\big[\h{1pt}U_n \h{1pt}\big] \big( \sigma_k\h{1pt}\widehat{\xi}\h{2pt} \big) +\dfrac{|\h{1pt}\xi \h{1pt}|- ( 1 - s )\h{1pt}\sigma_k}{s \h{0.5pt}\sigma_k}\h{1pt}U_n\big( \sigma_k\h{1pt}\widehat{\xi}\h{2pt} \big) \quad \quad\text{for } \xi \in D_{\sigma_k} \setminus D_{(1-s)\sigma_k}.\end{eqnarray*} Plugging this identity into the right--hand side of (\ref{identity of overline}) and using (1) in Lemma \ref{uniform convergence of u_n in S_sigma_k in large scale}, we obtain
\begin{eqnarray*}
  a_n \h{1pt}\left( \dfrac{\lambda_n r_n}{s_n}\right)^2\int_{D_{\sigma_k}} \Big(\h{1pt}\big|\h{1pt}\overline{v}_{n, s, R}\h{1pt}\big|^2 - 1 \Big)^2  \h{2pt}\lesssim\h{2pt}  a_n \h{1pt}\left( \dfrac{\lambda_n r_n}{s_n}\right)^2\int_{D_{\sigma_k} \setminus D_{(1 - s)\h{0.5pt}\sigma_k}} \left(\h{1pt}\left|\h{1pt}U_n\big( \sigma_k\h{1pt}\widehat{\xi}\h{2pt}\big) \h{1pt}\right|^2 - 1 \right)^2 \h{3pt}\lesssim\h{3pt}s \h{1pt}b_k.
\end{eqnarray*}Taking $s \rightarrow 0$ yields (\ref{convergence of ginzburg-landau term}).\vspace{0.2pc}

In the next, we assume that Case B in (\ref{definition of J_n}) holds. Case C in (\ref{definition of J_n}) can be considered by the same arguments as Case B. Still by (\ref{definition of J_n}) and (\ref{definition of v_n in large scale}), it turns out
\begin{eqnarray*}\overline{v}_{n, s, R}(\xi) =  \dfrac{\sigma_k - |\h{1pt}\xi \h{1pt}|}{s \h{0.5pt} \sigma_k} \h{2pt} \mathrm{Rot}_{\alpha_1} \Pi_{\mathbb{S}^2} \big[\h{1pt}U_n\h{1pt}\big] \big( \sigma_k\h{1pt}\widehat{\xi}\h{2pt} \big)  +\dfrac{|\h{1pt}\xi \h{1pt}|- ( 1 - s )\h{1pt}\sigma_k}{s \h{0.5pt}\sigma_k}\h{1pt}U_n\big( \sigma_k\h{1pt}\widehat{\xi}\h{2pt} \big) \quad \quad\text{for } \xi \in  D_{\sigma_k} \setminus D_{(1-s)\sigma_k}.\end{eqnarray*}
Utilizing the above representation and denoting by $t$ the quantity $\dfrac{\sigma_k - |\h{1pt}\xi \h{1pt}|}{s \h{0.5pt} \sigma_k}$, we obtain
\begin{align}\big|\h{1pt}\overline{v}_{n, s, R}(\xi)\h{1pt}\big|^2 - 1 &= \left(1 - t\right)^2\left[\h{2pt}\left|\h{1pt}U_n\big(\sigma_k \h{1pt}\widehat{\xi}\h{2pt}\big) \right|^2 - 1 \h{2pt}\right] + 2 \h{1pt} t \h{1pt} (1 - t) \h{1pt} \left[ \h{2pt} \Big<\mathrm{Rot}_{\alpha_1} \Pi_{\mathbb{S}^2} \big[\h{1pt}U_n\h{1pt}\big] \big( \sigma_k\h{1pt}\widehat{\xi}\h{2pt} \big), U_n\big(\sigma_k\h{1pt}\widehat{\xi}\h{2pt}\big)  \Big>  - 1 \right]\label{computation of v_n - 1}\nonumber\\[2mm]
&= \left(1 - t\right)^2\left[\h{2pt}\left|\h{1pt}U_n\big(\sigma_k \h{1pt}\widehat{\xi}\h{2pt}\big) \right|^2 - 1 \h{2pt}\right] +  2 \h{1pt} t \h{1pt} (1 - t) \h{1pt} \left[ \h{2pt} \left|\h{1pt}U_n\big(\sigma_k\h{1pt}\widehat{\xi}\h{2pt}\big) \h{1pt}\right|  - 1 \right]\nonumber\\[2mm]
&+ 2 \h{1pt} t \h{1pt} (1 - t)\h{2pt} \left<\Big[\h{1pt}\mathrm{Rot}_{\alpha_1} - \mathrm{I}_3\h{1pt}\Big] \Pi_{\mathbb{S}^2} \big[\h{1pt}U_n\h{1pt}\big] \big( \sigma_k\h{1pt}\widehat{\xi}\h{2pt} \big), U_n\big(\sigma_k\h{1pt}\widehat{\xi}\h{2pt}\big)  \right>.
\end{align}
The quantity in the last line of (\ref{computation of v_n - 1}) equals  $$  - 4 \h{1pt} t \h{1pt} (1 - t) \sin^2 \frac{\alpha_1}{2}\left<\begin{small}\left(\begin{array}{lcl}1 &0 &0\\
0 &1 &0\\
0 &0 &0 \end{array}\right)\end{small} \Pi_{\mathbb{S}^2} \big[\h{1pt}U_n\h{1pt}\big] \big( \sigma_k\h{1pt}\widehat{\xi}\h{2pt} \big), U_n\big(\sigma_k\h{1pt}\widehat{\xi}\h{2pt}\big)  \right>.$$ It then follows from (\ref{computation of v_n - 1}) that \begin{eqnarray*}\Big|\h{1pt}\big|\h{1pt}\overline{v}_{n, s, R}(\xi)\h{1pt}\big|^2 - 1 \h{1pt}\Big| \h{2pt}\lesssim\h{2pt} \sin^2 \dfrac{\alpha_1}{2} + \left|\h{2pt}\left|\h{1pt}U_n\big(\sigma_k \h{1pt}\widehat{\xi}\h{2pt}\big) \right|^2 - 1 \h{2pt}\right| \h{20pt}\text{for $\xi \in D_{\sigma_k} \setminus D_{(1 - s)\h{0.5pt}\sigma_k}$.}
\end{eqnarray*}Applying this estimate to the right--hand side of (\ref{identity of overline}), by  (1) in Lemma \ref{uniform convergence of u_n in S_sigma_k in large scale}, we then obtain \begin{eqnarray}a_n \h{1pt}\left( \dfrac{\lambda_n r_n}{s_n}\right)^2\int_{D_{\sigma_k}} \Big(\h{1pt}\big|\h{1pt}\overline{v}_{n, s, R}\h{1pt}\big|^2 - 1 \Big)^2 \h{2pt}\lesssim\h{2pt}s\h{1pt}b_k + a_n \h{1pt}\left( \dfrac{\lambda_n r_n}{s_n}\right)^2 s \h{2pt}\sin^4 \dfrac{\alpha_1}{2}.
\label{bounds of 1 - |v_n|^2}
\end{eqnarray}
Now we illustrate the relative positions of $X_n$, $U_n\big(\sigma_k, 0\big)$ and $\Pi_{\h{0.5pt}\mathbb{S}^2}\big[ \h{1.5pt}U_n\big(\sigma_k, 0 \big) \h{1.5pt} \big]$ in Case B as follows:  \\
\begin{center}
    \includegraphics[scale=0.08]{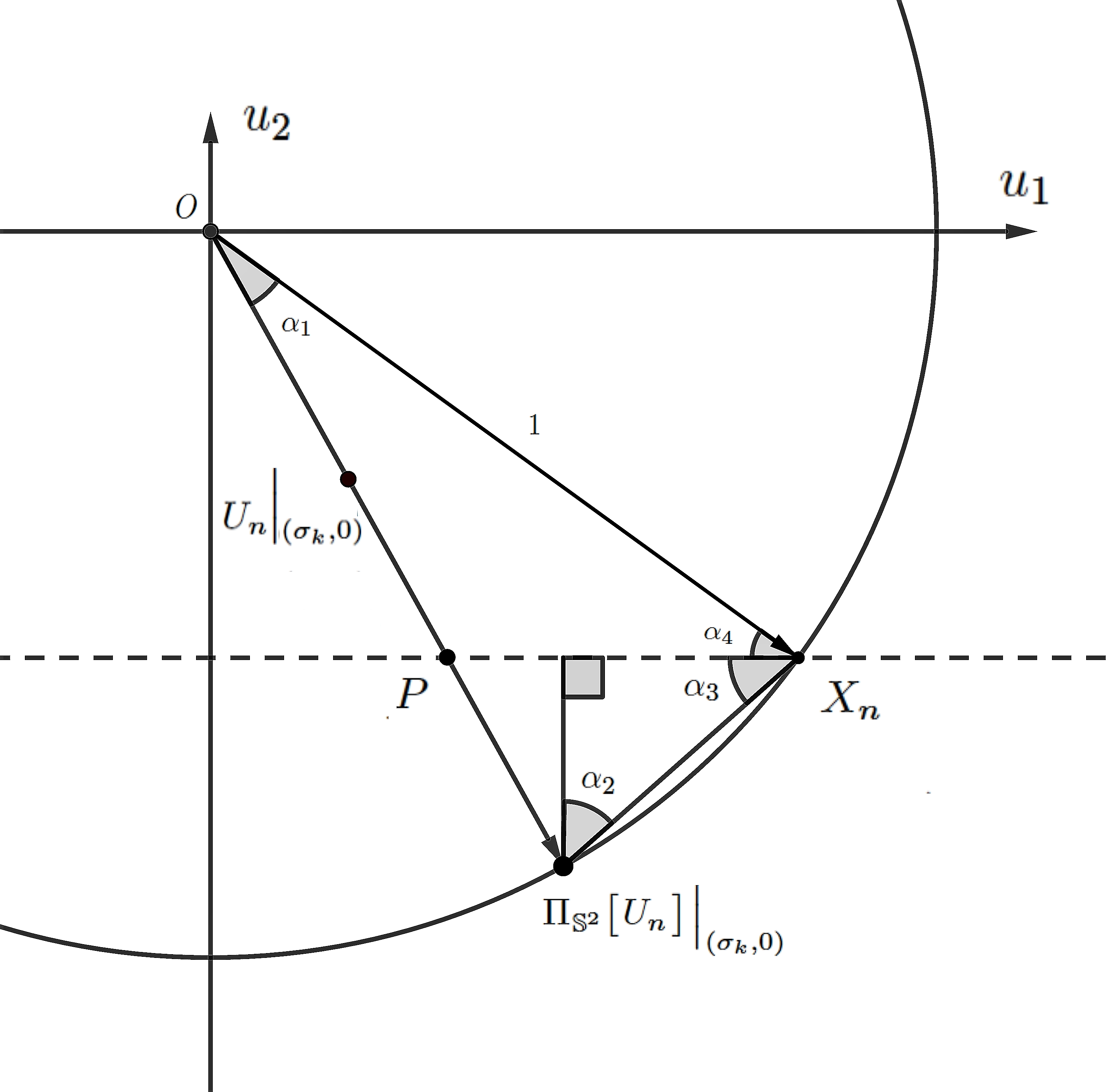}\
\\[2mm]
Figure 3
\end{center} It can be shown from Figure 3 that
\begin{equation}
     2 \sin^2 \dfrac{\alpha_1}{2}=\dfrac{1}{2}\h{2pt}\Big|\h{1pt}X_n - \Pi_{\mathbb{S}^2} \big[\h{1pt}U_n\big(\sigma_k, 0 \big)\big]\h{1pt}\Big|^2 =  \dfrac{1}{2\cos^2\alpha_2} \h{2pt}\Big|\h{1pt}X_{n,2}- \big[\h{1pt}\Pi_{\mathbb{S}^2} \big[\h{1pt}U_n\big(\sigma_k, 0 \big)\big] \h{1pt}\big]_2\h{1pt}\Big|^2.
    \label{sin law cos law}
\end{equation}
Still by Figure 3, $\alpha_2+\alpha_3 = \pi/2$ and $\alpha_3 + \alpha_4 \rightarrow \pi/2$ as $n \to \infty$. Hence for large $n$, the angle $\alpha_2$ almostly equals  $\alpha_4$. Since $\alpha_4$ converges to $\arcsin\h{2pt}(-b)$ as $n \to \infty$, we then can find a $\alpha_0 \in \big( 0, \pi/2\big)$ depending only on $b$ so that $\alpha_2 \in \big(0, \pi/2 - \alpha_0\big)$ for large $n$, which furthermore infers by (\ref{sin law cos law}) the following estimate:\begin{equation}\label{est of sin2 alpha}
2 \sin^2 \dfrac{\alpha_1}{2} \h{2pt}\leq\h{2pt}\dfrac{1}{2\sin^2\alpha_0} \h{2pt}\Big|\h{1pt}X_{n,2}-\big[\h{1pt}\Pi_{\mathbb{S}^2} \big[\h{1pt}U_n\big(\sigma_k, 0 \big)\big] \h{1pt}\big]_2\h{1pt}\Big|^2 \h{2pt}\leq\h{2pt} \dfrac{1}{2\sin^2\alpha_0}\h{2pt} \Big|\h{1pt}U_n\big(\sigma_k, 0\big) - \Pi_{\mathbb{S}^2}\big[\h{1pt}U_{n}\big(\sigma_k, 0\big)\h{1pt}\big]\Big|^2.
\end{equation}
Here we have used the item (3) in Lemma \ref{uniform convergence of u_n in S_sigma_k in large scale}. This estimate and  the item (6) in Lemma \ref{uniform convergence of u_n in S_sigma_k in large scale} yield \begin{eqnarray} a_n \h{1pt}\left( \dfrac{\lambda_n r_n}{s_n}\right)^2 \sin^4 \dfrac{\alpha_1}{2} \h{2pt}\lesssim\h{2pt} \dfrac{1}{\sin^4\alpha_0} a_n\left(\dfrac{\lambda_n\h{1pt}r_n}{s_n}\right)^2  \Big|\h{1pt} \big|\h{1.5pt}U_n\big(\sigma_k, 0\big) \h{1pt}\big| - 1 \h{1pt}\Big|^4 \h{2pt}\leq\h{2pt}  \dfrac{b_k}{\sin^4 \alpha_0}.
\label{sin alpha}
\end{eqnarray}
By applying this estimate to the right--hand side of (\ref{bounds of 1 - |v_n|^2}), (\ref{convergence of ginzburg-landau term}) then follows for Case B in (\ref{definition of J_n}). \\

\noindent \textbf{Step 4. Convergence of Dirichlet energy}\vspace{0.5pc}\\
\noindent The Dirichlet energy of $\overline{v}_{n, s, R}$ is computed as follows:
\begin{eqnarray}
\int_{D_{\sigma_k}}\Big|\h{.5pt} D_\xi \overline{v}_{n, s, R} \h{.5pt}\Big|^2
= \int_{D_{ \sigma_k }} \Big| D_{\xi} \h{1pt} \Pi_{\h{0.5pt}\mathbb{S}^2} \big[ \h{1pt}F^{\h{0.5pt}l}_{n, R}\h{0.5pt} [\h{0.5pt} v \h{0.5pt}] \h{1pt}\big] \h{1pt}\Big|^2 + \int_{D_{\sigma_k} \setminus D_{(1-s)\h{1pt}\sigma_k}} \Big| D_\xi \h{1pt}\Pi_{\h{0.5pt}\mathbb{S}^2} \big[ h_{n, s, R} \big] \Big|^2 +  \Big| D_\xi \overline{v}_{n, s, R} \Big|^2.
\label{decomp. of D_k in large scale}
\end{eqnarray}

\noindent \textbf{Step 4.1.} For the first term on the right--hand side above, it can be computed that
$$
\int_{D_{ \sigma_k }} \Big| D_{\xi} \h{0.5pt}  \Pi_{\h{0.5pt}\mathbb{S}^2} \big[ \h{1pt}F^{\h{0.5pt}l}_{n, R}\h{0.5pt} [\h{0.5pt} v \h{0.5pt}] \h{1pt}\big] \h{1pt}\Big|^2
= \int_{D_{ \sigma_k }} \Big|\h{1pt} F^{\h{0.5pt}l}_{n, R}\h{0.5pt} [\h{0.5pt} v \h{0.5pt}] \h{1.5pt} \Big|^{-2}\Big| D_{\xi}   \h{1pt}F^{\h{0.5pt}l}_{n, R}\h{0.5pt} [\h{0.5pt} v \h{0.5pt}]\h{1pt}\Big|^2-  \Big|\h{1pt} F^{\h{0.5pt}l}_{n, R}\h{0.5pt} [\h{0.5pt} v \h{0.5pt}] \h{1.5pt} \Big|^{-4}
\sum^2_{i=1}\Big<  F^{\h{0.5pt}l}_{n, R}\h{0.5pt} [\h{0.5pt} v \h{0.5pt}],  D_{\xi_i} F^{\h{0.5pt}l}_{n, R}\h{0.5pt} [\h{0.5pt} v \h{0.5pt}]  \Big>^2.
$$ Recall the definition of $F^{\h{0.5pt}l}_{n, R}\h{0.5pt} [\h{0.5pt} v \h{0.5pt}]$ in (\ref{definition of F_n^l}) and note that $F^{\h{0.5pt}l}_{n, R}\h{0.5pt} [\h{0.5pt} v \h{0.5pt}]$ converges to $y_*$  in $C^0\big( \overline{D_{\sigma_k}}\big) $ as $n \to \infty$. It then turns out from the last equality that
$$
\lim_{n\to \infty} s_n^{-2}\int_{D_{ \sigma_k }} \Big|\h{1pt} D_{\xi} \h{1pt}  \Pi_{\h{0.5pt}\mathbb{S}^2} \big[ \h{1pt}F^{\h{0.5pt}l}_{n, R}\h{0.5pt} [\h{0.5pt} v \h{0.5pt}] \h{1pt}\big] \h{1pt}\Big|^2
= \int_{D_{ \sigma_k }} R^2 \left| D_{\xi}   \h{1pt}\dfrac{v - Z_*}{|\h{1pt}v - Z_*\h{1pt}|\vee R}\right|^2
-  R^2 \h{1.5pt}\sum^2_{i=1} \Big<  y_*,  D_{\xi_i}  \h{1pt}\dfrac{v - Z_*}{|\h{1pt}v - Z_*\h{1pt}|\vee R}  \Big>^2.
$$
In light of the orthogonality of $y_*$ and $v-Z_*$, the inner product on the right--hand side above equals $0$. Hence, by  the fact that $v \in H^1\big(D_{\sigma_k}; \mathbb{R}^3\big)$, we obtain from the last equality that
\begin{eqnarray}
\lim_{R \to \infty} \lim_{n \rightarrow \infty}  s_n^{-2}\int_{ D_{ \sigma_k } } \Big| D_{ \xi }  \h{1pt} \Pi_{ \h{0.5pt} \mathbb{S}^2 } \big[ \h{1pt} F^{ \h{0.5pt} l }_{n, R} \h{0.5pt} [ \h{0.5pt} v \h{0.5pt}] \h{1pt} \big] \h{1pt} \Big|^2 =  \int_{D_{\sigma_k}}\big| D_\xi v \h{1pt}\big|^2.
\label{conv. main part in large scale}
\end{eqnarray}

\noindent \textbf{Step 4.2.} Let $\big(\tau, \varphi\big)$ be the polar coordinates of the $\big(\xi_1, \xi_2\big)$--space. Then \begin{eqnarray}
 \int_{D_{\sigma_k} \setminus D_{(1-s)\sigma_k}} \big| D_\xi \h{1pt}h_{n, s, R} \big|^2 = \int_{D_{\sigma_k} \h{1pt}\setminus \h{1pt}D_{(1-s)\h{1pt}\sigma_k}} \big|\partial_\tau h_{n, s, R} \big|^2   +   \tau^{-2} \big| \partial_\varphi h_{n, s, R} \big|^2.
 \label{decomp of D_k / D_(1-s), large}
\end{eqnarray}
\noindent \textbf{Step 4.2.1.} For the radial derivative of $h_{n, s, R}$, the definition of $h_{n, s, R}$ in (\ref{definition of h_n in large scale}) induces
\begin{align}\label{com of partial tau h}
\big(s  \h{1pt}\sigma_k\big)^{2}\int_{D_{\sigma_k} \setminus D_{(1-s)\h{1pt}\sigma_k}} \big|\h{1pt}\partial_\tau h_{n, s, R}\h{1pt}\big|^2
&=  \int_{D_{\sigma_k} \setminus D_{( 1-s )\h{1pt}\sigma_k}}\h{1pt}\left|   F^{\h{0.5pt}l}_{n, R} \h{1pt}\big[\h{1pt} U^{sc}_{\infty} \h{1pt}\big] - J_n\right|^2\big(\h{1pt}\sigma_k\h{1pt}\widehat{\xi}\h{2pt}\big)\\[2mm] 
&\lesssim \int_{ D_{\sigma_k} \setminus D_{( 1-s ) \h{1pt} \sigma_k} } \h{1pt} \Big| F^{\h{0.5pt}l}_{n, R} \h{1pt}\big[\h{1pt} U^{sc}_{\infty} \h{1pt}\big] - U_n \Big|^2 \big(\sigma_k \h{1pt}\widehat{\xi}\h{2pt}\big) +  \Big| U_n  - J_n \Big|^2 \big(\h{1pt}\sigma_k \widehat{\xi} \h{2pt}\big).\nonumber
\end{align}Now we estimate the integral in the last line of (\ref{com of partial tau h}).\vspace{0.1pc}

 Firstly, we claim \begin{align}\label{est of first int}   \lim_{n \rightarrow \infty} s_n^{-2} \int_{ D_{\sigma_k} \setminus D_{( 1-s ) \h{1pt} \sigma_k} } \h{1pt} \Big| F^{\h{0.5pt}l}_{n, R} \h{1pt}\big[\h{1pt} U^{sc}_{\infty} \h{1pt}\big] - U_n \Big|^2 \big(\sigma_k \h{1pt}\widehat{\xi}\h{2pt}\big) = 0 \h{15pt}\text{for any $R > \big\|\h{1pt} U^{sc}_{\infty} - Z_*  \big\|_{\infty, \p D_{\sigma_k}}$.} 
\end{align}Due to the definition of $F_{n, R}^l [\h{1pt}v\h{1pt}]$ in (\ref{definition of F_n^l}), the cases in (\ref{infinite case in large scale}) and (\ref{pi_L convergence in large scale}) should be treated separately. In the following arguments, we always take  $R >  \big\|\h{1pt} U^{sc}_{\infty} - Z_*  \big\|_{\infty, \p D_{\sigma_k}}$. \\[2mm]
\textbf{Case (\ref{infinite case in large scale}):} In this case, it holds \begin{align*}s_n^{-2} \int_{ D_{\sigma_k} \setminus D_{( 1-s ) \h{1pt} \sigma_k} } \h{1pt} \Big| F^{\h{0.5pt}l}_{n, R} \h{1pt}\big[\h{1pt} U^{sc}_{\infty} \h{1pt}\big] - U_n \Big|^2 \big(\sigma_k \h{1pt}\widehat{\xi}\h{2pt}\big) \h{2pt} \lesssim \h{2pt} \Big\|\h{1pt}U^{sc}_n - U^{sc}_\infty \h{1pt}\Big\|^2_{\infty; \h{1pt}\p D_{\sigma_k}} + \left| \h{1pt}\dfrac{\Pi_{\h{0.5pt}\Gamma}\big(y_n\big) - y_n}{s_n} - v_*\h{1pt}\right|^2.
\end{align*}In light of (\ref{pi_gamma - y_n}) and (2) in Lemma \ref{uniform convergence of u_n in S_sigma_k in large scale}, (\ref{est of first int}) follows by taking $n \to \infty$ on both sides above.\\[2mm]
\textbf{Case  (\ref{pi_L convergence in large scale}):} In this case, we have  
\begin{align*}
s_n^{-2} \int_{ D_{\sigma_k} \setminus D_{( 1-s ) \h{1pt} \sigma_k} } \h{1pt} &\Big| F^{\h{0.5pt}l}_{n, R} \h{1pt}\big[\h{1pt} U^{sc}_{\infty} \h{1pt}\big] - U_n \Big|^2 \big(\sigma_k \h{1pt}\widehat{\xi}\h{2pt}\big)\\[2mm]
&  \h{2pt}\lesssim\h{2pt}  \Big\|\h{1pt}U_n^{sc} - U_\infty^{sc}\h{1pt}\Big\|^2_{\infty; \h{1pt}\p D_{\sigma_k}} +  \left| \dfrac{\Pi_{\h{0.5pt}\Gamma}\big(y_n\big) - y_n}{s_n} + \dfrac{X_n - \Pi_{\h{0.5pt}\Gamma}\big(y_n\big) }{s_n} - Z_*   \right|^2.
\end{align*}Note that in the case (\ref{pi_L convergence in large scale}), $Z_* = v_* + \gamma_*\h{1pt}t_*$. In light of (\ref{pi_L convergence in large scale}), (\ref{pi_gamma - y_n}) and (2) in Lemma \ref{uniform convergence of u_n in S_sigma_k in large scale}, we then obtain (\ref{est of first int})  by taking $n \to \infty$  on both sides above.\vspace{0.1pc}

In the next, we study the integral of $\big| U_n - J_n \big|^2\big(\h{1pt}\sigma_k \h{1pt}\widehat{\xi}\h{2pt}\big)$ in the last line of (\ref{com of partial tau h}). We claim that \begin{align}\label{conv of Un - Jn, case A}
\lim_{n \to \infty} s_n^{-2}\int_{ D_{\sigma_k}  \setminus D_{ (1 - s) \h{1pt} \sigma_k}  }  \h{1pt} \big| \h{1pt}U_n - J_n \h{1pt}\big|^2\big(\h{1pt}\sigma_k\h{1pt}\widehat{\xi}\h{2pt}\big) = 0.\end{align}
If Case A in (\ref{definition of J_n}) holds, then by (1) in Lemma \ref{uniform convergence of u_n in S_sigma_k in large scale}, it turns out  
\begin{align*} 
 s_n^{-2}\int_{ D_{\sigma_k}  \setminus D_{ (1 - s) \h{1pt} \sigma_k}  }  \h{1pt} \big| \h{1pt}U_n - J_n \h{1pt}\big|^2\big(\h{1pt}\sigma_k\h{1pt}\widehat{\xi}\h{2pt}\big)  \h{2pt} \leq \h{2pt}  s_n^{-2}\int_{ D_{\sigma_k}  \setminus D_{ (1 - s) \h{1pt}  \sigma_k}} \Big| \h{1pt}\big| U_n \big|^2 - 1 \h{1pt}\Big|^2\big(\h{1pt}\sigma_k\h{1pt}\widehat{\xi}\h{2pt}\big) \h{2pt}\lesssim\h{2pt}\dfrac{s\h{1pt}b_k}{a_n \left( \lambda_n r_n \right)^2}. 
\end{align*}By this estimate, we then obtain (\ref{conv of Un - Jn, case A}) since $a_n\big(\lambda_n r_n\big)^2 \to \infty$ as $n \to \infty$. If Case B in (\ref{definition of J_n}) holds, then \begin{align*} \big|\h{1pt}U_n - J_n \h{1pt}\big|^2 \h{2pt}& = \Big|\h{1pt}U_n - \mathrm{Rot}_{\alpha_1} \Pi_{\h{0.5pt}\mathbb{S}^2} \big[ U_n \big] \h{1pt}\Big|^2 \h{2pt}= \h{2pt} \big|\h{1pt}U_n  \h{1pt}\big|^2 + 1  - 2 \h{1pt} \Big< \h{1pt}U_n, \h{1.5pt} \mathrm{Rot}_{\alpha_1}  \Pi_{\h{0.5pt}\mathbb{S}^2} \big[ U_n \big] \h{1pt}\Big> \\[2mm]
&\h{2pt} = \Big|\h{1pt}\big|\h{1pt}U_n  \h{1pt}\big| -  1\h{1.5pt}\Big|^2  - 2 \h{1pt} \Big< \h{1pt}U_n, \h{1.5pt} \big[\h{1pt}\mathrm{Rot}_{\alpha_1} - \mathrm{I}_3\h{1pt}\big] \Pi_{\h{0.5pt}\mathbb{S}^2} \big[ U_n \big] \h{1pt}\Big> \h{2pt}\lesssim\h{2pt}\Big|\h{1pt}\big|\h{1pt}U_n  \h{1pt}\big| -  1\h{1.5pt}\Big|^2  + \sin^2\dfrac{\alpha_1}{2}   \h{15pt}\text{on $\p D_{\sigma_k}$.} 
\end{align*}Recalling (\ref{est of sin2 alpha}) and applying fundamental theorem of calculus, we have, for any $\varphi \in \big(0, 2\pi\big)$, that \begin{align*}\sin  \dfrac{\alpha_1}{2} &\h{2pt}\lesssim_b\h{2pt} \left| \h{1pt}\big| U_n\left(\sigma_k, 0 \right)\big|^2 - 1 \h{1pt}\right|  \\[2mm]&\h{2pt}\lesssim\h{2pt}\left| \h{1pt}\big| U_n\big(\sigma_k \cos \varphi, \sigma_k \sin \varphi\big)\big|^2 - 1 \h{1pt}\right| + \left|\h{2pt} \int_0^\varphi \Big[  \h{1pt}\big< U_n, \h{1.5pt} \xi^{\bot}\cdot D_{\xi} U_n\big> \h{1pt}\Big]\big(\sigma_k \cos\phi, \h{1.5pt}\sigma_k \sin \phi\big) \h{2pt}\mathrm{d} \phi \h{2pt} \right|.
\end{align*}Here $\xi^{\bot} = (- \xi_2, \xi_1)^\top$. The last two estimates then yield  \begin{align*} s_n^{-2}\int_{\p D_{\sigma_k}} \big|\h{1pt}U_n - J_n \h{1pt}\big|^2 \h{2pt}\lesssim_b\h{2pt}& \h{2pt}s_n^{-2}\int_{\p D_{\sigma_k}}\Big|\h{1pt}\big|\h{1pt}U_n  \h{1pt}\big|^2 -  1\h{1.5pt}\Big|^2  \\[2mm]
&+ \int_0^{2 \pi} \left|\h{2pt} \int_0^\varphi \Big[  \h{1pt}\big< U_n, \h{1.5pt} \xi^{\bot}\cdot D_{\xi} U^{sc}_n\big> \h{1pt}\Big]\big(\sigma_k \cos\phi, \h{1.5pt}\sigma_k \sin \phi\big) \h{2pt}\mathrm{d} \phi \h{2pt} \right|^2 \h{2pt}\mathrm{d} \varphi.
\end{align*} In light of (1) in Lemma \ref{uniform convergence of u_n in S_sigma_k in large scale}  and the convergence of $U_n$ to $y_*$ in $C^0\big(\p D_{\sigma_k}\big)$, in the current large--scale regime, we can take $n \to \infty$ on both sides above and obtain \begin{align}\label{est, b and c,} \limsup_{n \to \infty} s_n^{-2}\int_{\p D_{\sigma_k}} \big|\h{1pt}U_n - J_n \h{1pt}\big|^2 &\h{2pt}\lesssim_b\h{2pt}\int_0^{2 \pi} \left|\h{2pt} \int_0^\varphi \Big[  \h{1pt}\big< y_*, \h{1.5pt} \xi^{\bot}\cdot D_{\xi} U^{sc}_\infty\big> \h{1pt}\Big]\big(\sigma_k \cos\phi, \h{1.5pt}\sigma_k \sin \phi\big) \h{2pt}\mathrm{d} \phi \h{2pt} \right|^2 \h{2pt}\mathrm{d} \varphi \\[2mm]
&\h{2pt}=\h{2pt}\int_0^{2 \pi} \Big|\h{2pt} \big< y_*, \h{1.5pt}  U^{sc}_\infty\big(\sigma_k \cos \varphi, \sigma_k \sin \varphi\big) - v_*\big> - \big< y_*, \h{1.5pt}  U^{sc}_\infty\big(\sigma_k, 0\big) - v_*\big>   \h{2pt} \Big|^2 \h{2pt}\mathrm{d} \varphi. \nonumber\end{align}Since $U_\infty^{sc}  \in H^1\big( D_{\sigma_k}; v_* + \mathrm{Tan}_{y_*}\mathbb{S}^2\big)$, then by trace theorem and the absolute continuity of $U_\infty^{sc}$ on $\p D_{\sigma_k}$, the integrand in the last integral of (\ref{est, b and c,}) equals  $0$. (\ref{conv of Un - Jn, case A}) therefore follows if Case B in (\ref{definition of J_n}) holds. The Case C in (\ref{definition of J_n}) can be treated by the same arguments as Case B.\vspace{0.1pc}

Now we apply (\ref{est of first int})--(\ref{conv of Un - Jn, case A}) to the right--hand side of (\ref{com of partial tau h}) and obtain 
\begin{eqnarray}
\lim_{R \to \infty} \lim_{n \rightarrow \infty} \h{2pt} s_n^{-2} \int_{D_{\sigma_k}  \setminus D_{(1-s)\h{1pt}\sigma_k} }   \big|\h{1pt}\partial_\tau \h{1pt} h_{n, s, R} \h{1pt}\big|^2   \h{2pt}=\h{2pt}0.
\label{conv. of radial, large}
\end{eqnarray}
\noindent \textbf{Step 4.2.2.} Still by the definition of $h_{n, s, R}$ in (\ref{definition of h_n in large scale}), the angular derivative of $h_{n, s, R}$ can be estimated by
\begin{eqnarray*}
\int_{D_{\sigma_k} \setminus D_{ ( 1-s ) \h{1pt} \sigma_k} } \tau^{-2}\h{1pt} \big| \h{1pt} \partial_\varphi h_{n, s, R} \h{1pt} \big|^2 \h{2pt}   \h{2pt}\lesssim  \h{2pt}  s_n^2 \int_{ (1-s) \h{1pt} \sigma_k }^{ \sigma_k } \tau^{-1} \mathrm{d} \tau \int_{ \p D_{\sigma_k} } \big| \h{1pt} D_\xi U^{sc}_\infty \h{1pt}\big|^2 + \big|\h{1pt}D_\xi U^{sc}_n \h{1pt}\big|^2.
\end{eqnarray*}
Here we take $R > \big\| U_\infty^{sc} - Z_*\big\|_{\infty; \p D_{\sigma_k}}$ and $n$ suitably large. Therefore, (1) in Lemma \ref{uniform convergence of u_n in S_sigma_k in large scale} induces \begin{eqnarray}\label{conv angular} s_n^{-2} \int_{D_{\sigma_k}  \setminus D_{(1-s)\h{1pt}\sigma_k}} \tau^{-2} \big| \partial_\varphi \h{1pt} h_{n, s, R} \big|^2
\h{2pt}\lesssim\h{2pt} b_k \h{1pt} \log \left(\dfrac{1}{1-s}\right).
\end{eqnarray}Note that (\ref{conv angular}) holds for all $R > \big\| U_\infty^{sc} - Z_*\big\|_{\infty; \p D_{\sigma_k}}$,  $n$ suitably large and $\sigma_k > 1/2$.\\

\noindent \textbf{Step 4.2.3.} Applying (\ref{conv. of radial, large})--(\ref{conv angular}) to the right--hand side of (\ref{decomp of D_k / D_(1-s), large}) and noticing (\ref{unif conv of hnsr}), we obtain
\begin{align}\label{conv, large, second} \lim_{s \to 0} \lim_{R \to \infty}\h{1pt} \lim_{n \rightarrow \infty} \h{2pt} s_n^{-2}  \int_{D_{\sigma_k} \setminus D_{(1-s)\sigma_k}} \Big| D_\xi \h{1pt}\Pi_{\h{0.5pt}\mathbb{S}^2} \big[ h_{n, s, R} \big] \Big|^2  \h{2pt}=\h{2pt} 0 \h{15pt}\text{for any $\sigma_k$ satisfying $\sigma_k > 1/2$.}\end{align}
\textbf{Step 4.3.} In this step, we consider the integral of $D_\xi \overline{v}_{n, s, R}$ in (\ref{decomp. of D_k in large scale}). By the definition of $\overline{v}_{n, s, R}$ in (\ref{definition of v_n in large scale}), \begin{eqnarray*}
\int_{D_{\sigma_k} \setminus D_{(1-s)\h{1pt}\sigma_k}} \big| \p_\tau \overline{v}_{n, s, R} \big|^2
= (s \h{1pt}\sigma_k)^{-2} \int_{D_{\sigma_k} \setminus D_{(1-s)\h{1pt}\sigma_k}} \big| \h{1pt}U_n - J_n \h{1pt}\big|^2\big(\sigma_k\h{1pt}\widehat{\xi}\h{2pt}\big). \end{eqnarray*}Utilizing (\ref{conv of Un - Jn, case A}) induces
\begin{eqnarray*}
\lim_{n \to \infty} s_n^{-2 } \int_{D_{\sigma_k} \setminus D_{(1-s)\h{1pt}\sigma_k}}  \big| \p_\tau \overline{v}_{n, s, R} \big|^2 = 0.
\end{eqnarray*}
As for the angular derivative, still by (1) in Lemma \ref{uniform convergence of u_n in S_sigma_k in large scale}, we have
\begin{eqnarray*}
s_n^{-2}\int_{D_{\sigma_k} \setminus D_{(1-s)\h{1pt}\sigma_k}} \tau^{-2} \big| \p_\varphi \overline{v}_{n, s, R} \big|^2
\h{2pt}\lesssim\h{2pt}   \int_{ (1-s) \h{1pt} \sigma_k }^{ \sigma_k } \tau^{-1} \mathrm{d} \tau \int_{ \p D_{\sigma_k} } \big|\h{1pt}D_\xi U^{sc}_n \h{1pt}\big|^2 \h{2pt}\leq\h{2pt}b_k\h{1pt}\log \left(\dfrac{1}{1 - s}\right).\end{eqnarray*}
Here we take $n$ large and assume $\sigma_k > 1/2$. By the last two estimates, it follows
\begin{align}\label{conv, third, large}
 \limsup_{n \rightarrow \infty} \h{2pt} s_n^{-2}  \int_{D_{\sigma_k} \setminus D_{(1-s)\sigma_k}} \big| D_\xi \overline{v}_{n, s, R} \big|^2  \h{10pt}\text{is independent of $R$ and converges to $0$ as $s \to 0$.}
\end{align}In the last estimate, $\sigma_k > 1/2$. \vspace{0.3pc}\\
\noindent \textbf{Step 4.4.} Applying (\ref{conv. main part in large scale}), (\ref{conv, large, second}) and (\ref{conv, third, large}) to the right--hand side of (\ref{decomp. of D_k in large scale}), we get \begin{align}\label{con of dir energy}\lim_{s \to 0} \lim_{R \rightarrow \infty} \lim_{n \to \infty} s_n^{-2}\int_{D_{\sigma_k}} \big|D_{\xi} \overline{v}_{n, s, R}\big|^2 = \int_{D_{\sigma_k}} \big|\h{1pt}D_\xi v\h{1pt}\big|^2 \h{10pt}\text{for any $\sigma_k$ satisfying $\sigma_k > 1/2$.}
\end{align}

\noindent \textbf{Step 5.} We complete the proof in this step. In light of the energy--minimizing property of $u_n$, it turns out
\begin{eqnarray}&&\int_{D_{\sigma_k}} \left\{ \big|D_{\xi} U^{sc}_n\big|^2   + \left(\dfrac{\lambda_n\h{1pt}r_n}{s_n}\right)^2 G_{a_n, \mu}\big(\rho_{x_n} + \lambda_n\h{1pt}r_n\h{1pt}\xi_1,  \h{2pt}U_n  \h{3pt} \big) \right\}   \left(1 + \dfrac{\lambda_n \h{1pt}r_n}{\rho_{x_n}} \h{1pt}\xi_1\right)\nonumber\\ [2mm]
 &&\leq \int_{D_{\sigma_k}}\left\{  s_n^{-2 } \h{2pt} \big|D_{\xi} \overline{v}_{n, s, R}\big|^2 + \left(\dfrac{\lambda_n\h{1pt}r_n}{s_n}\right)^2 G_{a_n, \mu}\big(\rho_{x_n} + \lambda_n\h{1pt}r_n\h{1pt}\xi_1,  \h{2pt}\overline{v}_{n, s, R}  \h{3pt} \big)    \right\} \left(1 + \dfrac{\lambda_n \h{1pt}r_n}{\rho_{x_n}} \h{1pt}\xi_1\right).
\label{large scale basic inequality}
\end{eqnarray}
Utilizing (\ref{decom of potential})--(\ref{limit of J_2^s}) and lower--semi continuity, we obtain \begin{align}\label{low bod}4 \pi \h{0.5pt}c_3\h{0.5pt}\sigma_k^2 \h{2pt} + &\int_{D_{\sigma_k}} \big| \h{1pt}D_\xi U_\infty^{sc}\h{1pt}\big|^2 \\[2mm]
&\h{2pt}\leq \h{2pt}\liminf_{n \to \infty} \int_{D_{\sigma_k}} \left\{ \big|D_{\xi} U^{sc}_n\big|^2   + \left(\dfrac{\lambda_n\h{1pt}r_n}{s_n}\right)^2 G_{a_n, \mu}\big(\rho_{x_n} + \lambda_n\h{1pt}r_n\h{1pt}\xi_1,   U_n  \h{3pt} \big) \right\}   \left(1 + \dfrac{\lambda_n \h{1pt}r_n}{\rho_{x_n}} \h{1pt}\xi_1\right).
\nonumber\end{align}Recall $h_{n, s, R}$ defined in (\ref{definition of h_n in large scale}). By (\ref{pi_gamma - y_n}), (\ref{pi_L convergence in large scale}), (\ref{est of first int})--(\ref{conv of Un - Jn, case A}) and (1) in Lemma \ref{uniform convergence of u_n in S_sigma_k in large scale}, the $L^2\big(D_{\sigma_k}\big)$--norm of the mapping $\frac{h_{n, s, R} - y_n}{s_n}$ is bounded with the upper bound independent of $n$. Due to this uniform boundedness, (\ref{unif conv of hnsr}) and (\ref{ration convergence }), the  $L^2\big(D_{\sigma_k}\big)$--norm of the mapping $\frac{\Pi_{\h{0.5pt}\mathbb{S}^2}\big[\h{1pt}h_{n, s, R}\h{1pt}\big] - y_n}{s_n}$ is bounded with the upper bound independent of $n$. Noticing the definition of $\overline{v}_{n, s, R}$ in (\ref{definition of v_n in large scale}), we then can use the boundedness of $\frac{\Pi_{\h{0.5pt}\mathbb{S}^2}\big[\h{1pt}h_{n, s, R}\h{1pt}\big] - y_n}{s_n}$ in $L^2\big(D_{\sigma_k}\big)$, (\ref{conv of Un - Jn, case A}) and (1) in Lemma \ref{uniform convergence of u_n in S_sigma_k in large scale} to obtain  the uniform boundedness of $\frac{\overline{v}_{n, s, R} - y_n}{s_n}$ in $L^2\big(D_{\sigma_k}\big)$. Here the upper bound of the  $L^2\big(D_{\sigma_k}\big)$--norm of $\frac{\overline{v}_{n, s, R} - y_n}{s_n}$ might depend on $s$ and $R$ but independent of $n$. Now we can apply the same derivation for (\ref{limit of J_2^s}) to show that \begin{align*} \left(\dfrac{\lambda_n r_n}{s_n}\right)^2 \int_{D_{\sigma_k}} \dfrac{4 \h{1pt}\big[\h{1pt}\overline{v}_{n, s, R}\h{1pt}\big]_1^2 + \big[\h{1pt}\overline{v}_{n, s, R}\h{1pt}\big]_3^2}{\big(\rho_{x_n} + \lambda_n r_n \xi_1\big)^2} \longrightarrow 4\pi\h{0.5pt}c_3\h{0.5pt}\sigma_k^2 \h{15pt}\text{as $n \to \infty$.}
\end{align*}By this limit, (\ref{con of dir energy}) and (\ref{convergence of ginzburg-landau term}), it follows \begin{align*}
 \lim_{s \to 0} \lim_{R \to \infty} \lim_{n \to \infty} &\int_{D_{\sigma_k} }\left\{  s_n^{-2 } \h{2pt} \big|D_{\xi} \overline{v}_{n, s, R}\big|^2 + \left(\dfrac{\lambda_n\h{1pt}r_n}{s_n}\right)^2 G_{a_n, \mu}\big(\rho_n + \lambda_n\h{1pt}r_n\h{1pt}\xi_1, \overline{v}_{n, s, R} \big) \right\} \left(1 + \dfrac{\lambda_n \h{1pt}r_n}{\rho_n} \h{1pt}\xi_1\right) \nonumber\\[2mm]
&= 4 \pi \h{0.5pt}c_3\h{0.5pt}\sigma_k^2 + \int_{D_{\sigma_k}}  \big|\h{1pt}D_\xi v\h{1pt}\big|^2.
\end{align*}
The proof is completed by the above limit and (\ref{large scale basic inequality})--(\ref{low bod}).
\end{proof}
\begin{proof}[\bf Proof of Proposition \ref{variant decay lemma} in large--scale regime]\
\\[3mm]
Any mapping $u$ in $N_k$ can be expressed as $u = v_* + f_1 \h{1pt}t_* + f_2 \h{1pt}e^*_3$ for some $f_1$ and $f_2$ in $H^1\big(D_{\sigma_k}\big)$. In light of Lemma \ref{strong H1 convergence in large case}, the proof follows by the similar arguments as the proof for the small--scale regime in Section 3.2. We omit it here.
\end{proof}

\section{Emptiness of the coincidence set \vspace{0.5pc}}

\noindent  Recalling the sequence $\big\{ w_{a_n, b}^+\big\}$ obtained in Step 1 of Section 1.4.2, in this section, we show \begin{prop}\label{prop 4.1} If $a_n$ is large, then there is no $x_n \in T$ satisfying (\ref{contradiction to violate}). 
\end{prop}
\noindent With this proposition, we obtain the existence of biaxial--ring solutions in Theorem \ref{biaxial--ring solution} for $b$ satisfying (\ref{cond of b parameter}) and  $a$ large. The existence of split--core solutions in Theorem \ref{split--core solution} can be obtained by similar arguments used here.
\begin{proof}[\bf Proof of Proposition \ref{prop 4.1}] We divide the proof into five steps. In the following, $w_{a_n, b}^+ = \mathscr{L}\big[ u_{a_n, b}^+\big]$.\vspace{0.4pc}

\noindent \textbf{Step 1.} Fix a small $\epsilon_0 > 0$ and $r \in (0, \epsilon_0)$. Moreover, we assume that $B_{r}(x)$ lies in $\mathscr{J}$ and satisfies \begin{align}\label{small energy condition for uanb} E_{a_n, \mu;\h{1pt}x, 2^{-2}r}\big[ u_{a_n, b}^+\big] < \epsilon_1.
\end{align}Here $\epsilon_1$ is given in Proposition \ref{variant decay lemma}.  Defining $\lambda_0 := \lambda\h{0.5pt}\theta_0$, we then can apply Proposition \ref{variant decay lemma} to obtain \begin{align*} E_{a_n, \mu;\h{1pt}x, 2^{-2}\lambda_0 \h{0.2pt} r}\big[ u_{a_n, b}^+\big] \h{2pt}\leq\h{2pt}  \frac{1}{2}\h{1pt} E_{a_n, \mu; \h{1pt}x, 2^{-2}\lambda \h{0.3pt}r}\big[ u_{a_n, b}^+\big] + \left(\dfrac{r}{4}\right)^{3/2}, \h{15pt}\text{for any $a_n > a_0$.}\end{align*} Now we take  $\epsilon_0$ sufficiently small (depending on $\epsilon_1$). The last energy--decay estimate and (\ref{small energy condition for uanb}) infer  \begin{align*}E_{a_n, \mu; \h{1pt}x, 2^{-2}\lambda_0 \h{0.2pt} r}\big[ u_{a_n, b}^+\big] \h{2pt}\leq\h{2pt} \dfrac{\epsilon_1}{2} + \left(\frac{\epsilon_0}{4}\right)^{3/2} \h{1.5pt}<\h{1.5pt}\epsilon_1.
\end{align*}Inductively we suppose that $E_{a_n, \mu; \h{1pt}x, 2^{-2}\lambda_0^k\h{0.2pt}r}\big[ u_{a_n, b}^+\big] < \epsilon_1$ for some $k \in \mathbb{N}$. By Proposition \ref{variant decay lemma}, it follows  \begin{align}\label{energy---decay induct} E_{a_n, \mu;\h{1pt}x, 2^{-2}\lambda_0^{k+1} \h{0.2pt} r}\big[ u_{a_n, b}^+\big] \h{2pt}\leq\h{2pt}  \frac{1}{2} \h{1pt}E_{a_n, \mu; \h{1pt}x, 2^{-2}\lambda_0^k \h{0.2pt}r}\big[ u_{a_n, b}^+\big] + \left(\dfrac{\lambda_0^k\h{0.5pt}r}{4}\right)^{3/2} \h{1.5pt}<\h{1.5pt} \frac{\epsilon_1}{2} + \left(\frac{\epsilon_0}{4}\right)^{3/2}\h{1.5pt}<\h{1.5pt} \epsilon_1. \end{align} Hence, the last estimate holds for any $k \in \big\{0\big\} \cup \mathbb{N}$. With a standard iteration argument, it  yields \begin{align}\label{energy decay of mathcal E, not origin} E_{a_n, \mu;\h{1pt}x, s}\big[ u_{a_n, b}^+\big] \h{1.5pt}\lesssim\h{1.5pt} \left(\dfrac{s}{r}\right)^{\alpha_0}, \h{15pt}\text{where $s \in \left(0, \frac{r}{4}\h{1.5pt}\right]$ and $\alpha_0 = - \dfrac{\ln 2}{\ln \lambda_0} \in (0, 1)$.} 
\end{align}In light of the definition of $\mathcal{E}_{a, \mu; \h{0.5pt}x, r}\big[ w_{a, b}^+\big]$ in (\ref{small energy cond, off 0}), we then obtain from (\ref{energy decay of mathcal E, not origin}) the estimate: \begin{align}\label{energy decay of e, not origin, x coordinate} \mathcal{E}_{a_n,\h{0.3pt}\mu; \h{0.3pt}x, \h{0.2pt}s} \big[ w_{a_n, b}^+\big] &\h{2pt}\leq\h{2pt} 2\left(1 + \dfrac{\rho_x}{s}\right) \arcsin \left(\dfrac{s}{\rho_x}\right) \int_{D_s\left(\rho_x, 0\right)} e_{a_n, \mu}\big[ u_{a_n, b}^+\big]  \h{1.5pt}\lesssim\h{1.5pt} \left(\dfrac{s}{r}\right)^{\alpha_0},\nonumber\\[2mm]
&\h{2pt} \h{15pt}\text{for any $a_n > a_0$, $B_r(x) \in \mathscr{J}$ satisfying (\ref{small energy condition for uanb}), $r \in (0, \epsilon_0)$ and $s \in \left(0, \frac{r}{4}\h{1.5pt}\right]$.} 
\end{align}

Recall $r_{\epsilon_0}$ given in Step 2 of Section 1.4.2. Now we take $r = r_{\epsilon_0}$ and $x = 0$ in  (\ref{energy density small approx mapping}). It then follows \begin{align*}\mathcal{E}_{a_n, \mu; \h{1pt}0,  r_{\epsilon_0}}\big[ w_{a_n, b}^+\big] < \epsilon_0 \h{15pt}\text{ for $n$ suitably large.}\end{align*} Let $\epsilon_0 < \epsilon_1$ where $\epsilon_1$ is as in Proposition \ref{small energy implies energy decay}. Using the similar derivations for (\ref{energy---decay induct}),  we obtain with an use of Proposition  \ref{small energy implies energy decay}  that \begin{align*} \mathcal{E}_{\h{0.2pt}a_n, \h{0.2pt}\mu;  \h{1pt}0, \nu_0^{k+1} r_{\epsilon_0}}\big[ w_{a_n, b}^+\big] \h{2pt}\leq\h{2pt} \dfrac{1}{2}\h{1pt}\mathcal{E}_{a_n,  \mu; \h{1pt}0, \nu_0^k r_{\epsilon_0}}\big[ w_{a_n, b}^+\big] + \left(\nu_0^k \h{0.6pt}r_{\epsilon_0}\right)^{3/2} \h{1.5pt}<\h{1.5pt}\epsilon_1  \h{15pt}\text{for any  $k \in \big\{0\big\} \cup \mathbb{N}$.} 
\end{align*}Here $a_n$ is  large and $\epsilon_0$ is taken small. This estimate and standard iteration argument then induce \begin{align}\label{energy decay of mathcal E} \mathcal{E}_{a_n, \mu;  \h{1pt}0, s}\big[ w_{a_n, b}^+\big] \h{1.5pt}\lesssim\h{1.5pt} \nu_0^{-1} \left(\dfrac{s}{r_{\epsilon_0}}\right)^{\alpha_1}, \h{10pt}\text{where $a_n$ is large, $s \in (0, r_{\epsilon_0}]$ and $\alpha_1 = - \dfrac{\ln 2}{\ln \nu_0} \in (0, 1)$.} 
\end{align}
\noindent \textbf{Step 2.} Let $B'_s\left(x\right) := B_s\left(x\right) \cap T$ and  $C_s\left(x\right) := B'_{s/2}\left(x \right) \times \big(\h{0.5pt}0,  \sqrt{3}\h{0.5pt} s \big/ 2  \h{0.5pt}\big)$. Since $C_s\left(x\right) \subset B_s\left(x\right)$, by  trace theorem, Poincar\'{e}'s inequality and (\ref{energy decay of e, not origin, x coordinate}), it turns out \begin{align}\label{trace campa with aver on ball}&s^{-2} \int_{B'_{s/2}\left(x\right)} \left|\h{2pt} w_{a_n, b}^+ - \dashint_{C_s\left(x\right)} w_{a_n, b}^+ \h{2pt}\right|^2 \h{2pt} \h{2pt}\lesssim\h{2pt}s^{-1} \int_{C_s\left(x\right)} \big| \nabla w_{a_n, b}^+ \big|^2 \h{2pt}\lesssim\h{2pt}\left(\dfrac{s}{r}\right)^{\alpha_0}, \nonumber\\[2mm] & \h{80pt}\text{for any $a_n > a_0$, $B_r(x) \in \mathscr{J}$ satisfying (\ref{small energy condition for uanb}),  $r \in (0, \epsilon_0)$ and $s \in \left(0, \frac{r}{4}\h{1.5pt}\right]$.}
\end{align}
For any $y \in C_s(x)$, \begin{align*} w_{a_n, b}^+\left(y\right) - w_{a_n, b}^+\left(y'\right) = \int_0^{y_3}\p_z w_{a_n, b}^+ \h{2pt}\mathrm{d}z. 
\end{align*}Here $y' = (y_1, y_2, 0)$. Integrating the above identity with respect to the $y$--variable over $C_s\left(x\right)$, we obtain \begin{align*} \dashint_{C_s\left(x\right)} w_{a_n, b}^+ - \dashint_{B'_{s/2}\left(x\right)}w_{a_n, b}^+ = \dashint_{C_s\left(x\right)} \int_0^{y_3} \p_z w_{a_n, b}^+\h{2pt}\mathrm{d} z.
\end{align*} By this equality and (\ref{energy decay of e, not origin, x coordinate}), it turns out \begin{align*} &\left|\h{2pt} \dashint_{C_s\left(x\right)} w_{a_n, b}^+ - \dashint_{B'_{s/2}\left(x\right)}w_{a_n, b}^+ \h{2pt}\right|^2\h{2pt}\lesssim\h{2pt}s^{-1} \int_{C_s\left(x\right)} \big| \p_z w_{a_n, b}^+ \big|^2 \h{2pt}\lesssim\h{2pt}\left(\dfrac{s}{r}\right)^{\alpha_0},\\[2mm]
&\h{40pt}\text{for any $a_n > a_0$, $B_r(x) \in \mathscr{J}$ satisfying  (\ref{small energy condition for uanb}),  $r \in (0, \epsilon_0)$ and $s \in \left(0, \frac{r}{4}\h{1.5pt}\right]$.}
\end{align*}Combining this estimate with (\ref{trace campa with aver on ball}), we obtain \begin{align}\label{trace campa with aver on flat }s^{-2} \int_{B'_{s/2}\left(x\right)} \left|\h{2pt} w_{a_n, b}^+ - \dashint_{B'_{s/2}\left(x\right)} w_{a_n, b}^+ \h{2pt}\right|^2  &\h{2pt}\lesssim  \h{2pt}\left(\dfrac{s}{r}\right)^{\alpha_0}, \nonumber\\[2mm]
&\h{-95pt}\text{for any $a_n > a_0$, $B_r(x) \in \mathscr{J}$ satisfying (\ref{small energy condition for uanb}), $r \in (0, \epsilon_0)$ and $s \in \left(0, \frac{r}{4}\h{1.5pt}\right]$.}
\end{align}Similarly, by (\ref{energy decay of mathcal E}), it satisfies \begin{align}\label{trace campa with aver on flat, origin}s^{-2} \int_{B'_{s/2}\left(0\right)} \left|\h{2pt} w_{a_n, b}^+ - \dashint_{B'_{s/2}\left(0\right)} w_{a_n, b}^+ \h{2pt}\right|^2  \h{2pt}\lesssim \h{2pt}\nu_0^{-1}\left(\dfrac{s}{r_{\epsilon_0}}\right)^{\alpha_1}, \h{15pt}\text{for any large $a_n$ and $s \in \big(0, r_{\epsilon_0} \big]$.}
\end{align}\

\noindent \textbf{Step 3.} In the following, we fix $x \in T$ and assume  $B_{r}(x) \subset B_1$. Here $r \leq r_{\epsilon_0} < \epsilon_0$.  Note that,  $B_{r}(x)$ may have non--empty intersection with $l_z$. Now we estimate the H\"{o}lder norm of $w_{a_n, b}^+$ on $B'_{\sigma\h{0.3pt}r}(x)$. Here $\sigma \in (0, 1)$ is a small constant. It will be determined during the course of the proof. Throughout the following arguments, we always take $a_n$ large enough when it is necessary. \vspace{0.3pc}

Letting $y \in B'_{\sigma\h{0.3pt}r}(x)$ and $\rho \in \big(0,  2 \sigma\h{0.3pt}r\big) $, we divide our proof into four cases. \vspace{0.3pc}\\
\noindent \textbf{Case 1.} If $y = 0$, then we take $s = 2 \rho$ in (\ref{trace campa with aver on flat, origin}). It holds \begin{align*}  \rho^{-2} \int_{B'_{\rho}\left(0\right)} \left|\h{2pt} w_{a_n, b}^+ - \dashint_{B'_{\rho}\left(0\right)} w_{a_n, b}^+ \h{2pt}\right|^2  \h{2pt}\lesssim_{\h{1pt}\nu_0} \h{2pt} \left(\frac{\rho}{r_{\epsilon_0}}\right)^{\alpha_1}.
\end{align*}Here, we need $\sigma < 1/4$ so that $2 \rho < 4 \sigma r_{\epsilon_0} < r_{\epsilon_0}$.\vspace{0.4pc}

\noindent \textbf{Case 2.} In this case, we assume $y \neq 0$ and $\rho_y > 2^4 \sigma\h{0.5pt}r_{\epsilon_0}$. Letting $\sigma < 2^{-4}$ and taking $r = 2^4 \sigma \h{0.3pt}r_{\epsilon_0}$ in (\ref{energy density small approx mapping}), for large $n$, we can derive from (\ref{energy density small approx mapping}) the following small--energy condition: \begin{align*}\mathcal{E}_{a_n, \mu; \h{1pt}y, 2^4 \sigma\h{0.3pt}r_{\epsilon_0}}\big[ w_{a_n, b}^+\big] \h{1.5pt}<\h{1.5pt} \epsilon_0.
\end{align*}Since $B_{2^4 \sigma\h{0.5pt}r_{\epsilon_0}}(y) \cap l_z = \emptyset$, we can replace $x$, $r$, $a$ in (\ref{small energy condition rho z coordinate, simplified}) with $y$, $2^4 \sigma\h{0.3pt}r_{\epsilon_0}$, $a_n$, respectively. It follows \begin{align*}\int_{D_{4 \sigma \h{0.3pt}r_{\epsilon_0}}(\rho_y, 0) \h{1pt}\cap\h{1pt}\mathbb{D}} e_{a_n, \mu} \big[ u_{a_n, b}^+\big]\h{1.5pt}\lesssim\h{1.5pt}\epsilon_0.
\end{align*}Note that $2^4 \sigma\h{0.3pt}r \leq 2^4 \sigma\h{0.3pt}r_{\epsilon_0}$ and $B_{2^4 \sigma\h{0.5pt}r}(y) \cap l_z = \emptyset$. In addition, it holds $B_{2^4 \sigma\h{0.5pt}r}(y) \subset B_{2^5 \sigma\h{0.5pt}r}(x) \subset B_r(x) \subset B_1$ if $\sigma < 2^{-5}$. We therefore can use this set relationship and the last estimate to get \begin{align}\label{cond eneded}B_{2^4 \sigma\h{0.5pt}r}(y) \in \mathscr{J}\h{15pt}\text{and}\h{15pt}E_{a_n, \mu; \h{1pt}y, 4 \sigma\h{0.3pt}r}\big[ u_{a_n, b}^+\big] \h{1.5pt}<\h{1.5pt}\epsilon_1.
\end{align}Here we have taken $\epsilon_0$ suitably small (depending on $\epsilon_1$). In light of (\ref{cond eneded}) and $\rho < 2 \sigma\h{0.2pt}r$,   we now replace $r$, $s$, $x$ in (\ref{trace campa with aver on flat }) with $2^4 \sigma \h{0.3pt}r$, $2 \rho$ and $y$, respectively. It then follows \begin{align*}\rho^{-2} \int_{B'_{\rho}\left(y\right)} \left|\h{2pt} w_{a_n, b}^+ - \dashint_{B'_{\rho}\left(y\right)} w_{a_n, b}^+ \h{2pt}\right|^2  &\h{2pt}\lesssim \h{2pt}\left(\dfrac{\rho}{\sigma r}\right)^{\alpha_0}.
\end{align*}

\noindent\textbf{Case 3.} In this case, we assume that  $B_{2^3\rho}(y) \cap l_z \neq \emptyset$. Firstly, \begin{align*}  \int_{B'_{\rho}\left(y\right)} \left|\h{2pt} w_{a_n, b}^+ - \dashint_{B'_{\rho}\left(y\right)} w_{a_n, b}^+ \h{2pt}\right|^2 \h{2pt}\leq\h{2pt}\int_{B'_{\rho}\left(y\right)} \left|\h{2pt} w_{a_n, b}^+ - \dashint_{B'_{2^4 \rho}(0)} w_{a_n, b}^+ \h{2pt}\right|^2 \h{2pt}\leq\h{2pt}\int_{B'_{2^3 \rho}\left(y\right)} \left|\h{2pt} w_{a_n, b}^+ - \dashint_{B'_{2^4 \rho}(0)} w_{a_n, b}^+ \h{2pt}\right|^2.
\end{align*}On the other hand,  $B_{2^3\rho}(y) \cap l_z \neq \emptyset$ yields $0 \in B_{2^3\rho}(y)$. Hence, $B'_{2^3\rho}(y) \subset B'_{2^4 \rho}(0)$. This set relationship and the last estimate then induce \begin{align*}  \int_{B'_{\rho}\left(y\right)} \left|\h{2pt} w_{a_n, b}^+ - \dashint_{B'_{\rho}\left(y\right)} w_{a_n, b}^+ \h{2pt}\right|^2  \h{2pt}\leq\h{2pt}\int_{B'_{2^4 \rho}\left(0\right)} \left|\h{2pt} w_{a_n, b}^+ - \dashint_{B'_{2^4 \rho}(0)} w_{a_n, b}^+ \h{2pt}\right|^2.
\end{align*}Note that $2^5 \rho < 2^6 \sigma \h{0.3pt}r_{\epsilon_0} < r_{\epsilon_0}$, provided that $\sigma < 2^{-6}$. We can take $s = 2^5\rho$ in (\ref{trace campa with aver on flat, origin}) to get \begin{align*} \rho^{-2}\int_{B'_{2^4 \rho}\left(0\right)} \left|\h{2pt} w_{a_n, b}^+ - \dashint_{B'_{2^4 \rho}\left(0\right)} w_{a_n, b}^+ \h{2pt}\right|^2  \h{2pt}\lesssim_{\h{1pt}\nu_0} \h{2pt}  \left(\frac{\rho}{r_{\epsilon_0}}\right)^{\alpha_1}.
\end{align*}Combining the last two estimates infers \begin{align*}   \rho^{-2}\int_{B'_{\rho}\left(y\right)} \left|\h{2pt} w_{a_n, b}^+ - \dashint_{B'_{\rho}\left(y\right)} w_{a_n, b}^+ \h{2pt}\right|^2  \h{2pt}\lesssim_{\h{1pt}\nu_0}\h{2pt}\left(\frac{\rho}{r_{\epsilon_0}}\right)^{\alpha_1}.
\end{align*}

\noindent\textbf{Case 4.} In this last case, we suppose that $y \neq 0$ and $\rho_y \leq 2^4 \sigma \h{0.3pt}r_{\epsilon_0}$. Moreover, it satisfies \begin{align*}B_{2^{k} \rho}(y) \cap l_z = \emptyset \h{15pt}\text{ and } \h{15pt} B_{2^{k + 1} \rho}(y) \cap l_z \neq \emptyset, \h{20pt}\text{for some natural number $k \geq 3$.} \end{align*}By   $B_{2^{k} \rho}(y) \cap l_z = \emptyset$, we have \begin{align}\label{assumption satisfied by k}2^{k}\rho \leq \rho_y \leq 2^4 \sigma\h{0.3pt}r_{\epsilon_0}.\end{align} With $B_{2^{k + 1} \rho}(y) \cap l_z \neq \emptyset$, it follows $ |\h{0.5pt} y \h{0.5pt}| < 2^{k + 1} \rho$. This estimate together with (\ref{assumption satisfied by k}) induces \begin{align*} | \h{1pt}\eta \h{1pt}| \h{0.5pt}\leq\h{0.5pt} |\h{1pt} \eta - y \h{1pt}| + |\h{1pt} y \h{1pt}| \h{0.5pt}\leq\h{0.5pt} 2^{k} \rho + 2^{k+1}\rho \h{0.5pt}\leq\h{0.5pt}  2^{k + 2} \rho   \h{0.5pt}\leq\h{0.5pt} 2^6 \sigma\h{0.5pt} r_{\epsilon_0}, \h{15pt}\text{for any $\eta \in B_{2^{k} \rho}(y)$.}
\end{align*}Taking $\epsilon_0$ small enough then infers $B_{2^{k}\rho}(y) \in \mathscr{J}$.  Suppose that \begin{align}\label{assumptio use}  \left(2^{k}\rho\right)^{-2} \int_{B'_{2^{k}\rho}\left(y\right)} \left|\h{2pt} w_{a_n, b}^+ - \dashint_{B'_{2^{k}\rho}\left(y\right)} w_{a_n, b}^+ \h{2pt}\right|^2  \h{2pt}\leq\h{2pt}2^{-2k} \left(\frac{\rho}{r_{\epsilon_0}}\right)^{\frac{\alpha_1}{2}}.
\end{align}In addition, we can show \begin{align*}  \int_{B'_{\rho}\left(y\right)} \left|\h{2pt} w_{a_n, b}^+ - \dashint_{B'_{\rho}\left(y\right)} w_{a_n, b}^+ \h{2pt}\right|^2 \h{1pt}\leq\h{1pt}\int_{B'_{\rho}\left(y\right)} \left|\h{2pt} w_{a_n, b}^+ - \dashint_{B'_{2^{k} \rho}(y)} w_{a_n, b}^+ \h{2pt}\right|^2  \h{1pt}\leq\h{1pt}\int_{B'_{2^{k}\rho}\left(y\right)} \left|\h{2pt} w_{a_n, b}^+ - \dashint_{B'_{2^{k} \rho}(y)} w_{a_n, b}^+ \h{2pt}\right|^2.
\end{align*}The last estimate and (\ref{assumptio use}) yield \begin{align*} \rho^{-2}\int_{B'_{\rho}\left(y\right)} \left|\h{2pt} w_{a_n, b}^+ - \dashint_{B'_{\rho}\left(y\right)} w_{a_n, b}^+ \h{2pt}\right|^2 \h{2pt}\leq\h{2pt}\dfrac{2^{2k}}{\left(2^{k}\rho\right)^2}\int_{B'_{2^{k}\rho}\left(y\right)} \left|\h{2pt} w_{a_n, b}^+ - \dashint_{B'_{2^{k} \rho}(y)} w_{a_n, b}^+ \h{2pt}\right|^2 \h{2pt}\leq\h{2pt}\left(\frac{\rho}{r_{\epsilon_0}}\right)^{\frac{\alpha_1}{2}}.
\end{align*} 

Now we assume (\ref{assumptio use}) fails. Note that $B_{2^{k + 1} \rho}(y) \cap l_z \neq \emptyset$. Hence, $B_{2^{k}\rho}(y) \subset B_{2^{k + 1} \rho}(y) \subset B_{2^{k+2}\rho}$.  By this set relationship, it follows \begin{align*}\int_{B'_{2^{k}\rho}\left(y\right)} \left|\h{2pt} w_{a_n, b}^+ - \dashint_{B'_{2^{k}\rho}\left(y\right)} w_{a_n, b}^+ \h{2pt}\right|^2 &\h{1pt}\leq\h{1pt}\int_{B'_{2^{k}\rho}\left(y\right)} \left|\h{2pt} w_{a_n, b}^+ - \dashint_{B'_{2^{k +2} \rho}(0)} w_{a_n, b}^+ \h{2pt}\right|^2 \\[1.5mm]
& \h{1pt}\leq\h{1pt}\int_{B'_{2^{k + 2}\rho}\left(0\right)} \left|\h{2pt} w_{a_n, b}^+ - \dashint_{B'_{2^{k+2} \rho}(0)} w_{a_n, b}^+ \h{2pt}\right|^2.
\end{align*}  Using (\ref{assumption satisfied by k}), we obtain $ B_{2^{k+3}\rho} \subset B_{r_{\epsilon_0}} $, provided that $\sigma < 2^{- 7}$. By taking $s = 2^{k+3}\rho$ in (\ref{trace campa with aver on flat, origin}), it holds \begin{align*}\left(2^{k+3}\rho\right)^{-2} \int_{B'_{2^{k+2}\rho}\left(0\right)} \left|\h{2pt} w_{a_n, b}^+ - \dashint_{B'_{2^{k+2}\rho}\left(0\right)} w_{a_n, b}^+ \h{2pt}\right|^2  \h{2pt}\lesssim_{\h{1pt}\nu_0} \h{2pt}\left(\dfrac{2^{k+3}\rho}{r_{\epsilon_0}}\right)^{\alpha_1}.
\end{align*} The last two estimates infer \begin{align*}\left(2^{k}\rho\right)^{-2} \int_{B'_{2^{k}\rho}\left(y\right)} \left|\h{2pt} w_{a_n, b}^+ - \dashint_{B'_{2^{k}\rho}\left(y\right)} w_{a_n, b}^+ \h{2pt}\right|^2 \h{2pt}\lesssim_{\h{1pt}\nu_0}\h{2pt}\left(\dfrac{2^{k+3}\rho}{r_{\epsilon_0}}\right)^{\alpha_1}.
\end{align*} If  (\ref{assumptio use}) fails, the above estimate gives us \begin{align}\label{low bd of k}   \left(\frac{r_{\epsilon_0}}{\rho}\right)^{\frac{\alpha_1}{2}} \h{2pt}\lesssim_{\h{1pt}\nu_0}\h{2pt}  2^{3k}.
\end{align}On the other hand, (\ref{assumption satisfied by k}) induces $2 \rho_y \leq 2^5 \sigma \h{0.3pt}r_{\epsilon_0} < r_{\epsilon_0}$. Hence, $B_{\rho_y}(y) \subset B_{2 \rho_y} \subset B_{r_{\epsilon_0}}$. We then can replace $s$ in (\ref{energy decay of mathcal E}) with $2 \rho_y$ and obtain \begin{align*} \mathcal{E}_{a_n, \mu; \h{1pt}y, \rho_y}\big[ w_{a_n, b}^+\big] \h{1.5pt}\leq\h{1.5pt}2 \mathcal{E}_{a_n, \mu; \h{1pt}0, 2 \rho_y} \big[ w_{a_n, b}^+ \big] \h{1.5pt}\lesssim_{\h{1pt}\nu_0}\h{1.5pt} \left( \frac{2 \rho_y}{r_{\epsilon_0}}\right)^{\alpha_1} \h{1.5pt}\lesssim_{\h{1pt}\nu_0}\h{1.5pt}\sigma^{\alpha_1}.
\end{align*}The last estimate above have used (\ref{assumption satisfied by k}) again. Taking $\sigma$ small enough (depending on $\nu_0$ and $\epsilon_1$), we have from the above estimate  and  (\ref{small energy condition rho z coordinate, simplified}) the small--energy condition: $$E_{a_n, \mu; \h{1pt}y, 2^{-2}\rho_y}\big[ u_{a_n, b}^+\big] < \epsilon_1.$$   Notice that $B_{\rho_y}(y) \in \mathscr{J}$. We then can replace $x$, $r$  and $s$ in (\ref{trace campa with aver on flat }) with $y$, $\rho_y$  and $2 \rho$, respectively. It follows \begin{align*} \rho^{-2} \int_{B'_{\rho}\left(y\right)} \left|\h{2pt} w_{a_n, b}^+ - \dashint_{B'_{\rho}\left(y\right)} w_{a_n, b}^+ \h{2pt}\right|^2  &\h{2pt}\lesssim \h{2pt} \left(\frac{2 \rho}{\rho_y}\right)^{\alpha_0} \h{1.5pt}\lesssim\h{1.5pt} 2^{-  k \h{0.5pt}\alpha_0}. 
\end{align*}In the above derivation, we also have used (\ref{assumption satisfied by k}) and $k \geq 3$. The last estimate and (\ref{low bd of k}) finally give us \begin{align*} \rho^{-2} \int_{B'_{\rho}\left(y\right)} \left|\h{2pt} w_{a_n, b}^+ - \dashint_{B'_{\rho}\left(y\right)} w_{a_n, b}^+ \h{2pt}\right|^2 \h{2pt}\lesssim_{\h{1pt}\lambda_0, \nu_0}\h{2pt} \left(\dfrac{\rho}{r_{\epsilon_0}}\right)^{  \frac{\alpha_0 \alpha_1}{6}}.
\end{align*}

Based on the arguments in the above four cases, we conclude that \begin{align*}\rho^{-2} \int_{B'_{\rho}\left(y\right)} \left|\h{2pt} w_{a_n, b}^+ - \dashint_{B'_{\rho}\left(y\right)} w_{a_n, b}^+ \h{2pt}\right|^2 \h{2pt}\lesssim_{\h{1pt}\lambda_0, \nu_0}\h{2pt} \left(\dfrac{\rho}{\sigma\h{0.3pt}r}\right)^{  \frac{\alpha_0 \alpha_1}{6}}, \h{15pt}\text{for any $y \in B_{\sigma\h{0.3pt}r}'(x)$ and $\rho \in (0, 2 \sigma\h{0.3pt}r)$.}
\end{align*}We then can apply Morrey--Campanato type estimate to get \begin{align}\label{loc unifor holder} \Big|\h{1.5pt}w_{a_n, b}^+\left(z_1\right) - w_{a_n, b}^+\left(z_2\right) \Big| \h{2pt}\lesssim_{\h{1pt}\lambda_0, \nu_0}\h{2pt}  \big(\sigma \h{0.3pt}r\big)^{- \beta} \h{1pt}\big| z_1 - z_2 \big|^\beta, \h{15pt}\text{for any $z_1$ and $z_2$ in $B'_{2^{-1}\sigma \h{0.3pt}r}(x)$.}
\end{align} Here we simply use $\beta$ to denote $\frac{\alpha_0\alpha_1}{12}$. Note that the above estimate holds for all $B_{r}(x) \subset B_1$ with $x \in T$ and $r \leq r_{\epsilon_0}$. \vspace{0.4pc}

\noindent \textbf{Step 4.} In this step, we show that there exist $\delta_0 > 0$ sufficiently small and $a_0 > 0$ sufficiently large so that $x_n$ satisfying (\ref{contradiction to violate}) must be in $B'_{1 - \delta_0}(0)$, for all $a_n > a_0$. Recalling   $w_{a_n, b}^+ = \mathscr{L}\big[ u_{a_n, b}^+\big]$, we restrict our study on the $(\rho, z)$--plane by considering $u_{a_n, b}^+$.  \vspace{0.3pc}

If we take $\delta < 2^{-10} \h{0.5pt}r_{\epsilon_0}$, then $B_{2^{-1}\delta}\big(1 - \delta, 0, 0\big) \subset B_1$ with $2^{-1} \delta  < r_{\epsilon_0}$. By (\ref{loc unifor holder}), it turns out
\begin{align*}
   & \Big|\h{1pt} u_{a_n, b}^+(\rho_1, 0) - u_{a_n, b}^+(\rho_2, 0) \h{1pt}\Big| \h{2pt}\lesssim_{\h{1pt}\lambda_0, \nu_0}\h{2pt} \big(\sigma\h{0.3pt}\delta\h{0.3pt}\big)^{-\beta} \h{1pt}\big|\h{1pt} \rho_1 - \rho_2 \h{1pt}|^{\beta}, \nonumber\\[2mm]
&\h{30pt}\text{for large $a_n$, $\delta < 2^{-10}\h{0.3pt} r_{\epsilon_0}$ and $\rho_1, \rho_2 \in \big(1 - \delta - 2^{-2}\sigma\h{0.3pt}\delta, 1 - \delta + 2^{-2} \h{0.3pt}\sigma\delta\h{0.3pt}\big)$.}   
\end{align*} Hence, for sufficiently small $\epsilon > 0$, the last estimate yields \begin{align}\label{path 1}\Big|\h{1pt} u_{a_n, b}^+(\rho, 0) - u_{a_n, b}^+(1 - \delta, 0) \h{1pt}\Big| \h{2pt}\lesssim_{\h{1pt}\lambda_0, \nu_0}\h{2pt} \epsilon^{\beta}, \h{15pt}\text{for any $\rho \in \big(1 - \delta - \epsilon \h{1pt}\sigma\h{0.5pt}\delta, 1 - \delta + \epsilon \h{1pt}\sigma\h{0.5pt} \delta \h{0.3pt}\big)$.}
\end{align}Here $a_n$ is large and $\delta < 2^{-10}\h{0.3pt} r_{\epsilon_0}$.\vspace{0.3pc}

Taking $r = r_{\epsilon_0}$ and $x = e_1^* = (1, 0, 0)$ in (\ref{energy density small approx mapping}), we obtain $\mathcal{E}_{a_n, \mu; \h{1pt}e_1^*, r_{\epsilon_0}}\big[ w_{a_n, b}^+\big] < \epsilon_0$. Utilizing (\ref{small energy condition rho z coordinate, simplified}) then induces  \begin{align*} \int_{\mathbb{D} \h{1pt}\cap\h{1pt} D_{2^{-2} r_{\epsilon_0}}\left(1, 0\right)} \big| Du_{a_n, b}^+\big|^2 \h{2pt}\lesssim\h{2pt}\epsilon_0.
\end{align*}Let $(r, \varphi)$ be the polar coordinates on $\mathbb{D}$ with respect to the center  $0$. Moreover, we assume $\varphi \in (- \pi/2, \pi/2)$. If $\epsilon$ is sufficiently small, then the subset in $\mathbb{D}$ whose points satisfy $r \in ( 1 - \delta - \epsilon\h{1pt}\sigma\h{0.3pt}\delta, 1 - \delta + \epsilon\h{1pt}\sigma\h{0.3pt}\delta)$ and $\varphi \in (- \epsilon \h{0.3pt}\sigma\h{0.3pt}\delta, \epsilon \h{0.3pt}\sigma\h{0.3pt}\delta)$ is contained in $\mathbb{D} \h{1pt}\cap D_{2^{-2}r_{\epsilon_0}}(1, 0)$. The last estimate then induces
\begin{align*} \int_{1 -\delta - \epsilon\h{1pt}\sigma\h{0.3pt} \delta}^{1 -\delta + \epsilon\h{1pt}\sigma\h{0.3pt}\delta}
\int_{- \epsilon\h{1pt}\sigma\h{0.3pt} \delta}
^{ \h{2pt}\epsilon\h{1pt}\sigma\h{0.3pt} \delta}
\big|\h{1pt} D u_{a_n, b}^+\h{1pt} \big|^2 \h{2pt}r \h{1pt} \mathrm{d} \h{1pt} r \h{1pt}\mathrm{d} \h{1pt} \varphi \h{2pt}\leq\h{2pt} \int_{\mathbb{D}\h{1pt}\cap \h{1pt} D_{2^{-2}r_{\epsilon_0}}(1, 0)} \big|\h{1pt} D u_{a_n, b}^+\h{1pt} \big|^2
\h{2pt}\lesssim\h{2pt} \epsilon_0.\end{align*}
We can find a $\rho_1 \in \big(1 -\delta  - \epsilon\h{1pt}\sigma\h{0.3pt} \delta, 1 -\delta  + \epsilon\h{1pt}\sigma\h{0.3pt} \delta \h{1pt}\big)$ such that \begin{align*}\int_{- \epsilon\h{1pt}\sigma\h{0.3pt} \delta}
^{ \h{2pt}\epsilon\h{1pt}\sigma\h{0.3pt} \delta}
\big|\h{1pt} D u_{a_n, b}^+\h{1pt} \big|^2\left(\rho_1, \varphi\right) \h{2pt}\mathrm{d} \varphi \h{2pt}\leq\h{2pt}\dfrac{1}{\epsilon \h{1pt}\sigma\h{0.3pt}\delta \h{1pt}\big(1 - \delta\big)}\int_{1 -\delta - \epsilon\h{1pt}\sigma\h{0.3pt} \delta}^{1 -\delta + \epsilon\h{1pt}\sigma\h{0.3pt}\delta}
\int_{- \epsilon\h{1pt}\sigma\h{0.3pt} \delta}
^{ \h{2pt}\epsilon\h{1pt}\sigma\h{0.3pt} \delta}
\big|\h{1pt} D u_{a_n, b}^+\h{1pt} \big|^2 \h{2pt}r \h{1pt} \mathrm{d} \h{1pt} r \h{1pt}\mathrm{d} \h{1pt} \varphi.
\end{align*}The last two estimates then give us
\begin{align*}  \int_{- \epsilon\h{1pt}\sigma\h{0.3pt} \delta}
^{\h{2pt}\epsilon\h{1pt}\sigma\h{0.3pt}\delta}
\big|\h{1pt} \p_\varphi u_{a_n, b}^+\h{1pt} \big|^2\big(\rho_1, \varphi\big) \h{2pt}\mathrm{d} \h{1pt} \varphi
\h{2pt}\lesssim\h{2pt}  \epsilon_0  \h{1pt}\big(\epsilon\h{1pt}\sigma\h{0.3pt} \delta\h{1pt}\big)^{-1}.\end{align*}
For any $\varphi_0 \in ( - \epsilon\h{1pt}\sigma\h{0.3pt}\delta, \epsilon\h{1pt}\sigma\h{0.3pt} \delta)$, we obtain
\begin{equation}
\Big|\h{1pt} u_{a_n, b}^+\big(\rho_1,\varphi_0\big) - u_{a_n, b}^+(\rho_1,0)  \h{1pt}\Big|
\h{2pt}\leq \h{2pt} \big|\h{1pt} \varphi_0 \h{1pt}\big|^{1/2} \left( \int^{\varphi_0}_0 \big|\h{1pt} \p_\varphi u_{a_n, b}^+ \h{1pt}\big|^2\big(\rho_1, \varphi\big)\h{2pt} \mathrm{d} \h{1pt} \varphi\right)^{1/2}
\h{2pt}\lesssim\h{2pt}\epsilon_0^{1/2}.
\label{path 2}
\end{equation}

In the next, we note that the subset in $\mathbb{D}$ whose points satisfy $r \in ( \rho_1, 1 )$ and $\varphi \in (- \epsilon\h{1pt}\sigma\h{0.3pt} \delta, \epsilon\h{1pt}\sigma\h{0.3pt} \delta)$ is also contained in $\mathbb{D} \h{1pt}\cap D_{2^{-2}r_{\epsilon_0}}(1, 0)$.
It then follows
\begin{align*} \int_{\rho_1}^1 \int_{- \epsilon\h{1pt}\sigma\h{0.3pt}\delta}^{\h{2pt} \epsilon\h{1pt}\sigma\h{0.3pt} \delta}
\big|\h{1pt} \p_r u_{a_n, b}^+\h{1pt} \big|^2 \h{2pt}\mathrm{d} \h{1pt} r \h{1pt} \mathrm{d} \h{1pt} \varphi
 \h{2pt}\lesssim\h{2pt} \int_{\mathbb{D}\h{1pt}\cap\h{1pt} D_{2^{-2}r_{\epsilon_0}}(1, 0)} \big|\h{1pt}D u_{a_n, b}^+\h{1pt}\big|^2
\h{2pt}\lesssim\h{2pt} \epsilon_0.\end{align*}
There is $\varphi_1 \in ( - \epsilon\h{1pt}\sigma\h{0.3pt} \delta,  \epsilon\h{1pt}\sigma\h{0.3pt} \delta)$ such that
\begin{align*}
\int_{\rho_1}^1 \big|\h{1pt}\p_r u_{a_n, b}^+\h{1pt}\big|^2(r, \varphi_1) \h{1pt} \mathrm{d} \h{1pt} r
\h{2pt}\leq\h{2pt} \big(\epsilon \h{1pt}\sigma\h{0.3pt}\delta\big)^{-1} \int_{\rho_1}^1 \int_{- \epsilon\h{1pt}\sigma\h{0.3pt}\delta}^{\h{2pt} \epsilon\h{1pt}\sigma\h{0.3pt} \delta}
\big|\h{1pt} \p_r u_{a_n, b}^+\h{1pt} \big|^2 \h{2pt}\mathrm{d} \h{1pt} r \h{1pt} \mathrm{d} \h{1pt} \varphi \h{2pt}\lesssim\h{2pt}\epsilon_0 \big(\epsilon\h{1pt}\sigma\h{0.3pt}\delta\h{1pt}\big)^{-1}.
\end{align*}
By fundamental theorem of calculus, we obtain
\begin{equation}
\Big|\h{1pt} u_{a_n, b}^+(1, \varphi_1) - u_{a_n, b}^+(\rho_1, \varphi_1) \h{1pt}\Big|
\h{2pt}\leq \h{2pt} \sqrt{1- \rho_1} \left( \int_{\rho_1}^1 \big|\h{1pt} \p_r u_{a_n, b}^+\big|^2(r,\varphi_1)  \h{2pt}\mathrm{d} \h{1pt} r   \right)^{1/2}
\h{2pt}\lesssim\h{2pt} \left(\frac{\epsilon_0}{\epsilon\h{1pt}\sigma}\right)^{1/2}.
\label{path 3}
\end{equation}

Finally, we recall the boundary condition  (\ref{bdy cond of u_a, b plus}). Choosing $\delta$ small enough, for large $a_n$, we then have
\begin{equation}
\Big|\h{1pt} u_{a_n, b}^+(1,\varphi_1) - u_{a_n, b}^+(1,0) \h{1pt} \Big|\h{2pt}\leq \h{2pt} \epsilon_0.
\label{path 4}
\end{equation}

Combining the estimates in (\ref{path 1})--(\ref{path 4}), for $a_n $ sufficiently large and $\delta$ sufficiently small, we get
\begin{align*} \Big|\h{1pt} u_{a_n, b}^+(1 - \delta, 0) - U^*_{a_n}(1, 0) \h{1pt}\Big| = \Big| u_{a_n, b}^+(1-\delta, 0) - u_{a_n, b}^+(1, 0)\h{1pt} \Big| \h{2pt}\lesssim_{\h{1pt}\lambda_0, \nu_0}\h{2pt} \epsilon^\beta + \epsilon^{1/2}_0 + \left(\frac{\epsilon_0}{\epsilon\h{1pt}\sigma}\right)^{1/2}. \end{align*}We now take $\epsilon = \epsilon^{1/2}_0$ and let $\epsilon_0$ sufficiently small. The above estimate induces \begin{align*} \Big|\h{1pt} u_{a_n, b}^+(1 - \delta, 0) - U^*_{a_n}(1, 0) \h{1pt}\Big|   \h{2pt}\leq\h{2pt} \epsilon_0^\gamma,\h{15pt}\text{for some $\gamma \in (0, 1)$.} \end{align*}Hence, it follows
\begin{align}\label{small on annulus} \big[ u_{a_n, b}^+ \big]_2\big(1-\delta,0\big) \h{2pt}\geq\h{2pt}- \dfrac{1}{2} H_{a_n} -  \epsilon_0^\gamma. \end{align}Recall $b_0$ and $b$ used in Step 1 of Section 1.4.2. Now  we take $a_n$ large and let $\epsilon_0$ and $\delta_0$ be sufficiently small. (\ref{small on annulus}) induces that $ \big[ u_{a_n, b}^+ \big]_2\big(1-\delta,0\big) \h{2pt}>\h{2pt}(b_0 + b)/2$, for any $\delta < \delta_0$ and $a_n$ large.  \vspace{0.4pc}

\noindent \textbf{Step 5.} In light of (\ref{loc unifor holder}), we also have equicontinuity of $w_{a_n, b}^+$ on the closure of $B'_{1 - \delta_0}(0)$. It then turns out that $w_{a_n, b}^+$ converges to $w_b^+$ uniformly on the closure of $B'_{1 - \delta_0}(0)$ as $n \to \infty$. Since $\big[ w_b^+\big]_3 \geq b_0 > b$ on $T$, then for  $a_n$ large enough,  $\big[ w_{a_n, b}^+\big]_3 > (b_0 + b)/2$ on the closure of $ B_{1 - \delta_0}'(0)$. The proof is completed. 
\end{proof}



\cftaddtitleline{toc}{section}{\normalsize Part II: Structure of director fields near disclinations \vspace{0.5pc}}{}

\section{Singularities and structure of phase mapping}

Since this section, we begin to study the structures of disclinations of $w_{a, b}^+$ and $w_{a, c}^-$. Here for some fixed $b \in \mathrm{I}_-$ and $c \in (0, 1)$, $w_{a, b}^+$ and $w_{a, c}^-$ are biaxial--ring and split--core solutions respectively obtained from the previous sections. Particularly in this section, we consider the asymptotic behavior of the phase mapping $\Pi_{\mathbb{S}^4}\big[ w_{a, b}^+\big]$ and $\Pi_{\mathbb{S}^4}\big[ w_{a, c}^- \big]$ near their singularities on $l_z$.  Now we summarize the main results in this section. \begin{prop}\label{num and asy of singu} There exists a $a_0 > 0$ so that when $a > a_0$, the biaxial--ring solution $w_{a, b}^+$ has even number (the number might be $0$) of zeros on $l^+_z$ and the split--core solution  $w_{a, c}^-$ has odd number of zeros on $l^+_z$. In addition, the followings hold for the phase mappings: \begin{enumerate}\item[$\mathrm{(1)}.$] Let $z_{a, 1}^+$, ..., $z_{a, k_a}^+$ be the zeros of $w_{a, b}^+$ on $l_z^+$, where $k_a$ is the total number of zeros of $w_{a, b}^+$ on $l_z^+$. Moreover, the zeros are ordered so that for fixed $a$, $z_{a, j; 3}^+$ are increasing with respect to $j$. Then  \begin{align}\label{asy of pha map}\lim_{(a^{-1}, \h{0.5pt}r) \h{0.5pt}\to\h{0.5pt} (0, 0)} \h{2pt}  \max_{k = 1, ..., k_a } \h{1.5pt}\sum^2_{j=0} \h{1.5pt}r^j \h{1pt} \Big\lVert \h{1pt} \nabla^j \Pi_{\mathbb{S}^4}\h{1pt}\big[\h{0.5pt}w_{a, b}^+ \h{0.5pt}\big]  - \nabla^j \left[\h{1pt}\Lambda_k  \big(\cdot-z_{a, k}^+\big)\h{1pt}\right]\h{1.5pt} \Big\rVert_{\infty;\h{1pt} \p B_{r}  (z_{a, k}^+ ) } = 0.
\end{align}Given $\phi$ and $\theta$, the polar and azimuthal angles, respectively, we define
\begin{equation}\label{def of lambda plus and minus}
\Lambda_+ := \big(0,0,\cos \phi,\sin \phi \cos \theta, \sin \phi \sin \theta\big)^\top \h{5pt} \text{ and } \h{5pt} \Lambda_-:= \big(0,0,-\cos \phi,\sin \phi \cos \theta, \sin \phi \sin \theta\big)^\top.
\end{equation}Then it satisfies $\Lambda_k = \Lambda_-$ if $k$ is odd and $\Lambda_k = \Lambda_+$ if $k$ is even;
\item[$\mathrm{(2).}$] Let $z_{a, 1}^+$, ..., $z_{a, s_a}^+$ be the zeros of $w_{a, c}^-$ on $l_z^+$, where $s_a$ is the total number of zeros of $w_{a, c}^-$ on $l_z^+$. Moreover, the zeros are ordered so that for fixed $a$, $z_{a, j; 3}^+$  are increasing with respect to $j$. Then  \begin{align}\label{asy of pha map}\lim_{(a^{-1}, \h{0.5pt}r) \h{0.5pt}\to\h{0.5pt} (0, 0)} \h{2pt} \max_{k = 1, ..., s_a}\h{1.5pt}  \sum^2_{j=0} \h{1.5pt}r^j \h{1pt} \Big\lVert \h{1pt} \nabla^j \Pi_{\mathbb{S}^4}\h{1pt}\big[\h{0.5pt}w_{a, c}^- \h{0.5pt}\big]  - \nabla^j \left[\h{1pt}\Lambda_k  \big(\cdot-z_{a, k}^+\big)\h{1pt}\right]\h{1.5pt} \Big\rVert_{\infty;\h{1pt} \p B_{r}  (z_{a, k}^+ ) } = 0.
\end{align}In this case, $\Lambda_k = \Lambda_+$ if $k$ is odd and $\Lambda_k = \Lambda_-$ if $k$ is even.
\end{enumerate} 
\end{prop}

To prove this proposition, we will frequently use some lemmas from Section A.3 in the appendix.


\subsection{Strictly isolated zeros and their non--degeneracy}

In this section we study some general results on the  zeros of $w_{a}$ for large $a$. Here and throughout the following, $\big\{ w_a\big\}$ denotes either the family $\big\{ w_{a, b}^+\big\}$ or $\big\{ w_{a, c}^-\big\}$. We focus on the  mutual distances and non--degeneracy of the zeros of $w_a$ on $l_z$.    The main result is \begin{prop}\label{strict isolation} Let $\big\{w_{a}\big\}$ denote either $\big\{ w_{a, b}^+\big\}$ or $\big\{w_{a, c}^-\big\}$. Then we have \begin{enumerate}
\item[$\mathrm{(1)}.$] There exist $a_0 > 0$ suitably large and  $\delta_0 \in \big(0, 1/4\big)$ so that all zeros of $w_{a}$ are contained in $l^+_{\delta_0} \h{1.5pt}\bigcup\h{1.5pt} l_{\delta_0}^-$, provided that $a > a_0$. Here $\delta_0$ and $a_0$ may depend on $b$, $c$ and $\mu$. $l^+_{\delta_0}$ is the closed segment connecting $\big(0, 0, \delta_0\big)$ and $\big(0, 0, 1 - \delta_0\big)$. $l_{\delta_0}^-$ is the symmetric segment of $l^+_{\delta_0}$ with respect to the origin;
\item[$\mathrm{(2)}.$] For any $R > 0$, it satisfies \begin{align}\label{non degeneracy} \liminf_{a \to \infty} \h{2pt}\min_{\big\{\h{.5pt} z_a \h{1pt}: \h{1pt}w_{a}(z_a) \h{0.5pt}= \h{0.5pt}0\h{.5pt}\big\}} \min_{\big\{\h{0.5pt}x \h{0.5pt}:\h{0.5pt} | x - z_a | \h{1pt}\leq\h{1pt} R\h{.5pt} a^{- \frac{1}{2}} \big\}} \dfrac{\big|\h{1pt} w_{a}(x) \h{1pt}\big|}{\sqrt{a} \h{1.5pt} | x - z_a |} \h{2pt}=\h{2pt}c_\mu(R\h{0.5pt}) \h{1pt}:=\h{1pt} \sqrt{\mu} \h{1pt}\min_{r \h{1pt}\leq\h{1pt}R\h{0.5pt}\mu^{\frac{1}{2}}} \dfrac{f(r)}{r} \h{2pt}>\h{2pt}0.
\end{align}In the above limit, $f \in C^2\big[ \h{0.5pt}0, \infty \h{0.5pt}\big)$ is the unique solution of \begin{align*} \left\{ \begin{array}{lcl} f'' + \dfrac{2}{r} \h{1pt}f' - \dfrac{2}{r^2} \h{1pt}f + f \big( 1 - f^2 \big) = 0 \h{15pt}&\text{on $\big(0, \infty\big)$;}\\[3mm]
f(0) = 0 \h{15pt}\text{and}\h{15pt}f(+ \infty) = 1. \end{array}\right.
\end{align*}Moreover, it holds
\begin{equation}
f'(r) > 0 \h{5pt} \text{in $\big[\h{0.5pt}0,\infty\h{0.5pt}\big)$} \h{15pt} \text{ and} \h{15pt} R f'(R) + \Big| 1 - R^2 \h{1pt}\big[\h{0.5pt} 1-f(R)\h{0.5pt}\big] \h{1pt} \Big| = o(1) \h{5pt} \text{ as $R \to \infty$};
\label{property of f}
\end{equation}Here $o(1)$ denotes a quantity so that it converges to $0$ as $R \to \infty$;
\item[$\mathrm{(3)}.$] There exist $a_0 > 0$ suitably large and $\delta_1 \in (0, 1)$ so that for all $a > a_0$, either $w_{a} $ has only one zero on $l_z^+$, or the distance between two different zeros of $w_{a}$ on $l_z^+$ is greater than $\delta_1$.
\end{enumerate} 
\end{prop}\begin{rmk}The (2) in Proposition \ref{strict isolation} is referred to as the non--degeneracy result of the zeros of $w_{a}$. The properties of the radial function $f$ have been obtained in \cite{AF97, FG00, G97, MP10}. The items (1) and (3) in Proposition \ref{strict isolation} infer that any zero of $w_{a}$ keeps strictly away from poles, the origin and other zeros of $w_{a}$, provided that $a$ is suitably large.
\end{rmk}We firstly show item (1) in Proposition \ref{strict isolation}.  \begin{proof}[\bf Proof of (1) in Proposition \ref{strict isolation}] Suppose that there are a sequence $\big\{ a_n\big\}$ which diverges to $\infty$ as $n \to \infty$ and a sequence $\big\{z_n\big\} \subset l_z \cap B_1^+$ such that they satisfy $w_{a_n}\left(z_n\right) = 0$ and $ z_n \to 0$ as $n \to \infty$. Up to a subsequence, we can assume that $w_{a_n}$ converges to some $w_b^+$ or $w_c^-$ strongly in $H^1\big(B_1\big)$ as $n \to \infty$. Notice that both $w_b^+$ and $w_c^-$ are smooth near $0$. We then can apply Lemma   \ref{|w_a|>1/2, away from sing.} in the appendix to obtain $\big| w_{a_n}\big| \geq 1/4$ in a neighborhood $\mathcal{O}$ of $0$, provided that $n$ is suitably large. However, $z_n \in \mathcal{O}$ when $n$ is large. We then obtain a contradiction since $w_{a_n}\big(z_n\big) = 0$.   Similarly,  there cannot have a sequence $\big\{ a_n\big\}$ which diverges to $\infty$ as $n \to \infty$ and a sequence $\big\{z_n\big\} \subset l_z \cap B_1^+$ such that they satisfy $w_{a_n}\left(z_n\right) = 0$ and $z_n \to e_3^* = (0, 0, 1)^\top$ as $n \to \infty$. Here one just needs Lemma  \ref{|w_a|>1/2, near poles.}. The proof for the case when $ z_n \to - e_3^*$ as $n \to \infty$ can be obtained by symmetry.
\end{proof}
In the following four  sections, we prove items (2) and (3) in Proposition \ref{strict isolation}.
\subsubsection{Accumulation of zeros} In this section, we  consider some accumulation properties of the zeros of $w_{a}$ up to a subsequence. \begin{lem}\label{accu of zero} Assume that $a_n \to \infty$ and $w_{a_n}$ converges to some $w_\infty$ strongly in $H^1(B_1; \mathbb{R}^5)$ as $n \to \infty$. Here $w_\infty$ equals  either $w_b^+$ or $w_c^-$. In addition, we suppose that \begin{align*} a_n \int_{B_1} \left[\h{1pt}\big| w_{a_n} \big|^2 - 1 \h{1pt}\right]^2 \longrightarrow 0\h{15pt}\text{as $n \to \infty$.}
\end{align*}Let $\mathcal{A}_\infty$ be the accumulation set of the zeros of all $w_{a_n}$. Then  \begin{align}\label{eq of accu and sing} \mathcal{A}_\infty =  \Big\{ \text{$\mathrm{Singularities \h{5pt}of}$  $w_\infty$} \Big\} \h{2pt}\subset\h{2pt} l_{\delta_0}^+ \h{1.5pt} \bigcup \h{1.5pt} l_{\delta_0}^-.
\end{align}Here $l_{\delta_0}^+$ and $l_{\delta_0}^-$ are defined in item (1) of Proposition \ref{strict isolation}. Moreover, it holds $\mathrm{Card}\big(\mathcal{A}_\infty\big) < \infty$.
\end{lem}
\begin{proof}[\bf Proof] By (1) in Proposition \ref{strict isolation}, it satisfies $\mathcal{A}_\infty \h{1pt}\subseteq\h{1pt}l_{\delta_0}^+ \cup l_{\delta_0}^-$. Now we prove the equality in (\ref{eq of accu and sing}). Given $z^* \in \mathcal{A}_\infty$, there exist a subsequence, still denoted by $\big\{a_n\big\}$,  and $z_n \in l_z$ so that $z_n \to z^*$ as $n \to \infty$. Moreover, $w_{a_n}\big(z_n\big) = 0$ for all $n$.  If $z^*$ is a smooth point of $w_\infty$, then by Lemma \ref{|w_a|>1/2, away from sing.}, there is an open neighborhood of $z^*$, denoted by $\mathcal{O}_{z^*}$, so that $\big| w_{a_n} \big| \geq 2^{-1}$ on $\mathcal{O}_{z^*}$ for large $n$. However, this is impossible since for large $n$, we have $z_n \in \mathcal{O}_{z^*}$ and $w_{a_n}\big(z_n\big) = 0$. Hence $z^*$ is a singularity of $w_\infty$. In the next, we assume that  $z^*$ is a singularity of  $w_\infty$.  Note that $w_\infty$ is smooth except at finitely many singularities on $l_z$. Meanwhile, the singularities of $w_\infty$ are different from the two poles and the origin. Fixing an $\epsilon > 0$ arbitrarily small, we know that $z^*_{\epsilon, +} := z^* + (0, 0, \epsilon)$ is a smooth point of $w_\infty$. Still by Lemma \ref{|w_a|>1/2, away from sing.}, we can take $n$ large enough so that $\big| w_{a_n} \big| \geq 2^{-1}$ on $B_{2^{-1}\epsilon}\big(z^*_{\epsilon, +}\big)$. In light of Lemma   \ref{energy density e^L_a[w_a] bdd}, $\nabla w_{a_n}$ is uniformly bounded on $B_{2^{-2}\epsilon'}\big(z^*_{\epsilon, +}\big)$ for large $n$ and small $\epsilon' < \epsilon$. Applying Arzel\`{a}--Ascoli theorem, we know that $w_{a_n}\big(z^*_{\epsilon, +}\big)$ converges to $w_\infty\big(z^*_{\epsilon, +}\big)$ as $n \to \infty$. Similarly if we define  $z^*_{\epsilon, -} := z^* - (0, 0, \epsilon)$, then $w_{a_n}\big(z^*_{\epsilon, -}\big)$ converges to $w_\infty\big(z^*_{\epsilon, -}\big)$ as well when $n \to \infty$. Recalling Theorems 1.6 and 1.7 in \cite{Y20}, particularly item (2) in these two theorems, we know that the third component of $w_\infty$ equals $\pm 1$ at $z^*_{\epsilon, +}$ and $z^*_{\epsilon, -}$. Meanwhile   $w_{\infty; 3} $ has different signs at $z^*_{\epsilon, +}$ and $z^*_{\epsilon, -}$. Hence, for large $n$, the third component of $w_{a_n}$ also has different signs at $z^*_{\epsilon, +}$ and $z^*_{\epsilon, -}$. By continuity of $w_{a_n}$ on $l_z$, for large $n$, there is a point on the segment connecting $z^*_{\epsilon, +}$ and $z^*_{\epsilon, -}$ so that the third component of $w_{a_n}$ equals  $0$ at this point. Hence, $w_{a_n}$ vanishes at this point. Since $\epsilon > 0$ is arbitrarily small, it follows $z^* \in \mathcal{A}_\infty$.
\end{proof}
\begin{lem}\label{Haus conv}Let $\big\{w_{a_n}\big\}$ be as in Lemma \ref{accu of zero} and $y_0 \in \mathcal{A}_\infty$. In addition, we assume that $\sigma_0$ is a positive constant suitably small so that $y_0$ is the only singularity of $w_\infty$ in $B_{2\sigma_0}\big(y_0\big)$. If   \begin{align*}\mathcal{V}_n := \overline{B_{\sigma_0}(y_0)} \h{2pt}\bigcap \h{2pt}\Big\{\h{1pt} \big| w_{a_n} \big| \leq 2^{-2}\h{1pt} \Big\} \neq \emptyset  \h{15pt}\text{ for all $n$,}\end{align*} then $\mathcal{V}_n$ converges to $y_0$ as $n \to \infty$ in the sense of Hausdorff. 
\end{lem}\begin{proof}[\bf Proof]  Suppose on the contrary that there are $\epsilon_0 > 0$ and a subsequence, still denoted by $\{ a_n \}$, such that \begin{align*} \mathrm{d}_H\big[\h{1pt}\mathcal{V}_n, y_0\h{1pt}\big]:=\max_{x \h{1pt}\in\h{1pt} \mathcal{V}_n} \big| \h{.5pt} x - y_0 \h{.5pt} \big| = \big| \h{.5pt} y_n - y_0 \h{.5pt}\big| \geq \epsilon_0 \h{15pt}\text{for any $n \in \mathbb{N}$.} \end{align*}Here $y_n \in \mathcal{V}_n$. Up to a subsequence, $\big\{y_n\big\}$ converges to a point, denoted by $y_0^* \in \overline{B_{\sigma_0}(y_0)}$. The point $y_0^*$ is not a singularity of $w_\infty$. Applying Lemma \ref{|w_a|>1/2, away from sing.} yields $\big| w_{a_n} \big| \geq 2^{-1}$ in an open neighborhood $\mathcal{O}_{y_0^*}$ of $y_0^*$, provided that $n$ is large. However, this is impossible since for $n$ large, $y_n \in \mathcal{O}_{y_0^*}$ and  $\big| w_{a_n} \big| \leq 2^{-2}$ at $y_n$.
\end{proof}
\subsubsection{Blow--up in the exterior core}

In this section, we consider the blow--up sequence of $w_{a_n}$ near one of its zeros.  We focus on the blow--up occurring in the exterior core.
\begin{lem}Let $y_0$, $\sigma_0$, $w_\infty$ and $\big\{w_{a_n}\big\}$  be  as in Lemma \ref{Haus conv}. Moreover, for any $n \in \mathbb{N}$, we assume that there is a point  $z_n \in B_{\sigma_0}(y_0)$ so that $w_{a_n}(z_n) = 0$. Then the set $\mathcal{V}_n$ defined in  Lemma \ref{Haus conv} is not empty. In addition, it holds \begin{align}\label{diam conv 0} \nu_n := \max_{z \h{1pt}\in\h{1pt}\mathcal{V}_n} \big|\h{1pt}z - z_n \h{1pt}\big| \longrightarrow 0 \h{15pt}\text{as $n \to \infty$.}
\end{align}For any sequence $\big\{r_n\big\}$ with $r_n \to 0$ as $n \to \infty$, if we have $r_n \geq \nu_n$ for any $n$ and $a_n r_n^2 \to \infty$ as $n \to \infty$, then there is a mapping $w^\infty \in H^1_{\mathrm{loc}}\big(\mathbb{R}^3; \mathbb{S}^4\big)$ so that up to a subsequence, the following convergences hold for  $w^{(n)}(\zeta) := w_{a_n}\left(z_n + r_n \zeta\right)$:\begin{align*}
w^{(n)} \longrightarrow w^\infty \h{5pt} \text{ strongly in } \h{2pt} H^1_\textup{loc}\big(\mathbb{R}^3; \mathbb{R}^5\big) \h{10pt} \text{ and }  \h{10pt}
a_n r_n^2 \Big(\big|\h{.5pt} w^{(n)} \h{.5pt}\big|^2 - 1\Big)^2 \longrightarrow 0 \h{5pt} \text{ strongly in } \h{2pt} L^1_\textup{loc}\big(\mathbb{R}^3\big).
\end{align*}The limiting map $w^\infty$ is a local minimizer in the sense that for any $B_r \subset \mathbb{R}^3$,  $w^\infty$ minimizes the Dirichlet energy in  $H \big(r,w^\infty\big)$. Here for any $r > 0$ and a $\mathbb{S}^4$--valued mapping $w_*$ on $\p B_r$,   \begin{align*}H \big(r, w_*\big) := \Big\{ w \in H^1\big(B_r ; \mathbb{S}^4\big)
\h{1pt} : \h{1pt} w=w_* \h{2pt} \text{ on $\p B_r$} \h{2pt} ,\h{2pt}  w = \mathscr{L}[\h{0.5pt}u\h{0.5pt}] \h{2pt} \text{ for some $3$--vector field } \h{2pt} u=u(\rho, z)\Big\}.\end{align*}

\label{stong conv. and min. property of w^infty}
\end{lem}
\begin{proof}[\bf Proof] Fixing an arbitrary $R > 0$  and recalling the $F_n$ defined in (\ref{def of F_n}), we have \begin{align*}\dfrac{1}{R} \int_{B_{R}} \big| \nabla_\zeta w^{(n)} \big|^2 + r_n^2 \h{1pt}F_n\big(w^{(n)}\big) = \dfrac{1}{r_n R} \int_{B_{r_n R}(z_n)} \left| \nabla w_{a_n} \right|^2 + F_n \big( w_{a_n}\big).
\end{align*}By (1) in Proposition \ref{strict isolation} and the monotonicity formula in Lemma \ref{Energy Monotonicity}, for a fixed $\delta \in \big(\h{0.5pt}0, 2^{-1}\delta_0\h{0.5pt}\big)$ and large $n$, the above equality yields \begin{align*}\dfrac{1}{R} \int_{B_{R}} \big| \nabla_\zeta w^{(n)} \big|^2 + r_n^2 \h{1pt}F_n\big(w^{(n)}\big) \h{2pt}\leq\h{2pt} \delta^{-1} \int_{B_{\delta}(z_n)} \left| \nabla w_{a_n} \right|^2 + F_n \big( w_{a_n}\big).
\end{align*} Here $\delta_0$ is given in item (1) of Proposition \ref{strict isolation}. Recall that $z_n \to y_0$ as $n \to \infty$. For any $\epsilon > 0$ small, we can take $n$ large enough and get from the last estimate that  \begin{align*} \dfrac{1}{R} \int_{B_{R}} \big| \nabla_\zeta w^{(n)} \big|^2 + r_n^2 \h{1pt}F_n\big(w^{(n)}\big) \h{2pt}\leq\h{2pt} \delta^{-1} \int_{B_{\delta + \epsilon}(y_0)} \left| \nabla w_{a_n} \right|^2 + F_n \big( w_{a_n}\big). 
\end{align*}Here we have used $B_\delta(z_n) \subset B_{\delta + \epsilon}(y_0)$ for large $n$. Now we take $n \to \infty$ in the above estimate. It turns out  \begin{align*} \sup_{R \h{1pt}>\h{1pt}0}\h{2pt}\limsup_{n \to \infty}\h{2pt} \dfrac{1}{R} \int_{B_{R}} \big| \nabla_\zeta w^{(n)} \big|^2 + r_n^2 \h{1pt}F_n\big(w^{(n)}\big)   \h{2pt}\leq\h{2pt}\delta^{-1} \int_{B_{\delta + \epsilon}(y_0)} \left| \nabla w_\infty \right|^2 + \sqrt{2}\h{0.5pt}\mu\h{0.5pt}\big( 1 - 3 \h{0.5pt}S\h{0.5pt}[ w_\infty ]\h{1pt}\big).
\end{align*}Utilizing the results in \cite{Y20} (see (4.4) and Proposition 4.4 there), we take $\epsilon \to 0$ and $\delta \to 0$ successively in the above estimate. It then follows \begin{align}\label{unifo of scaled energy} \sup_{R \h{1pt}>\h{1pt}0}\h{2pt}\limsup_{n \to \infty}\h{2pt} \dfrac{1}{R} \int_{B_{R}} \big| \nabla_\zeta w^{(n)} \big|^2 + r_n^2 \h{1pt}F_n\big(w^{(n)}\big)   \h{2pt}\leq\h{2pt}\int_{B_1} \big| \nabla \Lambda \big|^2 = 8\pi.
\end{align}Here $\Lambda$ equals  either $\Lambda_+$ or $\Lambda_-$ given in (\ref{def of lambda plus and minus}).\vspace{0.2pc}

With the assumption that $a_n r_n^2 \to \infty$ as $n \to \infty$, there then exists a $w^\infty$ with unit length so that up to a subsequence, $w^{(n)}$ converges to $w^\infty$ weakly in $H^1_{\mathrm{loc}}\big(\mathbb{R}^3; \mathbb{R}^5\big)$ and strongly in $L^2_{\mathrm{loc}}\big(\mathbb{R}^3; \mathbb{R}^5\big)$ as $n \to \infty$. Still by Fatou's lemma, for any $R > 1$, we can find a $\sigma \in \big(R, 2R\big)$ so that up to a subsequence, it holds \begin{equation}
\sup_{n \h{1pt}\in\h{1pt}\mathbb{N}} \int_{\p B_{\sigma}}  \big|\h{.5pt} \nabla_\zeta w^{(n)} \h{.5pt}\big|^2
+  a_n r_n^2 \Big(\big|\h{.5pt} w^{(n)} \h{.5pt}\big|^2 - 1\Big)^2
< \infty \h{5pt} \text{ and } \h{5pt}
 \int_{\p B_{\sigma}}  \big|\h{.5pt} w^{(n)} - w^\infty \h{.5pt}\big|^2 \to 0 \h{5pt} \text{ as } n \to \infty.
\label{fatou lemma for w^n}
\end{equation}
Since $\sigma r_n > \nu_n$, by the definition of $\nu_n$ in (\ref{diam conv 0}), it holds $\big| w_{a_n} \big| > 2^{-2}$ on $\p B_{\sigma r_n}(z_n)$. The normalized vector field $\widehat{w_{a_n}}$ is well--defined on $\p B_{\sigma r_n}(z_n)$. Equivalently, the vector field $\widehat{w^{(n)}}$ is also well--defined on $\p B_\sigma$. Let $W^{(n)}$ minimize the Dirichlet energy in  $H\Big(\sigma, \widehat{w^{(n)}}\Big)$. With this  $W^{(n)}$, we  define a comparison map as follows:\begin{eqnarray*}
v_{n, s} :=\left\{
\begin{aligned}
& W^{(n)}\left(\dfrac{\zeta}{1-s}\right) \quad &\text{ in }& B_{(1-s)\sigma};\\[2mm]
&    \dfrac{\sigma - |\h{1pt} \zeta \h{1pt}|}{s\h{1pt}\sigma} \h{1.5pt}\widehat{w^{(n)}} \big(\sigma \widehat{\zeta}\h{1pt}\big) +  \dfrac{|\h{1pt} \zeta \h{1pt}| - (1 - s)\sigma }{s\h{1pt}\sigma} \h{1.5pt}w^{(n)}  \big(\sigma \widehat{\zeta}\h{1pt}\big)  \quad &\text{ in }& B_{\sigma} \setminus B_{(1-s)\sigma}.
\end{aligned}
\right.
\end{eqnarray*}Here $s \in (0, 1)$ is arbitrary. Due to the energy  minimality of $w^{(n)}$, \begin{align}\label{bs es}  \int_{B_{\sigma}} \big| \nabla_\zeta w^{(n)} \big|^2 + r_n^2 \h{1pt}F_n\big(w^{(n)}\big) &\h{2pt}\leq\h{2pt}\int_{B_{\sigma}} \big| \nabla_\zeta v_{n, s} \big|^2 + r_n^2 \h{1pt}F_n\big(v_{n, s}\big)\\[1.5mm] &
\h{2pt}=\h{2pt}(1 - s) \int_{B_\sigma} \big| \nabla_\zeta W^{(n)}\big|^2 + \int_{B_{\sigma} \setminus B_{(1 - s)\sigma}} \big| \nabla_\zeta v_{n, s} \big|^2 \nonumber\\[1.5mm]
& \h{4pt}+\h{4pt} \mu \h{1pt}r_n^2 \int_{B_{\sigma}} D_{a_n} - 3 \sqrt{2} \h{1pt}S\h{1pt}\big[ v_{n, s}\big] + \dfrac{\mu}{2} \h{1pt}a_n\h{0.5pt} r_n^2 \h{1pt}\int_{B_{\sigma} \setminus B_{(1 - s)\sigma}}\left( \big| v_{n, s} \big|^2 - 1 \right)^2. \nonumber
\end{align}In light of (\ref{fatou lemma for w^n}) and $\big| w^{(n)} \big| > 2^{-2}$ on $\p B_{\sigma}$, we can apply Lemma \ref{strong H1 convergence in H(r, w) space} to obtain \begin{align}\label{conv of WN}\int_{B_\sigma} \big|\h{1pt} \nabla_\zeta W^{(n)}\h{1pt}\big|^2 \longrightarrow \int_{B_\sigma} \big|\h{1pt} \nabla_\zeta W^\infty\h{1pt}\big|^2 \h{15pt}\text{as $n \to \infty$.}
\end{align}Here $W^\infty$ minimizes the Dirichlet energy in $H\big(\sigma, w^\infty\big)$. By the uniform boundedness of $v_{n, s}$ on $B_\sigma$, we get \begin{align}\label{conv of S term}\mu \h{1pt}r_n^2 \int_{B_{ \sigma}} D_{a_n} - 3 \sqrt{2} \h{1pt}S\h{1pt}\big[ v_{n, s}\big] \longrightarrow 0 \h{15pt}\text{as $n \to \infty$.}
\end{align}Moreover, we have \begin{align*}a_n\h{0.5pt} r_n^2 \h{1pt}\int_{B_{\sigma} \setminus B_{(1 - s)\sigma}}\left( \big| v_{n, s} \big|^2 - 1 \right)^2 \h{2pt}\lesssim\h{2pt}s\h{0.5pt}\sigma\h{0.5pt} a_n  r_n^2 \h{1pt}\int_{\p B_\sigma} \left| \h{1pt}\big| w^{(n)} \big|^2 - 1 \h{1pt}\right|^2.
\end{align*}Hence, the first bound in (\ref{fatou lemma for w^n}) yields \begin{align}\label{potential on shell} \sup_{n\h{1pt}\in \h{1pt}\mathbb{N}} a_n\h{0.5pt} r_n^2 \h{1pt}\int_{B_{\sigma} \setminus B_{(1 - s)\sigma}}\left( \big| v_{n, s} \big|^2 - 1 \right)^2 \longrightarrow 0 \h{15pt}\text{as $s \to 0$.}
\end{align}  We are left to consider the Dirichlet energy of $v_{n, s}$ on $B_\sigma \setminus B_{(1 - s)\sigma}$. In the following arguments, we still use $\big(r, \varphi, \theta\big)$ to denote the spherical coordinates for the  $\zeta$--variable. Note that  the limit in (\ref{fatou lemma for w^n}) infers  \begin{align}\label{conv of radia}\int_{B_{\sigma} \setminus B_{(1 - s)\sigma}} \big| \p_r v_{n, s} \big|^2 \h{2pt}\leq\h{2pt}(s \h{0.5pt}\sigma)^{-1} \int_{\p B_{\sigma}} \left| 1 - \big| w^{(n)} \big| \h{1pt}\right|^2 \longrightarrow 0 \h{15pt}\text{as $n \to \infty$.}
\end{align}In addition, the definition of $v_{n, s}$ induces \begin{align*} \big|\h{1pt} \p_\varphi v_{n, s} \h{1pt}\big|^2     + \big|\h{1pt} \p_\theta v_{n, s} \h{1pt}\big|^2  \left(\sin \varphi\right)^{-2} \h{2pt}\lesssim\h{2pt}\sigma^2 \big| \nabla_\zeta w^{(n)} \big|^2\big(\sigma \widehat{\zeta}\h{1pt}\big) \h{15pt}\text{ on $B_\sigma \setminus B_{(1 - s)\sigma}$.}
\end{align*}Here we also have used $\big| w^{(n)} \big| > 2^{-2}$ on $\p B_\sigma$. Utilizing the last estimate, one can show that \begin{align*}\int_{B_\sigma \setminus B_{(1 - s) \sigma}} \dfrac{1}{r^2}\big|\h{1pt} \p_\varphi v_{n, s} \h{1pt}\big|^2     + \dfrac{1}{r^2 \sin^2 \varphi}\big|\h{1pt} \p_\theta v_{n, s} \h{1pt}\big|^2  \h{2pt}\lesssim\h{2pt} s\h{0.5pt}\sigma \h{1pt}\int_{\p B_\sigma} \big| \nabla_\zeta w^{(n)} \big|^2,
\end{align*}which furthermore implies by (\ref{fatou lemma for w^n}) that \begin{align*}\sup_{n \h{1pt}\in \h{1pt}\mathbb{N}}\int_{B_\sigma \setminus B_{(1 - s) \sigma}} \dfrac{1}{r^2}\big|\h{1pt} \p_\varphi v_{n, s} \h{1pt}\big|^2     + \dfrac{1}{r^2 \sin^2 \varphi}\big|\h{1pt} \p_\theta v_{n, s} \h{1pt}\big|^2 \longrightarrow 0 \h{15pt}\text{as $s \to 0$.} 
\end{align*}Applying this estimate together with (\ref{conv of WN})--(\ref{conv of radia}) to the right--hand side of (\ref{bs es}) induces \begin{align*}\limsup_{n \to \infty}\int_{B_{\sigma}} \big| \nabla_\zeta w^{(n)} \big|^2 + r_n^2 \h{1pt}F_n\big(w^{(n)}\big)  \h{2pt}\leq\h{2pt} \int_{B_\sigma} \big| \nabla_\zeta W^\infty \big|^2 \h{2pt}\leq\h{2pt} \int_{B_\sigma} \big| \nabla_\zeta w^\infty \big|^2.
\end{align*}On the other hand, lower--semi continuity infers \begin{align*}\liminf_{n \to \infty}\int_{B_{\sigma}} \big| \nabla_\zeta w^{(n)} \big|^2 + r_n^2 \h{1pt}F_n\big(w^{(n)}\big)  \h{2pt}\geq\h{2pt} \int_{B_\sigma} \big| \nabla_\zeta w^\infty \big|^2.
\end{align*}The proof is then completed by the last two estimates.
\end{proof}
Utilizing Lemma \ref{Energy Monotonicity} and the same arguments used in the proof of Proposition 4 in \cite{MZ10}, we obtain \begin{lem}\label{clearing lemma for scaling}Let $w^{(n)}$ and $w^\infty$ be as in Lemma \ref{stong conv. and min. property of w^infty}. Given $x_0 \in \mathbb{R}^3$, if for any $\epsilon > 0$, it satisfies \begin{align*} r_\epsilon^{-1} \int_{B_{r_\epsilon}(x_0)}\big| \nabla_\zeta w^\infty \big|^2 < \epsilon, \h{15pt}\text{for some $r_\epsilon > 0$ suitably small,}
\end{align*}then there is an open neighborhood of $x_0$, denoted by $\mathcal{O}_{x_0}$, so that  $\big| w^{(n)}  \h{1.5pt}\big| > 2^{-1}$ on $\mathcal{O}_{x_0}$ for large $n$.
\end{lem}With this lemma, we can characterize the limit $w^\infty$ obtained in Lemma \ref{stong conv. and min. property of w^infty} as follows.\begin{lem}
The limiting map $w^\infty$ obtained in Lemma \ref{stong conv. and min. property of w^infty}  equals either $\Lambda_+$ or $\Lambda_-$ in (\ref{def of lambda plus and minus}).
\label{w^infty is homogeneous zero}
\end{lem}
\begin{proof}[\bf Proof] Recall the monotonicity formula in Lemma \ref{Energy Monotonicity}. Letting $a = a_n$, $w_a = w_{a_n}$ and $y = z_n$ in this monotonicity formula, we then integrate the variable $R$ from $r\h{0.5pt}r_n$ to $\rho \h{0.5pt} r_n$. Here we take $\rho > r > 0$. Applying change of variables to the resulting equality, we obtain \begin{align*} \dfrac{1}{\rho} \int_{B_\rho} \big| \nabla_\zeta w^{(n)} \big|^2 + r_n^2 \h{1pt}F_n\big(w^{(n)}\big) &- \dfrac{1}{r} \int_{B_r} \big| \nabla_\zeta w^{(n)} \big|^2 + r_n^2 \h{1pt}F_n\big(w^{(n)}\big) \\[2mm]
&= \int_r^\rho \dfrac{2}{\sigma} \int_{\p B_\sigma} \left| \dfrac{\p \h{0.5pt}w^{(n)}}{\p \h{0.5pt}\vec{n}} \right|^2 \h{2pt}\mathrm{d} \sigma + \int_r^\rho \dfrac{2}{\sigma^2} \int_{B_\sigma} r_n^2\h{1pt}F_n\big(w^{(n)}\big) \h{2pt}\mathrm{d}\sigma.
\end{align*}  By the strong convergence in Lemma \ref{stong conv. and min. property of w^infty}, we can take $n \to \infty$ in the above equality and obtain  
\begin{equation}
\dfrac{1}{\rho} \int_{B_\rho} \big|\h{.5pt} \nabla_\zeta w^\infty \h{.5pt}\big|^2
- \dfrac{1}{r} \int_{B_r} \big|\h{.5pt} \nabla_\zeta w^\infty \h{.5pt}\big|^2
= \int_r^\rho \dfrac{2}{\sigma} \int_{\p B_\sigma} \left| \dfrac{\p \h{0.5pt} w^\infty}{\p \h{0.5pt} \vec{n}} \right|^2 \mathrm{d} \sigma.
\label{harmonic map energy mono.}
\end{equation}

Let $l_n$ be a positive sequence converging to $0$ as $n \to \infty$ and define $w^\infty_{l_n}(\zeta) := w^\infty\big(l_n \zeta\big)$. By (\ref{unifo of scaled energy}),
\begin{align*} \sup_{R \h{1pt}>\h{1pt}0, \h{2pt}n \h{1pt}\in\h{1pt}\mathbb{N}} \h{3pt}\dfrac{1}{R} \int_{B_R} \big|\h{.5pt} \nabla_\zeta w_{l_n}^\infty \h{.5pt}\big|^2 \h{2pt}\leq\h{2pt}8\pi.\end{align*}
Hence up to a subsequence, $w^\infty_{l_n}$ converges  to a limiting map $W_0$ as $n \to \infty$, weakly in $H^1_\textup{loc}(\mathbb{R}^3)$ and strongly in $L_{\mathrm{loc}}^2\big(\mathbb{R}^3\big)$. Still by Fatou's lemma, for any $R > 0$, we can find a $\sigma \in \big(R, 2R\big)$ so that up to a subsequence, it holds \begin{equation*}
\sup_{n \h{1pt}\in\h{1pt}\mathbb{N}} \int_{\p B_{\sigma}}  \big|\h{.5pt} \nabla_\zeta w_{l_n}^{\infty} \h{.5pt}\big|^2 
< \infty \h{5pt} \text{ and } \h{5pt}
 \int_{\p B_{\sigma}}  \big|\h{.5pt} w_{l_n}^{\infty} - W_0 \h{.5pt}\big|^2 \to 0 \h{5pt} \text{ as } n \to \infty.
\end{equation*} In light of Lemma \ref{strong H1 convergence in H(r, w) space},  $w_{l_n}^\infty$ converges to $W_0$ strongly in $H^1\big(B_\sigma\big)$ as $n \to \infty$. Moreover, $W_0$ minimizes the Dirichlet energy in $H\big(\sigma, W_0\big)$. In light that $R > 0$ is arbitrary, up to a subsequence, we can assume $w_{l_n}^\infty$ converges to $W_0$ strongly in $H_{\mathrm{loc}}^1\big(\mathbb{R}^3\big)$ as $n \to \infty$. Moreover, $W_0$ is a local minimizer in the sense that it minimizes Dirichlet energy in $H\big(R, W_0\big)$ for any $R > 0$. \vspace{0.2pc}

Notice the monotonicity formula of $w^\infty$ in  (\ref{harmonic map energy mono.}). The limit $\lim\limits_{r \to 0}r^{-1} \displaystyle\int_{B_r} \big|\h{.5pt} \nabla_\zeta w^\infty \h{.5pt}\big|^2 = \ell_*$ is well--defined. Here we also have used the uniform bound in (\ref{unifo of scaled energy}). Now, we take $r \to 0$ in (\ref{harmonic map energy mono.}) and get
$$
\dfrac{1}{\rho} \int_{B_\rho} \big|\h{.5pt} \nabla_\zeta w^\infty \h{.5pt}\big|^2 - \ell_* = \int_0^\rho \dfrac{2}{\sigma} \int_{\p B_\sigma} \left| \dfrac{\p \h{0.5pt} w^\infty}{\p \h{0.5pt}  \vec{n}} \right|^2 \mathrm{d} \sigma.
$$
By replacing $\rho$ in this equality with $\rho \h{1pt} l_n$ and changing variables, it turns out
$$
\dfrac{1}{\rho \h{1pt} l_n} \int_{B_{\rho\h{1pt}l_n}}  \big|\h{.5pt} \nabla_\zeta w^\infty \h{.5pt}\big|^2 - \ell_* = \int_0^\rho \dfrac{2}{\sigma} \int_{\p B_\sigma} \left| \dfrac{\p\h{0.5pt} w^\infty_{l_n}}{\p \h{0.5pt}  \vec{n}} \right|^2 \mathrm{d} \sigma.
$$Utilizing the strong convergence of $w_{l_n}^\infty$ in $H^1_{\mathrm{loc}}\big(\mathbb{R}^3\big)$, the limit of $r^{-1} \displaystyle\int_{B_r} \big|\h{.5pt} \nabla_\zeta w^\infty \h{.5pt}\big|^2$ as $r \to 0$ and Fatou's lemma, we can take $n \to \infty$ in the above equality and obtain $$
0 = \liminf_{n \to \infty} \dfrac{1}{\rho \h{1pt} l_n} \int_{B_{\rho\h{1pt}l_n}}  \big|\h{.5pt} \nabla_\zeta w^\infty \h{.5pt}\big|^2 - \ell_* = \liminf_{n \to \infty} \int_0^\rho \dfrac{2}{\sigma} \int_{\p B_\sigma} \left| \dfrac{\p\h{0.5pt} w^\infty_{l_n}}{\p \h{0.5pt}  \vec{n}} \right|^2 \mathrm{d} \sigma \h{2pt}\geq\h{2pt} \int_0^\rho \dfrac{2}{\sigma} \int_{\p B_\sigma} \left| \dfrac{\p\h{0.5pt} W_0}{\p \h{0.5pt}  \vec{n}} \right|^2 \mathrm{d} \sigma.
$$  Since $\rho > 0$ is arbitrary, $W_0$ is therefore homogeneous zero. Note that $W_0$ minimizes the Dirichlet energy in $H\big(B_1, W_0\big)$. The results in \cite{Y20} (see Lemma 4.3 and Proposition 4.4 there) induce that \begin{align*} W_0 = e_3, \h{3pt}- e_3,\h{3pt}\Lambda_+  \h{3pt}\text{or}\h{3pt}\Lambda_-.
\end{align*} Suppose that $W_0 = e_3$ or $- e_3$. By the strong convergence of $w_{l_n}^\infty$ to $W_0$ in $H^1_{\mathrm{loc}}\big(\mathbb{R}^3\big)$, for any $\epsilon > 0$, there is a $r_\epsilon > 0$ so that \begin{align*} \dfrac{1}{r_\epsilon}\int_{B_{r_\epsilon}} \big| \nabla w^\infty\big|^2 < \epsilon.
\end{align*}By Lemma \ref{clearing lemma for scaling}, it turns out $\big| w^{(n)}(0) \big| = \big| w_{a_n}\big(z_n\big) \big| > 2^{-1}$ for large $n$. However, this is impossible since $z_n$ is a zero of $w_{a_n}$. Therefore, $W_0$ equals  $\Lambda_+$ or $\Lambda_-$, which furthermore infers \begin{align}\label{determine l} \ell_* = \lim_{n \to \infty} l_n^{-1} \int_{B_{l_n}} \big| \nabla w^\infty \big|^2 = \lim_{n \to \infty}   \int_{B_1} \big| \nabla w_{l_n}^\infty \big|^2 = \int_{B_1} \big| \nabla W_0 \big|^2 = 8\pi.
\end{align}

Fix a $\rho > 0$ and let $r = l_n$ in (\ref{harmonic map energy mono.}). Then we take $n \to \infty$. By (\ref{unifo of scaled energy}) and (\ref{determine l}), it follows \begin{align*} \int_0^\rho \dfrac{2}{\sigma} \int_{\p B_\sigma} \left| \dfrac{\p \h{0.5pt} w^\infty}{\p \h{0.5pt} \vec{n}} \right|^2 \mathrm{d} \sigma \h{2pt}\leq\h{2pt} 0 \h{15pt}\text{for all $\rho > 0$.}
\end{align*}Hence, $w^\infty$ is also $0$--homogeneous. In light of the minimality of $w^\infty$ in Lemma \ref{stong conv. and min. property of w^infty}, we  conclude that $w^\infty$ equals  $\Lambda_+$ or $\Lambda_-$. Here we have used (\ref{determine l}), together with  Lemma 4.3 and Proposition 4.4 in \cite{Y20}.
\end{proof}

In the next, we give an upper bound of $\nu_n$ by utilizing Lemmas \ref{stong conv. and min. property of w^infty}--\ref{w^infty is homogeneous zero}. 
\begin{lem}\label{size of core} Let $\big\{a_n\big\}$ be as in Lemma \ref{stong conv. and min. property of w^infty}. Then $\nu_n$ defined in (\ref{diam conv 0}) satisfies \begin{align*} \limsup_{n \to \infty} \h{2pt}  a_n   \nu^2_n \h{2pt} < \h{2pt}\infty.
\end{align*}
\end{lem}\begin{proof}[\bf Proof] Suppose that up to a subsequence, $a_n \nu_n^2 \rightarrow \infty$ as $n \to \infty$. Then by Lemmas \ref{stong conv. and min. property of w^infty} and \ref{w^infty is homogeneous zero}, the mapping $w^{(n)}_*(\zeta) := w_{a_n}\big(z_n + \nu_n \zeta\big)$ converges strongly in $H^1_{\mathrm{loc}}\big(\mathbb{R}^3\big)$ to $\Lambda$ as $n \to \infty$. Here $\Lambda$ equals either $\Lambda_+$ or $\Lambda_-$. In light of the definition of $\nu_n$ in (\ref{diam conv 0}), there is a $p_n \in \mathcal{V}_n$ so that \begin{align}\label{contr to vio} \nu_n = \big| p_n - z_n \big| \h{10pt}\text{and}\h{10pt} \Big| w_{a_n} \h{0.5pt}(p_n) \h{1pt}\Big|  = 2^{-2}.
\end{align} The first equality in (\ref{contr to vio}) yields $\frac{p_n - z_n}{\nu_n} \in \p B_1$.  Up to a subsequence, we can assume $\frac{p_n - z_n}{\nu_n}$ converges to some $\zeta_0 \in \p B_1$ as $n \to \infty$. Since $\Lambda$ is smooth at $\zeta_0$, by Lemma \ref{clearing lemma for scaling}, there is an open neighborhood, denoted by $\mathcal{O}_{\zeta_0}$, so that $ \big| w^{(n)}_* \big| > 2^{-1}$ on $\mathcal{O}_{\zeta_0}$ for large $n$. We now take $n$ large. The above arguments   yield    $$\frac{p_n - z_n}{\nu_n} \in \mathcal{O}_{\zeta_0} \h{15pt}\text{and}\h{15pt} \Big| w_{a_n} (p_n)\h{1pt}\Big| = \left| w^{(n)}_*\left(\frac{p_n - z_n}{\nu_n}\right) \right| > 2^{-1}.$$ But this result violates the second equality in (\ref{contr to vio}). The proof is completed.
\end{proof}
\begin{rmk} In light of Lemma \ref{size of core}, the size of the core, i.e. $\nu_n$, is at most of order $\mathrm{O}\big(a_n^{- \frac{1}{2}}\big)$ as $n \to \infty$. Hence, for $r_n$ given in Lemma \ref{stong conv. and min. property of w^infty}, it satisfies $r_n>\!\!>\nu_n$. That is the reason why the blow--up considered in Lemmas \ref{stong conv. and min. property of w^infty} and \ref{w^infty is homogeneous zero} is referred to as the blow--up occurring in the exterior core.
\end{rmk}
\subsubsection{Blow--up in the interior core} Same as before, we assume that $a_n \to \infty$ and $r_n \to 0$ as $n \to \infty$. In this section, we additionally assume  \begin{align}\label{interior core size}a_n r_n^2 \to L  \h{10pt}\text{as $n \to \infty$, for some $L \in [\h{1pt}0, \infty\h{1pt})$.}\end{align} Moreover, we still use $w^{(n)}$ to denote the scaled mapping $w_{a_n}\big(z_n + r_n \zeta\big)$. \begin{lem}\label{result on wstar infty}Up to a subsequence, $w^{(n)}$ with $a_n$ and $r_n$ satisfying (\ref{interior core size}) converges   to some $w^\infty_\star$ in $C_{\mathrm{loc}}^2\big(\mathbb{R}^3\big)$. If $L = 0$, then $w^\infty_\star \equiv 0$ on $\mathbb{R}^3$. If $L > 0$, then for any $r > 0$, $w^\infty_\star$ minimizes the energy: \begin{equation}
F_{L\h{0.5pt}\mu}\big(w, B_r\big) := \int_{B_r} \big| \h{.5pt} \nabla_\zeta w \h{.5pt} \big|^2
 +\dfrac{L\h{0.5pt}\mu}{2}\Big(\big| \h{.5pt} w \h{.5pt} \big|^2-1\Big)^2 
 \label{GL energy}
\end{equation}
in   $$\overline{H}\big(r, w^\infty_\star\big) := \Big\{ w \in H^1\big(B_r;\mathbb{R}^5\big) \h{1pt}:\h{1pt}
w=w^\infty_\star  \h{2pt} \text{ on } \h{2pt} \p B_r, \h{2pt}
w=\mathscr{L}[\h{0.5pt}u\h{0.5pt}] \h{2pt} \text{ for some $3$--vector field} \h{2pt}
u=u(\rho,z) \h{1pt}
\Big\}.$$ It also holds \begin{align}\label{uni bou of winf star} \sup_{R \h{1pt}> \h{1pt}0} \h{2pt}\dfrac{1}{ R} \int_{B_{ R}} \big| \nabla_\zeta w^\infty_\star \big|^2 + \dfrac{L \h{0.5pt}\mu}{2}  \left( \big| w^\infty_\star \big|^2 - 1 \right)^2 \h{2pt}\leq\h{2pt}8\pi.
\end{align} Moreover, $w^\infty_\star$ satisfies $\big| w^\infty_\star \big| \leq 1$ on $\mathbb{R}^3$ and equals  $0$ at the origin.  \label{lim map in interior core}
\end{lem}\begin{proof}[\bf Proof] Note that (\ref{unifo of scaled energy}) still holds in the current case. Since $a_n r_n^2 \to L$ as $n \to \infty$, there exists a $w^\infty_\star$ so that up to a subsequence, $w^{(n)}$ converges to $w_\star^\infty$ weakly in $H^1_{\mathrm{loc}}\big(\mathbb{R}^3; \mathbb{R}^5\big)$ and strongly in $L^4_{\mathrm{loc}}\big(\mathbb{R}^3; \mathbb{R}^5\big)$. By lower--semi continuity, (\ref{unifo of scaled energy}) then yields (\ref{uni bou of winf star}). In light that $a_n r_n^2$ is bounded for all $n$, we can also obtain the convergence of $w^{(n)}$ to $w^\infty_\star$ in $C^2_{\mathrm{loc}}\big(\mathbb{R}^3\big)$ as $n \to \infty$. Here one just needs the elliptic equation satisfied by $w^{(n)}$, standard Schauder's estimate and Arzel\`{a}--Ascoli theorem. \vspace{0.2pc}

Now we let $w$ be an arbitrary mapping in $\overline{H}\big(r, w^\infty_\star\big)$. By the energy minimality of $w^{(n)}$, it turns out \begin{align*}   \int_{B_{r}} \big| \nabla_\zeta w^{(n)} \big|^2 + r_n^2 \h{1pt}F_n\big(w^{(n)}\big) &\h{2pt}\leq\h{2pt}\int_{B_{r}} \big| \nabla_\zeta \big( w + w^{(n)} - w^\infty_\star\big) \big|^2 + r_n^2 \h{1pt}F_n\big(w + w^{(n)} - w^\infty_\star\big).
\end{align*}Utilizing $C^2_{\mathrm{loc}}\big(\mathbb{R}^3\big)$--convergence of $w^{(n)}$, we can take $n \to \infty$ in the above estimate and get \begin{align*} \int_{B_{r}} \big| \nabla_\zeta w^\infty_\star \big|^2 + \dfrac{L \h{0.5pt}\mu}{2}  \left( \big| w^\infty_\star \big|^2 - 1 \right)^2 &\h{2pt}\leq\h{2pt}\int_{B_{r}} \big| \nabla_\zeta w \big|^2 + \dfrac{L \h{0.5pt}\mu}{2}  \left( \big| w \big|^2 - 1 \right)^2.
\end{align*}Hence, $w^\infty_\star$ is a minimizer of $F_{L\h{0.5pt}\mu}\left( \cdot, B_r \right)$ in $\overline{H}\big(r, w^\infty_\star\big)$. \vspace{0.2pc}

Notice that $\big| w^{(n)} \big| \leq H_{a_n}$ and $w^{(n)}(0) = 0$. As $n \to \infty$, it turns out $\big| w^\infty_\star \big| \leq 1$ on $\mathbb{R}^3$ and meanwhile $w^\infty_\star(0) = 0$. If $L = 0$, then $w^\infty_\star$ minimizes the standard Dirichlet energy in $\overline{H}\big(r, w^\infty_\star\big)$ for all $r > 0$. Hence, $w^\infty_\star$ is harmonic over $\mathbb{R}^3$. Using the uniform boundedness of $w^\infty_\star$ on $\mathbb{R}^3$ and Liouville's theorem, it follows $w^\infty_\star \equiv 0$ on $\mathbb{R}^3$. 
\end{proof} 

In the remaining of this section, we characterize the limiting map $w^\infty_\star$ with $L > 0$.  \begin{lem}\label{characterization of wstar infty} If $L > 0$, then $w^\infty_\star$ in Lemma \ref{lim map in interior core} equals  $f \left( \sqrt{L\h{0.5pt}\mu} \h{2pt}|\h{1pt}\zeta \h{1pt}|\h{1pt}\right) \Lambda$. Here $\Lambda = \Lambda_+$ or $\Lambda_-$  in (\ref{def of lambda plus and minus}). The function $f$ is the radial function introduced in item (2) of Proposition \ref{strict isolation}.
\end{lem}\begin{proof}[\bf Proof] Without loss of generality, we assume $L\h{0.5pt}\mu = 1$. The following proof is motivated by Millot--Pisante \cite{MP10}, which relies on the division trick of Mironescu  \cite{M96}  and the blow--down analysis of Lin--Wang \cite{LW02}. Now we define $v = w^\infty_\star \big/ f$. It follows that \begin{align*} \Delta_\zeta v + f^2 \big( 1 -  |\h{1pt} v \h{1pt} |^2 \big)v = - 2 \h{1pt}\dfrac{f'}{f} \h{1pt}\dfrac{\zeta}{|\h{0.5pt}\zeta\h{0.5pt}|} \cdot \nabla_\zeta v - \dfrac{2}{|\h{0.5pt}\zeta\h{0.5pt}|^2}\h{1pt}v \h{15pt}\text{in $\mathbb{R}^3 \setminus \big\{0\big\}$.}
\end{align*}Multiplying this  equation by $\p_r v = \left( \dfrac{\zeta}{ |\h{0.5pt}\zeta\h{0.5pt}| } \cdot \nabla_\zeta\right) v $, we get \begin{align}\label{identity of poho}\big|\h{1pt} \p_r v \h{1pt}\big|^2 \left( \dfrac{1}{r} + \dfrac{2f'}{f} \right) + \dfrac{\big(\h{1pt}   |\h{.5pt} v  \h{.5pt} |^2 - 1  \big)^2}{2} f^2\left(\dfrac{1}{r} + \frac{f'}{f} \right)  = \nabla_\zeta \cdot \h{1pt}\Phi(\zeta), 
\end{align}where \begin{align*} \Phi(\zeta) :=  -  ( \nabla_\zeta v )^\top  \p_r v   
+ \dfrac{\zeta}{ 2\h{0.5pt}r } \h{1pt}   \big|\h{1pt}\nabla_\zeta v  \h{1pt}\big|^2
+ \dfrac{\zeta }{ 4\h{0.5pt}r } \h{1pt} f^2 \h{1pt} \Big(   |\h{.5pt} v  \h{.5pt} |^2 - 1  \Big)^2
- \dfrac{\zeta }{r^3} \h{1pt} \left( \h{1pt}|\h{.5pt} v  \h{.5pt} |^2 - 1\right).
\end{align*}For $R > \rho  > 0$, we integrate (\ref{identity of poho}) on $B_R \setminus B_\rho$. Hence, \begin{align}\label{integral identity poho}\int_{B_R \setminus B_\rho} \big|\h{1pt} \p_r v \h{1pt}\big|^2 \left( \dfrac{1}{r} + \dfrac{2f'}{f} \right) + \dfrac{\big(\h{1pt}   |\h{.5pt} v  \h{.5pt} |^2 - 1  \big)^2}{2} f^2\left(\dfrac{1}{r} + \frac{f'}{f} \right)  = \int_{\p B_R}  \Phi(\zeta) \cdot \dfrac{\zeta}{|\h{1pt}\zeta\h{1pt}|} - \int_{\p B_\rho}  \Phi(\zeta) \cdot \dfrac{\zeta}{|\h{1pt}\zeta\h{1pt}|}.
\end{align}

Firstly, we consider the behavior of $\Phi$ near the origin. Notice that $f'(0) > 0$. It then turns out \begin{align}\label{lower bound of f}f(r) \geq \dfrac{f'(0)}{2} \h{0.5pt}r\h{15pt}\text{for $r$ sufficiently small.} \end{align} With an use of mean value theorem, we have \begin{align*}  \big| \h{0.5pt}v(\zeta) \h{0.5pt} \big| \h{2pt}\leq \h{2pt} \dfrac{2}{f'(0) } \h{1pt}\dfrac{ \big| \h{0.5pt}w_\star^\infty\big(\zeta\big)\h{0.5pt}\big|}{|\h{1pt}\zeta\h{1pt}|} \h{2pt}=\h{2pt}\dfrac{2}{f'(0) } \h{1pt}\dfrac{ \big| \h{0.5pt}\nabla_\zeta w_\star^\infty\big(\zeta_*\big) \cdot \zeta\h{0.5pt}\big|}{|\h{1pt}\zeta\h{1pt}|}.
\end{align*}Here $|\h{0.5pt}\zeta\h{0.5pt}|$ is sufficiently small. $\zeta_*$ is on the segment connecting $0$ and $\zeta$. Since $w^\infty_\star$ is smooth on $B_1$, it then turns out from the above estimate that \begin{align}\label{unf bd of v near 0} \big| \h{0.5pt}v(\zeta) \h{0.5pt} \big| \h{2pt}\leq \h{2pt} 2 \h{1pt}\dfrac{\big\| \h{1pt}\nabla_\zeta w^\infty_\star \big\|_{\infty;\h{1pt} B_1}}{f'(0)} \h{15pt}\text{for $\zeta$ sufficiently close to $0$.}
\end{align}In addition, by L'Hospital's rule, \begin{align}\label{pointwise conv of v} v\big(\rho\h{1pt}\zeta\h{1.5pt}\big) \longrightarrow B\h{1pt}\zeta \h{15pt}\text{as $\rho \to 0$, for all $\zeta \in \p B_1$. Here $B:= \dfrac{\nabla_\zeta w^\infty_\star (0)}{f'(0)}$.} 
\end{align}As for the first--order derivatives of $v$, simple computations yield \begin{align*} \rho  \h{1pt}\p_{\zeta_j} v   \h{2pt}\Big|_{\rho \h{1pt}\widehat{\zeta}} = \dfrac{\rho}{f(\rho)} \left[ \h{2pt} \p_{\zeta_j} w^\infty_\star \h{2pt}\Big|_{\rho\h{1pt}\widehat{\zeta}}     -  v\big(\rho \h{1pt}\widehat{\zeta}\h{1pt}\big) \h{1pt}  f'(\rho)   \h{1.5pt} \dfrac{\zeta_j}{|\h{1pt}\zeta\h{1pt}|} \h{2pt}\right].
\end{align*} In light of (\ref{lower bound of f})--(\ref{unf bd of v near 0}), \begin{align}\label{unifo bound of gradient} \left|\h{2pt}\rho  \h{1pt}\p_{\zeta_j} v   \h{2pt}\Big|_{\rho \h{1pt}\widehat{\zeta}} \h{2pt}\right|\h{2pt}\leq\h{2pt}10\h{1pt}\dfrac{\big\| \h{1pt}\nabla_\zeta w^\infty_\star \big\|_{\infty; \h{1pt} B_1}}{f'(0)}  \h{15pt}\text{for $\rho$ sufficiently close to $0$.}
\end{align}Moreover, \begin{align}\label{pointwise conv of gradient}\rho  \h{1pt}\p_{\zeta_j} v  \h{2pt}\Big|_{\rho \h{1pt}\widehat{\zeta}}\longrightarrow \dfrac{\p_{\zeta_j} w^\infty_\star (0)}{f'(0)} - \big(\h{0.5pt}B \widehat{\zeta}\h{1.5pt}\big) \dfrac{\zeta_j}{|\h{1pt}\zeta\h{1pt}|} \h{15pt}\text{as $\rho \to 0$ pointwisely.}
\end{align}Now we compute \begin{align*}\int_{\p B_\rho}  \Phi(\zeta) \cdot \dfrac{\zeta}{|\h{0.5pt}\zeta\h{0.5pt}|} = \int_{\p B_1}  1 -  \big|\h{1pt}v(\rho\h{1pt}\zeta)\h{1pt}\big|^2   - \big| \h{1pt}\rho\h{0.8pt}\p_r v\h{1pt}\big|^2(\rho\h{1pt} \zeta) + \dfrac{1}{2}\h{1pt}\big|\h{1pt} \rho\h{0.5pt}\nabla_\zeta v \h{1pt}\big|^2(\rho\h{1pt}\zeta) + \dfrac{\rho^2 f^2(\rho)}{4}\h{1pt}\left(\h{1pt} \big| \h{1pt}v(\rho\h{1pt}\zeta)\h{1pt} \big|^2 - 1 \right)^2.
\end{align*}In light of the uniform bounds given in (\ref{unf bd of v near 0}) and (\ref{unifo bound of gradient}), we can apply Lebesgue's dominated convergence theorem to the right--hand side of the above equality. By the convergence in (\ref{pointwise conv of v}) and (\ref{pointwise conv of gradient}), it then follows \begin{align}\label{limit in small ball} \int_{\p B_\rho}  \Phi(\zeta) \cdot \dfrac{\zeta}{|\h{0.5pt}\zeta\h{0.5pt}|} \longrightarrow 4 \pi + \int_{\p B_1}  \dfrac{1}{2} \left| \nabla_\zeta  \big(B\h{1pt} \widehat{\zeta}\h{1.5pt}\big) \right|^2 - \left|\h{1pt}B\h{1pt}\zeta\h{1pt}\right|^2 = 4\pi \h{15pt}\text{as $\rho \to 0$.}
\end{align}Hence, if we take $\rho \to 0$ in (\ref{integral identity poho}), then it holds \begin{align}\label{integral identity poho with rho equal 0}\int_{B_R } \big|\h{1pt} \p_r v \h{1pt}\big|^2 \left( \dfrac{1}{r} + \dfrac{2f'}{f} \right) + \dfrac{\big(\h{1pt}   |\h{.5pt} v  \h{.5pt} |^2 - 1  \big)^2}{2} f^2\left(\dfrac{1}{r} + \frac{f'}{f} \right)  = \int_{\p B_R}  \Phi(\zeta) \cdot \dfrac{\zeta}{|\h{1pt}\zeta\h{1pt}|} - 4 \pi.
\end{align}

In the next, we study the behavior of $\Phi(x)$ near $\infty$. Let $R_n$ be a sequence diverging to $\infty$ as $n \to \infty$.  In light of (\ref{uni bou of winf star}), it turns out \begin{align*} \sup_{n \h{1pt}\in\h{1pt}\mathbb{N}}\h{2pt}   \int_{B_{ 1}} \big| \nabla_\zeta w_n^{\star} \big|^2 + \dfrac{R_n^2}{2}  \left( \big| w_n^{\star} \big|^2 - 1 \right)^2 \h{2pt}\leq\h{2pt}8\pi, \h{15pt}\text{where $w^\star_n (\zeta) := w_\star^\infty\big(R_n \zeta \big)$.}
\end{align*}Therefore, by Fatou's lemma, there is a $\sigma_\star \in (0, 1)$ so that up to a subsequence, it satisfies\begin{align}\label{uniform on large sphere} \lim_{n \to \infty} \int_{\p B_{\sigma_\star}} \big| \nabla_\zeta w_n^{\star} \big|^2 + \dfrac{R_n^2}{2}  \left( \big| w_n^{\star} \big|^2 - 1 \right)^2 \h{2pt}=\h{2pt}\lim_{n \to \infty} \int_{\p B_{R_n \sigma_\star}} \big| \nabla_\zeta w_\star^\infty \big|^2 + \dfrac{1}{2}  \left( \big| w_\star^\infty \big|^2 - 1 \right)^2 \h{2pt}\leq\h{2pt}8\pi.
\end{align}Now we compute \begin{align*}\int_{\p B_{R_n \sigma_\star}}  \Phi(\zeta) \cdot \dfrac{\zeta}{|\h{1pt}\zeta\h{1pt}|} = \int_{\p B_{R_n \sigma_\star}} - \big| \h{1pt} \p_r v\h{1pt}\big|^2  + \dfrac{1}{2}\h{1pt}\big|\h{1pt}  \nabla_\zeta v \h{1pt}\big|^2  + \dfrac{ f^2 }{4}\h{1pt}\big(\h{1pt} | \h{1pt}v \h{1pt} |^2 - 1 \big)^2 - \dfrac{1}{r^2}\big(\h{1pt} |\h{1pt}v\h{1pt}|^2 - 1 \h{1pt}\big).
\end{align*}Owing to (\ref{uniform on large sphere}) and (\ref{property of f}), we can rewrite the above equality as follows: \begin{align*} \int_{\p B_{R_n \sigma_\star}}  \Phi(\zeta) \cdot \dfrac{\zeta}{|\h{0.5pt}\zeta\h{0.5pt}|} = \int_{\p B_{R_n \sigma_\star}} - \big| \h{1pt} \p_r w^\infty_\star\h{1pt}\big|^2  + \dfrac{1}{2}\h{1pt}\big|\h{1pt}  \nabla_\zeta w^\infty_\star \h{1pt}\big|^2  + \dfrac{ 1 }{4}\h{1pt}\big(\h{1pt} | \h{1pt}w^\infty_\star \h{1pt} |^2 - 1 \big)^2 \h{2pt}\mathrm{d} \mathscr{H}^2 + o_n(1).
\end{align*}Here $o_n(1)$ is a quantity which converges to $0$ as $n \to \infty$.  This equality together with (\ref{uniform on large sphere}) infer \begin{align*} \limsup_{n \to \infty} \int_{\p B_{R_n \sigma_\star}}  \Phi(\zeta) \cdot \dfrac{\zeta}{|\h{0.5pt}\zeta\h{0.5pt}|} \h{2pt}\leq\h{2pt}4 \pi.
\end{align*}Now we take $R = R_n\h{0.5pt}\sigma_\star$ in (\ref{integral identity poho with rho equal 0}) and take $n \to \infty$. The last estimate yields\begin{align*} \int_{\mathbb{R}^3} \big|\h{1pt} \p_r v \h{1pt}\big|^2 \left( \dfrac{1}{r} + \dfrac{2f'}{f} \right) + \dfrac{\big(\h{1pt}   |\h{.5pt} v  \h{.5pt} |^2 - 1  \big)^2}{2} f^2\left(\dfrac{1}{r} + \frac{f'}{f} \right)  =0.
\end{align*}
This equality induces that $\big|\h{1pt}w^\infty_\star \h{1pt}\big| = f$. Moreover, $v$ is $0$--homogeneous. Due to  (\ref{pointwise conv of v}), it holds $v(\zeta) = B \h{1pt} \widehat{\zeta}$ for any $\zeta \neq 0$. In light of the unit length of $v$, the $\mathscr{R}$--axial symmetry of $v$ and the fact that $$v_4 \cos \theta + v_5 \sin \theta \geq 0 \h{15pt}\text{for $\phi \in \big[\h{1pt}0, \pi / 2 \h{1pt}\big]$,}$$  $v$ equals  either $\Lambda_+$ or $\Lambda_-$. The proof is completed.
\end{proof}

\subsubsection{Proof of (2) and (3) in Proposition \ref{strict isolation}}

We firstly prove the non--degeneracy result in (\ref{non degeneracy}). Suppose that $\big\{a_n\big\}$ and $\big\{z_n\big\}$ satisfy \begin{align*}\liminf_{a \to \infty} \h{2pt}\min_{\big\{\h{.5pt} z_a \h{1pt}: \h{1pt}w_{a}(z_a) \h{0.5pt}= \h{0.5pt}0\h{.5pt}\big\}} \min_{\big\{\h{0.5pt}x \h{0.5pt}:\h{0.5pt} | x - z_a | \h{1pt}\leq\h{1pt} R\h{.5pt} a^{- \frac{1}{2}} \big\}} \dfrac{\big|\h{1pt} w_{a}(x) \h{1pt}\big|}{\sqrt{a} \h{1.5pt} | x - z_a |} = \lim_{n \to \infty}  \min_{\big\{\h{0.5pt}x \h{0.5pt}:\h{0.5pt} | x - z_n | \h{1pt}\leq\h{1pt} R\h{.5pt} a_n^{- \frac{1}{2}} \big\}} \dfrac{\big|\h{1pt} w_{a_n}(x) \h{1pt}\big|}{\sqrt{a_n} \h{1.5pt} | x - z_n |}, 
\end{align*} where $w_{a_n}(z_n) = 0$. Changing variables by letting $\overline{w}^{(n)}(\zeta) := w_{a_n}\big(z_n + a_n^{- \frac{1}{2}} \zeta\big)$, we can rewrite the above equality as follows: \begin{align}\label{lininf achieved}\liminf_{a \to \infty} \h{2pt}\min_{\big\{\h{.5pt} z_a \h{1pt}: \h{1pt}w_{a}(z_a) \h{0.5pt}= \h{0.5pt}0\h{.5pt}\big\}} \min_{\big\{\h{0.5pt}x \h{0.5pt}:\h{0.5pt} | x - z_a | \h{1pt}\leq\h{1pt} R\h{.5pt} a^{- \frac{1}{2}} \big\}} \dfrac{\big|\h{1pt} w_{a}(x) \h{1pt}\big|}{\sqrt{a} \h{1.5pt} | x - z_a |} = \lim_{n \to \infty}  \min_{\big\{\h{0.5pt}\zeta \h{0.5pt}:\h{0.5pt} |\h{0.5pt}\zeta\h{0.5pt} | \h{1pt}\leq\h{1pt} R  \big\}} \dfrac{\big|\h{1pt} \overline{w}^{(n)}(\zeta) \h{1pt}\big|}{ |\h{1pt}\zeta\h{1pt} |}. 
\end{align}By Lemma \ref{result on wstar infty}, up to a subsequence, $\overline{w}^{(n)}$ converges to $w_\star^\infty$ in $C^2\big(\overline{B_R}\big)$. Moreover, it holds \begin{align*} \dfrac{\big|\h{1pt} \overline{w}^{(n)}(\zeta) - w_\star^\infty(\zeta) \h{1pt}\big|}{|\h{1pt}\zeta \h{1pt}|} =   \left| \int_0^1 \widehat{\zeta} \cdot \left[ \nabla_\zeta  \overline{w}^{(n)} \h{1pt}\Big|_{t \h{0.5pt}\zeta}  - \nabla_\zeta w^\infty_\star \h{1pt}\Big|_{t \h{0.5pt}\zeta} \right] \h{2pt}\mathrm{d} t\h{1pt}\right| \h{2pt}\leq\h{2pt}\left\| \h{1pt}\nabla_\zeta \overline{w}^{(n)} - \nabla_\zeta w_\star^\infty \h{1pt}\right\|_{\infty; \h{1pt}\overline{B_R}}, \h{5pt}\text{for all $\zeta \in \overline{B_R}$.}
\end{align*}Therefore, $\frac{\big| \overline{w}^{(n)}(\zeta) \big|}{|\h{1pt}\zeta\h{1pt}|}$ uniformly converges to $\frac{\big| w_\star^\infty(\zeta)\big|}{|\h{1pt}\zeta\h{1pt}|}$ on $\overline{B_R}$ as $n \to \infty$. By this uniform convergence, (\ref{lininf achieved}) and the characterization of $w^\infty_\star$ in Lemma \ref{characterization of wstar infty},  (\ref{non degeneracy}) follows.\vspace{0.2pc}

We use a contradictory argument to prove the strict isolation of zeros. Suppose that   there exists a sequence $\big\{a^*_n\big\}$ tending to $\infty$ so that $w_{a^*_n}$ has at least two different zeros, denoted by $z_n^{(1)}$ and $z_n^{(2)}$, on $l_z^+$. In addition, these two zeros satisfy \begin{align}\label{close of zeros}\left|\h{1pt}z_n^{(1)} - z_n^{(2)}\h{1pt}\right| \longrightarrow 0 \h{15pt}\text{ as $n \to \infty$.}\end{align}Without loss of generality, we can assume that $w_{a^*_n}$ converges  to some $w_\infty$ strongly in $H^1(B_1; \mathbb{R}^5)$ and $a_n^* [\h{1pt} \big| w_{a_n^*} \big|^2 - 1 \h{1pt}]^2$ converges to $0$ strongly in $L^1\big(B_1\big)$ as $ n \to \infty$.   Moreover, by (\ref{close of zeros}), we let $z_n^{(1)}$ and $z_n^{(2)}$ converge to $y_0$ as $n \to \infty$. Here $y_0$ is a singularity of $w_\infty$. Taking $\sigma_0 > 0$ suitably small so that $y_0$ is the unique singularity of $w_\infty$ in the closure of $B_{\sigma_0}(y_0)$, we define \begin{align*}\mathcal{V}^*_n := \overline{B_{\sigma_0}(y_0)} \h{2pt}\bigcap \h{2pt}\left\{\h{1pt} \big| w_{a^*_n} \big| \leq \dfrac{1}{4}\h{1pt} \right\} \end{align*} and let $\displaystyle \nu_n^* := \max_{z\h{1pt}\in\h{1pt}\mathcal{V}^*_n}\h{2pt}\big|\h{1pt}z - z_n^{(1)}\h{1pt}\big|$. By Lemma \ref{size of core}, there is a $R_* > 0$ so that $\sqrt{a_n^*} \h{1pt}\nu_n^* \leq R_*$ for all $n$. Moreover,  (\ref{non degeneracy}) infers \begin{align*} \big|\h{1pt}w_{a^*_n}(x)\h{1pt}\big| \h{2pt}\geq\h{2pt}\dfrac{c_\mu\big(2 R_*\big) }{2} \h{1pt}\sqrt{a^*_n} \h{2pt} \left|\h{1pt} x - z_n^{(1)}\h{1pt}\right|, \h{15pt}\text{for large $n$ and  any $x$ satisfying $\left|\h{1pt} x - z_n^{(1)}\h{1pt}\right| \h{1pt}\leq\h{1pt} \frac{2 R_*}{\sqrt{a^*_n}}$.}
\end{align*}Since $\left|\h{1pt} z_n^{(2)} - z_n^{(1)} \h{1pt}\right| \h{1pt}\leq\h{1pt} \nu_n^* \h{1pt}\leq\h{1pt} \dfrac{R_*}{\sqrt{a_n^*}}$, the last estimate yields \begin{align*} \left|\h{1pt}w_{a^*_n}\big(z_n^{(2)}\big)\h{1pt}\right| \h{2pt}\geq\h{2pt}\dfrac{c_\mu\big(2 R_*\big) }{2} \h{1pt}\sqrt{a^*_n} \h{2pt} \left|\h{1pt} z_n^{(2)} - z_n^{(1)}\h{1pt}\right| \h{2pt}>\h{2pt}0.
\end{align*}However, this is impossible since $z_n^{(2)}$ is also a zero of $w_{a_n^*}$. The proof is completed.
\subsection{Asymptotic behavior of phase mapping near zeros}

Throughout the remaining arguments, the parameter $a$ is always assumed to be large enough. Due to the items (1) and (3) in Proposition \ref{strict isolation},   there exists a $\delta_2 > 0$ so that if $z_a$ is an arbitrary zero of $w_{a}$ on $l_z$, then it is the unique zero  of $w_{a}$ in $B_{\delta_2}(z_a)$.  Hence the total number of zeros of $w_{a}$ is uniformly bounded from above. Moreover,  the phase mapping $\widehat{w_{a}}$ is well--defined except at  finitely many zeros of $w_{a}$. This section is devoted to studying the asymptotic behavior of $\widehat{w_{a}}$ near each zero of $w_{a}$.  Our main result is read as follows:  \begin{prop}\label{auxialliary result for prop} Let $\big\{ z_{a, 1}, ..., z_{a, k_a}\big\}$  be the family of zeros of $w_{a}$ on $l_z^+$, where $k_a$ is  the total number of zeros of $w_a$ on $l_z^+$. Then it holds \begin{align*}\lim_{(a^{-1}, \h{0.5pt}r) \h{0.5pt}\to\h{0.5pt} (0, 0)} \h{1.5pt}\max_{k \h{0.5pt}=\h{0.5pt} 1, ..., k_a} \h{2pt}\min_{\Lambda \h{1pt}\in\h{1pt}\{\h{0.5pt} \Lambda_+, \Lambda_-\}} \h{1.5pt}\sum^2_{j=0} \h{1.5pt}r^j \h{1pt} \Big\lVert \h{1pt} \nabla^j \widehat{w_{a}}  - \nabla^j \big[\h{1pt}\Lambda  \big(\cdot-z_{a, k}\big)\h{1pt}\big]\h{1.5pt} \Big\rVert_{\infty;\h{1pt} \p B_{r}  (z_{a, k} ) } = 0.
\end{align*}
\end{prop} In light of the values of $w_a$ at the north pole and  the origin, the proof of Proposition \ref{num and asy of singu} follows easily from Proposition \ref{auxialliary result for prop}. \begin{proof}[\bf Proof of Proposition \ref{auxialliary result for prop}]  We assume on the contrary that Proposition \ref{auxialliary result for prop} fails. Then there are $\epsilon_0>0$, $\{ a_n \}$, $\{ r_n \}$ and $\{ z_n \}$ such that
\begin{equation}  \min_{\Lambda \h{1pt}\in\h{1pt}\{\h{0.5pt} \Lambda_+, \Lambda_-\}} \h{1.5pt}\sum^2_{j=0} \h{1.5pt}r_n^j \h{1pt} \Big\lVert \h{1pt} \nabla^j \widehat{w_{a_n}}  - \nabla^j \big[\h{1pt}\Lambda  \big(\cdot - z_n\big)\h{1pt}\big]\h{1.5pt} \Big\rVert_{\infty;\h{1pt} \p B_{r_n}  (z_n) } \geq \epsilon_0. 
\label{contradiction, convergence of hat w_a to tangent map, ori}
\end{equation}Here $a_n \to \infty$ and $r_n \to 0$ as $n \to \infty$. $z_n$ is a zero of $w_{a_n}$ on $l_z^+$. Without loss of generality, we can  assume that $w_{a_n}$ converges to some $w_\infty$ strongly in $H^1(B_1)$ and $a_n [\h{1pt}\big| w_{a_n} \big|^2 - 1 \h{1pt}]^2$ converges to $0$ strongly in $L^1(B_1)$ as $n \to \infty$. Moreover,  $z_n$  converges to some $y_0$  as $n \to \infty$, where $y_0$ is a singularity of  $w_\infty$. Still using  $w^{(n)}$ to denote the scaled mapping $w_{a_n}\big(z_n + r_n\h{1pt}\zeta\big)$, we then rewrite the assumption (\ref{contradiction, convergence of hat w_a to tangent map, ori}) as follows: \begin{align}\label{contradiction, convergence of hat w_a to tangent map} \min_{\Lambda \h{1pt}\in\h{1pt}\{\h{0.5pt} \Lambda_+, \Lambda_-\}} \h{1.5pt}\sum^2_{j=0} \h{1.5pt} \Big\lVert \h{1pt} \nabla_\zeta^j  \widehat{w^{(n)}}   - \nabla_\zeta^j  \Lambda \h{1.5pt} \Big\rVert_{\infty;\h{1pt} \p B_1 } \geq \epsilon_0 \h{15pt}\text{for all $n$.}
\end{align}

Owing to the equation (\ref{el-eq of u_a, intro}) satisfied by $w_a$ on $B_1$, for large $n$, we have \begin{align}
\Delta_\zeta   \widehat{w^{(n)}} + \Big| \nabla_\zeta  \widehat{w^{(n)}} \Big|^2 \widehat{w^{(n)}}
= &-\dfrac{2 \h{0.5pt}\nabla_\zeta  \big|w^{(n)}\big|}{\big|w^{(n)}\big|} \cdot \nabla_\zeta  \widehat{w^{(n)}}   \nonumber\\[3mm]
& - \dfrac{3 \h{1pt} \mu \h{1pt} r_n^2}{\sqrt{2} \h{1pt}} \left\{   \frac{\nabla_{w} S\h{1pt} [\h{1pt}w\h{1pt}] }{|\h{0.5pt}w\h{0.5pt}|} - \widehat{w} \h{1.5pt} \frac{ \widehat{w} \cdot \nabla_w S \h{1pt} [ \h{1pt} w  \h{1pt} ]}{|\h{0.5pt}w\h{0.5pt}|} \h{1pt}  \right\} \bigg|_{w = w^{(n)}} \h{15pt} \text{on $B_4 \setminus \{ 0 \}$.}
\label{eq of hat w^n}
\end{align}Now, we use this equation to show a contradiction to (\ref{contradiction, convergence of hat w_a to tangent map}) when $n$ is large. The arguments below are divided into two cases.\\[2mm] 
\noindent \textbf{Case I.} In this case, we suppose that  $a_n\h{0.5pt}r_n^2 \to \infty$ as $n \to \infty$. \\[4mm]
\textbf{I.1. $C^{1, \alpha}$--\h{1pt}estimate of $\widehat{w^{(n)}}$.} \vspace{0.6pc}

By  Lemmas \ref{size of core}, \ref{stong conv. and min. property of w^infty} and \ref{w^infty is homogeneous zero}, $w^{(n)}$ converges to $\Lambda$ strongly in $H^1(B_4)$ as $n \to \infty$. Here $\Lambda$ equals  either $\Lambda_+$ or $\Lambda_-$ in (\ref{def of lambda plus and minus}). Moreover, we also have $a_n\h{0.5pt}r_n^2 [ \h{1pt}| w^{(n)} |^2 - 1 \h{1pt}]^2$ converges to $0$ strongly in $L^1\big(B_4\big)$ as $n \to \infty$. Notice that $\Lambda$ is smooth on $B_3 \setminus B_{1/3}$.  We therefore can apply Lemma \ref{clearing lemma for scaling} to obtain \begin{align}\label{lo bd on strip}\big| w^{(n)} \big| > 1/2 \h{15pt}\text{ on $B_2 \setminus B_{1/2}$ for large $n$.}\end{align} In addition, the local gradient estimate in Lemma \ref{energy density e^L_a[w_a] bdd} infers that $\big| \nabla_\zeta w^{(n)} \big|$ is uniformly bounded in the thin shell $B_{1 + 4r_\star} \setminus B_{1 - 4r_\star}$, where $r_\star \in \big(0, 1/8\big)$ is sufficiently small. Hence, $\Delta_\zeta \widehat{w^{(n)}}$ is uniformly bounded on the thin shell $B_{1 + 4r_\star} \setminus B_{1 - 4r_\star}$ by the uniform boundedness of $\nabla_\zeta w^{(n)}$, (\ref{lo bd on strip}) and the equation (\ref{eq of hat w^n}).  Standard interior $L^p$--estimate for elliptic equations and Morrey's inequality then yield the uniform boundedness of $\widehat{w^{(n)}}$ in $C^{1, \alpha}\big(B_{1  + 3 r_\star} \setminus B_{1 - 3 r_\star}\big)$ for any $\alpha \in (0, 1)$.\\[4mm]
\noindent \textbf{I.2. Uniform bound of $a_n \h{0.5pt}r_n^2 \h{0.5pt}\left| \big| w^{(n)} \big|^2 - 1 \right|$ near $\p B_1$.} \vspace{0.6pc}

In the next, we  show the uniform boundedness result of $a_n \h{0.5pt}r_n^2 \h{0.5pt}\left| \big| w^{(n)} \big|^2 - 1 \right|$ near $\p B_1$ by following the idea of Bethuel--Brezis--H\'{e}lein \cite{BBH93}. Still by (\ref{el-eq of u_a, intro}), the partial differential equation satisfied by $\big| w^{(n)} \big|^2$ can be read as follows: \begin{align}
\dfrac{1}{2}\h{0.5pt}\Delta_\zeta \h{0.5pt}   \big|  w^{(n)} \big|^2 
= \big| \nabla_\zeta    w^{(n)} \big|^2
- \dfrac{9 \h{1pt} \mu \h{1pt} r_n^2}{\sqrt{2}}\h{1pt} S \h{1pt}\big[\h{1pt} w^{(n)}\h{1pt}\big]
+ a_n \h{0.5pt} r_n^2 \h{0.5pt} \mu \h{0.5pt}\Big(\big|  w^{(n)} \big|^2-1\Big)\h{0.5pt}
\big|  w^{(n)} \big|^2
\h{15pt} \text{ on $B_4$.}
\label{eq for |w^n|^2}
\end{align}Utilizing (\ref{lo bd on strip}) and the fact that $\big| w^{(n)} \big| \leq H_{a_n}$, we get from the above equation that \begin{align}\label{lo bd of max prin}\dfrac{1}{2}\h{1pt}\Delta_\zeta  \Big(H_{a_n}^2 - \big| w^{(n)} \big|^2 \Big)
&\h{1pt}\geq \h{1pt} \dfrac{a_n \h{1pt} r_n^2 \h{1pt} \mu}{4} \Big(H_{a_n}^2 - \big| w^{(n)} \big|^2 \Big) \nonumber\\[2mm]
&  -\big| \nabla_\zeta  w^{(n)} \big|^2 - a_n \h{0.5pt} r_n^2 \h{0.5pt} \mu\h{1pt}\Big(H_{a_n}^2-1\Big)\h{1pt}\big| w^{(n)} \big|^2 +  \dfrac{9\h{1pt} \mu \h{1pt} r_n^2}{\sqrt{2}}\h{1pt} S\h{1pt}\big[\h{0.5pt}w^{(n)}\h{0.5pt}\big] \h{10pt}\text{on $B_2 \setminus B_{1/2}$.}
\end{align}Recall that $\nabla_\zeta w^{(n)}$ is uniformly bounded on  $B_{1 + 4 r_\star}\setminus B_{1 - 4 r_\star}$. In light of this uniform boundedness, the upper bound of $a_n \big(H_{a_n}^2 - 1 \big)$ when $n$ is large and the uniform boundedness of $w^{(n)}$, there is a positive constant $c_{\star}$ depending only on $r_\star$ and $\mu$ so that \begin{align*} \left|\h{2pt}-\big| \nabla_\zeta  w^{(n)} \big|^2 - a_n \h{0.5pt} r_n^2 \h{0.5pt} \mu\h{1pt}\Big(H_{a_n}^2-1\Big)\h{1pt}\big| w^{(n)} \big|^2 +  \dfrac{9\h{1pt} \mu \h{1pt} r_n^2}{\sqrt{2}}\h{1pt} S\h{1pt}\big[\h{0.5pt}w^{(n)}\h{0.5pt}\big] \h{2pt}\right| \h{2pt}\leq\h{2pt}\dfrac{c_\star \mu}{4} \h{15pt}\text{on $B_{1 + 4 r_\star} \setminus B_{1 - 4 r_\star}$.} 
\end{align*}Applying the above estimate to the right--hand side of (\ref{lo bd of max prin}), we have 
\begin{align*}\Delta_\zeta  \Psi_n -  \dfrac{a_n \h{0.5pt} r_n^2 \h{0.5pt} \mu}{2}\h{1pt}\Psi_n \h{2pt}\geq\h{2pt} 0 \h{15pt} \text{in $B_{1 + 4 r_\star} \setminus B_{1 - 4 r_\star}$. \h{1.5pt}Here $\Psi_n := H_{a_n}^2 - \big| w^{(n)} \big|^2 -  c_\star \left(a_n \h{1pt} r_n^2\right)^{-1}$.}\end{align*}
Utilizing the lower bound in (\ref{lo bd on strip}) and the fact that $H_a \longrightarrow 1 $ as $a \to \infty$, we can choose $N  \in \mathbb{N}$ sufficiently large such that
\begin{align*} \Psi_n \leq 1 \h{10pt} \text{on} \h{5pt} \p B_{1 + 4 r_\star} \h{1.5pt}\bigcup \h{1.5pt}\p B_{1 - 4 r_\star}  \h{15pt} \text{  for any $n > N$}. \end{align*}
Now we pick up a comparison function $\eta_n(\zeta) := \exp\left\{\frac{\sqrt{a_n \h{.5pt} r^2_n \h{.5pt} \mu}}{4}\big(\h{1pt}r - (1 - 4 r_\star) \h{1pt}\big)\big(\h{1pt}r - (1 + 4 r_\star)\h{1pt}\big)\right\}$, where $r = |\h{0.5pt}\zeta\h{0.5pt}|$.  If we keep taking  $N$ large, then it satisfies 
\begin{align*}\Delta_\zeta \eta_n - \frac{a_n \h{0.5pt} r_n^2 \h{0.5pt} \mu}{2} \h{1.5pt}\eta_n \h{2pt}<\h{2pt}0 \h{10pt} \text{in} \h{5pt} B_{1 + 4 r_\star}\setminus B_{1 - 4 r_\star} \h{10pt} \text{and}  \h{10pt} \eta_n \equiv  1 \h{10pt} \text{on} \h{5pt} \p B_{1 + 4 r_\star} \h{1pt}\bigcup\h{1pt}\p B_{1 - 4 r_\star} \h{2pt},\h{10pt} \text{for $n > N$}.\end{align*}
Hence, for all $n > N$, 
\begin{align*} \Delta_\zeta  \big(\Psi_n - \eta_n\big) - \frac{a_n \h{0.5pt} r_n^2 \h{0.5pt} \mu}{2}\h{1.5pt}\big(\Psi_n - \eta_n \big) > 0 \h{10pt} \text{in} \h{5pt} B_{1 + 4 r_\star}\setminus B_{1 - 4 r_\star} \h{10pt} \text{and}  \h{10pt} \Psi_n - \eta_n\leq 0 \h{10pt} \text{on} \h{5pt} \p B_{1 + 4 r_\star}\h{1.5pt}\bigcup\h{1.5pt}\p B_{1 - 4 r_\star}.\end{align*}
	Due to the maximum principle, we obtain \begin{align*} \Psi_n(\zeta) \h{2pt}\leq\h{2pt}\eta_n(\zeta) \h{2pt}\leq\h{2pt}\exp\left\{-  \sqrt{a_n \h{.5pt} r^2_n \h{.5pt} \mu} \h{1.5pt} r_\star^2 \right\} \h{15pt}\text{on $B_{1 + 3 r_\star} \setminus B_{1 - 3 r_\star}$.}
\end{align*}This estimate then yields \begin{align}\label{uni bd of an rn hn}a_n\h{0.5pt}r_n^2 \h{0.5pt}\left|\h{1pt} \big| w^{(n)} \big|^2 - 1 \h{1pt}\right| \h{2pt}\leq\h{2pt}c_\star \h{15pt}\text{on $B_{1 + 3 r_\star}\setminus B_{1 - 3r_\star}$.}
\end{align}Here $c_\star$ depends on $r_\star$ and $\mu$.\\[4mm]
\noindent \textbf{I.3. $C^{1, \alpha}$--\h{1pt}estimate of $\big| w^{(n)}\big|$.} \vspace{0.6pc}

 In light of the uniform bound of $\nabla_\zeta w^{(n)}$ obtained in \text{I.1} and (\ref{uni bd of an rn hn}), the right--hand side of (\ref{eq for |w^n|^2}) is uniformly bounded when it is restricted on $B_{1 + 3 r_\star}\setminus B_{1 - 3 r_\star}$. Standard $L^p$--estimate for elliptic equations infers that $\big| w^{(n)}\big|^2$ is uniformly bounded in $C^{1, \alpha}\big(B_{1 + 2r_\star}\setminus B_{1 - 2 r_\star}\big)$ for any $\alpha \in (0, 1)$. Then by the lower bound in (\ref{lo bd on strip}), $\big| w^{(n)}\big|$ is uniformly bounded in $C^{1, \alpha}\big(B_{1 + 2r_\star}\setminus B_{1 - 2 r_\star}\big)$ for any $\alpha \in (0, 1)$.\\[4mm]
\noindent \textbf{I.4. Contradiction to (\ref{contradiction, convergence of hat w_a to tangent map}) in Case I.}\vspace{0.6pc}

 Using the $C^{1, \alpha}$--estimate of $\widehat{w^{(n)}}$ in I.1, the lower bound in (\ref{lo bd on strip}) and the $C^{1, \alpha}$--estimate of $\big| w^{(n)}\big|$ in I.3, we have the uniform boundedness of $\Delta_\zeta \widehat{w^{(n)}}$ in $C^\alpha\big(B_{1 + 2r_\star}\setminus B_{1 - 2 r_\star}\big)$. Here we have also used the equation (\ref{eq of hat w^n}). By Schauder's estimate and Arzel\`{a}--Ascoli theorem, up to a subsequence, $\widehat{w^{(n)}}$ converges in $C^2\big(B_{1 + r_\star}\setminus B_{1 - r_\star}\big)$ as $n \to \infty$. Since $w^{(n)}$ converges to $\Lambda$, then by (\ref{uni bd of an rn hn}), the limit of $\widehat{w^{(n)}}$ on $B_{1 + r_\star}\setminus B_{1 - r_\star}$ equals  $\Lambda$ as well. We therefore obtain a contradiction to (\ref{contradiction, convergence of hat w_a to tangent map}). \vspace{0.6pc}

\noindent \textbf{Case II.} In this case, we assume that $a_n\h{0.5pt}r_n^2 \to L$ as $n \to \infty$. Here $L$ is a finite non--negative constant.\\[4mm]
\noindent \textbf{II.1. $C^{1, \alpha}$--\h{1pt}estimate of $\widehat{w^{(n)}}$.}\vspace{0.6pc}

Define $\overline{w}^{(n)}\big(\h{0.5pt}\zeta'\h{0.5pt}\big) := w_{a_n}\big(z_n + a_n^{- \frac{1}{2}} \zeta'\h{1pt}\big)$. By Lemmas \ref{lim map in interior core} and \ref{characterization of wstar infty}, for any $R > 0$, it holds \begin{align*}\overline{w}^{(n)}(\zeta') \longrightarrow f\big(\sqrt{\mu}\h{1pt}|\h{1pt}\zeta'\h{1pt}|\h{1pt}\big) \Lambda(\zeta') \h{15pt}\text{ in $C^2\big(B_R\big)$ as $n \to \infty$.}\end{align*}Here $\Lambda$ still equals  $\Lambda_+$ or $\Lambda_-$ in (\ref{def of lambda plus and minus}). $f$ is the radial function defined in the item (2) of Proposition \ref{strict isolation}. Changing variables by letting $\zeta' = \sqrt{a_n} \h{0.5pt}r_n \h{0.5pt} \zeta$, we have from the last convergence that  \begin{equation*}
     \dfrac{1}{\sqrt{a_n}\h{0.5pt} r_n }\h{1.5pt}\left\lVert \nabla_\zeta w^{(n)} - \nabla_\zeta\Big[ f\Big(\sqrt{a_n\h{0.5pt}\mu}\h{1.5pt} r_n \h{1pt} |\h{1pt} \zeta \h{1pt}|\h{1pt}\Big) \Lambda(\zeta)\h{1pt}\Big]\h{1.5pt} \right\rVert_{\infty; \h{1pt}B_{\frac{R}{\sqrt{a_n}\h{0.5pt}r_n}}} \longrightarrow 0 \h{15pt}\text{as $n \to \infty$.}
\end{equation*}If $L = 0$, then we take $R = 1$. It follows $\frac{1}{\sqrt{a_n}\h{0.5pt}r_n} > 4$ for large $n$. If $L > 0$, then we take $R = 8 \sqrt{L}$. It turns out that $\frac{8\sqrt{L}}{\sqrt{a_n}\h{0.5pt}r_n} > 4$ for large $n$. Hence, for any $L \geq 0$, the above convergence induces \begin{align}\label{key gradient conv}\dfrac{1}{\sqrt{a_n}\h{0.5pt} r_n }\h{1.5pt}\left\lVert \nabla_\zeta w^{(n)} - \nabla_\zeta\Big[ f\Big(\sqrt{a_n\h{0.5pt}\mu}\h{1.5pt} r_n \h{1pt} |\h{1pt} \zeta \h{1pt}|\h{1pt}\Big) \Lambda(\zeta)\h{1pt}\Big] \h{1.5pt} \right\rVert_{\infty; \h{1pt}B_4} \longrightarrow 0 \h{15pt}\text{as $n \to \infty$.}
\end{align}In addition, the non--degeneracy result (\ref{non degeneracy}) yields \begin{align} \label{non--deg, app} \big|\h{1pt}w^{(n)}(\zeta)\h{1pt}\big| \h{2pt}\geq\h{2pt}c_{\mu, L}\h{1pt}\sqrt{a_n}\h{0.5pt}r_n\h{1pt}| \h{1pt}\zeta\h{1pt}| \h{15pt}\text{on $B_4$, where $n$ is large and $c_{\mu, L} = \dfrac{1}{2} c_\mu \left(\max \Big\{1, 8\sqrt{L}\Big\}\right)$.}
\end{align}Note that the constant $c_\mu(R)$  has been given in (\ref{non degeneracy}). In light of (\ref{key gradient conv}), it holds \begin{align}\label{key gradient conv, upper bound}\dfrac{1}{\sqrt{a_n}\h{0.5pt} r_n }\h{1.5pt}\left\lVert \nabla_\zeta w^{(n)} \right\rVert_{\infty; \h{1pt}B_4} &\h{2pt}\leq\h{2pt}1 + \dfrac{1}{\sqrt{a_n}\h{0.5pt} r_n }\h{1.5pt}\left\lVert  \nabla_\zeta\Big[ f\Big(\sqrt{a_n\h{0.5pt}\mu}\h{1.5pt} r_n \h{1pt} |\h{1pt} \zeta \h{1pt}|\h{1pt}\Big) \Lambda(\zeta)\h{1pt}\Big] \h{1.5pt} \right\rVert_{\infty; \h{1pt}B_4}\nonumber\\[2mm]
&\h{2pt}\lesssim\h{2pt}1 + \big\| \h{1pt}f'\h{1pt}\big\|_{\infty; \h{1pt}[\h{1pt}0, \infty)} + \sup_{r \h{1pt}\in\h{1pt}[\h{1pt}0, \infty)} \dfrac{f(r)}{r} \h{35pt}\text{for large $n$}.
\end{align}On the other hand, (\ref{non--deg, app}) induces \begin{align}\label{low boudd of wn} \big|\h{1pt}w^{(n)} \h{1pt}\big| \h{2pt}\geq\h{2pt}\frac{c_{\mu, L}}{4}\h{1pt}\sqrt{a_n}\h{0.5pt}r_n \h{15pt}\text{ on $B_4 \setminus B_{1/4}$ for large $n$.}\end{align} By this pointwise lower bound and (\ref{key gradient conv, upper bound}), it turns out \begin{align}\label{point control of gradient in terms of function} \left| \h{1pt}\nabla_\zeta w^{(n)}\h{1pt}\right|\h{2pt}\leq \overline{c}_{\mu, L}\h{1pt}\left|\h{1pt}w^{(n)}\h{1pt}\right| \h{15pt}\text{on $B_4 \setminus B_{1/4}$ pointwisely, where  $n$ is large.}
\end{align}In (\ref{point control of gradient in terms of function}), $\overline{c}_{\mu, L}$ is a positive constant depending on $\mu$ and $L$. By (\ref{point control of gradient in terms of function}) and the equation (\ref{eq of hat w^n}), $\Delta \widehat{w^{(n)}}$ is uniformly bounded on $B_4 \setminus B_{1/4}$. Standard interior $L^p$--estimate for elliptic equations and Morrey's inequality then yield the uniform boundedness of $\widehat{w^{(n)}}$ in $C^{1, \alpha}\big(B_3 \setminus B_{1/3}\big)$ for any $\alpha \in (0, 1)$. Here $n$ is taken large. \\[4mm]
\noindent \textbf{II.2. $C^{\alpha}$--\h{1pt}estimate of $\nabla_{\zeta}\log \big| w^{(n)}\big|$.} \vspace{0.6pc}

Recalling (\ref{key gradient conv, upper bound}), we can apply mean value theorem to obtain \begin{align*} \left|\h{1pt}w^{(n)}\h{1pt}\right| \h{2pt}\leq\h{2pt}4 \h{1pt}\left\| \nabla_\zeta w^{(n)}\right\|_{\infty;\h{1pt} B_4} \h{2pt}\lesssim\h{2pt}\sqrt{a_n}\h{0.5pt}r_n \h{15pt}\text{on $B_4$.}
\end{align*}This upper bound together with (\ref{low boudd of wn}) show that $\log \big|w^{(n)}\big| - \log \big(\sqrt{a_n}\h{0.5pt}r_n\big)$ is uniformly bounded from both above and below on $B_4 \setminus B_{1/4}$. The upper and lower bounds are independent of $n$.
By rewriting  (\ref{eq for |w^n|^2}), it turns out \begin{align*} \Delta_\zeta \log \big| w^{(n)}\big| = \frac{\big| \nabla_\zeta    w^{(n)} \big|^2 - 2 \h{1pt}\big| \h{1pt}\nabla_\zeta \big|\h{1pt} w^{(n)} \h{1pt}\big| \h{1pt}\big|^2}{\big| w^{(n)} \big|^2}
- \dfrac{9 \h{1pt} \mu \h{1pt} r_n^2}{\sqrt{2}}\h{1pt} \frac{S \h{1pt}\big[\h{1pt} w^{(n)}\h{1pt}\big]}{\big| w^{(n)}\big|^2}
+ a_n \h{0.5pt} r_n^2 \h{0.5pt} \mu \h{0.5pt}\Big(\big|  w^{(n)} \big|^2-1\Big) 
\h{15pt} \text{ on $B_4 \setminus \big\{0\big\}$.}
\end{align*}In light of (\ref{point control of gradient in terms of function}), the convergence of $a_n\h{0.5pt}r_n^2$ as $n \to \infty$ and the uniform boundedness of $w^{(n)}$,  the last equation infers that $\Delta_\zeta \big(\log \big|w^{(n)}\big| - \log \big(\sqrt{a_n}\h{0.5pt}r_n\big) \big)$ is uniformly bounded on $B_4 \setminus B_{1/4}$. Standard interior $L^p$--estimate for elliptic equations and Morrey's inequality then yield the uniform boundedness of $\log \big|w^{(n)}\big| - \log \big(\sqrt{a_n}\h{0.5pt}r_n\big) $ in $C^{1, \alpha}\big(B_3 \setminus B_{1/3}\big)$ for any $\alpha \in (0, 1)$. Therefore, $\nabla _{\zeta}\log \big| w^{(n)}\big|$ is uniformly bounded in $C^\alpha\big(B_3 \setminus B_{1/3}\big)$ for all $\alpha \in (0, 1)$.\\[4mm]
\noindent \textbf{II.3. Contradiction to (\ref{contradiction, convergence of hat w_a to tangent map}) in Case II.}\vspace{0.6pc}

By the uniform boundedness of $\widehat{w^{(n)}}$ in $C^{1, \alpha}\big(B_3 \setminus B_{1/3}\big)$ and the uniform boundedness of $\nabla _{\zeta}\log \big| w^{(n)}\big|$ in $C^\alpha\big(B_3 \setminus B_{1/3}\big)$, it can be shown from (\ref{eq of hat w^n}) that $\Delta \widehat{w^{(n)}}$ is uniformly bounded in $C^\alpha \big(B_3 \setminus B_{1/3}\big)$. By Schauder's estimate and Arzel\`{a}--Ascoli theorem, up to a subsequence, $\widehat{w^{(n)}}$ converges in $C^2\big(B_{2}\setminus B_{1/2}\big)$ as $n \to \infty$. Finally, we determine the limit of $\widehat{w^{(n)}}$. On $B_4 \setminus \big\{0\big\}$, it holds by triangle inequality that
\begin{align*} \Big| \widehat{w^{(n)}}  - \Lambda \Big|
\h{2pt}\leq\h{2pt} \left| \widehat{w^{(n)}}  - \dfrac{f\big(\sqrt{a_n\h{0.5pt}\mu} \h{1.5pt} r_n \h{1pt}|\h{1pt} \zeta \h{1pt}|\h{1pt}\big)}{\big| w^{(n)} \big|}\Lambda  \right|
+  \left| \dfrac{f\big(\sqrt{a_n\h{0.5pt}\mu} \h{1.5pt} r_n \h{1pt} |\h{1pt} \zeta \h{1pt}|\h{1pt}\big)}{\big| w^{(n)} \big|}\Lambda  - \Lambda \h{1pt}\right|
\h{2pt}\leq\h{2pt} 2\left| \dfrac{w^{(n)} - f\big(\sqrt{a_n\h{0.5pt}\mu} \h{1pt}r_n \h{1pt}|\h{1pt} \zeta \h{1pt}|\h{1pt}\big)\Lambda}{\big| w^{(n)} \big|} \right| .\end{align*}
Using (\ref{non--deg, app}), mean value theorem and (\ref{key gradient conv}), we obtain
\begin{align*} \Big| \widehat{w^{(n)}}  - \Lambda \h{1pt} \Big|
\h{2pt}\leq\h{2pt} \dfrac{2}{c_{\mu, L}} \dfrac{\left|w^{(n)} - f \big(\sqrt{a_n\h{0.5pt}\mu} \h{1pt} r_n \h{1pt} |\h{1pt} \zeta \h{1pt}|\h{1pt}\big)\Lambda  \h{1pt}\right| }{ \sqrt{a_n} \h{0.5pt} r_n \h{0.5pt} |\h{1pt} \zeta \h{1pt}|} \longrightarrow 0 \h{10pt} \text{ as $n \to \infty$ on $B_4  \setminus \big\{0\big\}$.} \end{align*}Therefore, $\widehat{w^{(n)}}$ also converges to $\Lambda$ in $C^2\big(B_2 \setminus B_{1/2}\big)$ as $n \to \infty$. This is a contradiction to (\ref{contradiction, convergence of hat w_a to tangent map}). The proof is completed.
\end{proof}

\section{Half--degree ring disclinations}

Recall the tensor field  matrix $\mathscr{Q}$ in (\ref{ansatz}) and let $\mathscr{Q}_{a, b}^+$ denote $a^{-1}\mathscr{Q}\big(Rx\big)$ with $v(y) = u_{a, b}^+\big(R^{-1}y\big)$. As the discussions in Section 1.4.3, there is a $\rho_* > 0$ so that $\mathscr{Q}_{a, b}^+$ is negative uniaxial on the circle $\mathscr{C}_* := \big\{ \big(x_1, x_2, 0 \big) : x_1^2 + x_2^2 = \rho_*^2\big\}$.  In addition, there is  an $\epsilon_* > 0$  suitably small  so that   
\begin{align}
\begin{aligned}
&u^+_{a, b; 1} -\sqrt{3} \h{1pt} u^+_{a,b; 2}  < 0  \h{10pt} \text{on } \Big\{\big(\rho,0\big) \in T \h{1pt}:\h{1pt} \rho \in \big[\h{1pt}\rho_* - \epsilon_* , \rho_*\h{1pt}\big)\Big\},
 \\[2mm]
&u^+_{a, b; 1} -\sqrt{3} \h{1pt} u^+_{a, b; 2}  > 0 \h{10pt} \text{on } \Big\{\big(\rho,0\big) \in T \h{1pt}:\h{1pt} \rho \in \big(\h{1pt}\rho_*, \rho_* + \epsilon_*\h{1pt}\big]\Big\}.
\end{aligned}
\label{+ve and -ve of u_a,1 -sqrt3 u_a,2 on T}
\end{align}Here the constants $\rho_*$ and $\epsilon_*$ may depend on $a$, $b$ and $\mu$. In the next section, we firstly consider the biaxial structure of $\mathscr{Q}_{a, b}^+$ on $\mathscr{T}_{a, \epsilon_*}\setminus \mathscr{C}_*$. Here $\mathscr{T}_{a, \epsilon_*}$ is the torus $\big\{ x \in \mathbb{R}^3: \mathrm{dist}\big(x, \mathscr{C}_* \big)\leq \epsilon_*\big\}$.

\subsection{Biaxial structure}
Let $\lambda_{a, b; j}^+$ be the three eigenvalues in (\ref{evalue of D with u= u^+_a,b}) computed in terms of $u_{a, b}^+$. It then follows
\begin{align*}\lambda_{a, b; 2}^+ - \lambda_{a, b; 1}^+ = \dfrac{3}{4} \left(u^+_{a, b; 1}+\dfrac{1}{\sqrt{3}} \h{1pt} u^+_{a, b; 2}\right) - \dfrac{1}{4}\sqrt{\big(u^+_{a, b; 1}-\sqrt{3} \h{1pt} u^+_{a, b; 2}\big)^2+4\big(u^+_{a, b; 3}\big)^2}. \end{align*}
Thus, $\lambda_{a, b; 2}^+ > \lambda_{a, b; 1}^+$ if and only if
\begin{align}
2\h{0.5pt}u^+_{a, b; 1}\left(u^+_{a, b; 1}+\sqrt{3} \h{1pt} u^+_{a, b; 2}\right) > \big(u^+_{a, b; 3}\big)^2.
\label{biaxial, lambda_2>lambda_1}
\end{align}
 By the sequential uniform convergence of $u_{a, b}^+$ on $T$ as $a \to \infty$,  the point $x_* := \big(\rho_*, 0\big)$ must be strictly away from the origin and $\p B_1$ for $a$ sufficiently large. Furthermore, Lemma \ref{energy density e^L_a[w_a] bdd} tells us that $u^+_{a, b}$ is also equi--continuous on $\big\{\big(\rho, z\big) : \rho \in \big[\h{1pt}\delta_0, 1 - \delta_0\h{1pt}\big] \h{3pt}\text{and}\h{3pt} |\h{0.5pt} z \h{0.5pt}| \leq \delta_0\big\}$ for all $\delta_0 \in (0, 1/2)$. Hence, for a fixed $\epsilon_1 \in \big(0,1/9\big)$, we can choose $\delta_1$ small enough and $a_0$ sufficiently large such that
\begin{align*}\big| u^+_{a, b; 1} -\sqrt{3} \h{1pt} u^+_{a, b; 2}  \big| + \big| u^+_{a, b; 3}  \big| < \epsilon_1
\h{10pt} \text{ and } \h{10pt}
\big|\h{0.5pt} u^+_{a, b} \h{0.5pt}\big| > 2/3
\h{15pt}
\text{on $D_{\delta_1} (x_*)$  for any $a > a_0$}.\end{align*}
 From these inequalities, we compute that
\begin{eqnarray*}
2u^+_{a, b; 1} \big(u^+_{a, b; 1}+\sqrt{3} \h{1pt} u^+_{a,b; 2}\big)-\big(u^+_{a, b; 3}\big)^2
&>& 2u^+_{a, b; 1}\big(2u^+_{a, b; 1}-\epsilon_1\big) - \epsilon_1^2 \\[1.5mm]
&=& 4\big|u^+_{a, b}\big|^2 - 4\big(u^+_{a, b; 2}\big)^2 - 4\big(u^+_{a, b; 3}\big)^2 - 2 \epsilon_1 u^+_{a, b; 1} - \epsilon_1^2 \\[1.5mm]
&>&\dfrac{16}{9} -\dfrac{4}{3}\big(u^+_{a, c; 1}+\epsilon_1\big)^2 - 2 \epsilon_1u^+_{a, b; 1}  - 5 \epsilon_1^2  \h{2pt}>\h{2pt} 0 
\end{eqnarray*}
on $D_{\delta_1}(x_*)$ for any $a>a_0$. Therefore, $\lambda^+_{a, b; 2} > \lambda^+_{a, b; 1}$ on $ D_{\delta_1}(x_*)$ for any $a>a_0$.\vspace{.2pc}

It can also be computed that
\begin{align*}\lambda^+_{a, b; 3} - \lambda^+_{a, b; 2} =  \dfrac{1}{2}\sqrt{\big(u^+_{a, b; 1}-\sqrt{3} \h{1pt} u^+_{a, b; 2}\big)^2+4\big(u^+_{a, b; 3}\big)^2}.\end{align*} Then $\lambda^+_{a, b; 3} > \lambda^+_{a, b; 2}$ if and only if
\begin{align}
 u^+_{a, b; 1}-\sqrt{3} \h{1pt} u^+_{a, b; 2} \neq 0 \quad \text{or }\quad  u^+_{a, b; 3} \neq 0.
\label{biaxial, lambda_3>lambda_2}
\end{align}In light of (1) in Remark \ref{positivity and gamma conv}, we have 
\begin{align}
u^+_{a, b; 3} > 0 \h{15pt} \text{in } \mathbb{D}^+ \h{10pt} \text{and} \h{10pt} u^+_{a, b; 3} < 0 \h{15pt} \text{in } \mathbb{D}^-.
\label{u_a,3 > 0 and u_a,3 < 0 }
\end{align} With the above arguments, we see that $\mathscr{Q}_{a, b}^+$ is biaxial with $\lambda^+_{a, b; 3} > \lambda^+_{a, b; 2} > \lambda^+_{a, b; 1}$ on $D_{\delta_1} (x_*) \setminus T$ for $a > a_0$. Combined this consequence with (\ref{+ve and -ve of u_a,1 -sqrt3 u_a,2 on T}), $\mathscr{Q}_{a, b}^+$ is biaxial on $\mathscr{T}_{a, \epsilon_*} \setminus \mathscr{C}_*$, provided that $a$ is sufficiently large and $\epsilon_*$ is small. Here $\epsilon_*$ depends on  $a$ and $b$.

\subsection{Variation of the director field near disclination ring}
Now we discuss  the topology of the director field near the  ring disclination of $\mathscr{Q}_{a, b}^+$. Note that $\lambda^+_{a, b; 3}$ is the largest eigenvalue in $\mathscr{T}_{a, \epsilon_*}\setminus \mathscr{C}_*$. The director field, i.\h{0.5pt}e. the normalized eigenvector of $\mathscr{Q}_{a, b}^+$ associated with the eigenvalue $\lambda^+_{a, b; 3}$, can be oriented and represented by $\kappa\h{1pt}\big[\h{1pt}u_{a, b}^+\h{1pt}\big]$. See the definition of $\kappa\h{0.5pt}[\h{0.5pt}u\h{0.5pt}]$ from (\ref{hat kappa_3, intro}). The coefficient of  $e_z$ in $\kappa\h{1pt}\big[\h{1pt}u_{a, b}^+\h{1pt}\big]$, i.\h{0.5pt}e. $\left<\kappa\h{1pt}\big[\h{1pt}u_{a, b}^+\h{1pt}\big], e_z\right>$, can be expressed by
\begin{align*}  \mbox{\fontsize{9}{8}\selectfont\(\dfrac{\sqrt{2}\,\text{sign}\big(u^+_{a, b; 3})}{2}\Big[\big(u^+_{a, b; 1}-\sqrt{3} \h{1pt} u^+_{a, b; 2}\big)^2 +4\big(u^+_{a, b; 3}\big)^2\Big]^{-\frac{1}{4}}\left[\sqrt{\big(u^+_{a, b; 1}-\sqrt{3} \h{1pt} u^+_{a, b; 2}\big)^2 +4\big(u^+_{a, b; 3}\big)^2}-\big(u^+_{a, b; 1}-\sqrt{3} \h{1pt} u^+_{a, b; 2}\big)\right]^{\frac{1}{2}}.\)}
\end{align*}Here $\big<\cdot, \cdot\big>$ is the standard inner product in $\mathbb{R}^3$. It turns out that  $\left<\kappa\h{1pt}\big[\h{1pt}u_{a, b}^+\h{1pt}\big], e_z\right> \to - 1$ as we approach the point $x_*^{r} := \big(\rho_* - r, 0 \big)$ along $\p^- D_{r} (x_*)$. Note that  $\p^- D_{r} (x_*)$ denotes the lower--half part of $\p D_{r} (x_*)$. This  convergence follows from the fact that $u^+_{a, b; 3} < 0 $ on $\mathbb{D}^-$ and $u^+_{a, b; 1} -\sqrt{3} \h{1pt} u^+_{a, b; 2} < 0$ at $x_{*}^r$. See   (\ref{u_a,3 > 0 and u_a,3 < 0 })  and  (\ref{+ve and -ve of u_a,1 -sqrt3 u_a,2 on T})  respectively.   We then  conclude that the director field $ \kappa\h{1pt}\big[\h{1pt}u_{a, b}^+\h{1pt}\big]$ converges to $-e_z$ when we approach the point $x_{*}^r$ along $\p^- D_{r} (x_*)$. Similarly, when we approach $x^r_{*}$ along $\p^+ D_{r} (x_*)$, the upper--half part of $\p D_r(x_*)$, the director field $ \kappa\h{1pt}\big[\h{1pt}u_{a, b}^+\h{1pt}\big]$ converges to $e_z$. Here, we just need the fact that $u^+_{a, b; 3} > 0 $ on $\mathbb{D}^+$. Meanwhile, due to (\ref{hat kappa_3, intro}),  (\ref{+ve and -ve of u_a,1 -sqrt3 u_a,2 on T}) and (\ref{u_a,3 > 0 and u_a,3 < 0 }), the coefficient of $e_\rho$ in $ \kappa\h{1pt}\big[\h{1pt}u_{a, b}^+\h{1pt}\big]$ keeps strictly positive on $ \p D_{r} (x_*) \setminus \big\{x^r_{*}\big\}$. Therefore, when we start from $x^r_{*}$ and rotate counter--clockwisely along $\p D_{r} (x_*)$ back to $x^r_{*}$, the director field  $ \kappa\h{1pt}\big[\h{1pt}u_{a, b}^+\h{1pt}\big]$ varies from $-e_z$ to $e_z$ continuously. During this process,   $ \kappa\h{1pt}\big[\h{1pt}u_{a, b}^+\h{1pt}\big]$ keeps strictly on the right--half part of $(\rho,z)$--plane except at $x^r_{*}$. The angle of  $ \kappa\h{1pt}\big[\h{1pt}u_{a, b}^+\h{1pt}\big]$ is totally changed by $\pi$. This verifies that $\mathscr{Q}_{a, b}^+$ admits a half--degree ring disclination at $\mathscr{C}_*$. \vspace{.2pc}

To end this section, we compute the tangent map of the director field  $ \kappa\h{1pt}\big[\h{1pt}u_{a, b}^+\h{1pt}\big]$ at $x_*$ for large $a$. Let $\varphi'$ be an angular variable ranging from $[\h{0.5pt}-\pi,\pi\h{0.5pt}]$. If $\varphi'=0$, then $u^+_{a, b; 1}\big( \rho_* + \epsilon, 0 \big) -\sqrt{3} \h{1pt} u^+_{a,b; 2}\big( \rho_* + \epsilon, 0\big)  > 0$ by (\ref{+ve and -ve of u_a,1 -sqrt3 u_a,2 on T}). Moreover, it satisfies $u^+_{a, b; 3}\big( \rho_a + \epsilon, 0 \big) = 0$. Hence,  $ \kappa\h{1pt}\big[\h{1pt}u_{a, b}^+\h{1pt}\big]  = e_\rho$ at $\varphi'=0$  for large $a$  and $\epsilon \in (0,\epsilon_*)$. If $\varphi' \in (0,\pi)$, then L'Hospital's rule infers
\begin{align*} \lim_{\epsilon \to 0^+} \dfrac{ u^+_{a, b; 1} -\sqrt{3} \h{1pt} u^+_{a, b; 2} }{ u^+_{a, b; 3} } \Bigg|_{x_* + \epsilon(\cos \varphi', \sin \varphi')}
= \dfrac{ \big(Du^+_{a, b; 1}(x_*) - \sqrt{3} \h{1pt}Du^+_{a, b; 2}(x_*)\big) \cdot (\cos \varphi', \sin \varphi')^\top  }{ D u^+_{a, b; 3}(x_*)\cdot (\cos \varphi', \sin \varphi')^\top }.\end{align*}
Since with respect to $z$--variable, $u^+_{a, b; 1}$ and $u^+_{a, b; 2}$ are even, and $u^+_{a, b; 3}$ is odd, we obtain
\begin{align*} \lim_{\epsilon \to 0^+} \dfrac{ u^+_{a, b; 1} -\sqrt{3} \h{1pt} u^+_{a,2} }{ u^+_{a, b; 3} } \Bigg|_{x_* + \epsilon(\cos \varphi', \sin \varphi')}
= \varkappa_*\h{1pt} \mathrm{ctan} \h{1pt} \varphi', \h{10pt}\text{where $\varkappa_* := \dfrac{  \p_\rho u^+_{a, b; 1}(x_*) - \sqrt{3}\h{1pt} \p_\rho u^+_{a, b; 2}(x_*) }{ \p_z u^+_{a, b; 3}(x_*) }$.} \end{align*}
We note that $\p_z u_{a,3}(x_*)$ is positive due to Hopf's Lemma. In addition, it holds $\varkappa_*  \geq 0$. Recall (\ref{hat kappa_3, intro}). Then for any $\varphi' \in (0, \pi)$ and large $a$, the above convergence result  yields the corresponding limit in (5) of Theorem \ref{biaxial--ring solution}.  Here we also use $u_{a, b; 3}^+ > 0$ on $\mathbb{D}^+$. We can obtain similar result if $\varphi' \in (-\pi, 0)$. The item (5) in Theorem \ref{biaxial--ring solution} is obtained.  

\section{Split--core solutions  with strength--one disclinations}

Denote by $\mathscr{Q}_{a, c}^-$ the tensor field $a^{-1}\mathscr{Q}\big(Rx\big)$ with $v(y) = u_{a, c}^-\big(R^{-1}y\big)$. For large $a$, we let  $z_a^+ = (0, 0, z_a) $ be the lowest point on $l^+_z$ at which $w^-_{a,c}$ vanishes. Near $z_a^+$, we use  $(r_*, \psi, \theta)$ to denote  the spherical coordinate system with respect to the center $z_a^+$.  Here $r_*$ is the radial variable, $\psi$ is the polar angle, while $\theta$ is still the azimuthal angle. Using  (2) in Proposition \ref{num and asy of singu}, for $\epsilon_\star > 0$, there are $a_0>0$ and $r_0>0$ so that 
\begin{align}
\sum^2_{j=0}\left\|\h{1pt} \p^{\h{1pt}j}_\psi  \h{1.5pt}\Pi_{\mathbb{S}^2} \big[ u_{a, c}^- \big]  - \p^{\h{1pt}j}_\psi\big(0,\cos \psi, \sin \psi \big) \h{1pt}\right\|_{\infty\h{.5pt};\h{1pt}\{\sigma\} \times[\h{0.5pt}0,\pi \h{0.5pt}]} < \epsilon_\star \h{15pt}\text{for  $a > a_0$ and  $\sigma \in \big(0,r_0\big)$.}
\label{C2 conv of u_a in psi}
\end{align}Without ambiguity, we still use $u_{a,c}^-$ in (\ref{C2 conv of u_a in psi}) to represent the mapping $u_{a, c}^- \big(r_* \sin \psi, z_a + r_* \cos \psi\big)$. It  depends on the variables $(r_*, \psi)$ and is the expression of $u_{a, c}^-$ under the spherical coordinates $(r_*, \psi, \theta)$. In the next, we firstly study the structures of $\mathscr{Q}_{a, c}^-$ on $l_z$.
\subsection{Uniaxial and isotropic structures on $l_z$}

 Let $\lambda^-_{a, c; j}$ ($j = 1, 2, 3$) be the three eigenvalues of $\mathscr{Q}_{a, c}^-$. They are the three eigenvalues in (\ref{evalue of D with u= u^+_a,b}) computed in terms of $u_{a, c}^-$.  Recall that $u^-_{a, c; 1} = u^-_{a, c; 3}=0$ on $B_1 \cap l_z$. Then by (\ref{evalue of D with u= u^+_a,b}), $\mathscr{Q}_{a, c}^-$ is uniaxial or isotropic on $B_1 \cap l_z$. More precisely, $\mathscr{Q}_{a, c}^-$ is isotropic at the points on $B_1 \cap l_z$ with $u_{a, c; 2}^- = 0$. For the points on $B_1 \cap l_z$ where $u^-_{a, c; 2}$ is positive, $\mathscr{Q}_{a, c}^-$ is positive uniaxial in the sense that  $\lambda^-_{a, c; 2} = \lambda^-_{a, c; 1} < \lambda^-_{a, c; 3}$. The eigenspace of the largest eigenvalue of $\mathscr{Q}_{a, c}^-$ is given by $\mathrm{span}\big\{ e_z\big\}$ at these positive uniaxial locations. For the points on  $B_1 \cap l_z$ where $u^-_{a, c; 2}$ is negative, $\mathscr{Q}_{a, c}^-$ is negative uniaxial in the sense that  $\lambda^-_{a, c; 2} < \lambda^-_{a, c; 1} = \lambda^-_{a, c; 3}$. The eigenspace of the largest eigenvalue of $\mathscr{Q}_{a, c}^-$ is given by $\mathrm{span}\big\{ e_\rho, e_{\theta}\big\}$ at these negative uniaxial locations. Here $e_\theta := \left(-  \frac{x_2}{\rho}, \frac{x_1}{\rho}, 0 \right)^\top$. 

\subsection{Biaxial structure}
In this section, we consider the biaxial structure of $\mathscr{Q}_{a, c}^-$ in the dumbbell $D_{r_0, r_1}\big(z_a^+, z_a^-\big)$. See the definition of dumbbell in Definition \ref{dumbbell}. The dumbbell size parameter $r_0$ is as in (\ref{C2 conv of u_a in psi}).  $r_1$ is a positive number less than $r_0\big/2$. We first compare
the three eigenvalues in $D_{r_0}\big(z_a^+\big)$. Due to the $\mathscr{R}$--axial symmetry of $w^-_{a, c}$, the case for $D_{r_0}(z_a^-)$ can be similarly studied. Suppose that $\sigma$ is an arbitrary number in
$(0,r_0)$. Using the polar angle $\psi$ in the spherical coordinates $(r_*,\psi,\theta)$ with respect to the center $z^+_a$, we have  
\begin{align*}\Pi_{\mathbb{S}^2}\big[ u_{a, c}^- \big] (\sigma, \psi)
= \Pi_{\mathbb{S}^2}\big[ u_{a, c}^- \big] (\sigma, 0)
+  \left(\p_\psi \h{1pt}\Pi_{\mathbb{S}^2}\big[ u_{a, c}^- \big] \h{1pt} \Big|_{(\sigma,0)} \right) \psi
+ \int_0^\psi \int_0^{\psi_1} \p_\psi^{\h{.5pt}2} \h{1pt} \Pi_{\mathbb{S}^2}\big[ u_{a, c}^- \big] \h{1pt}  \Big|_{(\sigma, \zeta)} \mathrm{d} \h{1pt} \zeta \h{1pt} \mathrm{d} \h{1pt} \psi_1. \end{align*}By the regularity of $w^-_{a, c}$ on $l_z$, it follows that $\p_\psi u^-_{a, c}\h{1pt} \Big|_{(\sigma,0)}=0$. The above equality then infers
\begin{align*}
\sqrt{3} \h{1pt} \left[\Pi_{\mathbb{S}^2}\big[ u_{a, c}^- \big]\right]_1 \Big|_{ (\sigma, \psi)} &+ \left[ \Pi_{\mathbb{S}^2}\big[ u_{a, c}^- \big]\right]_2 \Big|_{ (\sigma, \psi)}
\h{2pt}=\h{2pt} \left[ \Pi_{\mathbb{S}^2}\big[ u_{a, c}^- \big]\right]_2 \Big|_{(\sigma, 0)} - (1 - \cos \psi)\\[2mm]
&\h{2pt}+\h{2pt} \int_0^\psi \int_0^{\psi_1} \sqrt{3} \h{1.5pt} \p_\psi^{\h{.5pt}2} \h{1pt} \left[ \Pi_{\mathbb{S}^2}\big[ u_{a, c}^- \big]\right]_1 \Big|_{(\sigma, \zeta)} + \p_\psi^{\h{.5pt}2} \left(\left[ \Pi_{\mathbb{S}^2}\big[ u_{a, c}^- \big]\right]_2  -\cos \psi\right)\Big|_{(\sigma, \zeta)} \mathrm{d} \h{1pt} \zeta \h{1pt} \mathrm{d} \h{1pt} \psi_1.
\end{align*}
According to the estimate in (\ref{C2 conv of u_a in psi}), for any $a > a_0$ and $\sigma \in (0,r_0)$, we have from the last equality that 
\begin{align}
\sqrt{3} \h{1pt} \left[\Pi_{\mathbb{S}^2}\big[ u_{a, c}^- \big]\right]_1 \Big|_{ (\sigma, \psi)} &+ \left[ \Pi_{\mathbb{S}^2}\big[ u_{a, c}^- \big]\right]_2 \Big|_{ (\sigma, \psi)}
\h{2pt}\leq\h{2pt} \left[ \Pi_{\mathbb{S}^2}\big[ u_{a, c}^- \big]\right]_2 \Big|_{(\sigma, 0)} - (1 - \cos \psi)+  \epsilon_\star \h{0.2pt} \psi^2.
\label{ineq. of sqrt3 u_1+u_2 u=u^-_a,c}
\end{align}
Note that $ \left[ \Pi_{\mathbb{S}^2}\big[ u_{a, c}^- \big]\right]_2 \Big|_{(\sigma, 0)} =1$. Referring to (\ref{ineq. of sqrt3 u_1+u_2 u=u^-_a,c}), we can find an $\epsilon_\star$ small enough  such that
\begin{align}
\sqrt{3} \h{1pt} \left[\Pi_{\mathbb{S}^2}\big[ u_{a, c}^- \big]\right]_1 \Big|_{ (\sigma, \psi)} &+ \left[ \Pi_{\mathbb{S}^2}\big[ u_{a, c}^- \big]\right]_2 \Big|_{ (\sigma, \psi)}
\h{2pt}\leq\h{2pt} \cos \psi +  \epsilon_\star \h{.2pt} \psi^2  \h{2pt}<\h{2pt} 1 \h{15pt} \text{for any $\psi \in \left(0, \frac{\pi}{4}\h{1pt}\right)$.}
\label{sqrt3 u_1 + u_2 < 1- small term}
\end{align}
Moreover, we can keep taking $\epsilon_\star$ small and infer from (\ref{C2 conv of u_a in psi}) that
\begin{align}\sqrt{3} \h{1pt} \left[\Pi_{\mathbb{S}^2}\big[ u_{a, c}^- \big]\right]_1 \Big|_{(\sigma, \psi)} + \left[ \Pi_{\mathbb{S}^2}\big[ u_{a, c}^- \big]\right]_2 \Big|_{ (\sigma, \psi)}&\h{2pt}=\h{2pt} \cos \psi + \sqrt{3} \h{1pt}  \left[\Pi_{\mathbb{S}^2}\big[ u_{a, c}^- \big]\right]_1  \Big|_{(\sigma,\psi)} +  \left[ \Pi_{\mathbb{S}^2}\big[ u_{a, c}^- \big]\right]_2 \Big|_{(\sigma,\psi)}-\cos \psi  \nonumber\\[2mm]
&\h{2pt}\leq\h{2pt} \cos \frac{\pi}{4} + 2 \h{1pt} \epsilon_\star \h{2pt}<\h{2pt} 1 \h{25pt}\text{for any $\psi \in \left[\h{1pt}\frac{\pi}{4}, \pi \h{1pt}\right]$.}\label{large angule}\end{align}
Combining  (\ref{sqrt3 u_1 + u_2 < 1- small term}) and (\ref{large angule}), we obtain\begin{align*}\sqrt{3} \h{1pt} u^-_{a,c; 1} \Big|_{(\sigma,\psi)} + u^-_{a, c; 2} \Big|_{(\sigma,\psi)} \h{2pt}< \h{2pt} \big|\h{0.5pt} u^-_{a, c}\h{0.5pt}\big|\h{1.5pt} \Big|_{(\sigma,\psi) } \h{15pt}\text{for any $a > a_0$, $\sigma \in (0, r_0)$ and $\psi \in (0, \pi)$.} \end{align*}
If $\psi \in (0,\pi)$, then $u^-_{a, c; 1}$ is strictly positive. It turns out
\begin{align*}\sqrt{3} \h{1pt} u^-_{a, c; 1} \Big|_{(\sigma,\psi)} + u^-_{a, c; 2} \Big|_{(\sigma,\psi)}
\h{2pt}>\h{2pt} u^-_{a, c; 2} \Big|_{(\sigma,\psi)} \h{2pt}>\h{2pt} -\big|\h{0.5pt}u^-_{a, c}\h{0.5pt}\big|\h{1pt}\Big|_{(\sigma,\psi)} \h{15pt} \text{for any $\psi \in (0,\pi)$}.\end{align*}
The last two inequalities yield
\begin{align*}\left(\sqrt{3} \h{1pt} u^-_{a, c; 1} \Big|_{(\sigma,\psi)} + u^-_{a, c; 2} \Big|_{(\sigma,\psi)} \right)^2
\h{2pt}<\h{2pt} \big|\h{0.5pt} u^-_{a, c}\h{0.5pt}\big|^2\h{1pt}\Big|_{ (\sigma,\psi) } \h{15pt} \text{for any $a > a_0$, $\sigma \in (0,r_0)$ and $\psi \in (0,\pi)$},\end{align*}which furthermore induces
\begin{align}\label{comparison of components of u^-_a,c}
 \left(u^-_{a, c; 1}-\sqrt{3} \h{1pt} u^-_{a, c; 2}\right)^2 + 4 \left(u^-_{a, c; 3}\right)^2
&\h{2pt}=\h{2pt}4\h{1pt}\big|\h{0.5pt}u^-_{a, c} \h{0.5pt}\big|^2 - \left(\sqrt{3}\h{1pt}u^-_{a, c; 1} + u^-_{a, c; 2}\right)^2 \nonumber\\[1.5mm]
&\h{2pt}>\h{2pt} 3 \left(\sqrt{3} \h{1pt} u^-_{a, c; 1} + u^-_{a, c; 2}\right)^2 = 9\left( u^-_{a, c; 1} +\dfrac{1}{\sqrt{3}} \h{1pt} u^-_{a, c; 2} \right)^2.
\end{align}
Note that  (\ref{comparison of components of u^-_a,c}) is evaluated at $(\sigma, \psi)$ and holds for any $a > a_0$, $\sigma \in (0,r_0)$ and $\psi \in (0,\pi)$. By direct computations and (\ref{evalue of D with u= u^+_a,b}), it satisfies
\begin{align}\label{com of three eigenvalues negati}\lambda^-_{a, c; 2} - \lambda^-_{a, c; 1} = \dfrac{3}{4} \h{1pt}\left( u^-_{a, c; 1}+\dfrac{1}{\sqrt{3}} \h{1pt} u^-_{a, c; 2} \right) - \frac{1}{4}\h{1pt}\sqrt{\big(u^-_{a, c; 1}-\sqrt{3} \h{1pt} u^-_{a, c; 2}\big)^2+4\big(u^-_{a,c; 3}\big)^2},\nonumber\\[1.5mm]
 \lambda^-_{a, c; 3}-\lambda^-_{a, c; 1} = \dfrac{3}{4} \h{1pt}\left( u^-_{a, c; 1}+\dfrac{1}{\sqrt{3}} \h{1pt} u^-_{a, c; 2} \right) + \frac{1}{4}\h{1pt}\sqrt{\big(u^-_{a, c; 1}-\sqrt{3} \h{1pt} u^-_{a, c; 2}\big)^2+4\big(u^-_{a, c; 3}\big)^2}.\end{align}
Therefore, (\ref{comparison of components of u^-_a,c}) and the $\mathscr{R}$--axial symmetry of $w^-_{a, c}$ induce 
\begin{equation}
\lambda^-_{a, c; 3} > \lambda^-_{a, c; 1} > \lambda^-_{a, c; 2} \h{15pt} \text{on} \h{5pt} \left[ D_{r_0}\big(z^+_a \big) \h{1.5pt}\bigcup\h{1.5pt} D_{r_0}\big(z^-_a\big)  \right] \setminus l_z.
\label{bixial in two ball, u^-_a,c}
\end{equation}

Denote by $R^*$ the rectangle in $(x_1,z)$--plane with four vertices $x^+_1$, $x^-_1$, $x^+_2$ and $x^-_2$. Here \begin{align*}x_1^\pm := \left( r_1^{1/2}\h{0.5pt}\sqrt{2 r_0 - r_1}, \h{1pt}\pm \big(z_a - r_0 + r_1\big)\right), \h{15pt}x_2^\pm := \left(- r_1^{1/2}\h{0.5pt}\sqrt{2 r_0 - r_1}, \h{1pt}\pm \big(z_a - r_0 + r_1\big)\right).
\end{align*} We are left to compare the three eigenvalues on $R^*$. Fix $r_0$. Then we take $r_1$  sufficiently small and $a_0$ sufficiently large so that $u^-_{a, c}$ is close to $(0,-1,0)^\top$ uniformly on $R^*$ for any $a>a_0$.  By the first equality in (\ref{com of three eigenvalues negati}), it follows $\lambda^-_{a, c; 2} < \lambda^-_{a, c; 1}$ on $R^*$, provided that $r_1$ is small and  $a$ is large. To compare the eigenvalues $\lambda^-_{a, c; 1}$ and $\lambda^-_{a, c; 3}$, we first notice that  $u^-_{a, c; 1} $ is strictly positive on $R^* \setminus l_z$. Therefore,  
\begin{align*} \sqrt{3}\h{0.5pt}u^-_{a, c; 1}+u^-_{a, c; 2} \h{1.5pt}>\h{1.5pt} u^-_{a, c; 2} \h{1.5pt}>\h{1.5pt} - \big|\h{0.5pt}u^-_{a, c}\h{0.5pt}\big| \h{15pt} \text{on} \h{5pt} R^* \setminus l_z.\end{align*}
In addition, for any large $a$ and small $r_1$, the inequality $\sqrt{3}\h{0.5pt}u^-_{a, c; 1}+u^-_{a, c; 2} < \big|\h{0.5pt}u^-_{a, c}\h{0.5pt}\big|$ holds on $R^*$ in that $u^-_{a, c}$ is sufficiently close to $(0,-1,0)^\top$ on $R^*$ if $r_1$ is small  and $a$ is large. It turns out that $$\left( \sqrt{3}\h{0.5pt}u^-_{a, c; 1}+u^-_{a, c; 2}\right)^2 \h{1.5pt}<\h{1.5pt} \big|\h{0.5pt}u^-_{a, c}\h{0.5pt}\big|^2 \h{15pt}\text{on $R^* \setminus  l_z$.}$$Here and in what follows, we still take $a$ large and $r_1$ small. We then obtain similarly from (\ref{comparison of components of u^-_a,c}) that  
$$\big(u^-_{a, c; 1}-\sqrt{3} \h{1pt} u^-_{a, c; 2}\big)^2 + 4 \big(u^-_{a, c; 3}\big)^2\h{0.5pt}>\h{0.5pt}9\left( u^-_{a, c; 1} +\dfrac{1}{\sqrt{3}} \h{1pt} u^-_{a, c; 2} \right)^2 \h{15pt} \text{on} \h{5pt} R^* \setminus l_z.$$ It implies by the second equality in (\ref{com of three eigenvalues negati}) that $\lambda^-_{a, c; 3} > \lambda^-_{a, c; 1}$ on $R^* \setminus l_z$. Together with (\ref{bixial in two ball, u^-_a,c}), it follows
$$\lambda^-_{a, c; 3} > \lambda^-_{a, c; 1} > \lambda^-_{a, c; 2} \h{15pt} \text{on} \h{5pt} D_{r_0, r_1}\big(z^+_a, z^-_a\big) \setminus l_z.$$ Here $r_0$ is small and $a$ is large. Meanwhile,  $r_1 < r_0\big/2$ is also small.

\subsection{Variation of the director field along the contour $\mathscr{C}_{r_0, r_1}\big(z_a^+, z_a^-\big)$}
Note that the largest eigenvalue of $\mathscr{Q}_{a, c}^-$  is $\lambda^-_{a, c; 3}$ when $\mathscr{Q}_{a, c}^-$ is restricted on  $D_{r_0, r_1}\big(z^+_a, z_a^-\big) \setminus l_z$. The director field of $\mathscr{Q}_{a, c}^-$ on $D_{r_0, r_1}\big(z^+_a, z_a^-\big) \setminus l_z$ then equals  $\kappa\h{0.5pt}\big[\h{0.5pt}u_{a, c}^-\h{0.5pt}\big]$, where $\kappa\h{0.5pt}[u]$ is defined in (\ref{hat kappa_3, intro}).   Let $z \in \big(z_a , z_a  + r_0)$. Then $u^-_{a, c}$ converges to  $\big(0,\h{0.5pt}\big|\h{0.5pt}u^-_{a, c}(0,0,z)\h{0.5pt}\big|, \h{0.5pt}0\big)^\top$ when we approach $(0,0,z)$. As a consequence,  the coefficient of $e_\rho$ in $\kappa\h{0.5pt}\big[\h{0.5pt}u_{a, c}^-\h{0.5pt}\big]$ tends to $0$ as $x \to (0,0,z)$.  Since $u^-_{a, c; 3}$ is positive when $x$ is close to $(0,0,z)$ and not on $l_z$, we then conclude that $\kappa\h{0.5pt}\big[ \h{0.5pt}u_{a, c}^-\h{0.5pt}\big]$ converges to $e_z$  as $x$ approaches $(0,0,z)$ for any $z \in (z_a , z_a  + r_0)$. Moreover, by the $\mathscr{R}$--axial symmetry, it follows that $\kappa\h{0.5pt}\big[\h{0.5pt}u_{a, c}^-\h{0.5pt}\big]$ converges to $ -e_z$  as $x$ approaches $(0,0,z)$ for any $z \in (-z_a - r_0, - z_a)$.\vspace{0.2pc}

Note that for sufficiently large $a$, the mapping $u^-_{a, c}$ is close to $(0,-1,0)^\top$ when $x$ is close to the origin. Recall that $u^-_{a, c; 3}=0$ on $T$. These results consequently yield 
$$\kappa\h{0.5pt}\big[\h{0.5pt}u_{a, c}^- \h{0.5pt}\big] \equiv e_\rho \h{15pt} \text{on $D_{r_0, r_1  } \big(z^+_a, z_a^-\big) \h{1.5pt}\cap\h{1.5pt}\Big\{ \text{$x_1$--axis} \Big\}$,  \h{5pt} provided that $r_1$ is small enough.}$$ Moreover, we notice that for any point in $D_{r_0, r_1}\big(z_a^+, z_a^-\big) \setminus l_z$, the coefficient of $e_\rho$ in $\kappa\h{0.5pt}\big[\h{0.5pt}u_{a, c}^-\h{0.5pt}\big]$ keeps strictly positive. The item (4.1) in Theorem \ref{split--core solution} then follows.

\begin{alphasection}

\addcontentsline{toc}{section}{\normalsize Appendix}

\section*{Appendix}
\addcontentsline{toc}{subsection}{\normalsize A.1 \h{5pt}Proof of (1) in Lemma \ref{decay for large L}}
\setcounter{equation}{0}
\setcounter{thm}{0}
\subsection*{A.1 \h{5pt}Proof of (1) in Lemma \ref{decay for large L}}
In this section, we prove (\ref{two possible cases}) in Lemma \ref{decay for large L}. Supposing on the contrary that item (1) in Lemma \ref{decay for large L} fails, we can find $b_* > 0$, $\nu_n \rightarrow 1^-$, $L_n \rightarrow \infty$, $\big(h_n, \overline{w}_{*, n}\big) \in \mathbb{R}^2$ and a solution, denoted by $W_n$, to the Problem $S_{L_n, h_n, \overline{w}_{*, n}}$ on $D_{\sigma_k}$ so that the followings hold: \begin{align}&\mathrm{(i).} \h{2pt}E_{L_n, h_n}^\star\left[W_n\right] \leq 1; \h{25pt}\mathrm{(ii).}\h{2pt} \big\| W_n \big\|_{1, 2; \h{1pt}D_{\sigma_k}} \leq b_*;  \nonumber\\[2mm]&\mathrm{(iii).} \h{2pt} \int_{D_{1/8}}
 e_{L_n, h_n}^\star \left[ W_n\right] \h{2pt} > \nu_n  \int_{D_{1/4}}
 e_{L_n, h_n}^\star \left[ W_n \right] \h{2pt}>\h{2pt} \dfrac{\nu_n}{16}, \h{15pt}\text{for all $n \in \mathbb{N}$.}
\label{contradiction assumption}
\end{align} By (i) and (ii) in (\ref{contradiction assumption}), the sequence $\big\{h_n\big\}$ is  uniformly bounded.  There are $h_\infty \in \mathbb{R}$ and $W_\infty \in H^1\big(D_{\sigma_k}; \mathbb{R}^2\big)$ so that  $h_n \longrightarrow h_\infty$ and $W_n \longrightarrow W_\infty$ weakly in $H^1\big(D_{\sigma_k}; \mathbb{R}^2\big)$ as $n \rightarrow \infty$. By Sobolev embedding, we can also assume  $W_n \longrightarrow W_\infty$ strongly in $L^2\big(D_{\sigma_k}; \mathbb{R}^2\big)$ as $n \to \infty$. This convergence, the fact that  $L_n \rightarrow \infty$ and  (i) in (\ref{contradiction assumption}) then yield  \begin{eqnarray}\label{almost every Winfty}
y_* \cdot W_\infty = - h_\infty \h{20pt}\text{a.e. in $D_{\sigma_k}$.}\end{eqnarray}
In the above, we still use $y_*$ to denote the $2$--vector $\left(\sqrt{1 - b^2}, b \right)^\top$. On the other hand, by trace theorem, we can also assume $W_n \longrightarrow W_\infty$ strongly in $L^2\big(T_{\sigma_k}; \mathbb{R}^2\big)$ as $n \to \infty$. Here $T_{\sigma_k} := T \h{0.8pt}\cap\h{0.8pt} D_{\sigma_k}$. Therefore up to a subsequence, the Signorini lower bound $\overline{w}_{*, n}$ either converges to some finite number $\overline{w}_{*,\infty}$ or diverges to  $- \infty$ as $n \rightarrow \infty$. If $\overline{w}_{*, n} \rightarrow \overline{w}_{*, \infty}$, then we have $W_{\infty; 2} \geq \overline{w}_{*, \infty}$ on $T_{\sigma_k}$ in the sense of trace.  One can now apply Fatou's lemma to find a positive constant $B_k$ and a radius $r_k \in \big(1/2, \sigma_k\big)$ so that up to a subsequence, it holds \begin{eqnarray} \sup_{n \h{1pt}\in \h{1pt} \mathbb{N}\h{1pt}\cup\h{1pt} \{\infty\}}  \h{2pt}\big\lVert  W_n \big\rVert_{\infty;\h{1pt}\p D_{r_k}} +  \int_{\p D_{r_k}} \big|\h{1pt}D_\xi W_\infty\h{1pt}\big|^2 +  \sup_{n \h{1pt}\in \h{1pt} \mathbb{N}} \h{2pt}\int_{\p D_{r_k}} e_{L_n, h_n}^\star\left[W_n\right]  \h{2pt}\leq\h{2pt}B_k.
\label{bounds for W_n}
\end{eqnarray}
Moreover by (\ref{almost every Winfty}) and $W_{n; 2} \geq \overline{w}_{*, n}$ on $T_{\sigma_k}$, we can also assume $y_* \cdot W_\infty = - h_\infty$ and $W_{n; 2} \geq \overline{w}_{*, n}$ at $\big(\pm r_k, 0\big)$ for any $n \in \mathbb{N}$.  If $\overline{w}_{*, n} \rightarrow \overline{w}_{*, \infty}$ as $n \rightarrow \infty$, we can in addition assume $W_{\infty; 2} \geq \overline{w}_{*,\infty}$ at $\big(\pm r_k, 0\big)$.\vspace{0.2pc}

We now construct comparison mappings.  Denote by  $t_*$ the vector $(- y_{*; 2}, y_{*; 1})^\top$, where $y_{*; 1} = \sqrt{1 - b^2}$ and $y_{*; 2} = b$. Let $W$ be any vector field satisfying \begin{eqnarray} \label{rep of W}W = w \h{1pt}t_*  - h_\infty \h{1pt}y_* \h{10pt}\text{in $D_{r_k}$} \h{20pt}\text{and}\h{20pt}W = W_\infty \h{10pt}\text{on $\p D_{r_k}$.}
\end{eqnarray}
Here $w$ is a scalar function in $H^1\big(D_{r_k}\big)$. It is even with respect to the variable $\xi_2$. If $\overline{w}_{*, n} \rightarrow \overline{w}_{*, \infty}$, then we also assume  \begin{eqnarray}w \h{1pt}  \h{2pt}\geq \h{2pt}Z_* := \dfrac{\overline{w}_{*, \infty}}{y_{*; 1}} + h_{\infty}\h{1pt}\dfrac{y_{*; 2}}{y_{*; 1}} \h{20pt}\text{on $T_{r_k}$.}
\label{Z_*}\end{eqnarray} Associated with $W$, we define for any $R > 0$ the mapping $K_{n, R}[\h{0.5pt}W\h{0.5pt}]$ as follows: \begin{eqnarray}K_{n, R}[\h{0.5pt}W\h{0.5pt}] := \left\{\begin{array}{lcl}
- h_n \h{1pt} y_* + \left[ \h{1pt} h_n \h{1pt}\dfrac{y_{*; 2}}{y_{*; 1}} + R  \h{1pt}\dfrac{w - W_*}{| w - W_* | \vee R} \right] t_*, \h{20pt}&&\text{if $\overline{w}_{*, n} \rightarrow -\infty$;}\\[6mm]
- h_n \h{1pt} y_* + \left[ \h{1pt}\dfrac{\overline{w}_{*,n}}{y_{*; 1}} +  h_n \h{1pt}\dfrac{y_{*; 2}}{y_{*; 1}} + R  \h{1pt}\dfrac{w - W_*}{| w - W_*| \vee R} \right] t_*, \h{20pt}&&\text{if $\overline{w}_{*,n} \rightarrow \overline{w}_{*, \infty}$.} \end{array}\right.
\label{K_n}
\end{eqnarray}In the case where $\overline{w}_{*, n } \rightarrow - \infty$, we let $W_* = h_{\infty}\h{1pt} \dfrac{y_{*; 2}}{y_{*; 1}} $. If $\overline{w}_{*, n} \rightarrow \overline{w}_{*, \infty}$, then  $W_* := Z_*$. 
With $K_{n, R}[\h{2pt}W\h{2pt}]$, our comparison map is defined by \begin{eqnarray} \overline{W}_{n, s, R}\big(\xi\big) := \left\{\begin{array}{lcl} K_{n, R}[\h{0.5pt}W\h{0.5pt}] \left(\dfrac{\xi}{1 - s}\right) &&\text{if $\xi \in D_{(1 - s) \h{1pt}r_k} $;} \\[5mm]
\dfrac{r_k - |\h{1pt} \xi \h{1pt}|}{s \h{1pt}r_k} \h{2pt}K_{n, R}[\h{1pt}W_\infty\h{1pt}] \h{2pt}\bigg|_{r_k \widehat{\xi}}  + \dfrac{|\h{1pt} \xi \h{1pt}| - (1 - s) \h{1pt}r_k}{s \h{1pt}r_k} \h{2pt} W_n \h{2pt} \bigg|_{r_k \widehat{\xi}} &&\text{if $\xi \in D_{r_k} \setminus D_{(1 - s) r_k}$.} \end{array}\right.
\label{bar W_n}
\end{eqnarray}
Here $s \in (0,1)$.  It can be checked from (\ref{bar W_n}) that $\overline{W}_{n, s, R} = W_n$ on $\p D_{r_k}$. If $\overline{w}_{*, n} \rightarrow - \infty$, then by (\ref{K_n}),  $\left[\h{1pt} K_{n, R}\left[W\right] \h{1pt}\right]_2$  is uniformly bounded from below by $- R$. If we take $n$ large enough depending on $R$, it then follows $\left[\h{1pt} K_{n, R}\left[ W\right] \h{1pt}\right]_2 \geq \overline{w}_{*, n}$ for any $2$--vector field $W$ satisfying the first condition in (\ref{rep of W}). This lower bound together with $W_{n; 2} \geq \overline{w}_{*, n}$ at $\big(\pm r_k, 0\big)$ yield $\left[ \overline{W}_{n, s, R}\right]_2 \geq \overline{w}_{*, n}$ on $T_{r_k}$. See (\ref{bar W_n}).   If $\overline{w}_{*, n} \rightarrow \overline{w}_{*, \infty}$ as $n \rightarrow \infty$, then by (\ref{Z_*})--(\ref{K_n}), we also have  \begin{align}\label{low bou interior}\big[\h{1pt} K_{n, R}\left[W\right] \h{1pt}\big]_2 = \overline{w}_{*, n} + y_{*; 1} R\h{1pt}\dfrac{w - Z_*}{\big| w - Z_* \big| \vee R} \h{2pt}\geq\h{2pt} \overline{w}_{*, n} \h{15pt}\text{on $T_{r_k}$.}\end{align} Noticing that at $\big(\pm r_k, 0 \big)$, $W_{\infty; 2} \geq \overline{w}_{*, \infty}$ and $y_* \cdot W_\infty = - h_\infty$, we obtain \begin{align}\label{low boun W infty}\big[\h{1pt} K_{n, R}\left[W_\infty\right] \h{1pt}\big]_2 = \overline{w}_{*, n} + y_{*; 1} R\h{1pt}\dfrac{w_\infty - Z_*}{\big| w_\infty - Z_* \big| \vee R} \h{2pt}\geq\h{2pt} \overline{w}_{*, n} \h{15pt}\text{at $\big(\pm r_k, 0\big)$,} \end{align} where \begin{align*}w_\infty = \dfrac{W_{\infty; 2}}{y_{*; 1}} + h_\infty \dfrac{y_{*;2}}{y_{*; 1}} \h{2pt}\geq\h{2pt} \dfrac{\overline{w}_{*, \infty}}{y_{*; 1}} + h_\infty \dfrac{y_{*;2}}{y_{*; 1}} \h{15pt}\text{at $\big(\pm r_k, 0\big)$.}\end{align*}In light of (\ref{low bou interior})--(\ref{low boun W infty}) and the fact that $W_{n; 2} \geq \overline{w}_{*, n}$ at $\big(\pm r_k, 0 \big)$, it then follows by the definition of $\overline{W}_{n, s, R}$ in (\ref{bar W_n}) that $\left[ \overline{W}_{n, s, R}\right]_2 \geq \overline{w}_{*, n}$ on $T_{r_k}$. All the above arguments yield $\overline{W}_{n, s, R} \in H_{k, \overline{w}_{*,n}}$. Here we extend to define $\overline{W}_{n, s, R} = W_n$ on $D_{\sigma_k} \setminus D_{r_k}$. \vspace{0.2pc}

By trace theorem, $W_n$ converges to $W_\infty$ strongly in $L^2\big(\p D_{r_k}; \mathbb{R}^2\big)$ as $n \to \infty$. In light of this convergence and (\ref{bounds for W_n}), similar arguments for (\ref{conv. of grad v_n, origin}) can be applied to get \begin{eqnarray} \int_{D_{r_k}} \big|\h{1pt}D_\xi \overline{W}_{n, s, R}\h{1pt}\big|^2 \longrightarrow \int_{D_{r_k}} \big|\h{1pt}D_\xi W\h{1pt}\big|^2, \h{20pt}\text{as $n \rightarrow \infty$, $R \rightarrow \infty$ and $s \rightarrow 0$, successively.}
\label{convergence of Dirichlet}
\end{eqnarray}As for the potential term, by (\ref{K_n})--(\ref{bar W_n}), it turns out
\begin{eqnarray*} L_n \int_{D_{r_k}} \left(h_n + y_* \cdot \overline{W}_{n, s, R}\right)^2 &=& L_n \int_{D_{r_k} \setminus D_{(1 - s) r_k}} \left(\dfrac{|\h{1pt} \xi \h{1pt}| - (1 - s) \h{1pt}r_k}{s \h{1pt}r_k}\right)^2\left(h_n + y_* \cdot W_n \h{2pt} \Big|_{r_k \widehat{\xi}}\right)^2\\[3mm]
&=& L_n \int_{(1 - s)r_k}^{r_k} \left(\dfrac{\tau - (1 - s) \h{1pt}r_k}{s \h{1pt}r_k}\right)^2 \dfrac{\tau}{r_k}\h{2pt} \mathrm{d} \tau \int_{\p D_{r_k}} \Big(h_n + y_* \cdot W_n \Big)^2.
\end{eqnarray*}Applying (\ref{bounds for W_n}) to the last estimate yields \begin{eqnarray}L_n \int_{D_{r_k} } \left(h_n + y_* \cdot \overline{W}_{n, s, R}\right)^2 \h{2pt}\lesssim\h{2pt}s\h{1pt}B_k, \h{20pt}\text{for all $n$.}
\label{convergence of potential}
\end{eqnarray}By the minimality of $W_n$, it satisfies \begin{eqnarray*}\int_{D_{r_k}}
 \big|\h{1pt} D_\xi W_n \h{1pt}\big|^2 + 2 L_n \mu \h{1pt}\big(h_n + y_* \cdot W_n \big)^2 \leq \int_{D_{r_k}}
 \big|\h{1pt} D_\xi \overline{W}_{n, s, R} \h{1pt}\big|^2 + 2L_n \mu \h{1pt}\big(h_n + y_* \cdot \overline{W}_{n, s, R} \big)^2.
\end{eqnarray*}Applying (\ref{convergence of Dirichlet})--(\ref{convergence of potential}) to the last estimate infers \begin{eqnarray} \label{str cov}\int_{D_{r_k}} \big|\h{1pt}D_\xi W_\infty\h{1pt}\big|^2 &\leq& \liminf_{n \rightarrow \infty} \int_{D_{r_k}} e_{L_n, h_n}^\star\big[ W_n\big] \h{2pt}\leq\h{2pt} \limsup_{n \rightarrow \infty}\int_{D_{r_k}}
 e_{L_n, h_n}^\star \big[ W_n \big] \nonumber\\[3mm]
&\leq&\lim_{s \rightarrow 0}\h{2pt}\lim_{R \rightarrow \infty} \h{2pt} \limsup_{n \rightarrow \infty} \int_{D_{r_k}}
 e_{L_n, h_n}^\star\big[ \overline{W}_{n, s, R} \big] = \int_{D_{r_k}} \big|\h{1pt}D_\xi W\h{1pt}\big|^2. 
\end{eqnarray}Here the energy density $e_{L_n, h_n}^\star$ has been defined in  (\ref{E_L^star}). Taking $W = W_\infty$ in the above estimate, we have $W_n$ converging to $ W_\infty$ strongly in $H^1\big(D_{r_k}; \mathbb{R}^2\big)$ as $n \to \infty$. Moreover, the potential term $ L_n \big(h_n + y_* \cdot W_n\big)^2$ converges to $0$ strongly in $L^1\big(D_{r_k}\big)$ as $n \to \infty$. We then can take $n \rightarrow \infty$ in item (iii) of (\ref{contradiction assumption}) and obtain \begin{eqnarray}\label{ine for harmonicity} \int_{D_{1/8}} \big|\h{1pt}D_\xi W_\infty\h{1pt}\big|^2 = \int_{D_{1/4}} \big|\h{1pt}D_\xi W_\infty\h{1pt}\big|^2 \h{2pt} \geq \h{2pt} \dfrac{1}{16}.
\end{eqnarray}The first equality above infers that $W_\infty$ is a constant on $D_{1/4} \setminus D_{1/8}$. In light of (\ref{str cov}), $W_\infty$ is harmonic in the upper--half part of $D_{r_k}$. By the analyticity of harmonic functions and the symmetry of $W_\infty$ with respect to the $\xi_2$--variable, $W_\infty$ must be a constant map  throughout $D_{r_k}$. However, this is impossible due to the second inequality in (\ref{ine for harmonicity}). The proof finishes.

\addcontentsline{toc}{subsection}{\normalsize A.2 \h{5pt}Proof of (2) in Lemma \ref{decay for large L}}
\subsection*{A.2 \h{5pt} Proof of (2) in Lemma \ref{decay for large L}}

We prove (\ref{linear growth}) in Lemma \ref{decay for large L}. Without loss of generality, we assume $\overline{w}_* = 0$. The proof in this section is inspired by the penalization method used in \cite{PSU12} for the scalar Signorini obstacle problem. Some necessary modifications are made in order to estimate solutions of our vectorial Signorini obstacle problem. The following arguments are divided into three steps.\vspace{0.2pc}
\\
\textbf{Step 1.} Let $\beta_{\epsilon} = \beta_\epsilon\left(s\right)$ be a smooth real--valued function on $\mathbb{R}$ such that
\begin{equation}
\beta_{\epsilon} \leq 0\,, \,\,\h{10pt} \dfrac{ \mathrm{d} \h{1pt} \beta_{\epsilon}}{ \mathrm{d} \h{1pt} s } \geq 0 \, ,\,\, \h{10pt}\beta_{\epsilon}(s) = 0 \h{3pt}\text{ for } s \geq 0 \, , \,\,\h{10pt} \beta_{\epsilon}(s) = \epsilon + \dfrac{s}{\epsilon} \h{3pt}\text{ for } s \leq -2\h{1pt}\epsilon^2.
    \label{def of beta}
\end{equation}
In terms of $\beta_\epsilon$,
\begin{equation}
B_\epsilon(t) := 2\int_0^{t} \beta_{\epsilon}(s) \h{1pt} \mathrm{d} \h{1pt} s, \h{15pt}\text{for all $t \in \mathbb{R}$.}
\label{B_e}
\end{equation}By (\ref{def of beta}), the function $B_\epsilon \geq 0$ on $\mathbb{R}$. Using $B_\epsilon$, we define   $$E^\star_{\epsilon, L, h} [\h{0.5pt}v\h{0.5pt}] := E^\star_{L, h} [\h{0.5pt}v\h{0.5pt}] + \int_{T_{\sigma_k}} B_\epsilon(v_2), \h{15pt}\text{for any $v = \big(v_1, v_2\big) \in H^1\big(D_{\sigma_k}; \mathbb{R}^2\big)$.}$$
Direct method of calculus of variation infers that there is a unique minimizer, denoted by $u^\epsilon$, of the energy functional $E^\star_{\epsilon, L, h}$ in the configuration space: \begin{align*} \widetilde{H}_{k, u} := \Big\{ v \in H^1 \big(D_{\sigma_k}; \mathbb{R}^2\big) : \text{$v$ is even with respect to $\xi_2$--variable and} \h{2pt}v = u\h{3pt}\text{on $\p D_{\sigma_k}$} \Big\}.
\end{align*} Let $v$ be an arbitrary mapping in $H^1_0(D_{\sigma_k}; \mathbb{R}^2)$. Moreover, we assume that $v$ has even symmetry with respect to the variable $\xi_2$. Then the minimizing property of $u^\epsilon$ induces\begin{equation}
  \int_{D_{\sigma_k}}  D_\xi u^\epsilon : D_\xi v + 2\h{0.2pt}L\h{0.2pt}\mu \big(h + y_* \cdot u^\epsilon \big)\h{0.4pt}y_* \cdot v + \int_{T_{\sigma_k}} \beta_\epsilon \left(u^\epsilon_2\right) v_2 = 0,
\label{el eq for u_e}
\end{equation}
where $A:B=\sum_{i,j}A_{ij}B_{ij}$ denotes the matrix inner product. Since $u \in \widetilde{H}_{k, u}$ and $u_2 \geq 0 $ on $T_{\sigma_k}$, the minimizing property of $u^\epsilon$ can also infer
\begin{equation}
 E_{L, h}^\star\left[ u^\epsilon \right] \h{2pt}\leq\h{2pt} E^\star_{\epsilon, L, h}\left[ u^\epsilon \right] \h{2pt}\leq\h{2pt}E^\star_{\epsilon, L, h}\left[ u \right] \h{2pt}=\h{2pt}E^\star_{L, h} \left[ u \right] \h{2pt} \leq \h{2pt} 1.
\label{comparison of u^epsilon and u_L}
\end{equation}In addition, one can apply Poincar\'{e}'s inequality to get
$$  \int_{D_{\sigma_k}}  | \h{1pt} u^\epsilon \h{1pt} |^2
\h{2pt}\lesssim \h{2pt} \int_{D_{\sigma_k}}  | \h{1pt} u \h{1pt} |^2 + | \h{1pt} u^\epsilon - u \h{1pt} |^2
\h{2pt}\lesssim\h{2pt} \int_{D_{\sigma_k}}  | \h{1pt} u \h{1pt} |^2 + \int_{D_{\sigma_k}} \big| \h{1pt} D u^\epsilon - D u \h{1pt} \big|^2.$$
Therefore, $u^\epsilon$ is uniformly bounded in $H^1\big(D_{\sigma_k}; \mathbb{R}^2\big)$ by the above estimate and (\ref{comparison of u^epsilon and u_L}). We can extract a subsequence, still denoted by $\big\{u^\epsilon\big\}$, so that $u^\epsilon$ converges to some $u^\dagger \in \widetilde{H}_{k, u}$ weakly in $H^1\big(D_{\sigma_k}; \mathbb{R}^2\big)$ and strongly in $L^2\big(D_{\sigma_k}; \mathbb{R}^2\big)$ as $\epsilon \to 0$.\vspace{0.2pc}\\
\textbf{Step 2.} We claim that $u^\dagger = u$. Firstly we show $u^\dagger_2 \geq 0$ on $T_{\sigma_k}$. By the upper bound of $E^\star_{\epsilon, L, h}\left[ u^\epsilon \right]$ in (\ref{comparison of u^epsilon and u_L}), it holds \begin{align*} 
\int_{T_{\sigma_k}}B_\epsilon\left( u^\epsilon_2 \right) \h{2pt}\leq \h{2pt}1.
\end{align*}Now we fix a $\delta > 0$ and let $\epsilon > 0$ small enough so that $\delta > 2 \h{0.4pt}\epsilon^2$. The last estimate and (\ref{def of beta}) then yield \begin{align*} \int_{T_{\sigma_k} \h{0.4pt}\cap\h{0.8pt} \left\{ u^\epsilon_2 \h{1pt}\leq\h{1pt} - \delta\right\}}  \int_{- 2 \epsilon^2}^{u^\epsilon_2}  \beta_\epsilon\left(s\right) \mathrm{d} s \h{2pt}=\h{2pt}\int_{T_{\sigma_k} \h{0.4pt}\cap\h{0.8pt} \left\{ u^\epsilon_2 \h{1pt}\leq\h{1pt} - \delta\right\}} \epsilon \h{0.5pt}  u^\epsilon_2  + \dfrac{1}{2 \epsilon}  \big(u^\epsilon_2\big)^2    \h{2pt}\leq\h{2pt}\dfrac{1}{2},
\end{align*}which furthermore induces \begin{align*}\delta^2 \h{1pt}\mathscr{H}^1 \Big\{ T_{\sigma_k} \h{0.4pt}\cap\h{0.8pt} \big\{ u^\epsilon_2 \h{1pt}\leq\h{1pt} - \delta\big\}\Big\} \h{2pt}\leq\h{2pt} \int_{T_{\sigma_k} \h{0.4pt}\cap\h{0.8pt} \left\{ u^\epsilon_2 \h{1pt}\leq\h{1pt} - \delta\right\}}  \big(u^\epsilon_2\big)^2    \h{2pt}\leq\h{2pt}\epsilon  - 2 \epsilon^2 \int_{T_{\sigma_k} \h{0.4pt}\cap\h{0.8pt} \left\{ u^\epsilon_2 \h{1pt}\leq\h{1pt} - \delta\right\}}  u^\epsilon_2.
\end{align*}Here $\mathscr{H}^1$ is the one--dimensional Hausdorff measure. Utilizing the uniform boundedness of $u^\epsilon_2$ in $L^2\left(T_{\sigma_k}\right)$ and the almost everywhere convergence of $u^\epsilon_2$ to $u^\dagger_2$ on $T_{\sigma_k}$, we then can take $\epsilon \to 0$ in the above estimate and obtain $\mathscr{H}^1 \Big\{ T_{\sigma_k} \h{0.4pt}\cap\h{0.8pt} \big\{ u^\dagger_2 \h{1pt}\leq\h{1pt} - \delta\big\}\Big\} = 0$. Since $\delta > 0$ is arbitrary, it turns out $\mathscr{H}^1 \Big\{ T_{\sigma_k} \h{0.4pt}\cap\h{0.8pt} \big\{ u^\dagger_2 \h{1pt}\leq\h{1pt} 0 \big\}\Big\} = 0$. Equivalently $u^\dagger_2 \geq 0$ almost everywhere on $T_{\sigma_k}$. Therefore, $u^\dagger \in H_{k, 0, u}$. See (\ref{min. problem S_L}). Now we take $\epsilon \to 0$ in (\ref{comparison of u^epsilon and u_L}). By the lower semi--continuity, it holds $E^\star_{L, h}\left[ u^\dagger\right] \leq E^\star_{L, h}\left[ u \right]$. Notice the convexity of the energy functional $E^\star_{L, h}$. The minimizer of $E^\star_{L, h}$ in $H_{k, 0, u}$ is unique. Hence $u^\dagger = u$ in $D_{\sigma_k}$. \vspace{0.2pc}
\\
\textbf{Step 3.} Let $v \in H_0^1\big(D_{3/4}\big)$ and $\widetilde{e}_1 = \big(1, 0\big)^\top$. Moreover, $v$ is even with respect to the $\xi_2$--variable. Note that we have taken $\sigma_k \in (7/8, 1)$. Hence, if we let $\tau > 0$ be small enough,  then we can plug $v\big(\cdot - \tau \widetilde{e}_1\big)$ into (\ref{el eq for u_e}) and obtain \begin{align*} \int_{D_{\sigma_k}}  D_\xi \left( u^\epsilon\big(\xi + \tau \widetilde{e}_1\big)\right) : D_\xi v   + 2\h{0.2pt}L\h{0.2pt}\mu \left(h + y_* \cdot u^\epsilon\big(\xi + \tau \widetilde{e}_1\big) \right)\h{0.4pt}y_* \cdot v + \int_{T_{\sigma_k}} \beta_\epsilon \left(u^\epsilon_2\big(\xi + \tau \widetilde{e}_1\big)\right) v_2  = 0.
\end{align*}Subtracting the equation (\ref{el eq for u_e}) from the above, we get  \begin{align*}\int_{D_{\sigma_k}}  D_\xi \left( \dfrac{u^\epsilon\big(\xi + \tau \widetilde{e}_1\big) - u^\epsilon\left(\xi\right)}{\tau}\right) : D_\xi v  \h{2pt} &+  \h{2pt} 2\h{0.2pt}L\h{0.2pt}\mu \left( y_* \cdot \left(\dfrac{u^\epsilon\big(\xi + \tau \widetilde{e}_1\big) - u^\epsilon\left(\xi\right)}{\tau}\right) \right)\h{0.4pt}y_* \cdot v \\[2mm]
\h{2pt}& = \h{2pt}  - \int_{T_{\sigma_k}} \dfrac{\beta_\epsilon \left(u^\epsilon_2\big(\xi + \tau \widetilde{e}_1\big)\right) - \beta_\epsilon\left(u^\epsilon_2\right)}{\tau} v_2.
\end{align*}Now we let $\eta$ be a test function compactly supported in $D_{3/4}$. Moreover, $\eta \equiv 1$ in $D_{1/2}$ and is even with respect to the $\xi_2$--variable. Then we take $v = \frac{u^\epsilon\left(\xi + \tau \widetilde{e}_1\right) - u^\epsilon\left(\xi\right)}{\tau} \eta^2$ in the above estimate. In light of the monotonicity of $\beta_\epsilon$, it then turns out  \begin{align*}\int_{D_{\sigma_k}}& \left| D_\xi \left( \dfrac{u^\epsilon\big(\xi + \tau \widetilde{e}_1\big) - u^\epsilon\left(\xi\right)}{\tau}\right) \right|^2 \eta^2  \h{2pt} +  \h{2pt} 2\h{0.2pt}L\h{0.2pt}\mu \left( y_* \cdot \left(\dfrac{u^\epsilon\big(\xi + \tau \widetilde{e}_1\big) - u^\epsilon\left(\xi\right)}{\tau}\right) \right)^2 \eta^2 \\[3mm]
\h{2pt}&\h{60pt} \leq \h{2pt}  - 2 \int_{D_{\sigma_k}} \eta \h{1pt} D_\xi \eta \cdot D_\xi\left( \dfrac{u^\epsilon\left(\xi + \tau \widetilde{e}_1 \right) - u^\epsilon\left(\xi \right)}{\tau}\right) \cdot \dfrac{u^\epsilon \left(\xi + \tau \widetilde{e}_1\right) - u^\epsilon\left(\xi\right)}{\tau}.
\end{align*}By the above estimate and the uniform bounds in (\ref{comparison of u^epsilon and u_L}), it satisfies \begin{align*} \int_{D_{1/2}} \big| D_\xi D_{\xi_1} u^\epsilon \big|^2  \h{2pt} +  \h{2pt} 2\h{0.2pt}L\h{0.2pt}\mu \big( y_* \cdot D_{\xi_1} u^\epsilon\big)^2 \h{2pt}\lesssim\h{2pt} 1. \end{align*}Here we have taken $\tau \rightarrow 0^+$. Now we let $\epsilon \rightarrow 0$. By Step 2, the above estimate induces \begin{align}\label{bound of u first two est} \int_{D_{1/2}} \big| D_\xi D_{\xi_1} u  \big|^2  \h{2pt} +  \h{2pt} 2\h{0.2pt}L\h{0.2pt}\mu \big( y_* \cdot D_{\xi_1} u \big)^2 \h{2pt}\lesssim\h{2pt} 1.
\end{align} It can be shown that \begin{align*} D_{\xi_2}^2 u = - D_{\xi_1}^2 u + 2\h{0.2pt}L\h{0.2pt}\mu\big(h + y_*\cdot u\big)\h{0.3pt}y_*,\h{15pt}\text{in $D_{1/2}^+$.}
\end{align*}We can obtain by (\ref{comparison of u^epsilon and u_L})--(\ref{bound of u first two est}) that \begin{align*} \int_{D_{1/2}^+} \left| D_{\xi_2}^2 u \right|^2 \h{2pt}\lesssim\h{2pt}\int_{D_{1/2}^+} \left| D_{\xi_1}^2 u \right|^2  + L^2 \int_{D_{1/2}^+} \big(h + y_* \cdot u \big)^2 \h{2pt}\lesssim\h{2pt}1 + L.
\end{align*} The proof then finishes by the above estimate, (\ref{bound of u first two est}), the uniform bound of $E^\star_{L, h}\left[u\right]$ in (\ref{comparison of u^epsilon and u_L}), H\"{o}lder's inequality and Sobolev's inequality.


\addcontentsline{toc}{subsection}{\normalsize A.3 \h{5pt}Some lemmas used in Part II}
\subsection*{A.3 \h{5pt} Some lemmas used in Part II}

For the convenience of readers, we list some lemmas used in Part II. Except otherwise stated, the vector field $w_a$ in the following refers to either the biaxial solutions $w_{a, b}^+$ or the split--core solutions $w_{a, c}^-$.

\begin{lem}[\bf Monotonicity formula] 
For any $B_R (y) \subset B_1$, it satisfies
\begin{align*}\dfrac{\p}{\p R} \left(\dfrac{1}{R} \int_{B_R(y)} f_{a,\mu}(\h{.5pt}w_a\h{.5pt}) \right)
= \dfrac{2}{R}\int_{\p B_R(y)} \left|\dfrac{\p w_a}{\p \vec{n} }\right|^2
+ \dfrac{2}{R^2} \int_{B_R(y)}F_a[\h{.5pt}w_a\h{.5pt}] \h{2pt}\geq\h{2pt} 0.\end{align*} 
Here $F_a \left[w\right] := \mu\left[D_{a}-3\sqrt{2}\h{0.5pt}S\left[ w \right]+\dfrac{a}{2}\h{0.5pt}\big(\h{0.5pt}|\h{0.2pt}w\h{0.2pt}|^2-1\big)^2\right]$. The notion $\vec{n}$ is the outer--normal direction on $\p B_R(y)$.
\label{Energy Monotonicity}
\end{lem}
\begin{lem}[\bf Uniform convergence of $\big| w_a \big|$ away from singularities and poles] Suppose that there is a sequence $\big\{a_n\big\}$ tending to $\infty$ as $n \to \infty$. In addition, we assume that there is $w^\sharp \in H^1(B_1; \mathbb{S}^4)$ so that $w_{a_n}$ converges to $w^\sharp$ strongly in $H^1(B_1; \mathbb{R}^5)$ as $n \to \infty$. Then for any compact set $K \subset B_1$ on which $w^\sharp$ is smooth, the modulus $\big|w_{a_n}\big|$ converges to $1$ in $C^0(K)$ as $n \to \infty$.
\label{|w_a|>1/2, away from sing.}
\end{lem}
\begin{lem}[\bf Local gradient estimate] There exist two universal positive constants $\epsilon_\star $ and $r_\star $ such that the following holds.\vspace{0.2pc}

 Let $U$ be an open set in $B_1$ satisfying $U \subset \overline{U} \subset B_1$.  In addition, we assume that $1/2 \leq \big|w_a\big| \leq H_a$ on $U$.  If it satisfies 
\begin{equation*}
\dfrac{1}{r}\int_{B_r(y)} f_{a,\mu} (\h{.5pt}w_a\h{.5pt}) \leq \epsilon_\star, \h{15pt}\text{for some $y \in U$ and $0 < r < \min\Big\{r_\star, \h{3pt}\textup{dist}\big(y, \p U\big)\Big\}$, }
\end{equation*}
then we have
\begin{align*}
r^2 \sup_{B_{r/2}(y)} f_{a,\mu} (\h{.5pt}w_a\h{.5pt}) \leq 144.
\end{align*}
\label{energy density e^L_a[w_a] bdd}
\end{lem}\vspace{-0.3pc}
\noindent Lemmas \ref{Energy Monotonicity}--\ref{energy density e^L_a[w_a] bdd} can be proved by following the arguments in \cite{MZ10}. We omit their proofs. \vspace{0.2pc}

In the next, we provide a boundary monotonicity formula for $w_a$ near the north pole $N_0 = (0, 0, 1)^\top$.
\begin{lem}\label{boundary monotonicity formula} Define \begin{align*} \mathcal{E}_{a; \h{1pt}N_0, \h{0.3pt}r} := \dfrac{1}{r} \int_{B_r\left(N_0\right)\h{0.5pt}\cap \h{0.5pt} B_1} f_{a, \mu}\big(w_a\big).
\end{align*} Then it holds  \begin{align*}\dfrac{\mathrm{d}}{\mathrm{d} r} \h{2pt}\mathcal{E}_{a; \h{1pt}N_0, \h{0.3pt}r} \h{1pt}\geq\h{1pt} - 33 \pi H_a^2, \h{15pt}\text{for any $r \in (0, 1)$.}
\end{align*}
\end{lem}\begin{proof}[\bf Proof] In light of the Dirichlet boundary condition in (\ref{boundary condition of u_a, intro}), one can apply Pohozaev identity associated with the system (\ref{el-eq of u_a, intro}) to obtain \begin{align*}\dfrac{\mathrm{d}}{\mathrm{d} r} \h{2pt}\mathcal{E}_{a; \h{1pt}N_0, \h{0.3pt}r} &\h{1pt}=\h{1pt}\dfrac{2}{r^2} \int_{B_r\left(N_0\right)\h{0.5pt}\cap \h{0.5pt}B_1} F_a\left[w_a\right] + \dfrac{2}{r} \int_{B_1 \h{0.5pt}\cap\h{0.5pt}\p B_r\left(N_0\right)} \left|\h{1pt} \dfrac{\p \h{0.2pt} w_a}{\p \h{0.2pt} \vec{n}}\h{1pt}\right|^2 + \dfrac{1}{r^2} \int_{B_r\left(N_0\right) \h{1pt}\cap\h{1pt}\p B_1} \big(1 - z \big)  \left|\h{1pt} \dfrac{\p \h{0.2pt} w_a}{\p \h{0.2pt}\vec{n}}\h{1pt}\right|^2\\[2mm]
&- \dfrac{6 H_a^2}{r^2}\int_{B_r\left(N_0\right) \h{0.5pt}\cap\h{0.5pt}\p B_1} \big(1 - z \big) + \dfrac{2}{r^2} \int_{B_r\left(N_0\right)\h{1pt}\cap\h{1pt}\p B_1}\sin \phi \h{2pt} \dfrac{\p \h{0.2pt} w_a}{\p \h{0.2pt}\phi}\cdot \dfrac{\p \h{0.2pt} w_a}{\p \h{0.2pt}\vec{n}}.
\end{align*}To derive the above identity, we have also used $\big| \p_\phi w_a \big|^2 = 3 H_a^2$ and $\big| \p_\theta w_a \big|^2 = 3 H_a^2 \sin^2 \phi$ on $\p B_1$. Since the first three terms on the right--hand side above are non--negative and in addition we have $\big| \sin \phi\big| \leq r$ on $B_r\left(N_0\right)\cap \p B_1$, it then turns out from the last identity that \begin{align}\label{low bound bdry mono}\dfrac{\mathrm{d}}{\mathrm{d} r} \h{2pt}\mathcal{E}_{a; \h{1pt}N_0, \h{0.3pt}r} &\h{1pt}\geq\h{1pt}- 9 \pi H_a^2      -   \int_{B_r\left(N_0\right)\h{1pt}\cap\h{1pt}\p B_1} \left| \dfrac{\p \h{0.2pt} w_a}{\p \h{0.2pt}\vec{n}} \right|^2.
\end{align}

On the other hand, by Lemma \ref{Energy Monotonicity}, \begin{align*} \dfrac{1}{R} \int_{\p B_R}f_{a, \mu}\big(w_a\big) \h{1pt}\geq\h{1pt}\dfrac{\p}{\p R} \left(\dfrac{1}{R} \int_{B_R} f_{a,\mu}(\h{.5pt}w_a\h{.5pt}) \right)
\h{1pt}\geq\h{1pt} \dfrac{2}{R}\int_{\p B_R} \left|\dfrac{\p\h{0.5pt} w_a}{\p \h{0.5pt}\vec{n} }\right|^2, \h{15pt}\text{for any $R \in (0, 1)$.}
\end{align*}We then can take $R \rightarrow 1^-$ and obtain \begin{align}\label{boun normal deri}\int_{\p B_1}\big| \nabla w_a \big|^2 \h{1pt}\geq\h{1pt}2\int_{\p B_1} \left|\dfrac{\p\h{0.5pt} w_a}{\p \h{0.5pt} \vec{n} }\right|^2.
\end{align}Here we have also used $F_a \left[ w_a \right] \equiv 0$ on $\p B_1$. Notice that $$\displaystyle \big| \nabla w_a \big|^2 = 6H_a^2 + \left| \frac{\p\h{0.5pt}w_a}{\p \h{0.5pt}\vec{n}}\right|^2 \h{15pt}\text{on $\p B_1$.}$$ Applying this identity to (\ref{boun normal deri}) yields \begin{align*}\int_{\p B_1} \left|\dfrac{\p\h{0.5pt} w_a}{\p \h{0.5pt} \vec{n} }\right|^2 \h{1pt}\leq\h{1pt}24 \pi H_a^2. \end{align*}The proof is completed by the last estimate and (\ref{low bound bdry mono}).
\end{proof}
With Lemmas \ref{Energy Monotonicity} and \ref{boundary monotonicity formula}, we have the following result for clearing singularities near the north pole. \begin{lem}[\bf Uniform lower bound of $\big| w_a \big|$ near north pole] Suppose that there is a sequence $\big\{a_n\big\}$ tending to $\infty$ as $n \to \infty$. Moreover, we assume that there is $w^\sharp \in H^1(B_1; \mathbb{S}^4)$ so that $w_{a_n}$ converges to $w^\sharp$ strongly in $H^1(B_1; \mathbb{R}^5)$ as $n \to \infty$. If $w^\sharp$ is smooth on $\overline{B_{r_0}\left(N_0\right)} \h{0.5pt}\cap\h{0.5pt}\overline{B_1}$ for some $r_0 \in (0, 1)$, then there is a $r_1 \in (0, r_0)$ so that  it satisfies $\big|w_{a_n}\big| \geq 1/2$ on $B_{r_1}\left(N_0\right)\h{1pt}\cap\h{1pt}B_1\h{1pt}\cap\h{1pt}l_z$ for  large $n$. Here $N_0$ still denotes the north pole. $l_z$ is the $z$--axis.
\label{|w_a|>1/2, near poles.}
\end{lem}\noindent The proof of this lemma follows by slightly modifying the proof of Proposition 4 in \cite{MZ10}. Note that here we need Lemma A.2 in \cite{BBH93} in combination with our  Lemmas  \ref{Energy Monotonicity} and \ref{boundary monotonicity formula}.\vspace{0.4pc}

Letting $w_*$ be an $\mathbb{S}^4$--valued mapping on $\p B_r$, we define\begin{align*}H \big(r,w_*\big) := \Big\{ w \in H^1(B_r ; \mathbb{S}^4)
\h{1pt} : \h{1pt} w=w_* \h{2pt} \text{ on $\p B_r$} \h{2pt} ,\h{2pt}  w = \mathscr{L}[u] \h{2pt} \text{ for some $3$--vector field } \h{2pt} u=u(\rho, z)\Big\}.\end{align*}We have the following  convergence result for a sequence of minimizers.
\begin{lem}Let $r > 0$ be an arbitrary radius. For any $n \in \mathbb{N}$, we suppose that $W^{(n)}$ minimizes the Dirichlet energy in $H\big(r, W^{(n)}\big)$ and satisfies \begin{align*}\sup_{n \h{1pt}\in\h{1pt}\mathbb{N}} \int_{\p B_{r}}  \big|\h{.5pt} \nabla_{\mathrm{tan}} W^{(n)} \h{.5pt}\big|^2 
< \infty.\end{align*}Here $\nabla_{\mathrm{tan}}$ is the tangential derivative on $\p B_r$.  If for some $\mathbb{S}^4$--valued mapping $W^\infty$ on $\p B_r$, it holds\begin{align*} 
 \int_{\p B_{r}}  \big|\h{.5pt} W^{(n)} - W^\infty \h{.5pt}\big|^2 \to 0 \h{5pt} \text{ as } n \to \infty,
\end{align*}then we have \begin{align*} \int_{B_r} \big|\h{1pt}\nabla W^{(n)} \h{1pt}\big|^2 \longrightarrow \int_{B_r} \big|\h{1pt}\nabla W^\infty \h{1pt}\big|^2\h{15pt}\text{as $n \to \infty$.} 
\end{align*}Here we still use $W^\infty$ in the convergence above to denote a Dirichlet energy minimizer in $H\big(r, W^\infty\big)$.\label{strong H1 convergence in H(r, w) space}
\end{lem}We can prove this lemma by following the proof of Convergence Theorem 5.5 in \cite{HLP92}. We omit it here.
\end{alphasection}
\
\\
\\
\noindent\textbf{\large Acknowledgement:} Ho--Man Tai would like to express his sincere gratitude to his Ph.D supervisor Yong Yu, for the insightful guidance and the kind invitation to join the current project. Part of the results in this article are contained in the Ph.D thesis of H.--M. Tai. Meanwhile, Y. Yu would like to express great gratitude to Professor Fanghua Lin, for his continuous support and constant encouragement while preparing the work. Y. Yu is partially supported by RGC grant of Hong Kong, Grant No. 14302718.

\vspace{2pc}
Department of Mathematics, The Chinese University of Hong Kong, Hong Kong \\[1mm]
\textit{E-mail address}: hmtai@math.cuhk.edu.hk\\[2mm]
Department of Mathematics, The Chinese University of Hong Kong, Hong Kong\\[1mm]
\textit{E-mail address}: yongyu@math.cuhk.edu.hk

\end{document}